\documentclass[twoside]{report}%
\usepackage{fancyhdr}
\usepackage{amsmath}
\usepackage{amsfonts}
\usepackage[arrow, matrix, curve]{xy}
\usepackage{pstricks,pst-node, pst-tree, pst-plot,pstcol}
\usepackage{colortbl}
\usepackage{amssymb}
\usepackage{graphicx}%
\usepackage{makeidx}
\newtheorem{theorem}{Theorem}[section]

\newtheorem{conjecture}[theorem]{Conjecture}
\newtheorem{corollary}[theorem]{Corollary}

\newtheorem{definition}[theorem]{Definition}
\newtheorem{example}[theorem]{Example}

\newtheorem{lemma}[theorem]{Lemma}

\newtheorem{remark}[theorem]{Remark}

\newenvironment{proof}[1][Proof]{\noindent\textbf{#1.} }{\ \rule{0.5em}{0.5em}}

\usepackage{helvet}
\usepackage[avantgarde]{quotchap}
\renewcommand\chapterheadstartvskip {\vspace*{-5\baselineskip}}
\usepackage{helvet}
\renewcommand\sectfont{\sffamily\bfseries}
\usepackage{makeidx} 
\makeindex
\begin{document}
\fancyhead{}
\fancyfoot{}
\author{\\\\\\\\An den Naturwissenschaftlich-Technischen Fakult\"{a}ten\\
der Universit\"{a}t des Saarlandes\\
zur Erlangung des Grades eines\\
Doktors der Naturwissenschaften (Dr. rer. nat.)\\
vorgelegte Abhandlung\\\\\\
von\\\\\\
Michael Sagraloff\\
geboren am 14. M\"{a}rz 1978 in Regensburg}
\title{Special Linear Series and Syzygies of Canonical Curves of Genus 9}
\maketitle
\newpage
\begin{abstract}
In this thesis we give a complete description of the syzygies of irreducible, nonsingular, canoncial curves $C$ of genus 9. This includes a collection of all possible Betti tables for $C$. Moreover a direct correspondence between these Betti tables and the number and types of special linear series on $C$ is given. Especially for $\operatorname{Cliff}(C)=3$ the curve $C$ is contained in determinantal surface $Y$ on a 4-dimensional rational normal scroll $X\subset \mathbb{P}^8$  constructed from a base point free pencil of divisors of degree 5.   
\end{abstract}
\newpage
\newpage
\pagenumbering{roman}
\setcounter{page}{0}
\fancyhead[RE,LO]{\thepage}
\cfoot{}
\tableofcontents
\newpage
°
\newpage
\noindent \textbf{Deutsche Zusammenfassung.} In der vorliegenden Arbeit werden kanonisch eingebettete Kurven $C\subset \mathbb{P}^{g-1}$ vom Geschlecht $g$ behandelt. Wir betrachten die minimal freie Aufl\"{o}sung des Koordinatenrings $S_{C}$, wobei $S=\Bbbk[x_0,...,x_{g-1}]$ den homogenen Koordinatenring des $\mathbb{P}^{g}$ bezeichnet. Diese Aufl\"{o}sung nimmt dann folgende Gestalt an:
\[
\overset{g-1}{\overbrace{%
\begin{tabular}
[c]{|p{0.7cm}|p{0.7cm}|p{0.7cm}|p{0.7cm}|p{0.7cm}|p{0.7cm}|p{0.7cm}|p{0.7cm}|}%
\hline
\multicolumn{1}{|>{\columncolor{gray}}c|}{1}  &  &  &  &  &  &  & \\\hline
&  \multicolumn{1}{|>{\columncolor{gray}}c|}{$\beta_{12}$}  &  \multicolumn
{1}{|>{\columncolor{gray}}c|}{$\cdots$}  &  \multicolumn{1}{|>{\columncolor
{gray}}c|}{ }  &  \multicolumn{1}{|>{\columncolor{gray}}c|}{ }  &  &  &
\\\hline
&  &  &  \multicolumn{1}{|>{\columncolor{gray}}c|}{ }  &  \multicolumn
{1}{|>{\columncolor{gray}}c|}{ }  &  \multicolumn{1}{|>{\columncolor{gray}}%
c|}{$\cdots$}  &  \multicolumn{1}{|>{\columncolor{gray}}c|}{$\beta_{g-3,g-1}$}
& \\\hline
&  &  &  &  &  &  &  \multicolumn{1}{|>{\columncolor{gray}}c|}{1} \\\hline
\end{tabular}
}}%
\]
\vspace{-0.78cm}
\[
\hspace{-6.3cm}\underset{p}{\underbrace{\hspace{3.4cm}}}%
\]
wobei $p\leq\frac{g-3}{2}$ die Anzahl der f\"{u}hrenden Nulleintr\"{a}ge in der dritten Zeile bezeichnet. Green$^{\prime}$s Vermutung besagt nun, dass ein direkter Zusammenhang zwischen dieser Anzahl $p$ und dem Auftreten spezieller Linearsysteme auf $C$ besteht: Er vermutet, dass $C$ genau dann Clifford Index $p$ besitzt, wenn genau die ersten $p$ Eintr\"{a}ge in der dritten Zeile des obigen Tableaus Nulleintr\"{a}ge sind. Der Clifford Index eines Divisors $D$ auf einer Kurve $C$ ist definiert als     
\[
\operatorname*{Cliff}(D)=\deg D-2(h^{0}\mathcal{O}_{C}(D)-1)
\]
und der Clifford Index von $C$ als%
\[
\operatorname*{Cliff}(C)=\operatorname{min}\left\{\operatorname*{Cliff}(D):D\in\operatorname{Div}(C)\text{ effektiv mit }h^{0}\mathcal{O}_{C}(D),h^{1}\mathcal{O}_{C}%
(D)\geq2\right\}
\]
Green and Lazarsfeld bewiesen in \cite{GreenLazarsfeld}, dass aus der Existenz spezieller Linearsysteme auf $C$ die Existenz von Extrasyzygien folgt. Neuere Arbeiten von Hirschowitz und Ramanan in
\cite{hirschowitz} (1998) und von Voisin in \cite{voisin3} zeigen, dass Green's Vermutung f\^{u}r Kurven von ungeradem Geschlecht $g$ und maximalem Clifford Index $\frac{g-1}{2}.$ gilt. 
\noindent Dar\"{u}ber hinaus gilt Green's Vermutung auch allgemein f\"{u}r Kurven vom Geschlecht $g\leq9$, was durch Ergebnisse von Mukai in \cite{mukai} (1995) und Schreyer in
$\cite{schreyer2}$ belegt wird.

Man kann nun einen Schritt weiter gehen und fragen, welche Betti Tableaus auftreten k\"{o}nnen im Falle einer irreduziblen, glatten, kanonischen Kurve $C.$ F\"{u}r Kurven vom Geschlecht $g\leq8$ ist diese Fragestellung bereits in \cite{schreyer1} (1986) beantwortet worden. Basierend auf Computerberechnungen gibt Schreyer eine Liste von  Betti Tableaus f\"{u}r Kurven vom Geschlecht $9$, $10$ and $11$ in \cite{schreyer3} (2003) an, deren Vollst\"{a}ndigkeit zu \"{u}berpr\"{u}fen ist. So vermutet der Autor bespielsweise, dass f\"{u}r Kurven vom Geschlecht $9$ aus der Existenz von drei $g_{5}^{1}$$^{\prime}$en bereits die eines $g_{7}^{2}$ folgt, was sich im Allgemeinen als nicht zutreffend erweist. 

\bigskip

Das Ergebnis dieser Arbeit ist eine vollst\"{a}ndige Liste von Betti Tableaus f\"{ur} glatte, irreduzible, kanonische Kurven vom Geschlecht $g=9$:\\
\vspace{5mm}
\footnotesize
\noindent
\begin{tabular}{|c@{\hspace{0.4mm}}|c@{\hspace{0.4mm}}|c@{\hspace{0.4mm}}|}
\hline
general & $\exists\; g^1_5$ & $\exists \hbox{ two } g^1_5$ \cr \hline
\begin{tabular}{l@{\hspace{1.8mm}}|c@{\hspace{1mm}}c@{\hspace{1mm}}c@{\hspace{1mm}}c@{\hspace{1mm}}c@{\hspace{1mm}}c@{\hspace{1mm}}c@{\hspace{1mm}}c@{\hspace{1mm}}c@{\hspace{1mm}}c@{\hspace{1mm}}c@{\hspace{1mm}}}
& 0 & 1 & 2 & 3 & 4 & 5 & 6 & 7 \cr \hline
0 & 1 & - & - & - & - & - & - & - \cr
1 & - & 21 & 64 & 70 & - & - & - & - \cr
2 & - & - & -  & -  & 70  & 64  & 21 & -\cr
3 & - & - & - & - & - & - & - & 1 \cr
\end{tabular} &
\begin{tabular}{l@{\hspace{1.8mm}}|c@{\hspace{1mm}}c@{\hspace{1mm}}c@{\hspace{1mm}}c@{\hspace{1mm}}c@{\hspace{1mm}}c@{\hspace{1mm}}c@{\hspace{1mm}}c@{\hspace{1mm}}c@{\hspace{1mm}}c@{\hspace{1mm}}c@{\hspace{1mm}}}
& 0 & 1 & 2 & 3 & 4 & 5 & 6 & 7 \cr \hline
0 & 1 & - & - & - & - & - & - & - \cr
1 & - & 21 & 64 & 70 & 4 & - & - & - \cr
2 & - & - & - & 4 & 70 & 64  & 21 & -\cr
3 & - & - & - & - & - & - & - & 1 \cr
\end{tabular} &
\begin{tabular}{l@{\hspace{1.8mm}}|c@{\hspace{1mm}}c@{\hspace{1mm}}c@{\hspace{1mm}}c@{\hspace{1mm}}c@{\hspace{1mm}}c@{\hspace{1mm}}c@{\hspace{1mm}}c@{\hspace{1mm}}c@{\hspace{1mm}}c@{\hspace{1mm}}c@{\hspace{1mm}}}
& 0 & 1 & 2 & 3 & 4 & 5 & 6 & 7 \cr \hline
0 & 1 & - & - & - & - & - & - & - \cr
1 & - & 21 & 64 & 70 & 8 & - & - & - \cr
2 & - & - & -   & 8 & 70 & 64  & 21 & -\cr
3 & - & - & - & - & - & - & - & 1 \cr
\end{tabular} \cr \hline
$\exists \hbox{ three } g^1_5$ & $\exists \,g^2_7$ & $\exists\, g^1_4$ \cr \hline
\begin{tabular}{l@{\hspace{1.8mm}}|c@{\hspace{1mm}}c@{\hspace{1mm}}c@{\hspace{1mm}}c@{\hspace{1mm}}c@{\hspace{1mm}}c@{\hspace{1mm}}c@{\hspace{1mm}}c@{\hspace{1mm}}c@{\hspace{1mm}}c@{\hspace{1mm}}c@{\hspace{1mm}}}
& 0 & 1 & 2 & 3 & 4 & 5 & 6 & 7 \cr \hline
0 & 1 & - & - & - & - & - & - & - \cr
1 & - & 21 & 64 & 70 & 12 & - & - & - \cr
2 & - & - & -   & 12 & 70 & 64 & 21 & -\cr
3 & - & - & - & - & - & - & - & 1 \cr
\end{tabular} &
\begin{tabular}{l@{\hspace{1.8mm}}|c@{\hspace{1mm}}c@{\hspace{1mm}}c@{\hspace{1mm}}c@{\hspace{1mm}}c@{\hspace{1mm}}c@{\hspace{1mm}}c@{\hspace{1mm}}c@{\hspace{1mm}}c@{\hspace{1mm}}c@{\hspace{1mm}}c@{\hspace{1mm}}}
& 0 & 1 & 2 & 3 & 4 & 5 & 6 & 7 \cr \hline
0 & 1 & - & - & - & - & - & - & - \cr
1 & - & 21 & 64 & 70 & 24 & - & - & - \cr
2 & - & - & -   & 24 & 70 & 64 & 21 & -\cr
3 & - & - & - & - & - & - & - & 1 \cr
\end{tabular} &
\begin{tabular}{l@{\hspace{1.8mm}}|c@{\hspace{1mm}}c@{\hspace{1mm}}c@{\hspace{1mm}}c@{\hspace{1mm}}c@{\hspace{1mm}}c@{\hspace{1mm}}c@{\hspace{1mm}}c@{\hspace{1mm}}c@{\hspace{1mm}}c@{\hspace{1mm}}c@{\hspace{1mm}}}
& 0 & 1 & 2 & 3 & 4 & 5 & 6 & 7 \cr \hline
0 & 1 & - & - & - & - & - & - & - \cr
1 & - & 21 & 64 & 75 & 24 & 5 & - & - \cr
2 & - & - & 5  & 24 & 75 & 64 & 21 & -\cr
3 & - & - & - & - & - & - & - & 1 \cr
\end{tabular} \cr \hline
$\exists \, g^1_4\times g^1_5$ & $\exists \,g^2_6$ & $\exists\, g^1_3$ \cr \hline
\begin{tabular}{l@{\hspace{1.8mm}}|c@{\hspace{1mm}}c@{\hspace{1mm}}c@{\hspace{1mm}}c@{\hspace{1mm}}c@{\hspace{1mm}}c@{\hspace{1mm}}c@{\hspace{1mm}}c@{\hspace{1mm}}c@{\hspace{1mm}}c@{\hspace{1mm}}c@{\hspace{1mm}}}
& 0 & 1 & 2 & 3 & 4 & 5 & 6 & 7 \cr \hline
0 & 1 & - & - & - & - & - & - & - \cr
1 & - & 21 & 64 & 75 & 44 & 5 & - & - \cr
2 & - & - & 5  & 44 & 75 & 64 & 21 & -\cr
3 & - & - & - & - & - & - & - & 1 \cr
\end{tabular} &
\begin{tabular}{l@{\hspace{1.8mm}}|c@{\hspace{1mm}}c@{\hspace{1mm}}c@{\hspace{1mm}}c@{\hspace{1mm}}c@{\hspace{1mm}}c@{\hspace{1mm}}c@{\hspace{1mm}}c@{\hspace{1mm}}c@{\hspace{1mm}}c@{\hspace{1mm}}c@{\hspace{1mm}}}
& 0 & 1 & 2 & 3 & 4 & 5 & 6 & 7 \cr \hline
0 & 1 & - & - & - & - & - & - & - \cr
1 & - & 21 & 64 & 90 & 64 & 20 & - & - \cr
2 & - & - & 20  & 64 & 90 & 64 & 21 & -\cr
3 & - & - & - & - & - & - & - & 1 \cr
\end{tabular} &
\begin{tabular}{l@{\hspace{1.8mm}}|c@{\hspace{1mm}}c@{\hspace{1mm}}c@{\hspace{1mm}}c@{\hspace{1mm}}c@{\hspace{1mm}}c@{\hspace{1mm}}c@{\hspace{1mm}}c@{\hspace{1mm}}c@{\hspace{1mm}}c@{\hspace{1mm}}c@{\hspace{1mm}}}
& 0 & 1 & 2 & 3 & 4 & 5 & 6 & 7 \cr \hline
0 & 1 & - & - & - & - & - & - & - \cr
1 & - & 21 & 70 & 105 & 84 & 35 & 6 & - \cr
2 & - & 6 & 35 & 84 & 105 & 70 & 21 & -\cr
3 & - & - & - & - & - & - & - & 1 \cr
\end{tabular} \cr\hline
\end{tabular}
\vspace{5mm}
\\{\normalsize \noindent Die Existenz mehrerer
$g_{5}^{1\prime}$e muss dabei so interpretiert werden, dass einige der Linearsysteme auch doppelt oder dreifach zu z\"{a}hlen sind. Es stellt sich heraus, dass eine allgemeine Kurve im Stratum, welches durch diese Betti Zahlen gegeben ist, eine entsprechende Anzahl unterschiedlicher $g_{5}^{1\prime}e$ besitzt.\\
Bei der Berechnung der Betti Tableaus orientieren wir uns an der Vorgehensweise Schreyers in \cite{schreyer1}: Mit Hilfe eines $g_{d}^{1}$ auf $C$ konstruieren wir einen $(d-1)$-dimensionalen Scroll $X\subset \mathbb{P}^{g-1}$, welcher die Kurve $C$ enth\"{a}lt. Die Rulings auf dem Scroll schneiden dann auf $C$ das entsprechende $g_{d}^{1}$ aus. Wir erhalten eine freie Aufl\"{o}sung von $\mathcal{O}_{C}$ als $\mathcal{O}_{X}$-Modul. Eine entsprechende Abbildungszylinderkonstruktion liefert dann eine (nicht unbedingt minimale) freie Aufl\"{o}sung von $S_{C}$. Hierbei treten schiefsymmetrische Matrizen $\psi$ mit Eintr\"{a}gen aus globalen Schnitten von Vektorb\"{u}ndeln $\mathcal{O}_{X}(aH+bR)$, $a,$$b\in \mathbb{Z}$, auf, welche die freie Aufl\"{o}sung bestimmen. Schlie\ss lich sind die R\"{a}nge eventuell auftretender, nicht minimaler Abbildungen zu berechnen.\\
Von besonderem Interesse ist der Fall $\operatorname{Cliff}(C)=3$: Im Falle der Existenz eines $g_{7}^{2}$ liegt die Kurve auf einer sogenannten Bordigatypfl\"{a}che, deren Betti Tableau bereits das der Kurve festlegt. Existiert jedoch kein $g_{7}^{2}$, so kann man in jedem Fall eine determinantielle Fl\"{a}che $Y\subset X$ auf dem Scroll angeben, welche die Kurve enth\"{a}lt: Hat $C$ genau ein $g_{5}^{1}$ mit einfacher Multiplizit\"{a}t, so ist $Y$ das Bild einer Aufblasung des $\mathbb{P}^2$ und im Falle der Existenz mehrerer $g_{5}^{1}$$^{\prime}$e das Bild einer Aufblasung von $\mathbb{P}^{1}\times \mathbb{P}^{1}$. Im Fall, dass ein $g_{5}^{1}$ mit h\"{o}herer Multiplizit\"{a}t auftritt, erhalten wir $Y$ als Bild einer Aufblasung der zweiten Hirzebruchfl\"{a}che.\\\\
Im Kapitel "Summary" geben wird ein kurzen Ausblick auf die Behandlung der F\"{a}lle $g=10$ und $g=11$. 
\newpage
\fancyhead[RO]{\rightmark}
\fancyhead[LE]{\leftmark}
\pagenumbering{arabic}

\setcounter{page}{0}

\setcounter{chapter}{-1} \setcounter{section}{0}

\chapter{{\protect\normalsize Introduction}} \index{$\mathcal{O}_{X}$ sheaf of regular functions on a variety $X$} \index{$\mathcal{I}_{X}$ ideal sheaf of a subvariety}
In this thesis we discuss curves $C\subset\mathbb{P}^{g-1}$ of geometric genus
$g=9$, which are embedded by the complete linear
series $\left\vert \omega_{C}\right\vert $ associated to the canonical bundle \index{$\omega_{X}$     canonical bundle on a variety $X$}
$\omega_{C}.$ It is known that for $C\subset\mathbb{P}^{n}$ embedded by a
very ample, complete linear series $\left\vert \mathcal{L}\right\vert $,
properties of the homogeneous coordinate ring \index{$S_{X}$     homogeneous coordinate ring of a variety $X$} \index{$I_{X}$   vanishing ideal of a variety $X$} $S_{C}=S/I_{C}$, as its graded Betti numbers depend both
on the curve itself and on $\mathcal{L}$ (cf. \cite{eisenbud2}) . Here we denote by \index{$S$   homogeneous coordinate ring of $\mathbb{P}^{n}$} $S=\Bbbk\lbrack
x_{0},...,x_{n}]=\operatorname*{Sym}H^{0}(C,\mathcal{L})$ the homogenous
coordinate ring \index{$\operatorname{Sym}$ symmetric algebra} of $\mathbb{P}^{n}.$ \index{$\mathbb{P}^{n}$ n-dimensional projective space} In the case
$\mathcal{L}=\omega_{C}$ we get an embedding of $C$ in $\mathbb{P}^{g-1}$ if
$C$ is non-hyperelliptic and properties of $S_{C}$ are intrinsic properties of
$C$ alone. For this reason one could ask, if it is possible to deduce
geometric properties of $C$ directly from the algebraic properties of $S_{C}.$
Hilbert gave us a first answer to this problem by introducing the Hilbert
polynomial, which contains information about the dimension and degree of the
embedding. For more detailed information of the embedded curve $C$ we have to
go further: Hilbert showed that there exists a minimal free resolution of
$S_{C}$ as $S-$module. This resolution, especially its Betti table, contains
further information, so we can ask in general which geometric
properties are encoded in the Betti numbers. From the Castelnuovo-Mumford
regularity of $S_{C}$ and its Gorenstein property we know that the Betti table
of the minimal free resolution of $S_{C}$ has the following form
\[
\overset{g-1}{\overbrace{%
\begin{tabular}
[c]{|p{0.7cm}|p{0.7cm}|p{0.7cm}|p{0.7cm}|p{0.7cm}|p{0.7cm}|p{0.7cm}|p{0.7cm}|}%
\hline
\multicolumn{1}{|>{\columncolor{gray}}c|}{1}  &  &  &  &  &  &  & \\\hline
&  \multicolumn{1}{|>{\columncolor{gray}}c|}{$\beta_{12}$}  &  \multicolumn
{1}{|>{\columncolor{gray}}c|}{$\cdots$}  &  \multicolumn{1}{|>{\columncolor
{gray}}c|}{ }  &  \multicolumn{1}{|>{\columncolor{gray}}c|}{ }  &  &  &
\\\hline
&  &  &  \multicolumn{1}{|>{\columncolor{gray}}c|}{ }  &  \multicolumn
{1}{|>{\columncolor{gray}}c|}{ }  &  \multicolumn{1}{|>{\columncolor{gray}}%
c|}{$\cdots$}  &  \multicolumn{1}{|>{\columncolor{gray}}c|}{$\beta_{g-3,g-1}$}
& \\\hline
&  &  &  &  &  &  &  \multicolumn{1}{|>{\columncolor{gray}}c|}{1} \\\hline
\end{tabular}
}}%
\]
\vspace{-0.78cm}
\[
\hspace{-6.3cm}\underset{p}{\underbrace{\hspace{3.4cm}}}%
\]
with a positive integer $p\leq\frac{g-3}{2}$, the number of leading zeros in the third
row. Following Green we say that $C$ fullfills \index{$N_{p}$ property} $N_{p}$ in this situation if there are at least $p$ zeroes. Green conjectured that there exists a direct correspondance
between the geometry of $C$ and the $N_{p}$ property:

\begin{conjecture}
(Green, 1984) Let $C\subset\mathbb{P}^{8}$ be a smooth nonhyperelliptic curve
over a field of characteristic 0 in its canonical embedding. Then%
\[
\beta_{p,p+2}\neq0\Leftrightarrow C\text{ has Clifford index }%
\operatorname*{Cliff}(C)\leq p
\]

\end{conjecture}

\noindent The Clifford index of an effective divisor $D$ on $C$ is defined as%
\[
\operatorname*{Cliff}(D)=\deg D-2(h^{0}\mathcal{O}_{C}(D)-1)
\]
and the Clifford index \index{$\operatorname{Cliff}$ Clifford index} of $C$ is%
\[
\operatorname*{Cliff}(C)=\operatorname{min}\left\{\operatorname*{Cliff}(D):D\in\operatorname{Div}(C)\text{ effective with }h^{0}\mathcal{O}_{C}(D),h^{1}\mathcal{O}_{C}%
(D)\geq2\right\}
\]

\noindent For a better understanding of this definition, we state some results
of Brill-Noether Theory, that studies the question whether there exists
a certain special linear series \index{$\left\vert D\right\vert$ complete linear series} $\left\vert D\right\vert $ on a curve $C$. 
\index{$W_{d}^{r}$ complete linear systems of degree $d$ and dimension $\geq r$} $C_{d}^{r}\subset\operatorname*{Div}C$ \index{$g_{d}^{r}$ element of $W_{r}^{d}$ with $r(D)=r$} as usual denotes the \index{$r(D)$ projective dimension of a $\left\vert D\right\vert$} variety \index{$C_{d}^{r}$ space of all divisors $D$ of degree $d$ with $r(D)\geq r$} of all
divisors $D$ of degree $d$ that fulfill $r(D)=\dim\left\vert D\right\vert \geq r.$
Further $g_{d}^{r}$ denotes an element of $W_{d}^{r}=\left\{\left\vert D \right\vert: \operatorname{deg}D=d,r(D)\geq r \right\}\subset \operatorname{Pic}^{d}(C)$. Brill Noether Theory (cf. \cite{arbarello}) then gives us a lower bound of the
dimension of this variety which is sharp for a generic curve. Together with
the Riemann Roch Theorem which says \index{$\operatorname{Pic}$ Picard group} \index{$\operatorname{Pic}^{d}$ Picard group of degree $d$} that $r=\dim\left\vert D\right\vert \geq
d-g$ and Clifford's Theorem we obtain the following picture:
\\\\
\bigskip%

\[
\psset{unit=0.8} \text{ \ \ \ \ \ \ \ \ \ }\begin{pspicture}(0,-1)(14,5)
\pscustom[linewidth=2pt,fillstyle=solid,fillcolor=gray]{
\psline[linewidth=2pt](0,0)(5,0)(9,4)(8,4)
}
\pscustom[linewidth=2pt,fillstyle=solid,fillcolor=white]{
\pscurve[linewidth=2pt](0,0)(5,1.4)(8,4)
}
\psline[linewidth=2pt](9,4)(10,5)
\psline[linewidth=2pt](0,0)(8,4)
\psline[linewidth=1.5pt]{->}(0,-0.2)(0,5.2)
\psline[linewidth=1.5pt]{->}(-0.2,0)(10.2,0)
\psline[linewidth=1.5pt]{-}(0,4)(-0.2,4)
\psline[linewidth=1.5pt]{-}(5,0)(5,-0.2)
\psline[linewidth=1.5pt]{-}(8,0)(8,-0.2)
\rput(1.5,2.5){$d\leq 2r$}
\rput(1.5,2.1){(Clifford)}
\rput*[Br](-0.3,3.9){$g-1$}
\rput*[Bl](8,4.2){$K_C$}
\rput[B](8,1){$r\geq d-g$}
\rput[B](8.5,0.6){(Riemann-Roch)}
\rput(5,4.2){dim$\;W^{r}_d \geq 0$}
\rput(5,3.8){(Brill-Noether)}
\psdots[dotstyle=*,linewidth=2pt](8,4)
\pscurve[linewidth=2pt]{->}(5,3.5)(4.5,2)(5,0.5)
\pscurve[linewidth=2pt]{->}(1.5,1.9)(1.4,1.3)(2,1.2)
\pscurve[linewidth=2pt]{<-}(6.2,0.7)(7,0.2)(8,0.4)
\rput[t](0,5.4){$r$}
\rput[Bl](10.3,-0.1){$d$}
\rput[t](5,-0.3){$g$}
\rput[t](8,-0.3){$2g-2$}
\rput[t](0,-0.3){$0$}
\rput[Br](-0.3,-0.1){$0$}
\end{pspicture}
\]

\noindent The gray shaded area gives us bounds for which $\dim W_{d}^{r}\geq0$
for a generic curve $C$ of genus $g.$ Clifford's Theorem says that for an
effective Divisor $D$ on $C$ of degree $d\leq2g-1$ we must have $r=h^{0}%
\mathcal{O}_{C}(D)-1\leq\frac{d}{2}.$ If equality holds then either $D$ is
zero, or $D$ is a canonical divisor, or $C$ is hyperelliptic and $D$ is linearly
equivalent to a multiple of a hyperelliptic divisor. 
\\In the last situation where $C$ is hyperelliptic the canonical map%
\[
j:C\rightarrow\mathbb{P}(H^{0}(C,\omega_{C}))=\mathbb{P}^{g-1}%
\]
is a $2:1$ map $\pi$ onto a rational normal curve which is a $(g-1)-$uple embedding
of $\mathbb{P}^{1}$ in $\mathbb{P}^{g-1}.$
\[
\text{\begin{xy}
\xymatrix{
C \ar[r]^{j} \ar[d]_{\pi} & \mathbb{P}^{g-1} \\
\mathbb{P}^1 \ar@{^{(}->}[ru]_{v_{g-1}}& \\
}
\end{xy}}%
\]
The sheaf $\pi_{*}\mathcal{O}_{C}$ is a rank two vector bundle $\mathcal{O}_{\mathbb{P}^{1}}\oplus \mathcal{O}_{\mathbb{P}^{1}}(a)$ on $\mathbb{P}^{1}$ with an $a\in \mathbb{Z}^{-}$. As $\pi^{*}\mathcal{O}_{\mathbb{P}^{1}}(1)\cong \mathcal{O}_{C}(D)$ we get $\pi^{*}\mathcal{O}_{\mathbb{P}^{1}}(2g-2)\cong \mathcal{O}_{C}(D)^{\otimes 2g-2}=\omega_{C}^{\otimes 2}$ and thus from the projection formula $\pi_{*}\pi^{*}\mathcal{O}_{\mathbb{P}^{1}}(2g-2)\cong \pi_{*}\omega_{C}^{\otimes 2}\cong \mathcal{O}_{\mathbb{P}^{1}}(2g-2) \otimes \pi_{*}\mathcal{O}_{C}$. Comparing dimensions of the global sections it follows that $a=-g-1$. Then $\Omega=\sum_{n\geq0}H^{0}(C,\omega_{C}^{\otimes n})$ regarded as a
module over the homogenous coordinate ring $S=\operatorname*{Sym}%
H^{0}(C,\mathcal{\omega}_{C})$ of $\mathbb{P}^{g-1}$ is the module of global
sections of the rank 2 vector bundle $j_{\ast}\mathcal{O}_{C}\cong
\mathcal{O}_{\mathbb{P}^{1}}\oplus\mathcal{O}_{\mathbb{P}^{1}}(-g-1)$ on the
rational normal curve $\mathbb{P}^{1}\subset\mathbb{P}^{g-1}.$ The Betti table for this rational normal curve is given as follows:
\[
\begin{tabular}{c|ccccccccc}
& 0 & 1 & 2 & $\cdots$ & i & $\cdots$ & g-2 \cr \hline
0 & 1 & - & - & - & - & - & - \cr
1 & - & $\binom{g-1}{2}$ & $2\binom{g-1}{3}$ & $\cdots$ & $i\binom{g-1}{i+1}$ & $\cdots$ & $g-2$ \cr
\end{tabular} 
\]
and as $\mathcal{O}_{\mathbb{P}^{1}}(-2)\cong\ \omega_{\mathbb{P}^{1}}$ the curve $C$ has the following Betti table:
\[
\begin{tabular}{c|cccccccccc}
& 0 & 1 & 2 & $\cdots$ & g-4 & g-3 & g-2 \cr \hline
0 & 1 & - & - & - & - & - & - \cr
1 & - & $\binom{g-1}{2}$ & $2\binom{g-1}{3}$ & $\cdots$ & $\cdots$ & $\cdots$ & $g-2$ \cr
2 & $g-2$ & $\cdots$ & $\cdots$ & $\cdots$ & $2\binom{g-1}{3}$ & $\binom{g-1}{2}$ & - \cr
3 & - & - & - & - & - & - & 1 \cr 
\end{tabular} 
\]
\\For non hyperelliptic $C$ the canonical map $j$ is an embedding and $\Omega$ is the
homogenous coordinate ring of $C\subset\mathbb{P}^{g-1}$ (cf. \cite{noether}%
). In this situation the definition of the Clifford index gives a natural
measure of the speciality of a divisor $D$ on $C.$

Green and Lazarsfeld proved in \cite{GreenLazarsfeld} that from the
existence of special linear systems on $C$ it follows the existence of
exceptional syzygies. Due to recent works of Hirschowitz and Ramanan in
\cite{hirschowitz} and Voisin in \cite{voisin3}, we know that
Green's Conjecture holds for curves $C$ of odd genus $g$ with maximal Clifford
index $\frac{g-1}{2}.$

\begin{theorem}\label{HRV}
(Hirschowitz-Ramanan-Voisin) Let $C$ be a smooth curve $C$ of genus
$g=2k+1\geq5$ with Betti number $\beta_{k,k+1}\neq0,$ then there exists a
$g_{k+1}^{1}$ on $C.$
\end{theorem}

\noindent Moreover for curves of genus $g\leq9$ Green's conjecture holds which
follows from results of Mukai in \cite{mukai} and Schreyer in
$\cite{schreyer2}$. From further results of Max Noether, Petri, Voisin and
Schreyer we know that curves which fullfill $N_{p}$ for $p\leq2$ have special
linear series of Clifford index $p.$ A substantial progress has been done in
\cite{voisin2} 2001 and \cite{voisin3} 2005, where Voisin proved the conjecture for general $k-$gonal curves of
arbritrary genus.

Going one step further one might ask which Betti tables actually occur for
irreducible, nonsingular curves $C.$ For curves of genus $g\leq8$ this is
already done in \cite{schreyer1} (1986). Based on computational evidence,
Schreyer gives a conjectural collection of Betti tables for curves of genus
$9$, $10$ and $11$ in \cite{schreyer3} (2003).

\bigskip

The result of this thesis is a complete table for smooth curves of genus 9. It
turns out that the conjectural table in \cite{schreyer3} is correct with the
exception that a curve $C$ of genus $9$ can even admit three linear systems of
type $g_{5}^{1}$ (counted with multiplicities) and no $g_{7}^{2}.$

\vspace{5mm}
\footnotesize
\noindent
\begin{tabular}{|c@{\hspace{0.4mm}}|c@{\hspace{0.4mm}}|c@{\hspace{0.4mm}}|}
\hline
general & $\exists\; g^1_5$ & $\exists \hbox{ two } g^1_5$ \cr \hline
\begin{tabular}{l@{\hspace{1.8mm}}|c@{\hspace{1mm}}c@{\hspace{1mm}}c@{\hspace{1mm}}c@{\hspace{1mm}}c@{\hspace{1mm}}c@{\hspace{1mm}}c@{\hspace{1mm}}c@{\hspace{1mm}}c@{\hspace{1mm}}c@{\hspace{1mm}}c@{\hspace{1mm}}}
& 0 & 1 & 2 & 3 & 4 & 5 & 6 & 7 \cr \hline
0 & 1 & - & - & - & - & - & - & - \cr
1 & - & 21 & 64 & 70 & - & - & - & - \cr
2 & - & - & -  & -  & 70  & 64  & 21 & -\cr
3 & - & - & - & - & - & - & - & 1 \cr
\end{tabular} &
\begin{tabular}{l@{\hspace{1.8mm}}|c@{\hspace{1mm}}c@{\hspace{1mm}}c@{\hspace{1mm}}c@{\hspace{1mm}}c@{\hspace{1mm}}c@{\hspace{1mm}}c@{\hspace{1mm}}c@{\hspace{1mm}}c@{\hspace{1mm}}c@{\hspace{1mm}}c@{\hspace{1mm}}}
& 0 & 1 & 2 & 3 & 4 & 5 & 6 & 7 \cr \hline
0 & 1 & - & - & - & - & - & - & - \cr
1 & - & 21 & 64 & 70 & 4 & - & - & - \cr
2 & - & - & - & 4 & 70 & 64  & 21 & -\cr
3 & - & - & - & - & - & - & - & 1 \cr
\end{tabular} &
\begin{tabular}{l@{\hspace{1.8mm}}|c@{\hspace{1mm}}c@{\hspace{1mm}}c@{\hspace{1mm}}c@{\hspace{1mm}}c@{\hspace{1mm}}c@{\hspace{1mm}}c@{\hspace{1mm}}c@{\hspace{1mm}}c@{\hspace{1mm}}c@{\hspace{1mm}}c@{\hspace{1mm}}}
& 0 & 1 & 2 & 3 & 4 & 5 & 6 & 7 \cr \hline
0 & 1 & - & - & - & - & - & - & - \cr
1 & - & 21 & 64 & 70 & 8 & - & - & - \cr
2 & - & - & -   & 8 & 70 & 64  & 21 & -\cr
3 & - & - & - & - & - & - & - & 1 \cr
\end{tabular} \cr \hline
$\exists \hbox{ three } g^1_5$ & $\exists \,g^2_7$ & $\exists\, g^1_4$ \cr \hline
\begin{tabular}{l@{\hspace{1.8mm}}|c@{\hspace{1mm}}c@{\hspace{1mm}}c@{\hspace{1mm}}c@{\hspace{1mm}}c@{\hspace{1mm}}c@{\hspace{1mm}}c@{\hspace{1mm}}c@{\hspace{1mm}}c@{\hspace{1mm}}c@{\hspace{1mm}}c@{\hspace{1mm}}}
& 0 & 1 & 2 & 3 & 4 & 5 & 6 & 7 \cr \hline
0 & 1 & - & - & - & - & - & - & - \cr
1 & - & 21 & 64 & 70 & 12 & - & - & - \cr
2 & - & - & -   & 12 & 70 & 64 & 21 & -\cr
3 & - & - & - & - & - & - & - & 1 \cr
\end{tabular} &
\begin{tabular}{l@{\hspace{1.8mm}}|c@{\hspace{1mm}}c@{\hspace{1mm}}c@{\hspace{1mm}}c@{\hspace{1mm}}c@{\hspace{1mm}}c@{\hspace{1mm}}c@{\hspace{1mm}}c@{\hspace{1mm}}c@{\hspace{1mm}}c@{\hspace{1mm}}c@{\hspace{1mm}}}
& 0 & 1 & 2 & 3 & 4 & 5 & 6 & 7 \cr \hline
0 & 1 & - & - & - & - & - & - & - \cr
1 & - & 21 & 64 & 70 & 24 & - & - & - \cr
2 & - & - & -   & 24 & 70 & 64 & 21 & -\cr
3 & - & - & - & - & - & - & - & 1 \cr
\end{tabular} &
\begin{tabular}{l@{\hspace{1.8mm}}|c@{\hspace{1mm}}c@{\hspace{1mm}}c@{\hspace{1mm}}c@{\hspace{1mm}}c@{\hspace{1mm}}c@{\hspace{1mm}}c@{\hspace{1mm}}c@{\hspace{1mm}}c@{\hspace{1mm}}c@{\hspace{1mm}}c@{\hspace{1mm}}}
& 0 & 1 & 2 & 3 & 4 & 5 & 6 & 7 \cr \hline
0 & 1 & - & - & - & - & - & - & - \cr
1 & - & 21 & 64 & 75 & 24 & 5 & - & - \cr
2 & - & - & 5  & 24 & 75 & 64 & 21 & -\cr
3 & - & - & - & - & - & - & - & 1 \cr
\end{tabular} \cr \hline
$\exists \, g^1_4\times g^1_5$ & $\exists \,g^2_6$ & $\exists\, g^1_3$ \cr \hline
\begin{tabular}{l@{\hspace{1.8mm}}|c@{\hspace{1mm}}c@{\hspace{1mm}}c@{\hspace{1mm}}c@{\hspace{1mm}}c@{\hspace{1mm}}c@{\hspace{1mm}}c@{\hspace{1mm}}c@{\hspace{1mm}}c@{\hspace{1mm}}c@{\hspace{1mm}}c@{\hspace{1mm}}}
& 0 & 1 & 2 & 3 & 4 & 5 & 6 & 7 \cr \hline
0 & 1 & - & - & - & - & - & - & - \cr
1 & - & 21 & 64 & 75 & 44 & 5 & - & - \cr
2 & - & - & 5  & 44 & 75 & 64 & 21 & -\cr
3 & - & - & - & - & - & - & - & 1 \cr
\end{tabular} &
\begin{tabular}{l@{\hspace{1.8mm}}|c@{\hspace{1mm}}c@{\hspace{1mm}}c@{\hspace{1mm}}c@{\hspace{1mm}}c@{\hspace{1mm}}c@{\hspace{1mm}}c@{\hspace{1mm}}c@{\hspace{1mm}}c@{\hspace{1mm}}c@{\hspace{1mm}}c@{\hspace{1mm}}}
& 0 & 1 & 2 & 3 & 4 & 5 & 6 & 7 \cr \hline
0 & 1 & - & - & - & - & - & - & - \cr
1 & - & 21 & 64 & 90 & 64 & 20 & - & - \cr
2 & - & - & 20  & 64 & 90 & 64 & 21 & -\cr
3 & - & - & - & - & - & - & - & 1 \cr
\end{tabular} &
\begin{tabular}{l@{\hspace{1.8mm}}|c@{\hspace{1mm}}c@{\hspace{1mm}}c@{\hspace{1mm}}c@{\hspace{1mm}}c@{\hspace{1mm}}c@{\hspace{1mm}}c@{\hspace{1mm}}c@{\hspace{1mm}}c@{\hspace{1mm}}c@{\hspace{1mm}}c@{\hspace{1mm}}}
& 0 & 1 & 2 & 3 & 4 & 5 & 6 & 7 \cr \hline
0 & 1 & - & - & - & - & - & - & - \cr
1 & - & 21 & 70 & 105 & 84 & 35 & 6 & - \cr
2 & - & 6 & 35 & 84 & 105 & 70 & 21 & -\cr
3 & - & - & - & - & - & - & - & 1 \cr
\end{tabular} \cr\hline
\end{tabular}
\vspace{5mm}
\\{\normalsize \noindent The interpretation of the existence of several
$g_{5}^{1\prime}s$ has to be taken as a count with multiplicity. It is true
that a general curve in the strata defined by these Betti numbers has that
many $g_{5}^{1\prime}s$ precisely. \newline\newline To obtain the Betti
numbers from a curve $C$ that admits special divisors $D$ of degree $d$ as in
Green's Conjecture we follow the approach in \cite{schreyer1}. We first
construct a rational normal scroll $X\subset\mathbb{P}^{g-1}$ from a pencil
$(D_{\lambda})_{\lambda}$ of divisors $D_{\lambda}\sim D:$%
\[
X=%
{\textstyle\bigcup\nolimits_{\lambda\in\mathbb{P}^{1}}}
\bar{D}_{\lambda}\subset\mathbb{P}^{g-1}%
\]
where $\bar{D}_{\lambda}$ denotes the linear span \index{$\bar{D}$ linear span of a divisor} of $D_{\lambda}.$ }

{\normalsize
\[
\raisebox{-0pt}{\includegraphics[width=6cm,height=6.5cm]{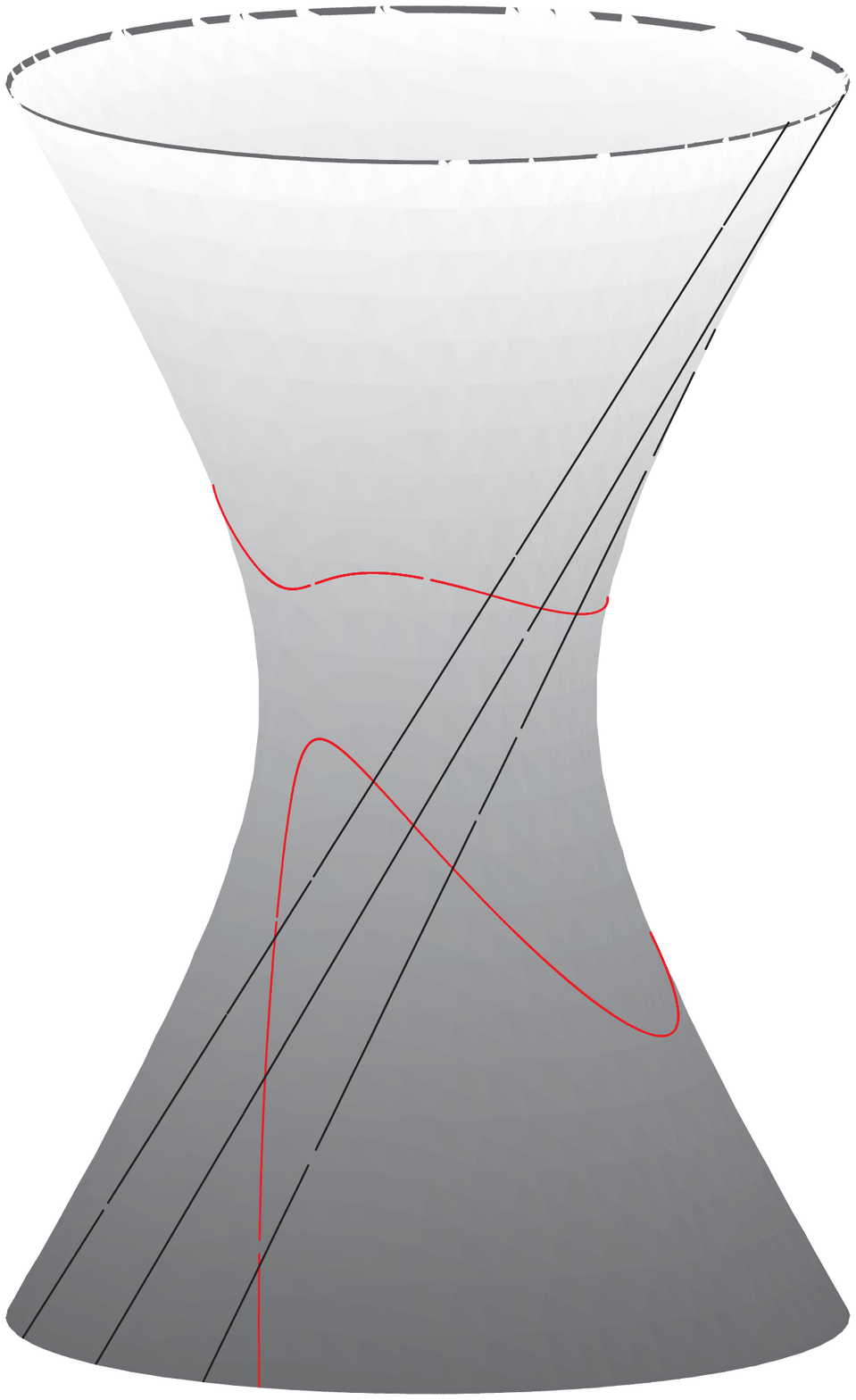}}
\]
}

\noindent\normalsize{The vanishing ideal $I_{X}$ of this scroll is given by
the $2\times2$ minors of a $2\times f$ matrix, $f=$ $h^{0}\mathcal{O}%
_{C}(K_{C}-D),$ with linear entries in $S$. The minimal free resolution of the
homogenous coordinate ring $S_{X}$ takes the following simple form:\newline%
\newline
\[
\underset{g-d+1}{\underbrace{%
\begin{tabular}
[c]{|p{0.7cm}|p{0.7cm}|p{0.7cm}|p{0.7cm}|p{0.7cm}|p{0.7cm}|}%
\hline
\multicolumn{1}{|>{\columncolor{gray}}c|}{1}  &  &  &  & \\\hline
&  \multicolumn{1}{|>{\columncolor{gray}}c|}{$\beta_{12}$}  &  \multicolumn
{1}{|>{\columncolor{gray}}c|}{$\cdots$}  &  \multicolumn{1}{|>{\columncolor
{gray}}c|}{ }  &  \multicolumn{1}{|>{\columncolor{gray}}c|}{ }
\\\hline
\end{tabular}
}}%
\]
\newline
\normalsize
and as $C$ is contained in $X$ the syzygies of $S_{X}$ are also
syzygies of $S_{C}.$ Further there exists a corresponding $\mathbb{P}^{d-2}%
$-bundle $\mathbb{P(\mathcal{E})}$ \index{$\mathbb{P(\mathcal{E})}$ projective space bundle} of degree $d-1$ over $\mathbb{P}^{1}$,
which is a desingularization of $X.$ With $\mathcal{E\cong O}_{\mathbb{P}^{1}%
}(e_{1})\oplus...\oplus\mathcal{O}_{\mathbb{P}^{1}}(e_{d-1})$ we say that $X$
is of \index{$S(e_{1},...,e_{d-1})$ type of a scroll} type $S(e_{1},...,e_{d-1})$. The type of a scroll constructed in this
way can be determined by the values $h^{0}\mathcal{O}_{C}(iD)$ for
$i\in\mathbb{N}$. Schreyer gives in \cite{schreyer1} a resolution of
$\mathcal{O}_{C}$ in terms of $\mathcal{O}_{X}$-modules. The minimal free resolution of these modules in $\mathbb{P}^{g-1}$ are given by complexes $\mathcal{C}^{i}$. Then Schreyer shows in \cite{schreyer1} that a mapping cone construction leads to a free resolution of $C\subset\mathbb{P}^{g-1}.$ Unfortunately this resolution
can contain non minimal maps, so it remains to determine the ranks of
them. }
\\\\
{\normalsize In Chapter \ref{cliff<=2} we repeat the results from
\cite{schreyer1} for the case of trigonal and tetragonal curves. In the
main part of the thesis we focus on pentagonal curves $C\subset\mathbb{P}^{8}$
of genus $9.$ We distinguish two cases with $\operatorname*{Cliff}(C)=3$: }

\bigskip

{\normalsize \noindent\textbf{I.} (\textit{C admits a }$g_{7}^{2}$) If there
exists a $g_{7}^{2}$ on $C$ then we get a plane model $C^{\prime}%
\subset\mathbb{P}^{2}$ of $C$ having six double points as only singularities.
Projection from each of them leads to a $g_{5}^{1}.$ The image
$S\subset\mathbb{P}^{8}$ of the blowup of $\mathbb{P}^{2}$ in these $6$ points
under the adjoint series is a Bordiga surface $S^{\prime}\subset$
$\mathbb{P}^{8}$ which is smooth iff $C^{\prime}$ has no infinitely near
double points. It has only isolated rational singularities that come from
contraction of a strict transform }$E_{i}^{\prime}$ {\normalsize of a point
}$p_{i}$ that admits one infinitely near double point $p_{j}${\normalsize .
The minimal free resolution for $S^{\prime}$ determines the Betti table
for $C$. }

\bigskip

{\normalsize \noindent\textbf{II.} (\textit{C admits no }$g_{7}^{2}$) Now we
assume that $C$ has no $g_{7}^{2}.$ Then the scroll $X$ constructed from one
of the existing $g_{5}^{1}=\left\vert D\right\vert $ turns out to be of type
$S(2,1,1,1),$ $S(2,2,1,0)$ or $S(3,1,1,0)$. We define the multiplicity
$m_{\left\vert D\right\vert }$ of the linear system $\left\vert D\right\vert $
to be equal to one, two or three depending on that type. Later on we will show
that there also exists a geometric interpretation that justifies this
definition: Let $k=%
{\textstyle\sum\nolimits_{\left\vert D\right\vert \sim g_{5}^{1}}}
m_{\left\vert D\right\vert }$ be the total number of all $g_{5}^{1}$ (counted
with multilplicities), then there exists a local one parameter family $(C_{\lambda
})_{\lambda}$ with $C_{0}=C$ and $C_{\lambda}$ a curve with Clifford index 3
that has exactly $k$ ordinary $g_{5}^{1}.$ Following the approach of Schreyer in \cite{schreyer1} we consider a representation of $C\subset\mathbb{P(\mathcal{E})}$ as vanishing locus of the Pfaffians of a $5\times5$
skew symmetric matrix $\psi$ with entries in $\mathbb{P(\mathcal{E})}$ obtained from the structure theorem for Gorenstein ideals of codimension $3$. Let
$H$ denote the class of a hyperplane section and $R$ that of a ruling on $X$,
then a closer examination of all possible types for $\psi$ leads to the
following complete table: \vspace{3mm}
\begin{align*}
&  Table\ for\ curves\ with\mathit{\ }a\ g_{5}^{1}\mathit{\ }but\mathit{\ }%
no\mathit{\ }g_{7}^{2}\\
& \\
&  {\scriptsize
\begin{array}
[c]{|l@{\hspace{1mm}}|l@{\hspace{1mm}}|l|}\hline
&  & \\
\text{\textbf{Special Linear Series}} & \text{\textbf{Determinantal Surface}
}\mathbf{Y} & \text{\textbf{\hspace{2mm} Type of matrix} }\psi\\
& C\subset Y\subset X\subset\mathbb{P}^{8} & \\\hline
&  & \\
\exists!\text{ }g_{5}^{1}\text{ with }m_{g_{5}^{1}}=1 &
\begin{array}
[c]{c}%
\mathbb{P}^{2}\text{ blown-up in}\\
\text{ }9\text{ doublepoints and}\\
\text{ }1\text{ triple point of a }g_{8}^{2}\\
\end{array}
&
\begin{array}
[c]{c}%
\left(
\begin{array}
[c]%
{c@{\hspace{1.15mm}}c@{\hspace{1.15mm}}c@{\hspace{1.15mm}}c@{\hspace{1.15mm}}c}%
0 & H & H & H & H\\
& 0 & H-R & H-R & H-R\\
&  & 0 & H-R & H-R\\
&  &  & 0 & H-R\\
&  &  &  & 0
\end{array}
\right) \\
\text{{on }}S(2,1,1,1)
\end{array}
\\\hline
&  & \\%
\begin{array}
[c]{c}%
\exists\text{ exactly two }g_{5}^{1}\text{ }\\
\text{with }m_{g_{5}^{1}}=1
\end{array}
&
\begin{array}
[c]{c}%
\mathbb{P}^{1}\times\mathbb{P}^{1}\text{ blown-up in}\\
7\text{ double points of a}\\
g_{5}^{1}\times g_{5}^{1}%
\end{array}
&
\begin{array}
[c]{c}%
\left(
\begin{array}
[c]%
{c@{\hspace{1.15mm}}c@{\hspace{1.15mm}}c@{\hspace{1.15mm}}c@{\hspace{1.15mm}}c}%
0 & H & H & H & H\\
& 0 & H-R & H-R & H-R\\
&  & 0 & H-R & H-R\\
&  &  & 0 & 0\\
&  &  &  & 0
\end{array}
\right) \\
\text{{on }}S(2,1,1,1)
\end{array}
\\\hline
&  & \\%
\begin{array}
[c]{c}%
\exists\text{ exactly three }g_{5}^{1}\text{ }\\
\text{with }m_{g_{5}^{1}}=1
\end{array}
&
\begin{array}
[c]{c}%
\mathbb{P}^{1}\times\mathbb{P}^{1}\text{ blown-up in}\\
7\text{ double points }\\
p_{1},...,p_{7}\text{ of a }g_{5}^{1}\times g_{5}^{1}\text{ }\\
\text{but }p\in C\diagdown\{\text{ }p_{1},...,p_{7}\}\\
\text{is base point}\\
\text{of }\left\vert (2,2)-p_{1}-...-p_{7}\right\vert
\end{array}
&
\begin{array}
[c]{c}%
\left(
\begin{array}
[c]%
{c@{\hspace{1.15mm}}c@{\hspace{1.15mm}}c@{\hspace{1.15mm}}c@{\hspace{1.15mm}}c}%
0 & H & H & H & H\\
& 0 & 0 & H-R & H-R\\
&  & 0 & H-R & H-R\\
&  &  & 0 & 0\\
&  &  &  & 0
\end{array}
\right) \\
\text{{on }}S(2,1,1,1)
\end{array}
\\\hline
&  & \\%
\begin{array}
[c]{c}%
\exists!\text{ }g_{5}^{1}\text{ with }\\
m_{g_{5}^{1}}=2
\end{array}
&
\begin{array}
[c]{c}%
P_{2}:=\mathbb{P}(\mathcal{O}_{\mathbb{P}^{1}}(2)\oplus \mathcal{O}_{\mathbb{P}^{1}})\\
\text{ blown-up in}\\
7\text{ double points } p_{1},...,p_{7}\\
\end{array}
&
\begin{array}
[c]{c}%
\left(
\begin{array}
[c]%
{c@{\hspace{1.15mm}}c@{\hspace{1.15mm}}c@{\hspace{1.15mm}}c@{\hspace{1.15mm}}c}%
0 & H+R & H & H & H\\
& 0 & H & H & H-R\\
&  & 0 & H-R & H-2R\\
&  &  & 0 & H-2R\\
&  &  &  & 0
\end{array}
\right) \\
\text{{on }}S(2,2,1,0)
\end{array}
\\\hline
&  & \\%
\begin{array}
[c]{c}%
\exists!\text{ }g_{5}^{1}\text{ with }\\
m_{g_{5}^{1}}=3
\end{array}
&
\begin{array}
[c]{c}%
P_{2}:=\mathbb{P}(\mathcal{O}_{\mathbb{P}^{1}}(2)\oplus \mathcal{O}_{\mathbb{P}^{1}})\\
\text{ blown-up in}\\
7\text{ double points } p_{1},...,p_{7}\\
\text{ lying on a rational}\\
\text{curve of class }A+B
\end{array}
&
\begin{array}
[c]{c}%
\left(
\begin{array}
[c]%
{c@{\hspace{1.15mm}}c@{\hspace{1.15mm}}c@{\hspace{1.15mm}}c@{\hspace{1.15mm}}c}%
0 & H+R & H & H & H\\
& 0 & H & H & H-R\\
&  & 0 & H-R & H-2R\\
&  &  & 0 & H-2R\\
&  &  &  & 0
\end{array}
\right) \\
\text{{on }}S(3,1,1,0)
\end{array}
\\\hline
&  & \\%
\begin{array}
[c]{c}%
\exists!\text{ }g_{5}^{1}\text{ with }m_{g_{5}^{1}}=1\\
\text{ and}\\
\exists!\text{ }g_{5}^{1}\text{ with }m_{g_{5}^{1}}=2
\end{array}
&
\begin{array}
[c]{c}%
\mathbb{P}^{1}\times\mathbb{P}^{1}\text{ blown-up in }\\
7\text{ double points }\\
p_{1},...,p_{7}\text{ of a }g_{5}^{1}\times g_{5}^{1}\\
\text{ lying on a rational}\\
\text{curve of type }(2,1)
\end{array}
&
\begin{array}
[c]{c}%
\left(
\begin{array}
[c]%
{c@{\hspace{1.15mm}}c@{\hspace{1.15mm}}c@{\hspace{1.15mm}}c@{\hspace{1.15mm}}c}%
0 & H & H & H & H\\
& 0 & H-R & \text{\textbf{F}} & H-R\\
&  & 0 & H-R & \text{\textbf{F}}\\
&  &  & 0 & 0\\
&  &  &  & 0
\end{array}
\right)  \text{ }\\
\text{{on }}S(2,1,1,1)\\
\text{{with an entry }\textbf{F}}\sim H-R
\end{array}
\\\hline
\end{array}
}%
\end{align*}
\newline\noindent Here $Y\subset X$ is a surface on the scroll $X$ given by
the $2\times2$ minors of a matrix%
\[
\omega\sim%
\begin{pmatrix}
H-a_{1}R & H-a_{2}R & H-a_{3}R\\
H-(a_{1}+k)R & H-(a_{2}+k)R & H-(a_{3}+k)R
\end{pmatrix}
\]
with entries in $\mathbb{P(\mathcal{E})}.$ }

{\normalsize The non minimal maps in the mapping cone construction,
representing a free resolution of $C\subset\mathbb{P}^{8}$, are given by
certain submatrices of $\psi$. Calculating ranks in the individual cases, it
turns out that the entries $\beta_{35}=\beta_{45}$ in the Betti table for
$C$ are given by $4k$ with $k$ the number of all $g_{5}^{1}$ counted with
multiplicities. }

{\normalsize \bigskip}

\noindent \textit{Acknowledgement}. I'm grateful to the DFG for financial support during my studies. I want to thank Dr. Stavros Papadakis who spent a lot of his time giving answer to my questions. I benefitted a lot from a multitude of discussions with him. I'm also grateful to Janko Boehm and Michael Kunte who read this text and pointed out a lot of things that needed to be fixed. Dominik Dietrich has put a lot of effort into the layout. If it looks well then it is due to him. Thank you.
\\Special thanks go to Prof. Dr. Frank-Olaf Schreyer who suggested the problem to me. I want to thank him for his guidance, encouragement and patience. This thesis is inspired by his ideas.   
\\\\
\noindent \textit{Notation}
Throughout the \index{$\Bbbk$ algebraically closed field of characteristic $0$} text $\Bbbk$ denotes an algebraically closed field of characteristic $0$. $S=\Bbbk[x_0,...,x_n]$ denotes the homogeneous coordinate ring of $\mathbb{P}^{n}$ and \index{$\mathfrak{m}$ maximal ideal in $S$} $\mathfrak{m}=(x_0,...,x_n)\subset S$ its maximal ideal. For two divisors $D$ and $D^{\prime}$ on a variety $V$ we write $D\sim D^{\prime}$ iff they are linear equivalent. \index{$D\sim D^{\prime}$ linear equivalence of divisors $D$ and $D^{\prime}$}

{\normalsize \newpage}
°
\newpage
\chapter{{\protect\normalsize Background}}

\section{Syzygies}

{\normalsize A projective variety $X\subset\mathbb{P}^{n}$ can be described by
its vanishing ideal $I_{X}\subset S:=\Bbbk\lbrack x_{0},...,x_{n}]$. We denote the corresponding \index{$S_{X}$     homogeneous coordinate ring of a variety $X$} \index{$I_{X}$   vanishing ideal of a variety $X$} homogeneous coordinate ring by $S_{X}:=S/I_{X}.$
Because of the Noetherian property of $S$, there exists a finite number of
generators of $I_{X}.$ These generators can also have certain relations which
can be described as a finitely generated module over $S.$ Then we consider the
relations of these relations and so on. Hilbert famous Syzygy Theorem says
that this process stops after finitely many steps: }

\begin{theorem}
{\normalsize (Hilbert Syzygy Theorem) Any finitely generated graded S-module M
has a finite graded free resolution%
\[
0\leftarrow M\overset{\varphi_{0}}{\leftarrow}F_{0}\overset{\varphi_{1}%
}{\leftarrow}F_{1}\leftarrow...F_{m-1}\overset{\varphi_{m}}{\leftarrow}%
F_{m}\leftarrow0
\]
with free S-modules $F_{i},$ $i=0,...,m\leq n+1.$ }
\end{theorem}

\begin{proof}
{\normalsize \cite{eisenbud2} sect.2B. }
\end{proof}

{\normalsize \bigskip}

{\normalsize Unfortunately a free resolution of $S_{X}$ has not to be minimal. But if we choose a minimal set of generators in each step above then we obtain a minimal free resolution, i.e. every map $\varphi_{i},$ $i=0,...,m$, has no degree
}zero {\normalsize part or equivalently it is not possible to seperate a
trivial subcomplex%
\[
0\leftarrow S(-d)\overset{\cong}{\leftarrow}S(-d)\leftarrow0
\]
To give a formal definition for this property we use the standard notation
$\mathfrak{m}$ to denote the homogeneous maximal ideal $(x_{0},...,x_{n}%
)\subset S:$ }

\begin{definition}
{\normalsize A complex of graded S-modules%
\[
...\leftarrow F_{i-1}\overset{\varphi_{i}}{\leftarrow}F_{i}\leftarrow...
\]
is called minimal if for each i the image of $\varphi_{i}$ is contained in
$\mathfrak{m}F_{i-1}.$ }
\end{definition}

{\normalsize \noindent Given a resolution of $S_{X},$ we obtain a minimal free
resolution by canceling trivial subcomplexes and in consequence we have the
following important property }

\begin{theorem}
{\normalsize Every minimal free resolution of a graded S-module M is unique
up to an isomorphism of complexes inducing the identity map on M. }
\end{theorem}

\begin{proof}
{\normalsize \cite{eisenbud2} Section 1.}
\end{proof}

{\normalsize \bigskip}

{\normalsize \noindent According to the last theorem, it follows the important
fact that for each minimal free resolution of a finitely generated graded
$S-$module the number of generators of each degree $j\in\mathbb{Z}$ required
for the free modules $F_{i}$ is the same in every minimal free resolution. We call these numbers $\beta_{ij}$
the Betti numbers of $M$. \index{$\beta_{ij}$ Betti numbers}}

\begin{definition}
{\normalsize (Betti numbers) Let M be a finitely generated, graded S-module
and
\[
0\leftarrow M\overset{\varphi_{0}}{\leftarrow}F_{0}\overset{\varphi_{1}%
}{\leftarrow}F_{1}\leftarrow...F_{m-1}\overset{\varphi_{m}}{\leftarrow}%
F_{m}\leftarrow0
\]
a minimal free resolution of M with free modules $F_{i}=%
{\textstyle\bigoplus\nolimits_{j}}
S(-j)^{\beta_{ij}}$ , then we call the numbers $\beta_{ij}$ the syzygy numbers
or graded Betti numbers of M. }
\end{definition}

{\normalsize For $X\subset\mathbb{P}^{n}$ a projective variety we call the
Betti numbers of $X$ those of the homogenous coordinate ring $S_{X}.$ Given a
set of generators of the vanishing ideal of $X,$ it is possible to determine a
minimal free resolution of $S_{X}$ by Gr\"{o}bner Basis algorithms
(implemented in Macaulay2, Singular,...) after finitely many steps. We use the
Macaulay notation to write down the Betti table of such a resolution:

\[
\begin{tabular}{c|ccccc}
& 0 & 1 & ... & m-1 & m \cr \hline
0 & $\beta_{00}$ & $\beta_{11}$ & .... & $\beta_{m-1,m-1}$ & $\beta_{m,m}$\cr
1 & $\beta_{01}$ & $\beta_{12}$ &  & $\beta_{m-1,m}$ & $\beta_{m,m+1}$\cr
$\overset{.}{:}$ & $\overset{.}{:}$ & $\overset{.}{:}$ &  & $\overset{.}{:}$ & $\overset{.}{:}$\cr
r & $\beta_{0,r}$ & $\beta_{1,r+1}$ & .... & $\beta_{m-1,m+r-1}$ & $\beta_{m,m+r}$\cr
\end{tabular}
\]
}

{\normalsize \bigskip}

\begin{example}{\label{compinter}(complete intersection in $\mathbb{P}^{3}$)}
{\textnormal{A complete intersection of two hypersurfaces of degree
2 in $\mathbb{P}^{3}$ is an elliptic curve $C\subset\mathbb{P}^{3}.$ $I_{C}$ is generated by exactly two quadratic forms $q_{1},q_{2}\in
S_{2}$ and the only relation between them is the Koszul relation $-q_{2}\cdot
q_{1}+q_{1}\cdot q_{2}=0$, hence the minimal free resolution of $S_{C}$ is
given by%
\[
0\leftarrow S_{C}\leftarrow S\leftarrow S(-2)\oplus S(-2)\leftarrow
S(-4)\leftarrow0
\]
and has the following Betti table%
\[
\begin{tabular}{c|ccc}
& 0 & 1 & 2 \cr \hline
0 & 1 & - & - \cr
1 & - & 2 & - \cr
2 & - & - & 1 \cr
\end{tabular}
\]}
}
\end{example}

{\normalsize \bigskip}

\begin{example}{\normalsize \label{koszul}(The Koszul complex in general) }
{\textnormal{Let $S$ be a ring and $N\cong S^{n+1}$ a free $S-$module
of rank $n+1.$ For $x=(x_{0},...,x_{n})\in N$ we define the Koszul complex to be
the complex%
\[
K(x):0\rightarrow R\rightarrow N\rightarrow\wedge^{2}N\rightarrow
...\rightarrow\wedge^{i}N\overset{d_{x}}{\rightarrow}\wedge^{i+1}%
N\rightarrow...
\]
where $d_{x}$ sends an element $a$ to $x\wedge a.$ Notice that
$\wedge^{0}N\cong S.$ In particular for $S=\Bbbk\lbrack x_{0},...,x_{n}]$ the
graded homogeneous polynomial ring in $\mathbb{P}^{n}$ and $x=(x_{0},...,x_{n})$ we get a graded free resolution of $\Bbbk$ as $S-$module:%
\[
K(x):\text{ }0\rightarrow S(-n-1)\rightarrow S^{n+1}(-n)\rightarrow S^{\binom{n+1}{2}(-n+1)%
}\rightarrow...\rightarrow S^{\binom{n+1}{i}}(-n+i-1)\overset{d_{x}%
}{\rightarrow}%
\]%
\[
\overset{d_{x}}{\rightarrow}S^{\binom{n+1}{i+1}}(-n+i)\rightarrow...\rightarrow
\wedge^{n+1}S^{n+1}\cong S\rightarrow\Bbbk\rightarrow0
\]
as the cokernel of $\wedge^{n}N\rightarrow\wedge^{n+1}N$ is isomorphic to
$S/(x_{0},...,x_{n})\cong\Bbbk$ (cf. \cite{eisenbud1} page 428)$.$}}
\end{example}
{\normalsize \bigskip}

{\normalsize We derive further properties of the Betti numbers: }

{\normalsize \bigskip}

\begin{theorem}
{\normalsize Let M be a finitely generated, graded S-module and%
\[
0\leftarrow M\overset{\varphi_{0}}{\leftarrow}F_{0}\overset{\varphi_{1}%
}{\leftarrow}F_{1}\leftarrow...F_{m-1}\overset{\varphi_{m}}{\leftarrow}%
F_{m}\leftarrow0
\]
a free resolution of M with free modules $F_{i}=%
{\textstyle\bigoplus\nolimits_{j}}
S(-j)^{\beta_{ij}}$. Then the \index{$H_{M}$ Hilbert function} Hilbert function $H_{M}$ is%
\[
H_{M}(d)=%
{\textstyle\sum\nolimits_{i=0}^{m}}
(-1)^{i}%
{\textstyle\sum\nolimits_{j}}
\beta_{ij}\binom{n+d-j}{n}%
\]
}
\end{theorem}

\begin{proof}
{\normalsize \cite{eisenbud2} Corollary 1.2 }
\end{proof}

{\normalsize \bigskip}

{\normalsize \noindent It is a direct consequence of the theorem above that
for sufficently large $d,$ the Hilbert function becomes polynomial: }

\begin{theorem}
{\normalsize There exists a \index{$P_{M}$ Hilbert polynomial} polynomial $P_{M}$ (called the Hilbert polynomial)
such that, if M has a free resolution as above, $P_{M}(d)=H_{M}(d)$ for
$d\geq\max_{i,j}\{\beta_{ij}-n\}.$ }
\end{theorem}

\begin{proof}
{\normalsize For $d\geq\max_{i,j}\{j-n:\beta_{ij}\neq0\}$ we get $n+d-j%
\geq0$ in every binomial coefficient $\binom{n+d-j}{n}$, thus it becomes polynomial in $d.$ }
\end{proof}

{\normalsize \bigskip}

{\normalsize \noindent In consequence, we are able to calculate the
Hilbert function and the Hilbert polynomial from the Betti numbers. In our
definition of the graded Betti numbers $\beta_{ij}$ of a minimal free
resolution of an $S-$module $M$, we used the important fact, that these
numbers are uniquely determined by $M$ and do not depend on the minimal
resolution. }

{\normalsize \noindent There is a further way of introducing the Betti numbers
$\beta_{ij\text{ }}$ and \index{$F_{i,j}$ space of $(i-1)$-th syzygies of degree
$j$} the space $F_{i,j}$ of $(i-1)$-th syzygies of degree
$j$ of the module $S_{X}$ as%
\[
F_{i,j}=\operatorname{Tor}_{i}^{S}(S_{X},\Bbbk)_{j}%
\]
\index{$\operatorname{Tor}_{i}^{S}(S_{X},\Bbbk)_{j}$ derived functors of $\otimes$}
These two definitions coincide as proved in \cite{eisenbud1} Exercise A 3.18
page 639. The advantage of this second definition of the Betti numbers is, that
$\operatorname{Tor}_{i}^{S}(S_{X},\Bbbk)_{j}$ can be computed as well by a free resolution of
$S_{X}$ as by one of $\Bbbk$ in terms of free $S-$modules, which is given by the Koszul komplex (cf. Example \ref{koszul}). But it is not
possible in general to obtain the complete Betti table only from the values of
the Hilbert function. Therefore we need more information of an $S-$module $M$
such as the Castelnuovo-Mumford regularity \index{$\operatorname{reg}M$ Castelnuovo-Mumford regularity}: }

\begin{definition}
{\normalsize (Castelnuovo-Mumford regularity) The Castelnuovo-Mumford
regularity of $M$ is given by%
\[
\operatorname*{reg}M=\max(i-j,\beta_{ij}\neq0)
\]
For $X$ a projective variety we define $\operatorname*{reg}%
X:=\operatorname*{reg}I_{X}$ with $I_{X}$ the vanishing ideal of $X.$ }
\end{definition}

{\normalsize \noindent The first property of the regularity of a module $M$ is
that it gives us a rather good lower boundary for the natural numbers $d$
\ such that the Hilbertfunction $H_{M}(d)$ agrees with the Hilbert polynomial
$P_{M}(d):$ }

\begin{theorem}
{\normalsize \label{hilbpoly}Let $M$ be a finitely generated graded module
over the polynomial ring $S=\Bbbk\lbrack x_{0},...,x_{n}],$ }then

{\normalsize 1) $H_{M}(d)=P_{M}(d)$ for all $d\geq\operatorname*{reg}%
M+\operatorname*{pd}M-n$ }

{\normalsize 2) For $M$ Cohen-Macaulay the boundary in 1) is sharp. }
\end{theorem}

\begin{proof}
{\normalsize \cite{eisenbud2} Theorem 4.2. }
\end{proof}

{\normalsize \bigskip}

{\normalsize \noindent The following theorem helps us to calculate the
regularity in terms of the vanishing of cohomology groups. It says that for a
projective variety $X,$ $\operatorname*{reg}I_{X}$ can be obtained as the
minimal value $r_{0}$, such that $I_{X}$ is $r-$regular for all $r\geq r_{0}$:
}

\begin{theorem}
{\normalsize Let $X\subset\mathbb{P}^{n}$ be \ projective variety and $r_{0}$
the minimal number with%
\[
H^{i}(\mathbb{P}^{n},\mathcal{I}_{X}(r-i))=0\text{ for all }i>0\text{ and all
}r\geq r_{0}%
\]
Then $\operatorname*{reg}X=\operatorname*{reg}I_{X}=r_{0}.$ }
\end{theorem}

\begin{proof}
{\normalsize \cite{eisenbud1} Exercise 20.20. }
\end{proof}

{\normalsize \bigskip}

{\normalsize \noindent According to the definition of $\operatorname*{reg}M$
the regularity gives the number of rows in the Betti table of $M$ with nonzero
entries. The number of columns with nonzero entries is determined by the
projective dimension of $M:$ }

\begin{definition} \index{$\operatorname{pd}M$ projective dimension}
{\normalsize (projective dimension) For $M$ an $S-$module, the
projective dimension $\operatorname*{pd}M$ is the minimal length of a
projective resolution of $M.$ }
\end{definition}
\index{$\operatorname{depth}(\mathfrak{m},M)$ length of maximal $\mathfrak{m}$-sequence}
{\normalsize \noindent There is a direct correspondence of the projective
dimension of a module $M$ and the length $\operatorname*{depth}(\mathfrak{m}%
,M)$ of a maximal $M$-sequence in $\mathfrak{m}$, $\mathfrak{m}$ the maximal ideal of $S$,
given by the the Auslander-Buchsbaum formula: }

\begin{theorem}
{\normalsize (Auslander-Buchsbaum formula) Let $S$ be a graded ring with
maximal ideal }$\mathfrak{m}$ {\normalsize and $M$ a finitely generated
$S-$module with $\operatorname*{pd}M<\infty$. Then
\[
\operatorname*{pd}M=\operatorname*{depth}(\mathfrak{m}%
,S)-\operatorname*{depth}(\mathfrak{m},M)
\]
}
\end{theorem}

\begin{proof}
{\normalsize \cite{eisenbud1} Theorem 19.9 and Exercise 19.8, page 479 and 489.
}
\end{proof}

{\normalsize \bigskip}

{\normalsize \noindent For $S=\Bbbk\lbrack
x_{0},...,x_{n}]$ we get $\operatorname*{depth}(\mathfrak{m},M)=n+1,$ hence
$\operatorname*{depth}(\mathfrak{m},M)$ and $\operatorname*{pd}M$
can directly be obtained from the Betti table of $M.$ }

{\normalsize \noindent For $X\subset\mathbb{P}^{n}$ a projective variety, the
definition of being arithmetically Cohen-Macaulay says that $\operatorname*{depth}(\mathfrak{m},S_{X})=\dim R_{X}$ , hence $X$ is
arithmetically Cohen-Macaulay if and only if%
\[
\max\{i|\beta_{ij}\neq0\text{ for at least one }j\}=\operatorname*{codim}X
\]
It will turn out that a canonical curve $C$ is always arithmetically Cohen
Macaulay. Further we will show that $C$ is even arithmetically Gorenstein. }

\begin{definition}
{\normalsize Let $X\subset\mathbb{P}^{n}$ be a projective, arithmetically
Cohen-Macaulay variety and $\operatorname*{codim}X=c.$ Then we call $X$
arithmetically Gorenstein if and only if there exists an $n\in\mathbb{Z}$ with%
\[
\operatorname*{Ext}{}^{c}(S_{X},S)=S_{X}(n)
\]
} \index{$\operatorname{Ext}$, extension functor}
\end{definition}

{\normalsize \noindent If $F$ is a minimal free resolution of $S_{X}$ then the
minimal resolution of $\operatorname*{Ext}{}^{c}(S_{X},S)$ can be obtained as
the dual $F^{\ast}$ of $F.$ Now it is an easy consequence that for $X$
arithmetically Gorenstein, $F$ has to be selfdual, i.e. for the Betti numbers we get%
\[
\beta_{ij}=\beta_{c-i,n-j}%
\]
with $c,$ $n$ as in the definition above. }

{\normalsize \noindent In Example \ref{compinter} we have seen that for an
elliptic curve $C\subset\mathbb{P}^{3}$ given as complete intersection of two
quadrics, its homogenous coordinate ring $S_{C}$ has the following symmetric
Betti table%
\[
\begin{tabular}{c|ccc}
& 0 & 1 & 2 \cr \hline
0 & 1 & - & - \cr
1 & - & 2 & - \cr
2 & - & - & 1 \cr
\end{tabular}
\]
Hence from the results above it follows that $C$ is arithmetically Gorenstein. }

\section{Canonical Curves and Green's Conjecture}

{\normalsize As we have remarked above we want to show that a canonical curve
$C\subset\mathbb{P}^{g-1}$ is arithmetically Gorenstein. Further we will see that
$\operatorname*{reg}C=3$ and as the Hilbert function of $S_{C}$ can easily be
calculated from the values $h^{0}(C,mK_{C})$, $m\in\mathbb{N}$ and $K_{C}$ a
canonical divisor \index{$K_{X}$ canonical divisor on a projective variety $X$} on $C$, we get a first approximation for the Betti table of
$S_{C}$. }

{\normalsize \noindent Let $C$ be a smooth, non hyperelliptic curve and
$K_{C}$ its canonical divisor, then we can embed $C$ in $\mathbb{P}^{g-1}$
canonically: \index{$\varphi_{\left\vert K_{C}\right\vert }$ canonical embedding}%
\[
\varphi_{\left\vert K_{C}\right\vert }:C\overset{\left\vert K_{C}\right\vert
}{\rightarrow}\mathbb{P}^{g-1}=\mathbb{P}H^{0}(C,\mathcal{O}_{C}(K_{C}))
\]
We denote by $I_{C}$ the vanishing ideal of $C\subset\mathbb{P}^{g-1}$ and
$S_{C}$ its homogeneous coordinate ring. The following theorem due to Max
Noether says that $C\subset\mathbb{P}^{g-1}$ is embedded projectively
normal and thus $C$ is arithmetically Cohen-Macaulay: }

\begin{theorem}
{\normalsize (Max Noether) If $C$ is non hyperelliptic then%
\[
\Omega=%
{\textstyle\sum\nolimits_{m\geq0}}
H^{0}(C,\mathcal{O}_{C}(K_{C})^{\otimes m})
\]
is the homogenous coordinate ring of $C\subset\mathbb{P}^{g-1}.$ It follows
that $H^{1}(C,\mathcal{I}_{C}(m))=0$ for all $m\geq0$.}
\end{theorem}

\begin{proof}
{\normalsize \cite{arbarello} page 117.}
\end{proof}

{\normalsize \bigskip}

{\normalsize \noindent Now as an easy consequence of Max Noether's Theorem and
Theorem \ref{hilbpoly} we obtain the following corollary: }

\begin{corollary}
{\normalsize Let $C\subset\mathbb{P}^{g-1}$ be the canonical model of a non
hyperelliptic curve of genus $g\geq3,$ then the Hilbert function }$H_{S_{C}}%
${\normalsize takes the following values:%
\[
H_{S_{C}}(d)=\left\{
\begin{array}
[c]{cc}%
0 & \text{if }d<0\\
1 & \text{if }d=0\\
g & \text{if }d=1\\
(2d-1)(g-1) & \text{if }d>0
\end{array}
\right.
\]
In particular $\beta_{1,2}(S_{C}),$ the number of quadratic generators of
$I_{C},$ is $\binom{g-1}{2}$ and $\operatorname*{reg}S_{C}=3.$ }
\end{corollary}

\begin{proof}
{\normalsize (\cite{eisenbud2} Corollary 9.4.) From Max Noether's Theorem we
already know that $S_{C}$ is arithmetically Cohen-Macaulay and $(S_{C})_{d}\cong H^{0}(C,\mathcal{O}_{C}(K_{C})^{\otimes d})=H^{0}(C,\mathcal{O}_{C}%
(dK_{C})).$ Therefore the values for $H_{S_{C}}$ follow from the Riemann-Roch
Theorem. From the Cohen Macaulay property we obtain the existence of a regular
sequence on $C$ consisting of $2$ linear forms $l_{1},l_{2}.$ The regularity
of $S_{C}$ is the same as that of $S_{C}/(l_{1},l_{2}).$ The Hilbert function
of this last module has values $(1,g-2,g-2,1)$, and thus applying Theorem
\ref{hilbpoly} we have $\operatorname*{reg}S_{C}=\operatorname*{reg}%
S_{C}/(l_{1},l_{2})=3.$ }
\end{proof}

{\normalsize \bigskip}

{\normalsize The Gorenstein property of $S_{C}$ is a direct consequence of
$\omega_{C}=\mathcal{O}_{C}(1).$ From the collected results we obtain that the
Betti table of $C$ has exactly $\operatorname*{codim}C+1=g-1$ columns and
$\operatorname*{reg}C+1=4$ rows. Further it is symmetric and $\beta_{ii}=0$
for all $i>0$ as $C$ is not contained in any hyperplane. Therefore the Betti
table of $C$ looks like: }
\[
\overset{g-1}{\overbrace{%
\begin{tabular}
[c]{|p{0.7cm}|p{0.7cm}|p{0.7cm}|p{0.7cm}|p{0.7cm}|p{0.7cm}|p{0.7cm}|p{0.7cm}|}%
\hline
\multicolumn{1}{|>{\columncolor{gray}}c|}{1}  &  &  &  &  &  &  & \\\hline
&  \multicolumn{1}{|>{\columncolor{gray}}c|}{$\beta_{12}$}  &  \multicolumn
{1}{|>{\columncolor{gray}}c|}{$\cdots$}  &  \multicolumn{1}{|>{\columncolor
{gray}}c|}{ }  &  \multicolumn{1}{|>{\columncolor{gray}}c|}{ }  &  &  &
\\\hline
&  &  &  \multicolumn{1}{|>{\columncolor{gray}}c|}{ }  &  \multicolumn
{1}{|>{\columncolor{gray}}c|}{ }  &  \multicolumn{1}{|>{\columncolor{gray}}%
c|}{$\cdots$}  &  \multicolumn{1}{|>{\columncolor{gray}}c|}{$\beta_{g-3,g-1}$}
& \\\hline
&  &  &  &  &  &  &  \multicolumn{1}{|>{\columncolor{gray}}c|}{1} \\\hline
\end{tabular}
}}%
\]
\\
\noindent with $\beta_{i,i+1}-\beta_{i-1,i+1}=i\cdot\binom{g-2}%
{i+1}-(g-1-i)\cdot\binom{g-2}{i-2}$, which we know from the values of the
Hilbert function $H_{S_{C}}.$ Due to recent results of Voisin in \cite{voisin2}
and \cite{voisin3}, for the generic curve over a field $\Bbbk$%
{\normalsize \ of characteristic }$0$, {\normalsize the Betti numbers
$\beta_{i-1,i+1}$ become zero for $i=0,...,\left\lfloor \frac{g-3}%
{2}\right\rfloor $. Especially in the case of our interest the Betti table of
a generic curve of genus }$9$ looks as follows:

\begin{example}
{\normalsize For $C\subset\mathbb{P}^{g-1}$ a generic curve of genus $g=9$ the
Betti table of $S_{C}$ has the following form:%
\textnormal{
\[
\begin{tabular}{c|cccccccccc}
& 0 & 1 & 2 & 3 & 4 & 5 & 6 & 7 \cr \hline
0 & 1 & - & - & - & - & - & - & - \cr
1 & - & 21 & 64 & 70 & - & - & - & - \cr
2 & - & - & -  & -  & 70  & 64  & 21 & -\cr
3 & - & - & - & - & - & - & - & 1 \cr
\end{tabular} 
\]}
}
\end{example}

{\normalsize \bigskip}

{\normalsize The question that arises is which Betti numbers $\beta_{ij}$ can be nonzero? A
classical result of Petri is the following theorem: }

\begin{theorem}
{\normalsize (Petri) Let $I_{C}$ be the homogenous ideal of a non
hyperelliptic canonical curve $C\subset\mathbb{P}^{g-1}$, then $I_{C}$ is
generated by $\binom{g-1}{2}$ quadrics except the two cases where }

{\normalsize 1) $C$ is trigonal, i.e. there exists a $g_{3}^{1}$ or }

{\normalsize 2) $C$ is isomorphic to a plane quintic $(g=6$ and $C$ has a
$g_{5}^{2})$ }

{\normalsize \noindent In these two exceptional cases the quadrics contained
in $I_{C}$ generate a rational normal scroll in 1) and the Veronese surface
$\mathbb{P}^{2}\hookrightarrow\mathbb{P}^{5}$ in 2). }
\end{theorem}

\begin{proof}
{\normalsize \cite{saint-donat} }
\end{proof}

{\normalsize \bigskip}

{\normalsize \noindent Thus we obtain $\beta_{1j}=0$ for $j>2$ except in the
two exceptional cases 1) and 2). This result suggests that the values
$\beta_{ij}$ that are nonzero correspond to special linear systems on the
curve $C.$ Green's conjecture, if shown to be true, would give an exact answer to this question$:$ }

\begin{conjecture}
{\normalsize (Green, 1984) Let $C\subset\mathbb{P}^{g-1}$ be a smooth
non hyperelliptic curve over a field of characteristic 0 in its canonical
embedding. Then%
\[
\beta_{p,p+2}\neq0\Leftrightarrow C\text{ has Clifford index }%
\operatorname*{Cliff}(C)\leq p
\]
}
\end{conjecture}

\noindent The Clifford index of an effective divisor $D$ on $C$
is defined as%
\[
\operatorname*{Cliff}(D)=\deg D-2r(D)
\]
and the Clifford index of $C$ is%
\[
\operatorname*{Cliff}(C)=\min{\{\operatorname*{Cliff}(D):D\in \operatorname{Div}(C) \text{ effective with }%
h^{i}\mathcal{O}_{C}(D)\geq2\text{ for }i=1,2\}}
\]
{\normalsize As we have remarked in the introduction the conjecture is already
shown to be true for the $"\Leftarrow"$ direction and also in several cases,
especially for genus $g=9,$ for the other direction. We will give a complete
list of all possible Betti tables for irreducible, nonsingular,
canonical curves $C$ of genus $g=9$. For $\operatorname*{Cliff}(C)\leq2$ the
results in \cite{schreyer1} can be applied to obtain these tables. This is
done in Chapter \ref{cliff<=2}. It remains to examine the case, where
$\operatorname*{Cliff}(C)=3$. Then there exists a }$g_{5}^{1},$ $g_{7}^{2},$
$g_{9}^{3}$ or $g_{11}^{4}$. The Brill Noether duals to $g_{9}^{3}$ and
$g_{11}^{4}$ are of type $g_{7}^{2}$ or $g_{5}^{1}$ correspondingly. From a
$g_{7}^{2}$ we get a plane model of $C$ of degree 7 that has exactly 6 double
points. Then $g_{5}^{1}$ on $C$ can be obtained from projection from one of the double points (cf. Theorem \ref{g27a}){\normalsize . Thus for a canonical curve }%
$C$ {\normalsize  of genus 9 to be pentagonal is equivalent to }%
$\operatorname*{Cliff}(C)=3${\normalsize .}

{\normalsize For }$\left(  D_{\lambda}\right)  _{\lambda}$ a base point
free pencil of divisors of degree, {\normalsize we consider the variety
swept out by the linear spans of these divisors:%
\[
X=%
{\textstyle\bigcup\nolimits_{D\in g_{5}^{1}}}
\bar{D}\subset\mathbb{P}^{8}%
\]
This is a rational normal scroll of dimension $4$. We apply the results of
\cite{schreyer1} to obtain a minimal free resolution of $\mathcal{O}_{C}$ as
$\mathcal{O}_{\mathbb{P}^{8}}$-module$:$ }

{\normalsize \bigskip}

{\normalsize \noindent1) Resolve $\mathcal{O}_{C}$ as an $\mathcal{O}_{X}%
$-module by direct sums of line bundles on $X:$%
\[
F_{\ast}:\ \ \ \ \ \ 0\rightarrow\mathcal{O}_{X}(-5H+3R)\rightarrow%
{\displaystyle\sum\limits_{i=1}^{5}}
\mathcal{O}_{X}(-3H+b_{i}R)\overset{\psi}{\rightarrow}%
\]
}
{\normalsize
\[
\overset{\psi}{\rightarrow}%
{\displaystyle\sum\limits_{i=1}^{5}}
\mathcal{O}_{X}(-2H+a_{i}R)\rightarrow\mathcal{O}_{X}\rightarrow
\mathcal{O}_{C}\rightarrow0
\]
2) Take the resolution $\mathcal{C}^{b}(a),$ $a,b\in\mathbb{Z},$ of each of
these line bundle as $\mathcal{O}_{\mathbb{P}^{8}}$-module. Then a mapping
cone construction leads to a (not necessarily minimal) resolution of
$\mathcal{O}_{C}$ as $\mathcal{O}_{\mathbb{P}^{8}}$-module:%
\[
\left[  \left[  \mathcal{C}^{3}(-5)\rightarrow%
{\displaystyle\sum\limits_{i}}
\mathcal{C}^{b_{i}}(-3)\right]  \rightarrow%
{\displaystyle\sum\limits_{i}}
\mathcal{C}^{a_{i}}(-2)\right]  \rightarrow\mathcal{C}^{0}%
\]
3) The non minimal parts of the mapping cone are related to submatrices of $\psi$. Hence to determine the ranks of these maps we have to study $\psi$. Then we obtain a minimal free resolution. }

{\normalsize \bigskip}

{\normalsize \noindent We will see that the matrix $\psi$ and the ranks of the
non minimal maps in the mapping cone are related to a certain number of
special linear systems of type $g_{5}^{1}$ or $g_{7}^{2}$. }
\\\\
{\normalsize We will provide some results on scrolls and give a description
how to manage the steps 1) and 2). }

\section{Scrolls in general\label{Scrollgeneral}}

{\normalsize Let $\mathcal{E=O(}e_{1})\oplus...\oplus\mathcal{O(}e_{d}),$
$e_{1}\geq...\geq e_{d}\geq0,$ be a globally generated, locally free sheaf of rank $d$ on
$\mathbb{P}^{1}$ and let
\[
\pi:\mathbb{P(\mathcal{E})\rightarrow P}^{1}%
\]
the corresponding $\mathbb{P}^{d-1}$-bundle. For $f=%
{\textstyle\sum\nolimits_{i=1}^{d}}
e_{i}\geq2$ consider the image of $\mathbb{P(\mathcal{E})}$ in $\mathbb{P}%
^{r}=\mathbb{P}H^{0}(\mathbb{P(\mathcal{E})},\mathcal{O}%
_{\mathbb{P(\mathcal{E})}}(1))$:%
\[
j:\mathbb{P(\mathcal{E})\rightarrow}X\subset\mathbb{P}^{r}\text{ \ , }r=f+d-1
\]
Then we call $X$ a rational normal \index{$R^{i}j_{\ast}$ higher direct image functor} scroll of type $S(e_{1},...,e_{d}).$ $X$ is
a non-degenerate, irreducible variety of minimal degree%
\[
\deg X=f=r-d+1=\operatorname*{codim}X+1
\]
in $\mathbb{P}^{r}.$ }

{\normalsize \noindent If all $e_{i}>0$ then $X$ is smooth and
$j:\mathbb{P(\mathcal{E})\rightarrow}X$ an isomorphism. Otherwise $X$ is
singular and $j:\mathbb{P(\mathcal{E})\rightarrow}X$ a resolution of
singularities. The singularities of $X$ are rational, i.e.
\[
j_{\ast}\mathcal{O}_{\mathbb{P(\mathcal{E})}}=\mathcal{O}_{X},\text{
\ \ }R^{i}j_{\ast}\mathcal{O}_{\mathbb{P(\mathcal{E})}}=0\text{ for }i>0
\]
Therefore there is no problem to replace $X$ by $\mathbb{P(\mathcal{E})}$ for
most cohomological considerations, even if $X$ is singular. }

\begin{remark}{\normalsize (Geometric description of scrolls) }
{\textnormal{A rational normal scroll $X$ of type $S(e_{1},...,e_{d})$ admits the following geometric description: Consider the
$e_{i}-th$ Veronese embedding of $\mathbb{P}^{1}:$%
\[
\gamma_{i}:\mathbb{P}^{1}\hookrightarrow\mathbb{P}^{e_{i}}%
\]
then the image $\Gamma_{e_{i}}$ is a rational normal curve of degree $e_{i}.$
Now we can embed the projective spaces $\mathbb{P}^{e_{i}}$ in $\mathbb{P}%
^{r},$ $r=d+f-1=%
{\textstyle\sum\nolimits_{i=1}^{d}}
(e_{i}+1)-1,$ as linearly independent subspaces of $\mathbb{P}^{r}.$
Identifying the curves $\Gamma_{e_{i}}$ with a common $\mathbb{P}^{1}$ we can
take the linear span of corresponding points which leads us to the scroll $X.$
If some of the $e_{i}=0$ the image $\Gamma_{e_{i}}$ under the mapping
$\gamma_{i}$ is a single point, thus $X$ becomes a cone with center spanned by all
points $\Gamma_{e_{i}}$ where $e_{i}=0.$ }}
\end{remark}

{\normalsize \bigskip}

{\normalsize Now let us examine the Picard group of $\mathbb{P(\mathcal{E})}$.
First we denote the \index{$H$ hyperplane class} hyperplane class $H=[j^{\ast}\mathcal{O}_{\mathbb{P}^{r}%
}(1)]$ and the \index{$R$ class of a ruling} 
ruling $R=[\pi_{\ast}\mathcal{O}_{\mathbb{P}^{1}}(1)].$ Then the
following theorem contains the needed results: }

\begin{theorem}
{\normalsize (Picard Group of $\mathbb{P(\mathcal{E})}$) }

{\normalsize \noindent1) $\operatorname*{Pic}\mathbb{P(\mathcal{E})=Z}%
H\oplus\mathbb{Z}R$ }

{\normalsize \noindent2) $H^{d}=f,$ $H^{d-1}.R=1$ and $R^{2}=0$ }

{\normalsize \noindent3) $K_{X}=-dH+(f-2)R$ }

{\normalsize \noindent4) There exist basic sections $\varphi_{i}\in
H^{0}(\mathbb{P(\mathcal{E})},\mathcal{O}_{\mathbb{P(\mathcal{E})}}%
(H-e_{i}R))$ and $s,t\in H^{0}(\mathbb{P(\mathcal{E})},\mathcal{O}%
_{\mathbb{P(\mathcal{E})}}(R))$, such that every section $\psi\in
H^{0}(\mathbb{P(\mathcal{E})},\mathcal{O}_{\mathbb{P(\mathcal{E})}}(aH+bR))$
can be identified with a homogeneous polynomial%
\[
\psi=%
{\textstyle\sum\nolimits_{\alpha}}
p_{\alpha}(s,t)\varphi_{1}^{\alpha_{1}}...\varphi_{d}^{\alpha_{d}}%
\]
of degree $a=\alpha_{1}+...+\alpha_{d}$ in the $\varphi_{i}^{\prime}s$ and
coefficients homogeneous polynomials $p_{\alpha}$ of degree%
\[
\deg p_{\alpha}=\alpha_{1}e_{1}+...+\alpha_{d}e_{d}+b
\]
}

{\normalsize \noindent5) For $b\geq-1$ the dimension $h^{0}%
(\mathbb{P(\mathcal{E})},\mathcal{O}_{\mathbb{P(\mathcal{E})}}(aH+bR))$ does
not depend on the type $S(e_{1},...,e_{d})$ of the scroll but only on its
degree $f:$%
\[
h^{0}(\mathbb{P(\mathcal{E})},\mathcal{O}_{\mathbb{P(\mathcal{E})}%
}(aH+bR))=f\binom{a+d-1}{d}+(b+1)\binom{a+d-1}{d-1}%
\]
Especially $h^{0}(\mathbb{P(\mathcal{E})},\mathcal{O}_{\mathbb{P(\mathcal{E}%
)}}(H-R))=f.$ }
\end{theorem}

\begin{proof}
{\normalsize For 1), 2) and 3) see \cite{schreyer1} 1.2.-1.7.. Applying the
Leray spectral sequence we can calculate the cohomology of a line bundle
$\mathcal{O}_{\mathbb{P(\mathcal{E})}}(aH+bR):$%
\[
H^{i}(\mathbb{P}^{1},R^{j}\pi_{\ast}\mathcal{O}_{\mathbb{P(\mathcal{E})}%
}(aH+bR))\Rightarrow H^{i+j}(\mathbb{P}(\mathcal{E}),\mathcal{O}%
_{\mathbb{P(\mathcal{E})}}(aH+bR))
\]
especially for $a\geq0$ and $i=j=0$ we obtain for the global section of
$\mathcal{O}_{\mathbb{P(\mathcal{E})}}(aH+bR)$:%
\[
H^{0}(\mathbb{P}(\mathcal{E}),\mathcal{O}_{\mathbb{P(\mathcal{E})}%
}(aH+bR))\cong H^{0}(\mathbb{P}^{1},(S_{a}\mathcal{E)(}b))
\]
Then there exist basic sections $\varphi_{i}\in H^{0}(\mathbb{P(\mathcal{E}%
)},\mathcal{O}_{\mathbb{P(\mathcal{E})}}(H-e_{i}R))$ obtained from the
inclusion of the $i^{th}-$summand%
\[
\mathcal{O}_{\mathbb{P}^{1}}\rightarrow\mathcal{E(-}e_{i})\cong\pi_{\ast
}\mathcal{O}_{\mathbb{P(\mathcal{E})}}(H-e_{i}R)
\]
Denote the global generators of $H^{0}(\mathbb{P(\mathcal{E})},\mathcal{O}%
_{\mathbb{P(\mathcal{E})}}(R))$ by $s$ and $t$, then there is a natural way to
identify a section $\psi\in H^{0}(\mathbb{P(\mathcal{E})},\mathcal{O}%
_{\mathbb{P(\mathcal{E})}}(aH+bR))$ with a polynomial in the way as claimed. The
formula for the dimension $h^{0}(\mathbb{P(\mathcal{E})},\mathcal{O}%
_{\mathbb{P(\mathcal{E})}}(aH+bR))$ is a direct consequence of this
representation. }
\end{proof}

{\normalsize \bigskip}

\begin{example}{\normalsize ($X$ a scroll of type $S(2,1,1,1)$ in $\mathbb{P}^{8}$) }
{\textnormal{Let $X$ be a scroll of type $S(2,1,1,1)$ and
$\mathbb{P}(\mathcal{E})$ the corresponding $\mathbb{P}^{3}-$bundle over
$\mathbb{P}^{1},$ then there exist basic sections $\varphi_{0}\in
H^{0}(\mathbb{P(\mathcal{E})},\mathcal{O}_{\mathbb{P(\mathcal{E})}}(H-2R))$
and $\varphi_{1},\varphi_{2},\varphi_{3}\in H^{0}(\mathbb{P(\mathcal{E}%
)},\mathcal{O}_{\mathbb{P(\mathcal{E})}}(H-R)).$ Especially the hyperplane sections $H^{0}%
(\mathbb{P(\mathcal{E})},\mathcal{O}_{\mathbb{P(\mathcal{E})}}(H))$ of the
scroll are generated by $x_{0}=s^{2}\varphi_{0},$ $x_{1}=st\varphi_{0},$
$x_{2}=t^{2}\varphi_{0},$ $x_{3}=s\varphi_{1},$ $x_{4}=t\varphi_{1}%
,...,x_{8}=t\varphi_{3}.$ It follows that the $2\times2$ minors of the
$2\times5$ matrix%
\[
\Phi=\left(
\begin{array}
[c]{ccccc}%
x_{0} & x_{1} & x_{3} & x_{5} & x_{7}\\
x_{1} & x_{2} & x_{4} & x_{6} & x_{8}%
\end{array}
\right)
\]
vanish on the scroll $X.$ Moreover the following theorem says that they even
generate its vanishing ideal. }}
\end{example}

\begin{theorem}
{\normalsize Let $X$ be a scroll of type $S(e_{1},...,e_{d})$ and
\[
\Phi=\left(
\begin{array}
[c]{ccc}%
x_{10}\text{ }x_{11}\text{ }..\text{ }x_{1e_{1}-1} & ... & x_{d0}\text{ }..\text{ }%
x_{de_{d}-1}\\
x_{11}\text{ }x_{12}\text{ }..\text{ }x_{1e_{1\text{ \ \ \ }}} & ... & x_{d1}\text{
}..\text{ }x_{de_{d\text{ \ \ \ }}}%
\end{array}
\right)
\]
be the $2\times f$ matrix given by the multiplication map%
\[
H^{0}(\mathbb{P(\mathcal{E})},\mathcal{O}_{\mathbb{P(\mathcal{E})}}(R))\otimes
H^{0}(\mathbb{P(\mathcal{E})},\mathcal{O}_{\mathbb{P(\mathcal{E})}%
}(H-R))\rightarrow H^{0}(\mathbb{P(\mathcal{E})},\mathcal{O}%
_{\mathbb{P(\mathcal{E})}}(H))
\]
i.e. $x_{ij}=t^{j}s^{e_{i}-j}\varphi_{i}$, then the vanishing ideal $I_{X}$ of
$X$ is generated by the $2\times2$ minors of $\Phi.$ }
\end{theorem}

\begin{proof}
{\normalsize \cite{schreyer1} 1.6. }
\end{proof}

{\normalsize \bigskip}
\label{ident}
{\normalsize As we have already remarked at the end of the last section, we
want to resolve the line bundles $\mathcal{O}_{\mathbb{P(\mathcal{E})}%
}(aH+bR)$ in terms of $\mathcal{O}_{\mathbb{P}^{r}}-$modules. From our
definition of a scroll $X$ we have
\[
H^{0}(\mathbb{P(\mathcal{E})},\mathcal{O}_{\mathbb{P(\mathcal{E})}}(H))\cong
H^{0}(X,\mathcal{O}_{X}(H))\cong H^{0}(\mathbb{P}^{r},\mathcal{O}%
_{\mathbb{P}^{r}}(1))
\]
If we further denote \index{$F=H^{0}\mathcal{O}(H-R)\otimes\mathcal{O}_{\mathbb{P}^{r}}\cong\mathcal{O}_{\mathbb{P}^{r}}^{f}$}
\index{$G=H^{0}\mathcal{O}(R)\otimes\mathcal{O}_{\mathbb{P}^{r}}\cong\mathcal{O}_{\mathbb{P}^{r}}^{2}$}
\[
F=H^{0}\mathcal{O}(H-R)\otimes\mathcal{O}_{\mathbb{P}^{r}}\cong\mathcal{O}_{\mathbb{P}^{r}%
}^{f}\text{ \ \ and \ \ }G=H^{0}\mathcal{O}(R)\otimes\mathcal{O}_{\mathbb{P}^{r}}\cong\mathcal{O}_{\mathbb{P}^{r}}^{2}%
\]
then the multiplication map%
\[
G\otimes F\rightarrow\mathcal{O}_{\mathbb{P}^{r}}(1)
\]
in the theorem above induces a map%
\[
\Phi:F\otimes\mathcal{O}_{\mathbb{P}^{r}}(-1)\rightarrow G^{\ast}%
\otimes\mathcal{O}_{\mathbb{P}^{r}}\cong G
\]
}

{\normalsize \noindent Now define the complexes \index{$\mathcal{C}^{b}$ complex} $\mathcal{C}^{b}$ by%
\[
\mathcal{C}_{j}^{b}=\left\{
\begin{array}
[c]{cccc}%
{\textstyle\bigwedge\nolimits^{j}}
F\otimes S_{b-j}G\otimes\mathcal{O}_{\mathbb{P}^{r}}(-j)\text{
\ \ \ \ \ \ \ \ \ \ \ \ \ \ } &  & \text{for} & 0\leq j\leq b\\%
{\textstyle\bigwedge\nolimits^{j+1}}
F\otimes D_{j-b-1}G^{\ast}\otimes\mathcal{O}_{\mathbb{P}^{r}}(-j-1) &  &
\text{for} & j\geq b+1
\end{array}
\right.
\]
}

{\normalsize \noindent and the differential map
\[
\mathcal{C}_{j}^{b}\rightarrow\mathcal{C}_{j-1}^{b}%
\]
by the multiplication with $\Phi\in H^{0}(F^{\ast
}\otimes G\otimes\mathcal{O}_{\mathbb{P}^{r}}(1))$ for $j\neq b+1$ and
$\wedge^{2}\Phi\in H^{0}(\wedge^{2}F^{\ast}%
\otimes\wedge^{2}G\otimes\mathcal{O}_{\mathbb{P}^{r}}(2))$ for $j=b+1.$ \\
E.g. For $j\leq b$ the differential of a term
\[
f_{1}\wedge...\wedge f_{j} \otimes g \in H^{0}(\mathcal{C}_{j}^{b}),
\]
with $f_{i}\in H^{0}(F\otimes \mathcal{O}_{\mathbb{P}^{r}}(-1))$, $g\in H^{0}(S_{b-j}G)$ is given by%
\[
f_{1}\wedge...\wedge f_{j}\otimes g\rightarrow%
{\textstyle\sum\nolimits_{i=1}^{j}}
(-1)^{i}f_{1}\wedge...\wedge\hat{f}_{i}\wedge...\wedge f_{j}\otimes\Phi
(f_{i})\cdot g
\]
}

{\normalsize \bigskip}

\begin{theorem}
{\normalsize \label{Cba}$\mathcal{C}^{b}(a)$ for $b\geq-1$ is the minimal free
resolution of $\mathcal{O}_{X}(aH+bR)$ as an $\mathcal{O}_{\mathbb{P}^{r}}%
-$module. }
\end{theorem}

\begin{proof}
{\normalsize \cite{schreyer1} 2.2. }
\end{proof}

{\normalsize \bigskip}

{\normalsize \noindent It follows that the resolution of $\mathcal{O}_{X}$ is
given by $\mathcal{C}^{0}.$ This is the well known Eagon-Northcott complex
with Betti table as follows:}

\[
\underset{f}{\underbrace{%
\begin{tabular}
[c]{|p{0.7cm}|p{0.7cm}|p{0.7cm}|p{0.7cm}|p{0.7cm}|p{0.7cm}|}%
\hline
\multicolumn{1}{|>{\columncolor{gray}}c|}{1}  &  &  &  & \\\hline
&  \multicolumn{1}{|>{\columncolor{gray}}c|}{$\binom{f}{2}$}  &  \multicolumn
{1}{|>{\columncolor{gray}}c|}{$2\binom{f}{3}$}  &  \multicolumn{1}{|>{\columncolor
{gray}}c|}{$\cdots$ }  &  \multicolumn{1}{|>{\columncolor{gray}}c|}{$(f-1)\binom{f}{f}$ }
\\\hline
\end{tabular}
}}%
\]

{\normalsize \bigskip}

\section{\label{ScrollVariety}Scrolls constructed from
varieties}

{\normalsize In the next step, for a linearly normal embedded smooth variety
$V\subset\mathbb{P}^{r}$ and a pencil of divisors $(D_{\lambda})_{\lambda}$ on
$V$, we construct a scroll $X\subset\mathbb{P}^{r}$ such that $V\subset X$ and
further the pencil $(D_{\lambda})_{\lambda}$ is cut out on $V$ by the class of
a ruling $R$ on $X.$ Assume that $D$ is a divisor on $V$ with $h^{0}\left(
V,\mathcal{O}_{V}(D)\right)  \geq2$ and $h^{0}\left(  V,\mathcal{O}%
_{V}(H-D)\right)  =f\geq2$. Furthermore let $G\subset H^{0}\left(  V,\mathcal{O}%
_{V}(D)\right)  $ be the $2$ dimensional subspace that defines the pencil of
divisors $(D_{\lambda})_{\lambda}$, then from the multiplication map%
\[
G\otimes H^{0}\left(  V,\mathcal{O}_{V}(H-D)\right)  \rightarrow H^{0}\left(
V,\mathcal{O}_{V}(H)\right)
\]
we obtain a $2\times f$ matrix $\Phi$ with linear entries whose $2\times2$
minors vanish on $V.$ It turns out that the variety $X\subset\mathbb{P}^{r}$
defined by these minors is a scroll of degree $f,$ such that $(D_{\lambda
})_{\lambda}$ is cut out by the class of a ruling $R$ on $X.$ Geometrically
$X$ can be obtained as the union of the linear spans of $D_{\lambda}$:%
\[
X=%
{\textstyle\bigcup\nolimits_{\lambda\in\mathbb{P}^{1}}}
\bar{D}_{\lambda}\subset\mathbb{P}^{r}%
\]
}

\begin{remark}{\normalsize (Scroll $X$ constructed from a $g_{d}^{1}$ on a canonical curve
$C$) }{\textnormal{If $\left\vert D\right\vert $ is a base point free complete
linear system of type $g_{d}^{1}$ on $C$, then from the geometric version of
Riemann-Roch we get%
\[
\dim\bar{D}=\deg D-\dim\left\vert D\right\vert -1=d-2
\]
and therefore the scroll $X$ constructed from $\left\vert D\right\vert $ is
$(d-1)$-dimensional$.$ Let $\varphi_{K_{C}}$ denote the canonical map%
\[
\varphi_{K_{C}}:C\rightarrow\mathbb{P}^{g-1}%
\]
and%
\[
\varphi_{\left\vert D\right\vert }:C\rightarrow\mathbb{P}^{1}%
\]
the map which corresponds to $\left\vert D\right\vert .$ Then
\[
\mathcal{E}:=(\varphi_{\left\vert D\right\vert })_{\ast}\mathcal{O}_{C}(1)
\]
is a locally free sheaf $\mathcal{O}_{\mathbb{P}^{1}}(e_{1})\oplus
...\oplus\mathcal{O}_{\mathbb{P}^{1}}(e_{d-1})$ and $\mathbb{P(\mathcal{E})}$
the $\mathbb{P}^{d-2}$-bundle corresponding to the scroll $X.$
$\mathbb{P(\mathcal{E})}$ is a desingularisation of $X$ and we obtain the
following diagram:%
\[
\text{\begin{xy}
\xymatrix{
\mathbb{P(\mathcal{E})} \ar[r]^{\sigma} \ar[d] & X \ar[r] & \mathbb{P}^{g-1} \\
\mathbb{P}^1 & C \ar[l]_{\varphi_{\left\vert D\right\vert }} \ar[u] \ar@{^{(}->}[ru]_{\varphi_{K_{C}}} & \\
}
\end{xy}}%
\]}}
\end{remark}

{\normalsize \bigskip}

{\normalsize \label{scrolldetermination}The type $S(e_{1}%
,...,e_{d})$ of a scroll $X$ constructed from a basepoint free pencil of
divisors $D_{\lambda}$ on a variety $V$ can be determined by considering the following partition of $r+1:$%
\begin{align*}
d_{0} &  =h^{0}(V,\mathcal{O}_{V}(H))-h^{0}(V,\mathcal{O}_{V}(H-D))\\
d_{1} &  =h^{0}(V,\mathcal{O}_{V}(H-D))-h^{0}(V,\mathcal{O}_{V}(H-2D))\\
&  \overset{.}{:}\\
d_{i} &  =h^{0}(V,\mathcal{O}_{V}(H-iD))-h^{0}(V,\mathcal{O}_{V}(H-(i+1)D))\\
&  \overset{.}{:}%
\end{align*}
In \cite{schreyer1} 2.5. the author shows that $X$ is a $d_{0}-1$ ($=d-2$ for
$V=C\subset\mathbb{P}^{g-1}$ and $(D_{\lambda})_{\lambda}=g_{d}^{1}$ )
dimensional scroll of type $S(e_{1},...,e_{d_{0}})$ and the numbers $e_{i}$
are given by the dual partition:%
\[
e_{i}=\#\{j|d_{j}\geq i\}-1
\]
}

{\normalsize \bigskip}

\begin{example}{\normalsize (Canonical Curve $C\subset\mathbb{P}^{3}$) }
{\textnormal{Let $C\subset\mathbb{P}^{3}$ be a canonical curve of
degree $3$ and genus $4.$ Then $C$ admits a divisor $D$ of
degree $3$ with $\dim\left\vert D\right\vert =1$. From the geometric version
of Riemann-Roch $\dim\bar{D}=1,$ thus the pencil of divisors $(D_{\lambda
})_{\lambda}$ is cut out by trisecants on $C.$ We denote the class of
hyperplane sections in $\mathbb{P}^{3}$ by $H$, that cuts out the canonical divisors on $C.$ Then we get%
\[
d_{0}=h^{0}(C,\mathcal{O}_{C}(K_{C}))-h^{0}(C,\mathcal{O}_{C}(K_{C}-D))=4-2=2
\]%
\[
d_{1}=h^{0}(C,\mathcal{O}_{C}(K_{C}-D))-h^{0}(C,\mathcal{O}_{C}(K_{C}-2D))=2-h^{0}(C,\mathcal{O}_{C}(H|_{C}-2D))
\]%
\[
d_{2}=h^{0}(C,\mathcal{O}_{C}(K_{C}-2D))-h^{0}(C,\mathcal{O}_{C}(K_{C}-3D))=h^{0}(C,\mathcal{O}_{C}(K_{C}-2D))
\]
and%
\[
d_{i}=0\text{ for }i>2\text{ as }\deg(K_{C}-iD)=6-3i<0
\]
The value $h^{0}(C,\mathcal{O}_{C}(K_{C}-D))$ gives the number of hyperplane sections passing through
the points of $D$. As they all lie on a line, it follows that there exist
exactly two of them, and they intersect exactly in this line. If one of the
hyperplanes $H_{0}$ is tangent to $C$ in all the points of $D$, then
$h^{0}(C,\mathcal{O}_{C}(K_{C}-D))=h^{0}(C,\mathcal{O}_{C}(H|_{C}-2D))=1$, otherwise we get $h^{0}%
(C,\mathcal{O}_{C}(K_{C}-2D))=0.$ In summary we obtain a $2$-dimensional scroll $X$ of type
$S(1,1)$ or $S(2,0)$ depending on whether there exists such a special
hyperplane $H_{0}$ or not. In the general case the corresponding
$\mathbb{P}^{1}$-bundle $\mathbb{P(\mathcal{E})}$ is isomorphic to
$\mathbb{P}^{1}\times\mathbb{P}^{1}$ and $C$ lies on a smooth quadric surface.
\[
{\includegraphics[height=6.5cm,width=6.5cm]{plot1}%
}%
\]
Whereas in the more special case we have $\mathbb{P(\mathcal{E})\cong
}\mathbb{P}(\mathcal{O}_{\mathbb{P}^{1}}\mathcal{(}2)\oplus\mathcal{O}_{\mathbb{P}^{1}})$
$=P_{2}$, the second Hirzebruch surface, where the scroll $X$ becomes
singular. It is the cone over a conic in $\mathbb{P}^{2}.$%
\[
\includegraphics[height=7.5cm,width=6.5cm]{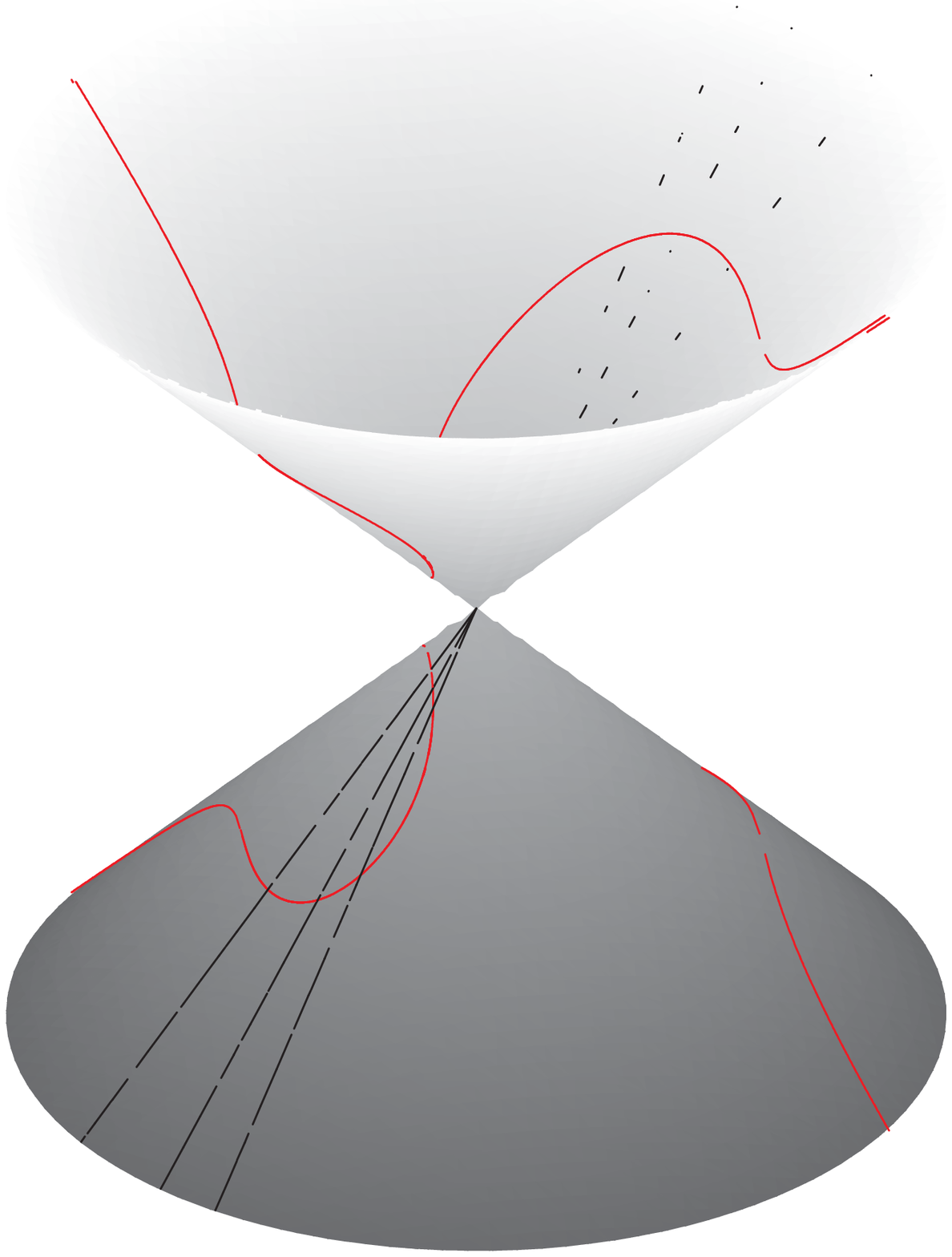}%
%
\]}}
\end{example}

\bigskip

{\normalsize Our aim was to determine a minimal free resolution of
$\mathcal{O}_{C}$ as $\mathcal{O}_{\mathbb{P}^{g-1}}-$module. We have
mentioned that we want to do this in two steps. First we consider a resolution
of $\mathcal{O}_{C}$ as $\mathcal{O}_{\mathbb{P(\mathcal{E})}}$-module. Then
we apply Theorem \ref{Cba} that leads to a mapping cone construction, which
gives us a (not necessarily) minimal resolution of $\mathcal{O}_{C}$ as
$\mathcal{O}_{\mathbb{P}^{g-1}}-$module: }

\begin{theorem}
{\normalsize \label{ResC}Let $C$ be a canonical curve $C\subset\mathbb{P}^{g-1}$
of genus $g$ that admits a base point free $g_{d}^{1}.$ Further let $X$ be the scroll
constructed from the $g_{d}^{1}$ as above and $\mathbb{P(\mathcal{E})}$ the
corresponding $\mathbb{P}^{d-2}-$bundle. Then }

\noindent{\normalsize 1) $C\subset\mathbb{P}(\mathcal{E})$ has a resolution
$F_{\ast}$ of type%
\[
0\rightarrow\mathcal{O}_{\mathbb{P(\mathcal{E})}}(-dH+(f-2)R)\rightarrow%
{\textstyle\sum\nolimits_{k=1}^{\beta_{d-2}}}
\mathcal{O}_{\mathbb{P(\mathcal{E})}}((-d+2)H+b_{k}R)\rightarrow
\]
}%
\[
{\normalsize \rightarrow%
{\textstyle\sum\nolimits_{k=1}^{\beta_{1}}}
\mathcal{O}_{\mathbb{P(\mathcal{E})}}(-2H+a_{k}R)\rightarrow\mathcal{O}%
_{\mathbb{P(\mathcal{E})}}\rightarrow\mathcal{O}_{C}\rightarrow0\ }%
\]

{\normalsize with $\beta_{i}=\frac{i(d-2-i)}{d-1}\binom{d}{i+1}$ and
$a_{k},b_{k}\in\mathbb{Z}$, $a_{k}+b_{k}=f-2,$ $a_{1}+...+a_{\beta_{d-2}%
}=2g-12.$ }

\noindent{\normalsize 2) $F_{\ast}$ is selfdual :%
\[
Hom(F_{\ast},\mathcal{O}_{\mathbb{P(\mathcal{E})}}(-dH+(f-2)R))\cong F_{\ast}%
\]
}

\noindent{\normalsize 3) If all $b_{k}\geq-1$ then an iterated mapping cone}%
\[
{\normalsize \left[  ...\left[  \mathcal{C}^{(f-2)}(-d)\rightarrow%
{\displaystyle\sum\limits_{k=1}^{\beta_{d-2}}}
\mathcal{C}^{b_{k}}(-d+2)\right]  \rightarrow...\right]  \rightarrow
\mathcal{C}^{0}}%
\]
{\normalsize  }

{\normalsize is a (not necessarily minimal) resolution of $\mathcal{O}_{C}$ as
an $\mathcal{O}_{\mathbb{P}^{g-1}}-$module. }
\end{theorem}

\begin{proof}
{\normalsize \cite{schreyer1} Corollary 4.4 and 6.7.}
\end{proof}

{\normalsize \bigskip}

From the existence of certain linear series on {\normalsize a canonical curve
$C$ one derives a nonsingular model $C^{\prime}\subset S$ on a smooth rational
surface $S$. Then the adjoint series $\left\vert K_{S}+C^{\prime}\right\vert $
embedds $C$ canonically. If $\left\vert K_{S}+C^{\prime}\right\vert $ is even
base point free on $S$, we can consider the image $S^{\prime}$ under the map
given by $\left\vert K_{S}+C^{\prime}\right\vert .$ In some cases we have an
analogue of the theorem above for the surface $S^{\prime}.$ Then from a
minimal free resolution of $S^{\prime}$ it is possible to deduce information
about the Betti numbers for $C.$ We write down the most important results as can
be found in \cite{schreyer1} Chapter 5: }

{\normalsize \bigskip}

{\normalsize \label{PkonScroll}For $H$ a sufficiently positive divisor on a
rational surface $S$ the image $S^{\prime}$ of $j:S\rightarrow\mathbb{P}%
H^{0}(S,\mathcal{O}_{S}(H))=\mathbb{P}^{r}$ can be described as a subvariety
of a scroll $X.$ It will turn out that in the cases of our interest $S^{\prime}$ is a determinantal surface on the scroll
$X,$ i.e. its vanishing ideal is given by the $2\times2$ minors of a $2\times
d$ matrix%
\[
\omega\sim%
\begin{pmatrix}
H-a_{1}R & ... & H-a_{d}R\\
H-(a_{1}+k)R & ... & H-(a_{d}+k)R
\end{pmatrix}
\]
with entries in $X.$ From this representation it is possible to obtain a free
resolution of $\mathcal{O}_{S^{\prime}}$ as $\mathcal{O}_{X}$-module. \\
\\
This can be seen as follows: Let us start with a rational ruled surface \index{$P_{k}$ k-th Hirzebruch surface}
\[
\pi:P_{k}:=\mathbb{P(\mathcal{O}}_{\mathbb{P}^{1}}\mathbb{\mathcal{(}%
}k\mathbb{\mathcal{)\oplus}\mathbb{\mathcal{O}}}_{\mathbb{P}^{1}%
}\mathbb{)\rightarrow P}^{1}\text{ , }k\geq0
\]
and consider $V:=S$ the surface obtained from $P_{k}$ via a sequence of
blowups:%
\[
\sigma:S\rightarrow P_{k}%
\]
By abuse of notation we denote the hyperplane class and the ruling of $P_{k}$
by $A$ and $B$ and also their pullbacks to $S.$ Further $E=%
{\textstyle\bigcup\nolimits_{i}}
E_{i}$ denotes the exceptional divisor of $\sigma.$ Let%
\[
H\sim dA+eB-%
{\textstyle\sum\nolimits_{i}}
\lambda_{i}E_{i}%
\]
be a divisor on $S$ with base point free complete linear series $\left\vert
H\right\vert $ and consider the map%
\[
j:S\rightarrow\mathbb{P}H^{0}(S,\mathcal{O}_{S}(H))=\mathbb{P}^{r}%
\]
with image $S^{\prime}\subset\mathbb{P}^{r}.$ }

{\normalsize \noindent Now suppose that }

\begin{enumerate}
\item {\normalsize $h^{0}(\mathcal{O}_{S}(H-B))\geq2$ }

\item {\normalsize $H^{1}(\mathcal{O}_{S}(kH-B))=0$ for $k\geq1$ and }

\item {\normalsize the map $S_{k}H^{0}\mathcal{O}_{S}(H)\rightarrow
H^{0}\mathcal{O}_{S}(kH)$ is surjective. }
\end{enumerate}

{\normalsize \noindent Then we apply our construction from above to obtain a
$(d+1)$-dimensional rational normal scroll%
\[
X=%
{\textstyle\bigcup\nolimits_{B_{\lambda}\in\left\vert B\right\vert }}
\bar{B}_{\lambda}%
\]
Let $\pi:\mathbb{P(\mathcal{E})\rightarrow P}^{1}$ denote the corresponding
$\mathbb{P}^{d}$-bundle and $S^{\prime\prime}$ the strict transform of
$S^{\prime}$ in $\mathbb{P(\mathcal{E})}$. Blowing up $S$ further we may
assume that $S\rightarrow S^{\prime}$ factors through $S^{\prime\prime}:$%
\[
\text{\begin{xy}
\xymatrix{
& S \ar[ld]\ar[r] \ar[rd] \ar[rdd]_{\pi_s} & S' \ar@{^{(}->}[r] & X \ar@{^{(}->}[r]  & \mathbb{P}^r \\
P_k \ar[rrd]_{\pi_{P_k}}&   & S'' \ar[u] \ar[d] \ar[r]& \mathbb{P}(\mathcal{E}) \ar[ld]^{\pi_{\mathbb{P}(\mathcal{E})}} \ar[u]  & \\
&     & \mathbb{P}^1                 & \\
}
\end{xy}}%
\]
If the conditions $1.-3.$ are fulfilled, then we can describe the syzygies of
$\mathcal{O}_{S^{\prime\prime}}$ in terms of $\mathcal{O}%
_{\mathbb{P(\mathcal{E})}}$-modules: }

\begin{theorem}
{\normalsize \label{Resolutiononscroll}$\mathcal{O}_{S^{\prime\prime}}$ has an
$\mathcal{O}_{\mathbb{P(\mathcal{E})}}$-module resolution of type%
\[
0\rightarrow%
{\textstyle\sum\nolimits_{k=1}^{\beta_{d-1}}}
\mathcal{O}_{\mathbb{P(\mathcal{E})}}(-dH+b_{k}^{(d-1)}R)\rightarrow...
\]%
\[
...\rightarrow%
{\textstyle\sum\nolimits_{k=1}^{\beta_{1}}}
\mathcal{O}_{\mathbb{P(\mathcal{E})}}(-2H+b_{k}^{(1)}R)\rightarrow
\mathcal{O}_{\mathbb{P(\mathcal{E})}}\rightarrow\mathcal{O}_{S^{\prime\prime}%
}\rightarrow0
\]
with $\beta_{i}=i\binom{d}{i+1}.$ }
\end{theorem}

\begin{proof}
{\normalsize \cite{schreyer1} 5.1-5.5.}
\end{proof}

{\normalsize \bigskip}

{\normalsize Especially in the cases of our interest, where we have $d=3,$ the
resolution above is of the following type:%
\[
0\rightarrow\mathcal{O}_{\mathbb{P(\mathcal{E})}}(-3H+b_{1}^{(2)}%
R)\oplus\mathcal{O}_{\mathbb{P(\mathcal{E})}}(-3H+b_{2}^{(2)}R)\rightarrow
\]%
\[
\overset{\omega}{\rightarrow}\mathcal{O}_{\mathbb{P(\mathcal{E})}}%
(-2H+b_{1}^{(1)}R)\oplus...\oplus\mathcal{O}_{\mathbb{P(\mathcal{E})}%
}(-2H+b_{3}^{(1)}R)\rightarrow\mathcal{O}_{\mathbb{P(\mathcal{E})}}%
\rightarrow\mathcal{O}_{S^{\prime\prime}}\rightarrow0
\]
}

{\normalsize \noindent Further the map $\omega$ is given by a $2\times3$
matrix%
\[
\omega\sim%
\begin{pmatrix}
H-a_{1}R & ... & H-a_{3}R\\
H-(a_{1}+k)R & ... & H-(a_{3}+k)R
\end{pmatrix}
\]
with entries in $\mathbb{P(\mathcal{E})}$ and the $2\times2$ minors of
$\omega$ generate the vanishing ideal of $S^{\prime\prime}\subset
\mathbb{P(\mathcal{E})}$ (cf. \cite{schreyer1} 5.5.). }
\\\\
{\normalsize The following important theorem due to Schreyer (\cite{schreyer1}
Theorem 5.7) states that there also exists a partial converse of this result: Let
$X\subset\mathbb{P}^{r}$ be a scroll of dimension $d+1$ and $\pi
:\mathbb{P(\mathcal{E})\rightarrow P}^{1}$ the corresponding $\mathbb{P}^{d}%
$-bundle. Further let $S^{\prime\prime}\subset\mathbb{P(\mathcal{E})}$ denote
the irreducible surface defined by the $2\times2$ minors of a matrix $\omega$
on $\mathbb{P(\mathcal{E})}$ as above, then the image $S^{\prime}$ of
$S^{\prime\prime}$ in $\mathbb{P}^{r}$ can be obtained as the image of a
blowup of $P_{k}=\mathcal{O}_{\mathbb{P}^{1}}(k)\oplus\mathcal{O}%
_{\mathbb{P}^{1}}$ defined by a complete linear system $H$ as above. We denote
by $A$ and $B$ the class of a hyperplane and a ruling on $P_{k}.$ Then with
\[
a=a_{1}+...+a_{d}\text{ \ and \ }f=\deg X
\]
we have: }

\begin{theorem}
{\normalsize \label{ConPk}$S^{\prime}\subset\mathbb{P}^{r}$ is the image of
$P_{k}$ under a rational map defined by a subseries of
\[
H^{0}(P_{k},\mathcal{O}_{P_{k}}(dA+(f-dk-a)B)
\]
which has%
\[
\delta=df-\frac{d(d+1)}{2}k-(d+1)a
\]
assigned base points. Furthermore, if $S^{\prime}\subset X\subset
\mathbb{P}^{r}$ contains a canonical curve $C$ of genus $r+1$, then the ruling
of $X$ cuts on $C$ a $g_{d+2}^{1}$ and the strict transform $C^{\prime}$ of
$C$ in $P_{k}$ is a divisor of class%
\[
C^{\prime}\sim(d+2)A+(f-(d+1)k-a+2)B
\]
and arithmetic genus
\[
p_{a}C^{\prime}=r+1+\delta.
\]
}
\end{theorem}

\begin{proof}
{\normalsize \cite{schreyer1} Theorem 5.7. }
\end{proof}

{\normalsize \bigskip}

\noindent{\normalsize We will use this important theorem to show that a curve
$C$ has a certain model on a blowup of $P_{k}$ and therefore that there exist
further special linear systems on $C.$ We have already mentioned above that we
are also interested in the converse, i.e. we start with a model $C^{\prime
}\subset$ $P_{k}$ of $C$ and want to give a description of the image
$S^{\prime}$ of $P_{k}$ under the mapping defined by the adjoint series. For
this reason we have to show that this series is base point free and fulfills
the conditions $1.-3.$ from above. We treat this problem in the next
chapter.\newpage}

\chapter{{\protect\normalsize Ampleness of the Adjoint series}}
{\normalsize \setcounter{theorem}{0} Let $C^{\prime}\subset X$ be a curve on a
surface $X$, that has only singularities in the points $p_{1},...,p_{s}$ with
multiplicity $2.$ After blowing up these singularities, we get a smooth curve
$C\subset S=\tilde{X}(p_{1},...,p_{m}).$ \index{$\tilde{X}(p_{1},...,p_{m})$ blowup of a variety in the points $p_{1},...,p_{m}$} Now we are interested if certain
linear systems on $S,$ especially the adjoint series }$\left\vert
K_{S}+C\right\vert ${\normalsize , are $i-$very ample on $S$ for $i=0,1.$
There we call a linear system $\left\vert L\right\vert $ on a smooth
projective surface $0-$very ample if it is base point free and $1-$very ample
if it is very ample. We follow the approach of Roland Weinfurtner in his
PhD thesis \cite{weinfurtner} that makes use of Reider's Theorem (cf.
\cite{Reider}) in a modified version: }

\begin{theorem}
{\normalsize (Modified version of Reider's Theorem) Let $L$ be a line bundle
on a projective surface $X$ and $L^{2}\geq5+4i.$ If $|K_{X}+L|$ is not $i-$very
ample then there exists an effective divisor $D$ on $X$ with $L-2D$
$\mathbb{Q-}$effective ($\exists n\in\mathbb{Z}^{+}:n(L-2D)$ is effective) and
a $0-$cycle $Z$ of degree $\leq i+1$, where $|K_{X}+L|$ is not $i-$very ample
and the following inequality holds:%
\[
D.(L-D)\leq i+1
\]
}
\end{theorem}

\begin{proof}
{\normalsize \cite{BFS} Theorem 2.1. and \cite{weinfurtner} Theorem 1.2. }
\end{proof}

{\normalsize \bigskip}

{\normalsize If $C^{\prime}$ has singularities in the points $p_{1}%
,...,p_{s},$ an iterated blowup of $X$ in $p_{1},...,p_{s}$ gives a
desingularisation of $C^{\prime}:$%
\[
S=\tilde{X}(p_{1},...,p_{s})\overset{\sigma_{s}}{\rightarrow}...\overset
{\sigma_{2}}{\rightarrow}\tilde{X}(p_{1})\overset{\sigma_{1}}{\rightarrow
}X%
\]
}

{\normalsize \noindent Let $\sigma$ be the composition of the blowups
$\sigma_{i}$ and $E_{j}:=\sigma_{s}^{\ast}(\sigma_{s-1}^{\ast}...(\sigma
_{j}^{-1}(p_{j}))...)$ the total transform of the exceptional divisor of the
blowup of the point $p_{j}$ on $S$. Further $C$ denotes the strict transform
of $C^{\prime}.$ Then from the adjunction formula we get%
\[
\omega_{S}(C)\otimes\mathcal{O}_{C}\cong\omega_{C}%
\]
and therefore%
\[
K_{C}\sim (K_{S}+C)|_{C}
\]
}

\noindent {\normalsize For $X=\mathbb{P}^{2}$ or $P_{2}$ let $H$ denote the class of a
hyperplane section and by abuse of notation also its pullback to $S$. $R$ denotes the class of a ruling on $P_{2}$. In the case $X=\mathbb{P}%
^{1}\times\mathbb{P}^{1}$ we have $\operatorname{Pic}(\mathbb{P}^{1}\times\mathbb{P}%
^{1})=\mathbb{Z\cdot(}0,1)\oplus\mathbb{Z\cdot}(1,0)$ with factor classes $(0,1)$ and $(1,0).$
Again by abuse of notation we also denote its pullbacks to $S$ by
$(0,1),(1,0)$. Then intersection theory gives the following theorem: }

\begin{theorem}
{\normalsize \label{totaltransform}The Picard group $Pic(S)$ is generated by
the pullback $\sigma^{\ast}Pic(X)$ of the Picard group of $X$ and the
exceptional divisors $E_{j}$ as defined above. Further }

{\normalsize 1) $E_{j}^{2}=-1$ }

{\normalsize 2) $E_{i}.E_{j}=0$ for $i\neq j$ }

{\normalsize 3) $\Gamma.E_{j}=0$ for $\Gamma\in\sigma^{\ast}Pic(X)$ }

{\normalsize 4) a) $H^{2}=1$ for $X=\mathbb{P}^{2}$ }

{\normalsize \ \ \ \ b) $H^{2}=2,$ $H.R=0$ and $R.R=0$ for $X=P_{2}$ }

{\normalsize \ \ \ \ c) $(0,1).(0,1)=(1,0).(1,0)=0$ and $(1,0).(0,1)=1$ for
$X=\mathbb{P}^{1}\times\mathbb{P}^{1}.$ }
\end{theorem}

\begin{proof}
{\normalsize \cite{hartshorne} V.3.2. }
\end{proof}

{\normalsize \bigskip}

{\normalsize If $C^{\prime}$ has only ordinary nodes, then all exceptional curves
$E_{j}$ are irreducible and intersect the strict transform $C$ of
$C^{\prime}$ in exactly two points. For }$\operatorname*{char}%
\Bbbk\neq2$, {\normalsize a singularity of multiplicity 2 is always of
analytical type $y^{2}-x^{k}=0$ with $k\geq2$. If there exists such a
singularity with $k\geq4$, for example if $C$ has a tacnode:%
\[
\psset{unit=0.6}\begin{pspicture}(0,-1)(4,5)
\pscurve[showpoints=false](0.5,0.5)(2.25,2)(4,0.5)
\pscurve[showpoints=false](0.5,3.5)(2.25,2)(4,3.5)
\rput[bl](4,3.6){$C$}
\end{pspicture}
\]
then blowing up $X$ we get the following picture:%
\[
\psset{unit=0.8}\begin{pspicture}(0,-1)(14,5)
\pscurve[showpoints=false](0.5,0.5)(2.25,2)(4,0.5)
\pscurve[showpoints=false](0.5,3.5)(2.25,2)(4,3.5)
\psdots(2.25,2)
\pscurve(6,0.5)(7,2)(9,3)
\pscurve(6,3.5)(7,2)(9,1)
\psline(7,0.5)(7,3.5)
\psdots(7,2)
\psdots(12.41,1.96)
\psdots(13.25,2.54)
\psline(11.5,0.5)(11.5,3.5)
\psline(11,1)(14,3)
\pscurve(11.75,0.5)(13,2.5)(14,2.25)
\psline[linewidth=3pt]{<-}(4.5,2)(5.5,2)
\psline[linewidth=3pt]{<-}(9.5,2)(10.5,2)
\rput[bl](4,3.6){$C$}
\rput[bl](9,3.1){$C$}
\rput[t](7,0.5){$E_1'=E_1$}
\rput[t](13.8,3.5){$E_2'=E_2$}
\rput[t](11.5,4){$E_1'$}
\rput[t](12,0){$E_1=E_1'+E_2'$}
\rput[t](14.3,2.2){$C$}
\end{pspicture}
\]
where the strict transform $\tilde{C}$ of $C$ under the
first blowup $\sigma_{1}$ still has a double point lying on the exceptional
divisor $E_{1}.$ Blowing up in this double point again, we get a smooth curve
$C$ that intersects the total transforms $E_{1}$ and $E_{2}$ in exactly two
points. In this situation $E_{1}$ decomposes into the strict transforms
$E_{1}^{\prime}$ and $E_{2}^{\prime}$ and $E_{2}=E_{2}^{\prime}$ is
irreducible. Therefore it is useful to state the following definition: }

\begin{definition}
{\normalsize (infinitely near points) A point $p_{j}$ is called
infinitely near to another point $p_{i}$ if and only if $p_{j}$ lies on the strict
transform of the exceptional divisor in the blowup of $p_{i}.$ }
\end{definition}

We assumed that $C^{\prime}$ has only double points as singularities. Then for
the strict transform $C^{(i)}\subset\tilde{X}(p_{1},...,p_{i})$ of $C^{\prime
}$ under the blowup $\sigma^{(i)}:=\sigma_{i}\circ...\circ\sigma_{1}$ there
are several possibilities: One is that $C^{(i)}$ meets the exceptional divisor
$E^{(i)}:=\sigma_{i}^{-1}(p_{i})$ transversally in two distinct points. The
second is that $C^{(i)}$ meets $E^{(i)}$ in one point $P$, with $C^{(i)}$
nonsingular in $P$, but $C^{(i)}$ and $E^{(i)}$ having intersection
multiplicity $2$ at $P.$ The third possibility is that $C^{(i)}$ has a double
point $P$ on $E^{(i)}$ (see the example above)$.$ In this case $E^{(i)}$ must
pass through $P$ in a direction not equal to any tangent direction of $P$,
since $C^{(i)}.E^{(i)}=2.$ It follows that for each point $p_{i}$ there exists
at most one further point $p_{j}$, $j>i,$ that lies infinitely near to
$p_{i}.$ From our definition of the total transforms $E_{i}$ above, it turns
out that $E_{i}=${\normalsize $\sigma_{s}^{-1}(\sigma_{s-1}^{-1}...(\sigma
_{j}^{-1}(p_{j}))...).$}

\noindent{\normalsize At any time the total transforms contain exactly the
strict transforms as components and they are all irreducible if and only if none of the
points $p_{j}$ is infinitely near to another point $p_{i}$. To be more
precisely the strict transforms $E_{i}^{\prime}$ are inductively given by%
\begin{align*}
E_{s}^{\prime}  &  =E_{s}\\
E_{s-1}^{\prime}  &  =E_{s-1}-\delta_{s-1,s}E_{s}\\
E_{s-2}^{\prime}  &  =E_{s-2}-\delta_{s-2,s-1}E_{s-1}-\delta_{s-2,s}E_{s-2}\\
&  :\\
&  :\\
E_{1}^{\prime}  &  =E_{1}-%
{\textstyle\sum\nolimits_{i=2}^{s}}
\delta_{1,i}E_{i}%
\end{align*}
}

\noindent{\normalsize with
\[
\mathcal{\delta}_{i,s}=\left\{
\begin{array}
[c]{l}%
1\text{ if }p_{s}\text{ is lying on }E_{i}^{\prime}\text{ in }\tilde{X}%
_{s-1}\\
0\text{ else}%
\end{array}
\right.
\]
Therefore we obtain a decomposition $E_{i}=%
{\textstyle\sum\nolimits_{k=i}^{s}}
\delta_{k}E_{k}^{\prime}$, $\delta_{k}=0$ or $1$, of $E_{i}$ as sum of
irreducible components. In our situation, this decomposition is given as follows}: If there exists a maximal chain of points
$p_{i},...,p_{i+k_{i}},$ such that $p_{i+k+1}$ lies infinitely near to
$p_{i+k}$ for $k=0,...,k_{i}-1$ and $p_{i}$ does not lie infinitely near to any other double point
$p_{j}$, then
\begin{align*}
E_{i}  &  =E_{i}^{\prime}+...+E_{i+k}^{\prime}\\
E_{i+1}  &  =E_{i+1}^{\prime}+...+E_{i+k}^{\prime}\\
&  \overset{.}{:}\\
E_{i+k}  &  =E_{i+k}^{\prime}%
\end{align*}
Especially in the case where $p_{i}$ is an ordinary double point of $C^{\prime
}$ we have $E_{i}=E_{i}^{\prime}.$ Further, we see that each strict
transform $E_{i}^{\prime}$ can be written in terms of the total transforms,
and $E_{i}^{\prime}=E_{i}$ iff $E_{i}$ is irreducible and $E_{i}^{\prime
}=E_{i}-E_{i+1}$ iff $p_{i+1}$ lies infinitely near to $p_{i}$.

\bigskip

{\normalsize The strict transform $C$ of }$C^{\prime}$ {\normalsize is a
divisor of class
\[
C\sim\sigma^{\ast}C^{\prime}-%
{\textstyle\sum\nolimits_{i=1}^{s}}
2E_{i}%
\]
and therefore the adjoint series is given by an effective divisor linear
equivalent to $K_{S}+C\sim\sigma^{\ast}\Delta-%
{\textstyle\sum\nolimits_{i=1}^{s}}
E_{i},$ $\Delta\in \operatorname{Pic(X)}.$ Let $L^{\prime}\sim\sigma^{\ast}\Sigma-%
{\textstyle\sum\nolimits_{i=1}^{s}}
\lambda_{i}E_{i},$ $\Sigma\in \operatorname{Pic(X)}$ be an effective divisor on $S$, then
assuming that $L^{2}\geq5+4i$, $L:=L^{\prime}-K_{S},$ Reider's Theorem says
that $\left\vert L^{\prime}\right\vert $ is $i-$very ample if we cannot find
an effective divisor $D\sim\sigma^{\ast}\Gamma-%
{\textstyle\sum\nolimits_{i=1}^{s}}
\alpha_{i}E_{i},$ $\Gamma\in \operatorname{Pic(X)}$ on $S$, such that $L-2D$ is $\mathbb{Q}$-effective and $D.(L-D)\geq2+i.$ Under the assumption $0\leq\lambda_{1}%
\leq...\leq\lambda_{s}\leq2,$ it is possible to restrict to a very manageable
set of divisors }$D${\normalsize . We start with a divisor $D$ as above that
fullfills $D.(L-D)<2+i$, then after several modification steps we get a
divisor $\tilde{D}$ with $\alpha_{i}=0$ or $1$, such that $L-2\tilde{D}$ is
also $\mathbb{Q}$-effective and $\tilde{D}.(L-\tilde{D})\leq D.(L-D)<2+i:$}

\begin{lemma}
\bigskip{\normalsize \label{ampleness}Let $X\mathbb{\ }$be a surface,
$L^{\prime}\sim\sigma^{\ast}\Sigma-%
{\textstyle\sum\nolimits_{i=1}^{s}}
\lambda_{i}E_{i},$ $\Sigma\in \operatorname{Pic(X)}$, $0\leq\lambda_{1}\leq...\leq\lambda
_{s}\leq2,$ an effective divisor on the iterated blowup $S=\tilde{X}%
(p_{1},...,p_{s})$ as above and $L^{2}=(L^{\prime}-K_{S})^{2}\geq5+4i.$ If
there exists no effective divisor $D\sim\sigma^{\ast}\Gamma-%
{\textstyle\sum\nolimits_{i=1}^{s}}
\alpha_{i}E_{i},$ $\Gamma\in \operatorname{Pic(X)}\backslash\{0\},$ $\alpha_{i}\in\{0,1\},$
or $D\sim E_{i}^{\prime}$\ on $S$, such that $L-2D$ is $\mathbb{Q}$-effective
and $D.(L-D)\leq 1+i,$ then $\left\vert L^{\prime}\right\vert $ is $i$-very
ample on }$S${\normalsize ($i=0,1$).}
\end{lemma}

\begin{proof}
According to Reider's Theorem we have to assure that there exists no divisor
$D${\normalsize $\sim\sigma^{\ast}\Gamma-%
{\textstyle\sum\nolimits_{i=1}^{s}}
\alpha_{i}E_{i},$ $\Gamma\in \operatorname{Pic(X)}$ on $S$, such that $L-2D$ is $\mathbb{Q}$-effective and $D.(L-D)\leq 1+i.$ For }$\Gamma\sim0$ the only effective divisors
are sums of the strict transforms $E_{i}^{\prime}.$ The calculation in (i)
and (ii) below shows that we can restrict to a single $E_{i}^{\prime}.$

For $\Gamma\not \sim 0,$ we start with an effective divisor
{\normalsize $D\sim\sigma^{\ast}\Gamma-%
{\textstyle\sum\nolimits_{i=1}^{s}}
\alpha_{i}E_{i},$ $\Gamma\in \operatorname{Pic(X)}$ on $S$, such that $L-2D$ is $\mathbb{Q}$
effective and $D.(L-D)\leq1+i.$ In the first step we will show that we can
restrict to those divisors }$D$ with $\alpha_{i}\geq0.$
\\\\
\textit{1st Modification step}: I{\normalsize f some of the coefficients
$\alpha_{i}$ are negativ, then we aim to modify $D$ to an effective divisor
$\tilde{D}\sim\sigma^{\ast}\Gamma-%
{\textstyle\sum\nolimits_{i=1}^{s}}
\tilde{\alpha}_{i}E_{i}$ with $\tilde{\alpha}_{i}\geq0$ for all $i.$ In each
step we have to assure that $\tilde{D}.(L-\tilde{D})\leq D.(L-D).$ This is
done in the following way: If none of the $\alpha_{i}$ is negative, there is
nothing to do. Otherwise let $j$ be the maximal index with $\alpha_{j}<0.$ If
$E_{j}$ is irreducible we have $E_{j}=E_{j}^{\prime}$ and in this situation we
set $\tilde{D}:=D-E_{j}$. In the case where $E_{j}$ is reducible we have
$E_{j}^{\prime}=E_{j}-E_{j+1}$. Here we consider $\tilde{D}:=D-(E_{j}-E_{j+1}).$
Because of }$E_{j}^{\prime}.D=\alpha_{j}<0$ the strict transform
$E_{j}^{\prime}$ is a component of $D,$ hence {\normalsize $\tilde{D}$
stays effective. The second condition is satisfied if }
\\\\
{\normalsize (i)%
\[
\tilde{D}.(L-\tilde{D})\leq D.(L-D)
\]%
\begin{align*}
&  \Leftrightarrow(D-E_{j}).(L-D+E_{j})\leq D.(L-D)\\
&  \Leftrightarrow2D.E_{j}+1-L.E_{j}\leq0\\
&  \Leftrightarrow2\alpha_{j}-(\lambda_{j}+1)+1=2\alpha_{j}-\lambda_{j}\leq0
\end{align*}
and}
\\\\
{\normalsize (ii)
\[
\tilde{D}.(L-\tilde{D})\leq D.(L-D)
\]%
\begin{align*}
&  \Leftrightarrow(D-(E_{j}-E_{j+1})).(L-D+(E_{j}-E_{j+1}))\leq D.(L-D)\\
&  \Leftrightarrow2D.(E_{j}-E_{j+1})+2-L.(E_{j}-E_{j+1})\leq0\\
&  \Leftrightarrow2(\alpha_{j}-\alpha_{j+1})+(\lambda_{j}-\lambda_{j+1}%
)+2\leq0.
\end{align*}
\noindent As $\lambda_{i}\geq0$ the condition (i) is always fulfilled and
because for $\lambda_{j}\leq\lambda_{j+1}$ also condition (i) is given. Now
the following sequence of modification steps leads us to an effective divisor
$\tilde{D}\sim\sigma^{\ast}\Gamma-%
{\textstyle\sum\nolimits_{i=1}^{s}}
\tilde{\alpha}_{i}E_{i}$ with $\tilde{\alpha}_{i}\geq0$ for all $i:$ If
$E_{j}=E_{j}^{\prime}$ we can substract the divisor $E_{j}^{\prime}$
$\left\vert \alpha_{j}\right\vert $ times to get a divisor $\tilde{D}%
\sim\sigma^{\ast}\Gamma-%
{\textstyle\sum\nolimits_{i=1}^{s}}
\tilde{\alpha}_{i}E_{i}$ with $\alpha_{i}\geq0$ for $i\leq j.$ }

{\normalsize \noindent In the case $E_{j}^{\prime}=E_{j}-E_{j+1}$ the
situation is a little more complicated: For $\alpha_{j+1}\geq\left\vert
\alpha_{j}\right\vert $ we can substract $E_{j}^{\prime}$ $\alpha_{j+1}+1$
times to get a divisor $\tilde{D}$ with $\alpha_{i}\geq0$ for $i\leq j.$
Otherwise $\tilde{D}\sim D-(\alpha_{j+1}+1)E_{j}^{\prime}$ has a negative
coefficient $\tilde{\alpha}_{j+1}=\alpha_{j+1}-(\alpha_{j+1}+1)=-1$. Now we
have to apply some correction steps to assure that $\tilde{D}$ has non negative
coefficients $\alpha_{i}$ for $i>j:$ }
\\\\
{\normalsize \textit{a)} $E_{j+1}$ irreducible: We substract
$E_{j+1}$ to get an effective divisor $\tilde{D}$ with $\alpha_{i}\geq0$ for
all $i>j$, such that we can go on with our procedure$.$ }
\\\\
\noindent
{\normalsize \textit{b)} $E_{j+1}^{\prime}=E_{j+2}-E_{j+1}$ reducible: We
substract $E_{j+1}^{\prime}$ first to get $\alpha_{j+1}\geq0.$ If
$\alpha_{j+2}$ becomes negative in this step, we apply one further correction
step of type a) or b) and so on. This chain of steps stops with the last
exceptional divisor $E_{s}=E_{s}^{\prime}$, which is irreducible. }
\\\\
{\normalsize Now we continue with the divisor $D=\tilde{D}$ which has non
negative coefficients $\alpha_{i}$ for $i\geq j-1.$ After finitely many steps
we obtain an effective divisor $\tilde{D}$ with only non negative $\alpha
_{i}.$}
\\\\
\textit{2nd Modification step}: {\normalsize The following consideration shows
that we can even restrict to those divisors $D$ with $\alpha_{i}=0$ or $1:$
Starting with a divisor $D$ with non negative $\alpha_{i}$, we consider
$\tilde{D}$ with $\tilde{\alpha}_{i}=\min(\alpha_{i},1).$ $\tilde{D}$ is again
effective, as we have only added multiples of the exceptional divisors
$E_{i}.$ It remains to show that $(D+E_{i}).(L-(D+E_{i}))\leq D.(L-D)$ for all
$i$ with $\alpha_{i}\geq2,$ which is true for $\lambda_{i}\leq2$ because of%
\begin{align*}
(D+E_{i}).(L-(D+E_{i}))  &  =D.(L-D)+1+L.E_{i}-2E_{i}.D\\
&  =D.(L-D)+2+\lambda_{i}-2\alpha_{i}=\\
&  =D.(L-D)+\lambda_{i}-2(\alpha_{i}-1)\leq D.(L-D)
\end{align*}
}
\end{proof}

{\normalsize \bigskip}

Now {\normalsize assume that $\left\vert L^{\prime}\right\vert $ is
base point free and $d=\dim\left\vert L^{\prime}\right\vert \geq3$, then we can
consider its image $S^{\prime}$ under the morphism defined by this complete
linear system }

{\normalsize
\[
\varphi:S\overset{\left\vert L\right\vert }{\rightarrow}S^{\prime}%
\subset\mathbb{P}^{d-1}%
\]
}

{\normalsize \noindent Our aim is to show that for certain linear systems the
image $S^{\prime}$ is an arithmetically Cohen-Macaulay surface. In the case
where $X=\mathbb{P}^{2}$ and $L^{\prime}\sim dH-%
{\textstyle\sum\nolimits_{i=1}^{s}}
E_{i}$ these surfaces are called Bordiga surfaces. Weinfurtner gives in his
PhD thesis \cite{weinfurtner} a complete description of them. We are
interested in the following linear systems: }

\section{Adjoint linear series on blowups of
$\mathbb{P}^{2}$}

{\normalsize As we have mentioned before we have to examine the case where $C$
is an irreducible, smooth, canonical curve of genus 9 that admits a plane
model $C^{\prime}\subset\mathbb{P}^{2}=X$ of degree $d=7.$ $C^{\prime}$ has
singularities of multiplicity 2 in exactly $s=\binom{d-1}{2}-9=6$ points
$p_{1},...,p_{6}.$ Blowing up these singularities%
\[
\sigma:S=\tilde{X}(p_{1},...,p_{6})\overset{\sigma_{6}}{\rightarrow
}...\overset{\sigma_{2}}{\rightarrow}\tilde{X}(p_{1})\overset{\sigma_{1}%
}{\rightarrow}\mathbb{P}^{2}%
\]
we can assume that $C$ is the strict transform of $C^{\prime}$ and
\[
C\sim7H-%
{\textstyle\sum\nolimits_{i=1}^{6}}
2E_{i}%
\]
its canonical system is cut out by the complete linear series $\left\vert
L^{\prime}\right\vert \sim\left\vert C+K_{S}\right\vert =\left\vert 4H-%
{\textstyle\sum\nolimits_{i=1}^{s}}
E_{i}\right\vert $}

\begin{theorem}
{\normalsize \label{ALSP2}Let $S$ be the iterated blowup of }%
$6$ {\normalsize points $p_{1},...,p_{6}$ on $\mathbb{P}^{2}$ and $C\sim 7H-%
{\textstyle\sum\nolimits_{i=1}^{6}}
2E_{i}$ an irreducible, nonsingular curve on $S$, then the adjoint linear
series $\left\vert C+K_{S}\right\vert =\left\vert 4H-%
{\textstyle\sum\nolimits_{i=1}^{s}}
E_{i}\right\vert $ is base point free. It is very ample on $S$ if and only if none
of the points $p_{i}$ lies infinitely near to another one. The image
$S^{\prime}\subset\mathbb{P}^{8}$ of $S$ under the morphism defined by
$\left\vert C+K_{S}\right\vert $ is arithmetically Cohen-Macaulay. Furthmermore if two of the points $p_{i}$ are infinitely near, then $S^{\prime}$ has only isolated singularities that are contractions of strict transforms }$E_{i}^{\prime}\sim
E_{i}-E_{i+1}.$
\end{theorem}

\begin{proof}
{\normalsize For $L:=L^{\prime}-K_{S}\sim C$ we get $L^{2}\geq49-24=25\geq
5+2\cdot1$, hence the first condition of Lemma \ref{ampleness} for $\left\vert
L^{\prime}\right\vert $ to be $i$-very ample ($i=0,1$) is fulfilled. We call
an effective divisor $D$ $i-$\textit{critical} if $L-2D$ is $\mathbb{Q-}%
$effective and $D.(C-D)<2+i.$ Then according to Lemma \ref{ampleness} it
remains to consider $D$ to be one of the following three types:}
\\\\
{\normalsize \textit{(a)} $D\sim E_{k}$ for a $k\in\{1,...,6\}.$}

\noindent
{\normalsize \textit{(b)} $D\sim E_{k}-E_{k+1}$ $k\in
\{1,...,6\}.$}

\noindent
{\normalsize \textit{(c)} $D\sim eH-%
{\textstyle\sum\nolimits_{i=1}^{6}}
\alpha_{i}E_{i}$ with $e\geq1$ and $\alpha_{i}=0$ or $1$. }
\\\\
{\normalsize \noindent As $E_{k}.(C-E_{k})=3\geq2+i$ there exists no
$i-$critical divisor of type (a). If two points $p_{k}$ and $p_{k+1}$ are
infinitely near, then there exists an effective divisor of type (b) and because
of $(E_{k}-E_{k+1}).(C-(E_{k}-E_{k+1}))=2$, this divisor is $1-$critical. For
a divisor of type (c) to be $i-$critical we can restrict to $e\leq3$ as
}$L-2D$ has to be {\normalsize $\mathbb{Q-}$effective. Let $\delta
=\#\{\alpha_{i}:\alpha_{i}=1\},$ then
\[
D.(C-D)=e(7-e)-\delta<3
\]
if and only if $e=1$ and $\delta>3$. In this situation we would get $D.C\leq-1$ and therefore $D$ and
$C$ have a common component, which contradicts to our assumption that $C$
is irreducible. In summary, there exists at most $1$-critical effective
divisors $D$, which is exactly the case if two points $p_{k}$ and $p_{l}$ are
infinitely near. According to Lemma \ref{ampleness}, the
map%
\[
\varphi:S\overset{\left\vert C+K_{S}\right\vert }{\rightarrow}S^{\prime
}\subset\mathbb{P}^{8}%
\]
defined by the adjoint series is base point free. The following
arguments show that $\varphi$ is even very ample outside the finite union
$\Delta=%
{\textstyle\bigcup\nolimits_{\lambda}}
D_{\lambda}$ of all $1-$critical critical divisors (being of type
$E_{k}-E_{k+1}$). With $U:=S\backslash\Delta$, this is exactly the case iff
$\left\vert L^{\prime}-E_{p}\right\vert $ is base point free for all $p\in U$,
where $E_{p}$ denotes the exceptional divisor of the blowup of $p\in S.$ We
want to apply Lemma \ref{ampleness} to $\left\vert L^{\prime}-E_{p}\right\vert
.$ This is possible as $(7H-%
{\textstyle\sum\nolimits_{i=1}^{6}}
2E_{i}-2E_{p})^{2}=49-28\geq5$. A $0$-critical divisors $D=D^{\prime}-\alpha_{p}E_{p}$} {\normalsize for
$\left\vert L^{\prime}-E_{p}\right\vert $, $\alpha_{p}=0$ or $1$, has to satisfy the following
inequality:%
\begin{align*}
D.(L-D-2E_{p})\leq1\Rightarrow (D^{\prime}-\alpha_{p}E_{p}).(L-D^{\prime}-(2-\alpha_{p})E_{p})\leq1\\
\Rightarrow D^{\prime}.(L-D^{\prime})-\alpha_{p}(2-\alpha_{p})\leq1
\end{align*}
Because of $D^{\prime}.(L-D^{\prime})\geq2$ this is only possible if
$\alpha_{p}=1$ and $D.(L^{\prime}-D)=2.$ But then $p\in D$ and $D$ is a divisor of type
$E_{k}-E_{k+1}$, which was excluded. }

{\normalsize Thus the effective divisors of type $E_{k}-E_{k+1}$ are
contracted by the morphism $\varphi$, whereas $\varphi$ is very ample on $U.$
The surface $S^{\prime}$ is called a Bordiga surface. It is smooth iff none of
the points $p_{i}$ is infinitely near to another one and it has isolated
singularities otherwise. }

It remains to show that $S^{\prime}$ is arithmetically Cohen-Macaulay.
{\normalsize For this purpose we aim to find a smooth hyperplane section on
$S^{\prime}$, that is arithmetically Cohen-Macaulay. As a direct
consequence $S^{\prime}$ is then arithmetically Cohen-Macaulay, too (see Lemma
\ref{aCM} below). $S^{\prime}$ has only isolated singularities and therefore
we can find a smooth hyperplane section $\Gamma$ applying Bertini's Theorem.
$\Gamma$ is obtained from an element of the linear system $\left\vert
L^{\prime}\right\vert =\left\vert 4H-%
{\textstyle\sum\nolimits_{j=1}^{6}}
E_{j}\right\vert $ on $S=\mathbb{\tilde{P}}^{2}(p_{1},...,p_{6}).$ Hence we
get $\deg(\Gamma)=(4H-%
{\textstyle\sum\nolimits_{j=1}^{6}}
E_{j})^{2}=10$ and $g(\Gamma)=\binom{4-1}{2}=3.$ It follows that
$\mathcal{L}:=\mathcal{O}_{\Gamma}(1)$ is very ample on $\Gamma$. From the
exact sequence }%
\[
0\rightarrow H^{0}(S^{\prime},\mathcal{O}_{S^{\prime}})\rightarrow
H^{0}(S^{\prime},\mathcal{O}_{S^{\prime}}(1))\rightarrow H^{0}(\Gamma
,\mathcal{O}_{\Gamma}(1))\rightarrow H^{1}(S^{\prime},\mathcal{O}_{S^{\prime}%
})=0
\]
($S^{\prime}$ is rational) it follows that $\Gamma$ is embedded projectively
normal by {\normalsize $\mathcal{L}$}. Then the two following lemmas show
that $\Gamma$ and therefore also $S^{\prime}$ is arithmetically Cohen-Macaulay.
\end{proof}

\bigskip

\begin{lemma}
Let $X\subset\mathbb{P}^{n}$ be a $d$-dimensional connected
variety, then $X$ is arithmetically Cohen Macaulay if and only if 

{\normalsize 1) $H^{1}(\mathbb{P}^{n},\mathcal{I}_{X}(m))=0$ for all
$m\in\mathbb{Z}$ and }

{\normalsize 2) $H^{i}(\mathbb{P}^{n},\mathcal{O}_{X}(m))=0$ for all
$i\neq0,d$ and $m\in\mathbb{Z}$ }
\end{lemma}

\begin{proof}
{\normalsize \cite{eisenbud1} page 472 Ex. 18.16 }
\end{proof}

{\normalsize \bigskip}

{\normalsize Then the following well known lemma can be obtained from this properties: }

\begin{lemma}
{\normalsize \label{aCM}Let $X^{d}\subset\mathbb{P}^{n},$ $d\geq2,$ be an
irreducible variety of dimension }$d,$ {\normalsize then $X$ is arithmetically
Cohen-Macaulay iff there exists a hyperplane section $\Gamma=$ $\mathbb{P}%
^{n-1}\cap X\subset\mathbb{P}^{n}$ that is arithmetically Cohen-Macaulay. }
\end{lemma}

\begin{proof}
{\normalsize $"\Leftarrow"$ Consider the exact sequence $(S1^{\ast})$%
\[
0\rightarrow\mathcal{I}_{X}(m-1)\rightarrow\mathcal{I}_{X}(m)\rightarrow
\mathcal{I}_{\Gamma}(m)\rightarrow0
\]
then $H^{1}(\Gamma,\mathcal{I}_{\Gamma}(m))=0$ for all $m\in\mathbb{Z}$ as
}$\Gamma$ {\normalsize is arithmetically Cohen-Macaulay. It follows that%
\[
h^{1}(X,\mathcal{I}_{X}(m-1))\geq h^{1}(X,\mathcal{I}_{X}(m))
\]
From the exact sequence $(S2^{\ast})$%
\[
0\rightarrow\mathcal{I}_{X}(m)\rightarrow\mathcal{O}_{\mathbb{P}^{n}%
}(m)\rightarrow\mathcal{O}_{X}(m)\rightarrow0
\]
we obtain the corresponding long exact sequence $(LS2^{\ast})$%
\begin{align*}
...  &  \rightarrow H^{0}(X,\mathcal{O}_{X}(m))\rightarrow H^{1}%
(X,\mathcal{I}_{X}(m))\rightarrow H^{1}(\mathbb{P}^{n},\mathcal{O}%
_{\mathbb{P}^{n}}(m))=0\rightarrow\\
&  \rightarrow H^{1}(X,\mathcal{O}_{X}(m))\rightarrow H^{2}(X,\mathcal{I}%
_{X}(m))\rightarrow...
\end{align*}
and from%
\[
0\rightarrow\mathcal{O}_{X}(-1)\rightarrow\mathcal{O}_{X}\rightarrow
\mathcal{O}_{\Gamma}\rightarrow0
\]
the long exact sequence%
\[
0\rightarrow H^{0}(X,\mathcal{O}_{X}(-1))\rightarrow H^{0}(X,\mathcal{O}%
_{X})\rightarrow H^{0}(\Gamma,\mathcal{O}_{\Gamma})
\]
As $X$ is irreducible and reduced we get $H^{0}(X,\mathcal{O}_{X}%
)\cong\mathbb{\Bbbk}$, and thus $H^{0}(X,\mathcal{O}_{X})\hookrightarrow
H^{0}(\Gamma,\mathcal{O}_{\Gamma})$ is the restriction of constant functions
on $X$ to $C.$ Further we have $H^{0}(X,\mathcal{O}_{X}(-1))=0$ and
$H^{0}(X,\mathcal{O}_{X}(m))=0$ for all negative values $m.$ Together with
$h^{1}(X,\mathcal{I}_{X}(m-1))\geq h^{1}(X,\mathcal{I}_{X}(m))$ we conclude
that $H^{1}(X,\mathcal{I}_{X}(m))=0$ for all $m.$ This shows 1)\ in the lemma
above. For 2) we look at the sequence $(LS2^{\ast})$ again: To determine the
group $H^{2}(X,\mathcal{I}_{X}(m))$ consider the long exact cohomology
sequence%
\[
...\rightarrow0=H^{1}(\Gamma,\mathcal{I}_{\Gamma}(m))\rightarrow
H^{2}(X,\mathcal{I}_{X}(m-1))\rightarrow H^{2}(X,\mathcal{I}_{X}%
(m))\rightarrow...
\]
and use Serre's vanishing theorem that says that $H^{2}(X,\mathcal{I}%
_{X}(m))=0$ for $m\gg0$ to see that $H^{2}(X,\mathcal{I}_{X}(m))=0$ for all
$m.$ Applying this result to $LS2^{\ast}$ leads to $H^{1}(X,\mathcal{O}%
_{X}(m))=0$ for all $m\in\mathbb{Z}.$ }

{\normalsize $"\Rightarrow"$ 2)\ is trivially fulfilled, so that it remains
to show that $H^{1}(\Gamma,\mathcal{I}_{\Gamma}(m))=0$ for all $m\in
\mathbb{Z}$. Consider the long exact sequence corresponding to $(S1^{\ast}):$%
\[
...\rightarrow H^{1}(X,\mathcal{I}_{X}(m))\rightarrow H^{1}(\Gamma
,\mathcal{I}_{\Gamma}(m))\rightarrow H^{2}(\Gamma,\mathcal{I}_{\Gamma
}(m-1)\rightarrow...
\]
The first group is $0$ as $X$ is arithmetically Cohen-Macaulay and $H^{2}(\Gamma,\mathcal{I}_{\Gamma
}(m-1)=0$, which follows from%
\[
0=H^{1}(X,\mathcal{O}_{X}(m-1))\rightarrow H^{2}(X,\mathcal{I}_{X}%
(m-1))\rightarrow H^{2}(\mathbb{P}^{n},\mathcal{O}_{\mathbb{P}^{n}}(m-1))=0
\]
thus $H^{1}(\Gamma,\mathcal{I}_{\Gamma}(m))=0.$ }
\end{proof}

\begin{remark}
\textnormal{For $C^{\prime}$ an irreducible plane curve of degree $7$ having exactly one
triple point $p_{0}$ and three doublepoints $p_{1},...,p_{3}$ as only
singularities, we consider the blowup in these singularities and the strict
transform $C\sim7H-3E_{0}-%
{\textstyle\sum\nolimits_{i=1}^{3}}
2E_{i}$ of $C^{\prime}.$ Then the adjoint linear series $\left\vert
L^{\prime}\right\vert =\left\vert 4H-2E_{0}-%
{\textstyle\sum\nolimits_{j=1}^{3}}
E_{j}\right\vert $ is base point free (see \cite{weinfurtner} page 35/36 with
the same methods as above) on $S=\tilde{P}^{2}(p_{0},...,p_{3})$ if there
exists no effective divisor $D\sim H-E_{0}-\sum_{i=1}^{3}E_{i}.$} 
\noindent \textnormal{But this follows easily from $D.C<0$ and the irreducibility of $C$. Therefore we obtain a
morphism}%
\[
\varphi:S\overset{\left\vert C+K_{S}\right\vert }{\rightarrow}S^{\prime
}\subset\mathbb{P}^{8}%
\]
\textnormal{from $\left\vert L^{\prime}\right\vert $ and $\varphi$ is even very ample
outside the union of all effective divisors $D$ of the following types:}
\\\\
$H-E_{0}-%
{\textstyle\sum\nolimits_{j\in\Delta}}
\alpha_{j}E_{j}$ \textnormal{with} $a_{j}=0$ \textnormal{or $1$ and} $|\Delta|=\#\{\alpha_{i}:\alpha_{i}=1\}=2.$

\noindent
$E_{l}-E_{l+1}$ \textnormal{for} $l\in\{1,2,3\}.$
\\\\
\noindent \textnormal{The existence of effective divisors $D$ of these types are given in
the situation where one triple point and two of the double points are lying on
a line or two of the double points are infinitely near. Outside the finite
union of these critical divisors $\varphi$ is very ample and thus the image
$S^{\prime}$ of $S$ under the morphism $\varphi$ is a surface with only
isolated singularities, such that we can find a smooth hyperplane section
$\Gamma$ applying Bertini's Theorem. $\Gamma$ is obtained from an element of
the linear system $\left\vert L^{\prime}\right\vert =\left\vert 4H-2E_{0}-%
{\textstyle\sum\nolimits_{j=1}^{3}}
E_{j}\right\vert $ on $S=\tilde{P}^{2}(p_{0},...,p_{3}).$ Thus we get
$\deg(\Gamma)=(4H-2E_{0}-%
{\textstyle\sum\nolimits_{j=1}^{3}}
E_{j})^{2}=9$ and $g(\Gamma)=\binom{4-1}{2}=3.$ It follows that $L:=\mathcal{O}_{\Gamma
}(1)$ is very ample on $\Gamma$ and $\Gamma$ is embedded by $L$ projectively
normal, hence $\Gamma$ is arithmetically Cohen-Macaulay. Then the surface $S^{\prime}$, which is called
a Castelnuovo surface is arithmetically Cohen-Macaulay, too.
\\Furthermore it is possible to generalize our above results, especially in the case where $C$ is an irreducible plane curve that has a singular point $p_{1}$ of multiplicity $m\geq3$ and further double points as only singularities. We remark that in the situation where $C$ has a quadruple point, it is possible to obtain an effective divisor of class $E_{1}-E_{2}-E_{3}$. This happens exactly in the case where two double points are lying infinitely near to the quadruple point. Then one has to check if this divisor is $i-critical$.}  
\end{remark}

{\normalsize \bigskip}

\section{Adjoint linear series on blowups of
$\mathbb{P}^{1}\times\mathbb{P}^{1}$}

{\normalsize In the section about pentagonal curves we consider canonical
curves $C$ of genus $9$ that admit two distinct linear systems of type
$g_{5}^{1}.$ From these linear systems we get a model $C^{\prime}%
\subset\mathbb{P}^{1}\times\mathbb{P}^{1}=:X$ of $C$, which is a divisor of
class $(5,5)$ on $\mathbb{P}^{1}\times\mathbb{P}^{1}.$ $C^{\prime}$ has
$s:=p_{a}(C^{\prime})-g(C)=4\cdot4-9=7$ double points }$p_{1},...,p_{7}%
$ {\normalsize as only singularities. Blowing up $X$ in these singular points
we can assume that $C$ is the strict transform of $C^{\prime}$ on the blowup%
\[
\sigma:S=\tilde{X}(p_{1},...,p_{7})\overset{\sigma_{7}}{\rightarrow
}...\overset{\sigma_{2}}{\rightarrow}\tilde{X}(p_{1})\overset{\sigma_{1}%
}{\rightarrow}\mathbb{P}^{1}\times\mathbb{P}^{1}%
\]
Then $C\sim(5,5)-%
{\textstyle\sum\nolimits_{i=1}^{7}}
2E_{i}$ and its canonical series is cut out by $\left\vert L^{\prime
}\right\vert =\left\vert K_{S}+C\right\vert =\left\vert (3,3)-%
{\textstyle\sum\nolimits_{i=1}^{7}}
E_{i}\right\vert .$ As in the previous section we want to apply Lemma
\ref{ampleness} again to show that $\left\vert L^{\prime}\right\vert $ is
}${\normalsize i-}$very ample:

\bigskip

\begin{theorem}
{\normalsize \label{ALSP1P1}Let $S$ be the iterated blowup of 7 points
$p_{1},...,p_{7}$ on $\mathbb{P}^{1}\times\mathbb{P}^{1}$ and $C\sim(5,5)-%
{\textstyle\sum\nolimits_{i=1}^{7}}
2E_{i}$ an irreducible, nonsingular curve on $S$, then the adjoint linear
series $\left\vert C+K_{S}\right\vert =\left\vert (3,3)-%
{\textstyle\sum\nolimits_{i=1}^{7}}
E_{i}\right\vert $ is base point free. It is very
ample on $S$ if and only if none of the points $p_{i}$ lies infinitely near to another
one. Furthermore the image $S^{\prime}\subset\mathbb{P}^{8}$ of $S$ under the morphism
defined by $\left\vert C+K_{S}\right\vert $ is arithmetically Cohen-Macaulay. If there exist two infinitely near double points then $S^{\prime}$ has isolated singularities that are contractions of strict transforms }%
$E_{i}^{\prime}\sim E_{i}-E_{i+1}$.
\end{theorem}

\begin{proof}
For $L:=L^{\prime}-K_{S}\sim C$ we get $L^{2}=${\normalsize $50-28>5+4i$ for
$i=0,1$ and from Theorem \ref{ampleness} it follows the $i-$very ampleness of
$\left\vert L^{\prime}\right\vert $ iff there exists no $i$-critical effective
divisor $D$, such that $L-2D$ is $\mathbb{Q-}$effective and $D.(C-D)<2+i.$ As
above for $X=\mathbb{P}^{2}$, we have to distinguish three types for $D:$ }
\\\\
{\normalsize \textit{(a)} $D\sim E_{k}$ for a $k\in\{1,...,7\}.$}

\noindent
{\normalsize \textit{(b)} $D\sim E_{k}-E_{k+1}$ for $k\in\{1,...,6\}.$}

\noindent
{\normalsize \textit{(c)} $D\sim(a,b)-%
{\textstyle\sum\nolimits_{i=1}^{7}}
\alpha_{i}E_{i}$ with $\alpha_{i}=0$ or $1$ and $(a,b)\in \operatorname{Pic}(\mathbb{P}%
^{1}\times\mathbb{P}^{1})^{\ast}$ }
\\\\
{\normalsize \noindent A divisor of type (a) cannot be $i-$critical as
$E_{k}.(C-E_{k})=3\geq2+i.$ Further there exists an effective divisor of type
(b) exactly if two of the points $p_{1},...,p_{7}$ are infinitely near. In
this situation every such divisor, which can be written as difference of two
total transforms $E_{k}$ and $E_{k+1}$ is $1$-critical but not $0-$critical
because of $(E_{k}-E_{k+1}).(C-(E_{k}-E_{k+1}))=2.$ Now let $D\sim(a,b)-%
{\textstyle\sum\nolimits_{i=1}^{7}}
\alpha_{i}E_{i}$ be a $1-$critical divisor, then from the condition that
$C-2D$ is $\mathbb{Q}$-effective we can restrict to $a,b\leq2. $ As above we
denote $\delta=\#\{\alpha_{i}:\alpha_{i}=1\},$ then the condition
$D.(C-2D)\leq2$ transforms into:%
\begin{align*}
2  &  \geq((a,b)-%
{\textstyle\sum\nolimits_{i=1}^{7}}
\alpha_{i}E_{i}).((5-a,5-b)-%
{\textstyle\sum\nolimits_{i=1}^{7}}
(2-\alpha_{i})E_{i})=\\
&  =a(5-b)+b(5-a)-%
{\textstyle\sum\nolimits_{i=1}^{7}}
\alpha_{i}(2-\alpha_{i})=5a+5b-2ab-\delta=\\
&  =D.C+\delta-2ab\geq0+\delta-2ab\\
&  \Rightarrow\delta\leq2ab+2
\end{align*}
We remark that the inequality $D.C\geq0$ is a consequence of the irreducibility of
$C.$ Therefore we get
\[
D.(C-D)=5a+5b-2ab-\delta\geq5a+5b-4ab-2
\]
In the case $a=0$ or $b=0$ we obtain $D.(C-D)\geq5-2=3.$ For $a=1$ we have
$D.(C-D)\geq5+b-2\geq3$ and the same for $b=1.$ If $a=2$ then $D.(C-D)\geq
10-3b-2\geq5$ for $b=0,1.$ In the case $a=b=2$ we use the inequality
$D.(C-D)\geq5a+5b-2ab-\delta\geq12-\delta\geq5.$ The same holds for $b=2.$ In
summary we have seen that there exists no $0$-critical divisor and the only
$1-$critical are exactly the divisors of type (b). Now with the same arguments
as in the "$\mathbb{P}^{2}$ case" it follows that $\left\vert L^{\prime
}\right\vert $ is base point free and even very ample outside the finite union
of divisors of type (b). Further the image $S^{\prime}\subset\mathbb{P}^{8}$
of $S$ under the morphism $\varphi$ defined by $\left\vert L^{\prime
}\right\vert $ has only isolated singularities. Due to Bertini's Theorem we
can find a smooth hyperplane section $\Gamma\in\left\vert (3,3)-%
{\textstyle\sum\nolimits_{i=1}^{7}}
E_{i}\right\vert $ with $\deg(\Gamma)=((3,3)-%
{\textstyle\sum\nolimits_{j=1}^{7}}
E_{i})^{2}=11$ and $g(\Gamma)=2\cdot2=4.$ It follows that $\mathcal{L}%
:=\mathcal{O}_{\Gamma}(1)$ is very ample on $\Gamma$. }Then the same arguments
as above {\normalsize show that $\Gamma$ is embedded by $\mathcal{L}$
projectively normal, hence $\Gamma$ and therefore also $S^{\prime}$ is arithmetically Cohen-Macaulay. }
\end{proof}

\bigskip

\section{Adjoint linear series on blowups of $P_{2}$}

{\normalsize Later on we will see that from a $g_{5}^{1}$ of higher
multiplicity (for the definition we refer to Chapter 4) on a canonical
curve $C$ of genus $9,$ we get a model $C^{\prime}\subset P_{2}=\mathbb{P}%
(\mathcal{O(}2\mathcal{)\oplus O})=:X$ of $C$, that is a divisor of class $5H$
on $P_{2}.$ $C^{\prime}$ has $p_{a}(C^{\prime})-g(C)=4\cdot4-9=7$ double
points }$p_{1},...,p_{7}$ {\normalsize as only singularities. In the case of our interest $C^{\prime}$ has no intersection with the exceptional divisor $E\sim H-2R$ on $P_{2}$.  Blowing up $X$
in the singular points we can assume that $C$ is the strict transform of
$C^{\prime}$ in the blowup%
\[
\sigma:S=\tilde{X}(p_{1},...,p_{7})\overset{\sigma_{7}}{\rightarrow
}...\overset{\sigma_{2}}{\rightarrow}\tilde{X}(p_{1})\overset{\sigma_{1}%
}{\rightarrow}P_{2}%
\]
Then $C\sim5H-%
{\textstyle\sum\nolimits_{i=1}^{7}}
2E_{i}$ \ and its canonical series is cut out by $\left\vert L^{\prime
}\right\vert =\left\vert K_{S}+C\right\vert =\left\vert 3A-%
{\textstyle\sum\nolimits_{i=1}^{7}}
E_{i}\right\vert .$ The following theorem can be obtained in analogous manner
as the main theorems in the last two sections:}

\bigskip

\begin{theorem}
{\normalsize \label{ALSF2}Let $S$ be the iterated blowup of 7 points
$p_{1},...,p_{7}$ on $P_{2}$ and $C\sim5H-%
{\textstyle\sum\nolimits_{i=1}^{7}}
2E_{i}$ an irreducible, nonsingular curve on $S$, then the adjoint linear
series $\left\vert C+K_{S}\right\vert =\left\vert 3H-%
{\textstyle\sum\nolimits_{i=1}^{7}}
E_{i}\right\vert $ is base point free. The image $S^{\prime}\subset\mathbb{P}^{8}$ of $S$ under the morphism defined
by $\left\vert C+K_{S}\right\vert $ is arithmetically Cohen-Macaulay and the
only singularities are contractions of the exceptional divisor $E\sim H-2R$ and strict transforms }$E_{i}^{\prime}\sim
E_{i}-E_{i+1}$ in the case where $p_{i+1}$ lies infinitely near $p_{i}.$
\end{theorem}

\begin{proof}
With $L:=L^{\prime}-K_{S}\sim C$ we have {\normalsize $L^{2}=50-28>5+4i$ for
$i=0,1.$ From Theorem \ref{ampleness} it follows the $i-$very ampleness of
$\left\vert L^{\prime}\right\vert $ iff there exists no $i$-critical effective
divisor $D$, such that $L-2D$ is $\mathbb{Q-}$effective and $D.(C-D)<2+i.$ As
above for $X=\mathbb{P}^{2}$ or $\mathbb{P}^{1}\times\mathbb{P}^{1}$, we have
to distinguish three types for $D:$}
\\\\
{\normalsize \textit{(a)} $D\sim E_{k}$ for a $k\in\{1,...,7\}.$}

\noindent
{\normalsize \textit{(b)} $D\sim E_{k}-E_{k+1}$ for $k\in\{1,...,6\}$}

\noindent
{\normalsize \textit{(c)} $D\sim aH+bR-%
{\textstyle\sum\nolimits_{i=1}^{7}}
\alpha_{i}E_{i}$ with $\alpha_{i}=0$ or $1$ and $aH+bR\in \operatorname{Pic}(F_{2})^{\ast},$
thus $a\geq0$ and $b\geq-2a.$ }
\\\\
{\normalsize \noindent In complete analogue to our arguments above, a divisor
of type (a) cannot be $i-$critical and an effective divisor of type (b) exists
if and only if two of the points $p_{1},...,p_{7}$ are infinitely near. In this
situation $E_{k}-E_{k+1}$ is $1$-critical but not $0-$critical$.$ Now let
$D\sim aH+bR-%
{\textstyle\sum\nolimits_{i=1}^{7}}
\alpha_{i}E_{i}$ be a $1-$critical divisor, then from the condition that
$C-2D$ is $\mathbb{Q}$-effective we can restrict to $a\leq2$ and
$b\leq\left\lfloor \frac{5-a}{2}\right\rfloor .$ Let $\delta=\#\{\alpha
_{i}:\alpha_{i}=1\},$ then from $D.(C-D)\leq2$ we conclude:%
\begin{align*}
2  &  \geq(aH+bR-%
{\textstyle\sum\nolimits_{i=1}^{7}}
\alpha_{i}E_{i}).((5-a)H-bR-%
{\textstyle\sum\nolimits_{i=1}^{7}}
(2-\alpha_{i})E_{i})=\\
&  =2a(5-a)+b(5-a)-ab-%
{\textstyle\sum\nolimits_{i=1}^{7}}
\alpha_{i}(2-\alpha_{i})=\\
&  =10a+5b-\delta-2a(a+b)=D.C+\delta-2a(a+b)\geq\delta-2a(a+b)\\
&  \Rightarrow\delta\leq2a(a+b)+2
\end{align*}
Therefore it follows that
\[
D.(C-D)=10a+5b-2a(a+b)-\delta\geq10a+5b-4a(a+b)-2
\]
In the case $a=0$ we obtain $D.(C-D)\geq5b-2\geq3$ and for $a=1$ we have
$D.(C-D)\geq4+b\geq2$ and $D.(C-D)=2\Leftrightarrow D\sim H-2R.$ If $a=2$ then $D.(C-D)\geq2-3b\geq5$ for $b<0.$ It
remains to discuss the cases $a=2$ and $b=0$ or $1:$ Here we use the
inequality $D.(C-D)=10a+5b-2a(a+b)-\delta\geq12+b-\delta\geq5.$ In summary we
have seen that there exists no $0$-critical divisor and the only $1-$critical
are exactly the divisors of type (b) and the exceptional divisor $E\sim H-2R$ on $P_{2}$. From the same arguments as in the
"$\mathbb{P}^{2}$ case" it follows that $\left\vert L^{\prime}\right\vert $ is
base point free and even very ample outside the finite union of divisors of
type (b) and the exceptional divisor $E$. Further the image $S^{\prime}\subset\mathbb{P}^{8}$ of $S$ under the
morphism $\varphi$ defined by $\left\vert L^{\prime}\right\vert $ has only
isolated singularities. Due to Bertini's Theorem we can find a smooth
hyperplane section $\Gamma\in\left\vert 3H-%
{\textstyle\sum\nolimits_{i=1}^{7}}
E_{i}\right\vert $ with $\deg(\Gamma)=(3H-%
{\textstyle\sum\nolimits_{i=1}^{7}}
E_{i})^{2}=11$ and $g(\Gamma)=2\cdot2=4.$ It follows that $\mathcal{L}%
:=\mathcal{O}_{\Gamma}(1)$ is very ample on $\Gamma$ and that $\Gamma$ is
embedded by $\mathcal{L}$ projectively normal, hence $\Gamma$ and therefore
also $S^{\prime}$ is arithmetically Cohen-Macaulay.}
\end{proof}

{\normalsize \bigskip}

{\normalsize In the situation where $C^{\prime}$ has infinitely near double
points we will show in Section \ref{deformation}\textbf{ }that it is possible
to separate these points, i.e. there exists a one paramter family of curves
$C_{\lambda}^{\prime}$ having only ordinary nodes. For this purpose we will
provide a further result here. We recall that $E_{i}^{\prime}$ denotes the
strict transform of the point $p_{i}$ in the iterated blowup of all singular
points$.$ Now let $p_{7}$ be infinitely near to $p_{6}$ and }%
\[
\sigma^{\prime}:S^{(6)}=\tilde{X}(p_{1},...,p_{6})\overset{\sigma_{6}%
}{\rightarrow}...\overset{\sigma_{2}}{\rightarrow}\tilde{X}(p_{1}%
)\overset{\sigma_{1}}{\rightarrow}P_{2}%
\]
{\normalsize the iterated blowup in the points $p_{1},...,p_{6}.$\ The
following theorem then states that the linear system $\left\vert L^{\prime
}\right\vert =\left\vert 5H-%
{\textstyle\sum\nolimits_{i=1}^{6}}
2E_{i}\right\vert $ is very ample on $S^{(6)}\backslash$}$%
({{\textstyle\bigcup\nolimits_{i=1}^{5}}
E_{i}^{\prime}\cup E})${\normalsize . With $L:=L^{\prime}-K_{S}\sim7H-%
{\textstyle\sum\nolimits_{i=1}^{6}}
3E_{i}$ we get $L^{2}=98-54-1>5$ and therefore the first condition to apply
Lemma \ref{ampleness} to $\left\vert L^{\prime}\right\vert $ is fulfilled.
}

\bigskip

\begin{theorem}
{\normalsize Let }$S$ be a surface as in Theorem \ref{ALSF2} above. {\normalsize Then the complete linear
series $\left\vert L^{\prime}\right\vert :=\left\vert 5H-%
{\textstyle\sum\nolimits_{i=1}^{6}}
2E_{i}\right\vert $ is very ample on $S\backslash%
({\textstyle\bigcup\nolimits_{i=1}^{5}}
E_{i}^{\prime}\cup E).$ }
\end{theorem}
\begin{proof}
{\normalsize The base point freeness of $\left\vert L^{\prime}\right\vert $ on
}$S^{(6)}$ {\normalsize follows if there exists no $0$-critical effective
divisor $D$. As
above we have to distinguish three types for $D:$ }
\\\\
{\normalsize \textit{(a)} $D\sim E_{k}$ for a $k\in\{1,...,6\}.$}

\noindent
{\normalsize \textit{(b)} $D\sim E_{k}-E_{k+1}$ for $k\in\{1,...,5\}$}

\noindent
{\normalsize \textit{(c)} $D\sim aH+bR-%
{\textstyle\sum\nolimits_{i=1}^{6}}
\alpha_{i}E_{i}$ with $\alpha_{i}=0$ or $1$ and $aH+bR\in \operatorname{Pic}(P_{2})^{\ast},$
thus $a\geq0$ and $b\geq-2a.$}
\\\\
{\normalsize \noindent It is easy to check that for every divisor $D$ of type (a) or (b) we
get $D.(L-D)\geq2$. Assume that $D\sim aH+bR-%
{\textstyle\sum\nolimits_{i=1}^{7}}
\alpha_{i}E_{i}$ is a $0-$critical divisor, then from the condition that
$L-2D$ is $\mathbb{Q}$-effective we can restrict to $a\leq3$ and
$b\leq\left\lfloor \frac{7-a}{2}\right\rfloor \leq3.$ On }$S^{(6)}$ the curve
$C$ is a divisor of class $5H-%
{\textstyle\sum\nolimits_{i=1}^{6}}
2E_{i}$ that has exactly one double point as only singularity. {\normalsize Let
$\delta=\#\{\alpha_{i}:\alpha_{i}=1\},$ then we conclude:%
\begin{align*}
D.(L-D) &  =2a(7-a)+b(7-a)-ab-%
{\textstyle\sum\nolimits_{i=1}^{6}}
\alpha_{i}(3-\alpha_{i})=\\
&  =14a+7b-2\delta-2a(a+b)=D.C+4a+2b-2a(a+b)\geq\\
&  \geq2(2a+b)-2a(a+b)
\end{align*}
For $a=0$ we get $D.(L-D)\geq2b\geq2$ and for $a=1:D.(L-D)\geq2.$ In the cases
$a=2,3$ we must have $-6\leq-2a\leq b\leq2.$ Thus if $a=2$ and $b<0$ then
$D.(L-2D)\geq-2b\geq2$ and for $a=2$, $b\geq0$ we get
$D.(L-2D)=D.C-2b=10a+3b-2\delta\geq20-12=8.$ In the situation $a=3,b\leq-2$ we
obtain $D.(L-2D)\geq-6-4b\geq2$ and for $a=3,b\geq-1$ this inequality also
holds because of $D.(L-2D)=10a+3b-2\delta\geq30-3-2\delta\geq15.$ In summary
we have seen that there exists no $0$-critical divisor, hence $\left\vert
L^{\prime}\right\vert $ is base point free on $S^{(6)}$.}

Now let $p\in${\normalsize $S^{(6)}\backslash$}$%
({\textstyle\bigcup\nolimits_{i=1}^{5}}
E_{i}^{\prime}\cup E)$ be an arbritrary point. The very ampleness of
{\normalsize $\left\vert L^{\prime}\right\vert $ on $S^{(6)}\backslash%
({\textstyle\bigcup\nolimits_{i=1}^{5}}
E_{i}^{\prime}\cup E)$ follows if $\left\vert L^{\prime}-E_{p}\right\vert =\left\vert
5H-%
{\textstyle\sum\nolimits_{i=1}^{6}}
2E_{i}-E_{p}\right\vert $ is base point free on the blowup $S_{p}:=\tilde
{S}^{(6)}(p)$} with exceptional divisor $E_{p}${\normalsize . With
$L_{p}:=L^{\prime}-K_{S_{p}}\sim7H-%
{\textstyle\sum\nolimits_{i=1}^{6}}
3E_{i}-2E_{p}$ we still have $L_{p}^{2}=98-54-8>5$, such that we can use
Lemma \ref{ampleness} again. A $0$-critical divisors $D$} {\normalsize for
$\left\vert L^{\prime}-E_{p}\right\vert $ has to satisfy the following
inequality:%
\[
D.(L^{\prime}-D-E_{p})\leq1\Rightarrow D.(L^{\prime}-D)-E_{p}.D\leq1
\]
Because of $E_{p}.D\leq1$ and $D.(L^{\prime}-D)\geq2$ this is only possible if
$E_{p}.D=1$ and $D.(L^{\prime}-D)=2.$ Then $p\in D$ with $D$ a strict transform 
$E_{k}^{\prime}=E_{k}-E_{k+1}$ or $D$ the exceptional divisor $E$, which was excluded.}
\end{proof}

{\normalsize \newpage}

\chapter{\label{cliff<=2}{\protect\normalsize Curves C with Clifford index Cliff(C)$\leq$2}}

\section{Trigonal Curves}

If $C$ is a canonical curve of genus $9$ and $\operatorname*{Cliff}(C)=1$, then
it follows the existence of a $g_{2r+1}^{r}$ with $r\in\{1,...,6\}$. As the
Brill Noether dual of a $g_{2r+1}^{r}$ is of type $g_{15-2r}^{7-r}$,
we can restrict to the existence of a $g_{3}^{1},$ $g_{5}^{2}$ or $g_{7}^{3}$.
The existence of a $g_{5}^{2}$ is not possible as a plane curve of degree $5$
has geometric genus less or equal than $\binom{4}{2}=6<9.$ Further Castelnuovo inequality
gives a boundary for the genus of space curves of degree $d\geq3$ (cf. Theorem
IV 6.4 in \cite{hartshorne}): According to this theorem we must have
$g=9\leq\frac{1}{4}(d^{2}-1)-d+1$ for odd degrees $d.$ As $d=7$ does not
fullfill this inequality we can omit the possibility of the existence of a
$g_{7}^{3}$. Hence, it follows: 
\begin{corollary}
{\normalsize Let $C\subset\mathbb{P}^{8}$ be a canonical curve of genus 9} that has Clifford index {\normalsize $\operatorname*{Cliff}(C)=1$, then} $C$ is trigonal.
\end{corollary}
{\normalsize We will repeat the main
results as can be found in \cite{schreyer1} Section 6.1.: A trigonal canonical
curve of genus $g$ is contained in a $2-$dimensional scroll $X\subset
\mathbb{P}^{g-1}$:%
\[
X=%
{\textstyle\bigcup\nolimits_{D\in g_{3}^{1}}}
\bar{D}\subset\mathbb{P}^{g-1}%
\]
of type $S(e_{1},e_{2})$ and degree $f=e_{1}+e_{2}=g-2$ (see Section
\ref{ScrollVariety}). Because of $\deg(K_{C}-nD)=2g-2-3n<0$ for $n>\frac
{2g-2}{3},$ for the values $e_{i}$ we get the following bounds:%
\[
\frac{2g-2}{3}\geq e_{1}\geq e_{2}\geq\frac{g-4}{3}%
\]
Further $C$ is a divisor of class $3H-(f-2)R$ on $X$ (cf. Theorem
\ref{ResC}). The mapping cone%
\[
\mathcal{C}^{f-2}(-3)\rightarrow\mathcal{C}^{0}%
\]
is a minimal resolution of $\mathcal{O}_{C}$ as an $\mathcal{O}_{\mathbb{P}%
^{g-1}}$-module and therefore the scroll $X$ and hence $g_{3}^{1}$ is uniquely
determined by $C$ for a trigonal curve. For $g(C)=9$ we get the minimal free
resolution of $\mathcal{O}_{C}$ from the following mapping cone construction
with $F=\mathcal{O}_{\mathbb{P(\mathcal{E})}}(H-R)\cong\mathcal{O}%
_{\mathbb{P}^{8}}^{f}$ and $G=\mathcal{O}_{\mathbb{P(\mathcal{E})}}%
(R)\cong\mathcal{O}_{\mathbb{P}^{8}}^{2}$ (cf. Section \ref{Scrollgeneral}): }

{\normalsize
\[%
\begin{array}
[c]{ccccccccc}%
0 & \rightarrow & \mathcal{O}_{\mathbb{P(\mathcal{E})}}(-3H+5R) & \rightarrow
& \mathcal{O}_{\mathbb{P(\mathcal{E})}} & \rightarrow & \mathcal{O}_{C} &
\rightarrow & 0\\
&  & \uparrow &  & \uparrow &  &  &  & \\
&  & S_{5}G\mathcal{(-}3) &  & \mathcal{O}_{\mathbb{P}^{8}} &  &  &  & \\
&  & \uparrow &  & \uparrow &  &  &  & \\
&  & F\otimes S_{4}G(-4) &  & \wedge^{2}F(-2) &  &  &  & \\
&  & \uparrow &  & \uparrow &  &  &  & \\
&  & \wedge^{2}F\otimes S_{3}G(-5) &  & \wedge^{3}F\otimes D_{1}G^{\ast}(-3) &
&  &  & \\
&  & \uparrow &  & \uparrow &  &  &  & \\
&  & \wedge^{3}F\otimes S_{2}G(-6) &  & \wedge^{4}F\otimes D_{2}G^{\ast}(-4) &
&  &  & \\
&  & \uparrow &  & \uparrow &  &  &  & \\
&  & \wedge^{4}F\otimes G(-7) &  & \wedge^{5}F\otimes D_{3}G^{\ast}(-5) &  &
&  & \\
&  & \uparrow &  & \uparrow &  &  &  & \\
&  & \wedge^{5}F(-8) &  & \wedge^{6}F\otimes D_{4}G^{\ast}(-6) &  &  &  & \\
&  & \uparrow &  & \uparrow &  &  &  & \\
&  & \wedge^{7}F(-10) &  & \wedge^{7}F\otimes D_{5}G^{\ast}(-7) &  &  &  & \\
&  & \uparrow &  & \uparrow &  &  &  & \\
&  & 0 &  & 0 &  &  &  &
\end{array}
\]
There exists no non-minimal map, thus it follows }

\begin{theorem}
{\normalsize \label{cliff1}Let $C\subset\mathbb{P}^{8}$ be an irreducible,
canonical and non hyperelliptic curve of genus 9 that admits a $g_{3}^{1}.$
Then the Betti table for }$C$ {\normalsize  is given as follows:
\textnormal{\[
\begin{tabular}{c|cccccccccc}
& 0 & 1 & 2 & 3 & 4 & 5 & 6 & 7 \cr \hline
0 & 1 & - & - & - & - & - & - & - \cr
1 & - & 21 & 70 & 105 & 84 & 35 & 6 & - \cr
2 & - & 6 & 35  & 84  & 105  & 70  & 21 & -\cr
3 & - & - & - & - & - & - & - & 1 \cr
\end{tabular} 
\]}}
\end{theorem}

\bigskip
\section{Tetragonal Curves}

For a canonical curve $C$ of genus 9, from the property of having Clifford index $2$ it follows the existence of a $g_{4}^{1}$. This can be seen as follows: From the
existence of a $g_{6}^{2}$ we get a plane model of $C$ that has exactly
$\binom{5}{2}-9=1$ doublepoint as only singularity. Projection from this point leads to a $g_{4}^{1}.$ {\normalsize In the case where $C$ admits a
$g_{8}^{3},$ we get a space model $C^{\prime}$ for $C$ of degree $l=8.$ If
$C^{\prime}$ has a singular point, projection from this point leads to a
$g_{d}^{2}$ with $d\leq6$. Otherwise according to a theorem of Castelnuovo
(cf. \cite{hartshorne} IV Theorem 6.4), it follows from $g(C^{\prime}%
)=9=\frac{1}{4}l^{2}-l+1$ that $C^{\prime}$ is contained in a quadric surface
$Q\subset\mathbb{P}^{3}.$ In the situation where $Q\cong\mathbb{P}^{1}%
\times\mathbb{P}^{1}$ is nonsingular, $C^{\prime}$ must be a divisor of type
$(4,4)$ on $\mathbb{P}^{1}\times\mathbb{P}^{1}.$ Projection along each
factor of $\mathbb{P}^{1}\times\mathbb{P}^{1}$ leads to two different
$g_{4}^{1}$$^{\prime}s$ on $C$. For $Q$ a singular quadric cone over an elliptic curve
$E$, $C^{\prime}$ is of class $C^{\prime}\sim4A$, $A$ the class of a
hyperplane section on $Q.$ The rulings $\left\vert R\right\vert $ on $Q$
cut out a $g_{4}^{1}$ on $C$. It remains to remark that for a }$g_{10}^{4}$
or a $g_{12}^{5}$ the Brill Noether dual is of type $g_{6}^{2}$
or $g_{4}^{1}$ respectively.
\begin{corollary}
{\normalsize Let $C\subset\mathbb{P}^{8}$ be a canonical curve of genus 9 with} Clifford index {\normalsize $\operatorname*{Cliff}(C)=2$}, then $C$ is tetragonal.
\end{corollary}
Now applying {\normalsize the main ideas, which can be found in
\cite{schreyer1} Section 6.2-6.6, to a canonical, tetragonal curve $C$ of
genus $g=9$ leads to the following results: Constructing $X$ as above from the
linear system of type $g_{4}^{1}$ on $C$, we see that $C$ is contained in a
$3-$dimensional rational normal scroll of type $S(e_{1},e_{2},e_{3})$ with%
\[
4=\frac{2g-2}{4}\geq e_{1}\geq e_{2}\geq e_{3}\geq0
\]
and degree $f=e_{1}+e_{2}+e_{3}=g-3=6.$ We denote by $\pi
:\mathbb{P(\mathcal{E})\rightarrow P}^{1}$ the corresponding $\mathbb{P}^{2}%
-$bundle. According to Theorem \ref{ResC} $C$ is given as complete
intersection of two divisors%
\[
Y\sim2H-b_{1}R\text{ \ , \ }Z\sim2H-b_{2}R
\]
on $X$ with%
\[
b_{1}+b_{2}=f-2=4
\]
and $b_{1}\geq b_{2}.$ In \cite{schreyer1} Section 6.3. the author verifies
that
\[
5=f-1\geq b_{1}\geq b_{2}\geq-1
\]
and even $b_{1}\leq f-2=4$ for genus $g\neq0.$ Therefore we can apply Theorem
\ref{ResC} to get a minimal free resolution of $\mathcal{O}_{C}$ as an
$\mathcal{O}_{\mathbb{P}^{g-1}}$-module via an iterated mapping cone:%
\[
\lbrack\mathcal{C}^{f-2}(-4)\rightarrow\mathcal{C}^{b_{1}}(-2)\oplus
\mathcal{C}^{b_{2}}(-2)]\rightarrow\mathcal{C}^{0}%
\]
The Betti table for the minimal free resolution is determined by the values
for $b_{1}$ and $b_{2}$ and vice versa$.$ $X$ and hence the $g_{4}^{1}$ is
uniquely determined by $C$ unless $b_{1}\geq f-2=4:$ Consider the map%
\[
\mathcal{O}^{\beta_{g-4,g-3}}(-g+3)\rightarrow\mathcal{O}^{\beta_{g-5,g-4}%
}(-g+4)
\]
in the resolution of $C$, then $X$ is given as support of the cokernel of its
dual. }

{\normalsize \bigskip}

\begin{theorem}
{\normalsize \label{cliff2}Let $C\subset\mathbb{P}^{8}$ be an irreducible,
nonsingular, canonical curve of genus 9 with $\operatorname*{Cliff}(C)=2$,
then }

\noindent{\normalsize a) $C$ admits a $g_{6}^{2}$ or a $g_{8}^{3}$ exactly if
$(b_{1},b_{2})=(4,0)$. The Betti table for $C$ then takes the following form%
\textnormal{\[
\begin{tabular}{c|cccccccccc}
& 0 & 1 & 2 & 3 & 4 & 5 & 6 & 7 \cr \hline
0 & 1 & - & - & - & - & - & - & - \cr
1 & - & 21 & 64 & 90 & 64 & 20 & - & - \cr
2 & - & - & 20  & 64  & 90  & 64  & 21 & -\cr
3 & - & - & - & - & - & - & - & 1 \cr
\end{tabular} 
\]}
}

\noindent{\normalsize b) $C$ admits a $g_{4}^{1}$, a base point free linear system of Clifford index 3 and no $g_{6}^{2}$ or $g_{8}^{3}$ exactly if $(b_{1},b_{2})=(3,1)$ with Betti table for $C$ as
follows%
\textnormal{\[
\begin{tabular}{c|cccccccccc}
& 0 & 1 & 2 & 3 & 4 & 5 & 6 & 7 \cr \hline
0 & 1 & - & - & - & - & - & - & - \cr
1 & - & 21 & 64 & 75 & 44 & 5 & - & - \cr
2 & - & - & 5 & 44  & 75  & 64  & 21 & -\cr
3 & - & - & - & - & - & - & - & 1 \cr
\end{tabular} 
\]}
}

\noindent{\normalsize c) $C$ admits a $g_{4}^{1}$ and no $g_{6}^{2}$, $g_{8}^{3}$ and
no base point free linear system of Clifford index 3 iff $(b_{1},b_{2})=(2,2)$. The Betti
table for $C$ is then given by%
\textnormal{\[
\begin{tabular}{c|cccccccccc}
& 0 & 1 & 2 & 3 & 4 & 5 & 6 & 7 \cr \hline
0 & 1 & - & - & - & - & - & - & - \cr
1 & - & 21 & 64 & 75 & 24 & 5 & - & - \cr
2 & - & - & 5  & 24  & 75  & 64  & 21 & -\cr
3 & - & - & - & - & - & - & - & 1 \cr
\end{tabular} 
\]}
}
\end{theorem}

\begin{proof}
{\normalsize From the conditions for $(b_{1},b_{2})$ we have to distinguish
three different cases for $g=9:$ }

{\normalsize \noindent\textbf{1. }$(b_{1},b_{2})=(4,0):$ From the mapping cone
construction above we get the following Betti table for $C\subset
\mathbb{P}^{8}:$%
\textnormal{\[
\begin{tabular}{c|cccccccccc}
& 0 & 1 & 2 & 3 & 4 & 5 & 6 & 7 \cr \hline
0 & 1 & - & - & - & - & - & - & - \cr
1 & - & 21 & 64 & 90 & 64 & 20 & - & - \cr
2 & - & - & 20  & 64  & 90  & 64  & 21 & -\cr
3 & - & - & - & - & - & - & - & 1 \cr
\end{tabular} 
\]}
\noindent The fibres of $Y^{\prime}\subset\mathbb{P(\mathcal{E})}$ over $\mathbb{P}^{1}$
are conics. If all these fibres are degenerate the $g_{4}^{1}$ is composed by
an elliptic or hyperelliptic involution%
\[
C\overset{2:1}{\rightarrow}E\overset{2:1}{\rightarrow}\mathbb{P}^{1}%
\]
and $Y$ is a birational ruled surface over $E$ with a rational curve
$\tilde{E}$ of double points. As $C$ is assumed to be nonsingular we must have
$\tilde{E}\cap Z=\varnothing$ and therefore
\[
0=\tilde{E}.Z=\tilde{E}.2H=2\deg\tilde{E}\Rightarrow\deg\tilde{E}=0
\]
A general divisor $\Gamma$ of class $H-\frac{b_{2}}{2}R$ intersects $Y$ in a
smooth curve isomorphic to $E,$ so the geometric genus is given by%
\[
2p_{a}E-2=\Gamma.Y.(f-2-b_{1}-\frac{b_{2}}{2})R=b_{2}=0
\]
(cf. \cite{schreyer1} Example 3.6). Thus $E$ is an elliptic curve that can be
embedded as a plane cubic. Then the composition $C\overset{2:1}{\rightarrow
}E\hookrightarrow\mathbb{P}^{2}$ gives a $g_{6}^{2}$ on $C.$ }

{\normalsize \noindent In the situation where a general fibre of $Y$ is a
nonsingular conic, the number of singular fibres is given by
\[
\delta=2f-3b_{1}=0
\]
hence there exists no singular fibre. Further $Y$ admits a determinantal presentation. It is obtained from the
determinant of a matrix $\psi$ with entries on $X$ as indicated below%
\[
\left(
\begin{array}
[c]{cc}%
H-a_{1}R & H-a_{2}R\\
H-(a_{1}+k)R & H-(a_{2}+k)R
\end{array}
\right)
\]
with $a_{1},a_{2}\in\mathbb{Z}$, $k\in\mathbb{N}$ and $a_{1}+a_{2}+k=b_{1}=4$.
Then $Y$ may be identified with the image of $P_{k}=\mathbb{P}(\mathcal{O}_{\mathbb{P}^{1}}%
(k)\oplus\mathcal{O}_{\mathbb{P}^{1}})$ under a rational map defined by a base point free
linear series. The composition%
\[
C\subset Y\rightarrow P_{k}\rightarrow H^{0}(P_{k},\mathcal{O}_{P_{k}%
}(A)))\cong\mathbb{P}^{k+1}%
\]
defines a linear series of degree $b_{2}+2+2(k+1)=A.C^{\prime}$ (cf. Theorem
\ref{ConPk}), hence of Clifford index 2. For $k=0$ we get a further $g_{4}%
^{1}$ and therefore a $g_{8}^{3}$ from $g_{4}^{1}\times g_{4}^{1}.$ In the
case $k=1$, $Y$ is a Del-Pezzo surface, precisely $\mathbb{P}^{2}$ blown up in one
point $p$, and $C$ the strict transform of a plane model $C^{\prime}$ of
degree 6 with a doublepoint in $p.$ We easily check that for $k\geq3$ the
determinant of the matrix above becomes reducible, which can be excluded as
$Y$ is irreducible. It remains to consider the case $k=2$, where $C^{\prime}$ is a
divisor of class $4A$ on the quadric cone $Y,$ $A$ the hyperplane class of
$Y.$ Then the $g_{4}^{1}$ is cut out by the class of a ruling $R$ on $Y$ and
$\left\vert (A|_{C^{\prime}})\right\vert =\left\vert (2R|_{C^{\prime}})\right\vert
$ is a linear series of type $g_{8}^{3}.$ }

{\normalsize Conversely in the case where $C$ admits a $g_{8}^{3},$ we get a
space model $C^{\prime}$ of $C$ of degree $l=8.$ If $C^{\prime}$ has a
singular point, projection from this point leads to a $g_{d}^{2}$ with $d\leq6$. Otherwise $C^{\prime}$ is contained in a quadric surface $Q\subset\mathbb{P}^{3}.$ In
the situation where $Q\cong\mathbb{P}^{1}\times\mathbb{P}^{1}$ is nonsingular,
$C^{\prime}$ must be a divisor of type $(4,4)$ on $\mathbb{P}^{1}%
\times\mathbb{P}^{1}.$ Projection along each factor of $\mathbb{P}^{1}%
\times\mathbb{P}^{1}$ leads to two different $g_{4}^{1}$ on $C,$ such that we
get $(b_{1},b_{2})=(4,0)$ as otherwise the $g_{4}^{1}$ is uniquely determined.}

{\normalsize \noindent For $Q$ a singular quadric cone, $C^{\prime}$ is of
class $C^{\prime}\sim4A$, $A$ the class of a hyperplane section on $Q.$ Then
the rulings $\left\vert R\right\vert $ on $Q$ cut out a $g_{4}^{1}$ on $C$ and
because of $(A-2R)|_{C^{\prime}}\sim0,$ we get
\[
K_{C^{\prime}}\sim2A|_{C^{\prime}}\sim4R|_{C^{\prime}}\Rightarrow\left\vert
K_{C}\right\vert =4g_{4}^{1}%
\]
The adjoint series $\left\vert 2A\right\vert $ defines a mapping%
\[
\varphi:Q\rightarrow\mathbb{P}^{8}%
\]
with image $Q^{\prime}.$ This surface is given on $X$ by the determinant of a
matrix $\psi$ of type%
\[
\left(
\begin{array}
[c]{cc}%
H-a_{1}R & H-a_{2}R\\
H-(a_{1}+2)R & H-(a_{2}+2)R
\end{array}
\right)
\]
on $X$ with $a_{1},a_{2}\in\mathbb{Z}$, $a_{1}+a_{2}=f-4=2$ (cf. Theorem
\ref{ConPk} and Section 6.4. in \cite{schreyer1}). The scroll $X$ is of type
$S(4,2,0)$ and therefore $a_{1},a_{2}\leq2$ as $Q^{\prime}$ is irreducible.
The case $a_{1}=a_{2}=1$ can also be excluded as in this situation it is
possible to make one of the $H-3R$ entries zero, thus $Q^{\prime}$ would
become reducible. It follows that $a_{1}=2$, $a_{2}=0$ and therefore
$Q^{\prime}=Y\sim2H-4R\Rightarrow(b_{1},b_{2})=(4,0).$ }

{\normalsize It remains to show that in the situation where $C$ admits a plane model $C^{\prime}$
of degree $6,$ we also get $(b_{1},b_{2})=(4,0)$: $C^{\prime}$ has exactly
$\binom{5}{2}-9=1$ doublepoint $p$ as only singularity. Blowing up
}$\mathbb{P}^{2}$ in $p$ we get a surface $S$ with exceptional divisor
$E_{p}.$ Then the image under the adjoint series $\left\vert 3H-E_{p}%
\right\vert $ {\normalsize  is a Del-Pezzo surface }$S^{\prime}\subset
\mathbb{P}^{8}$ {\normalsize of degree }$g-3=6,$ that {\normalsize  is arithmetically Cohen-Macaulay. Its
projective dimension $R_{S^{\prime}}$ is $6=8-\dim S^{\prime}$ and
further$\ \operatorname*{reg}R_{S^{\prime}}=2.$ According to Theorem
\ref{hilbpoly} the Hilbert function $H_{R_{S^{\prime}}}$ of $R_{S^{\prime}}$ is
given by the Hilbert polynomial $P_{S^{\prime}}$ of $S^{\prime}$ calculated in
\cite{hartshorne} Chapter V, Exercise 1.2.:%
\[
P_{S^{\prime}}(n)=\frac{1}{2}an^{2}+bn+c
\]
with $a=(K_{S}+C)^{2}=(3H-E_{p})^{2}=9-1=8$, $b=\frac{1}{2}(K_{S}%
+C)^{2}+1-g(K_{S}+C)=4+1-1=4$ and $c=1.$ Therefore we obtain%
\[
H_{R_{S^{\prime}}}(n)=P_{S^{\prime}}(n)=4n^{2}+4n+1\text{ for all }%
n\in\mathbb{N}%
\]
The minimal free resolution $F$ of $R_{S^{\prime}}$ over $R=\Bbbk[x_{0},...,x_{n}]$ reduces modulo
$(y_{1},y_{2},y_{3})$ to a minimal free resolution of $R_{S^{\prime}}^{\prime
}:=R_{S^{\prime}}/(y_{1},y_{2},y_{3})R_{S^{\prime}}$ over $R/(y_{1}%
,y_{2},y_{3})R\cong R^{\prime}:=\Bbbk\lbrack x_{0}^{\prime},...,x_{5}^{\prime
}]$ with $(y_{1},y_{2},y_{3})$ being a $R_{S^{\prime}}$ sequence of linear
polynomials in $x_{0},...,x_{8}$. The corresponding Betti numbers $\beta_{ij}$
stay the same. Thus the Hilbert function of $R_{S^{\prime}}^{\prime}$ can be
obtained by succesively dividing out $y_{1},y_{2}$ and $y_{3},$ hence it has
values $(1,6,1)$ and $H_{R_{S^{\prime}}^{\prime}}(n)=0$ for $n\geq3.$ We
consider a free resolution of $\Bbbk$ with free $R^{\prime}$ modules which is
given by the Koszul complex of length 6:%
\[
0\leftarrow\Bbbk\leftarrow R^{\prime}\leftarrow R^{\prime6}\leftarrow
R^{\prime15}\leftarrow R^{\prime20}\leftarrow R^{\prime15}\leftarrow
R^{\prime6}\leftarrow R^{\prime}\leftarrow0
\]
The Betti numbers $\beta_{ij}=\dim Tor_{i}^{R^{\prime}}(\Bbbk,R_{S^{\prime}%
}^{\prime})_{j}$ can be calculated by tensoring the complex above with
$R_{S^{\prime}}^{\prime}:$%
\[
0\leftarrow\underset{M^{(0)}}{\underbrace{R_{S^{\prime}}^{\prime}}}%
\overset{\varphi_{1}}{\leftarrow}\underset{M^{(1)}}{\underbrace{R^{\prime
6}\otimes R_{S^{\prime}}^{\prime}}}\overset{\varphi_{2}}{\leftarrow}%
\underset{M^{(2)}}{\underbrace{R^{\prime15}\otimes R_{S^{\prime}}^{\prime}}%
}\leftarrow...\overset{\varphi_{6}}{\leftarrow}\underset{M^{(6)}}%
{\underbrace{R_{S^{\prime}}^{\prime}}}\leftarrow0
\]
Taking into acount the graduation we get the following format:%
\[
\text{\begin{xy}
\xymatrix{
M^{(0)} & M^{(1)} \ar[l] & M^{(2)} \ar[l] & M^{(3)} \ar[l] & M^{(4)} \ar[l] & M^{(5)} \ar[l] & M^{(6)} \ar[l] \\
1 & 6 \ar[ld] & 15 \ar[ld] & 20 \ar[ld] & 15 \ar[ld] & 6 \ar[ld] & 1 \ar[ld] \\
6 & 36 \ar[ld] & 90 \ar[ld] & 120 \ar[ld] & 90 \ar[ld] & 36 \ar[ld] & 6 \ar[ld] \\
1 & 6 & 15 & 20 & 15 & 6 & 1 \\
}
\end{xy}}%
\]

\noindent where the arrows stand for the maps $\varphi_{k}^{(l)}:M_{l}%
^{(k)}\rightarrow M_{l}^{(k-1)},$ which give a decomposition of $\varphi_{k}$
in the parts of degree $l$ and the numbers in the format are given by
$c_{kl}:=\dim M_{l}^{(k)}.$ As $S^{\prime}$ is not contained in any hyperplane
we get $\beta_{11}=0$ and therefore $\beta_{m,m}=0$ for $m=1,...,6.$ The dual
$F^{\ast}$ of $F$ is a free resolution of $\omega_{S^{\prime}}$ (up to a shift
of degrees) and therefore we obtain from Green's Linear Syzygy Theorem (see
\cite{eisenbud2} Theorem 7.1) that the length $n$ of the linear strand of
$F^{\ast}$ satisfies $n\leq\beta_{68}-1=1-1=0,$ hence $\beta_{m,m+2}=0$ for
$m=0,...,5.$ It follows that the Betti table }for $S^{\prime}$ takes the
following form{\normalsize :%

\[
\begin{tabular}{c|ccccccccc}
& 0 & 1 & 2 & 3 & 4 & 5 & 6 \cr \hline
0 & 1 & - & - & - & - & - & - \cr
1 & - & 20 & 64 & 90 & 64 & 20 & - \cr
2 & - & - & - & - & - & - & 1 \cr
\end{tabular} 
\]
From the exact sequence
\[
0\rightarrow\mathcal{O}_{S^{\prime}}(-C)\rightarrow\mathcal{O}_{S^{\prime}%
}\rightarrow\mathcal{O}_{C}\rightarrow0
\]
and $\mathcal{O}_{S^{\prime}}(-C)\cong\omega_{S^{\prime}}$ the minimal free
resolution of $\mathcal{O}_{C}$ as $\mathcal{O}_{\mathbb{P}^{8}}$-module is
given as the mapping cone of the minimal free resolutions of $\mathcal{O}%
_{S^{\prime}}$ and $\omega_{S^{\prime}}$ as every map is minimal:%
\[%
\begin{array}
[c]{ccccccccc}%
0 & \rightarrow & \omega_{S^{\prime}} & \rightarrow & \mathcal{O}_{S^{\prime}}
& \rightarrow & \mathcal{O}_{C} & \rightarrow & 0\\
&  & \uparrow &  & \uparrow &  &  &  & \\
&  & \mathcal{O}_{\mathbb{P}^{8}}\mathcal{(-}2)^{1} &  & \mathcal{O}%
_{\mathbb{P}^{8}} &  &  &  & \\
&  & \uparrow &  & \uparrow &  &  &  & \\
&  & \mathcal{O}_{\mathbb{P}^{8}}\mathcal{(-}4)^{20} &  & \mathcal{O}%
_{\mathbb{P}^{8}}\mathcal{(-}2)^{20} &  &  &  & \\
&  & \uparrow &  & \uparrow &  &  &  & \\
&  & \mathcal{O}_{\mathbb{P}^{8}}\mathcal{(-}5)^{64} &  & \mathcal{O}%
_{\mathbb{P}^{8}}\mathcal{(-}3)^{64} &  &  &  & \\
&  & \uparrow &  & \uparrow &  &  &  & \\
&  & \mathcal{O}_{\mathbb{P}^{8}}\mathcal{(-}6)^{90} &  & \mathcal{O}%
_{\mathbb{P}^{8}}\mathcal{(-}4)^{90} &  &  &  & \\
&  & \uparrow &  & \uparrow &  &  &  & \\
&  & \mathcal{O}_{\mathbb{P}^{8}}\mathcal{(-}7)^{64} &  & \mathcal{O}%
_{\mathbb{P}^{8}}\mathcal{(-}5)^{64} &  &  &  & \\
&  & \uparrow &  & \uparrow &  &  &  & \\
&  & \mathcal{O}_{\mathbb{P}^{8}}\mathcal{(-}8)^{20} &  & \mathcal{O}%
_{\mathbb{P}^{8}}\mathcal{(-}6)^{20} &  &  &  & \\
&  & \uparrow &  & \uparrow &  &  &  & \\
&  & \mathcal{O}_{\mathbb{P}^{8}}\mathcal{(-}10) &  & \mathcal{O}%
_{\mathbb{P}^{8}}\mathcal{(-}8) &  &  &  & \\
&  & \uparrow &  & \uparrow &  &  &  & \\
&  & 0 &  & 0 &  &  &  &
\end{array}
\]
Therefore we obtain the Betti table for }$C:${\normalsize
\textnormal{\[
\begin{tabular}{c|cccccccccc}
& 0 & 1 & 2 & 3 & 4 & 5 & 6 & 7 \cr \hline
0 & 1 & - & - & - & - & - & - & - \cr
1 & - & 21 & 64 & 90 & 64 & 20 & - & - \cr
2 & - & - & 20  & 64  & 90  & 64  & 21 & -\cr
3 & - & - & - & - & - & - & - & 1 \cr
\end{tabular} 
\]}
\noindent We have shown that a tetragonal curve $C$ of genus $9$ admits a $g_{6}^{2}$ or
a $g_{8}^{3}$ exactly if its betti table looks as above, which is exactly the
case for $(b_{1},b_{2})=(4,0)$. In the next step we will give the Betti table
for a tetragonal curve }$C$ that admits no $g_{6}^{2},$ $g_{8}^{3}$ or a base
point free $g_{d}^{r}$ of Clifford index $3:$
\\\\
\textbf{2. }Let $g_{4}^{1}=\left\vert D\right\vert $ with an effective divisor
$D$ of degree $4$ on $C$. We can omit the case $h^{0}(C,\mathcal{O}%
_{C}(2D))=4$ as this would give a $g_{8}^{3}=\left\vert 2D\right\vert .$ Then
$F\sim K_{C}-2D$ gives a $g_{8}^{2}$. If $\left\vert F\right\vert $ has base
points we can deduce a $g_{d}^{2}$ with $d\leq7$, that we have excluded. In
the situation where $F\sim2D$ we must have $h^{0}(C,\mathcal{O}_{C}%
(K_{C}-3D))=2$ and $h^{0}(C,\mathcal{O}_{C}(K_{C}-4D))=1$, hence the scroll $X$
is of type $S(4,1,1).$ Then $(b_{1},b_{2})\neq(3,1)$ as the defining equation
of $Y\sim2H-3R$ would contain a factor $\varphi_{0}\in H^{0}(X,\mathcal{O}%
_{X}(H-4R)),$ which contradicts that $C\subset Y$ is irreducible. It follows
that $(b_{1},b_{2})=(2,2)$ in this case, hence $C$ has Betti table as follows%
\[
\begin{tabular}{c|cccccccccc}
& 0 & 1 & 2 & 3 & 4 & 5 & 6 & 7 \cr \hline
0 & 1 & - & - & - & - & - & - & - \cr
1 & - & 21 & 64 & 70 & 24 & 5 & - & - \cr
2 & - & - & 5 & 24  & 70  & 64  & 21 & -\cr
3 & - & - & - & - & - & - & - & 1 \cr
\end{tabular} 
\]
Now we assume that $F\not \sim 2D$, then from $\left\vert F\right\vert $ we
get a plane model $C^{\prime}\subset\mathbb{P}^{2}$ of degree 8 for $C$.
$C^{\prime}$ does not have any triple points as projection from such a point
would give a base point free $g_{5}^{1}.$ If $C^{\prime}$ has double points as
only singularities then from a base point free $g_{6}^{1}$ obtained from
projection from one of these double points and the $g_{4}^{1}$ we get a space
model $C^{\prime\prime}\subset\mathbb{P}^{1}\times\mathbb{P}^{1}%
\subset\mathbb{P}^{3}$ of $C.$ $C^{\prime\prime}$ is a divisor of class
$(4,6)$ on $\mathbb{P}^{1}\times\mathbb{P}^{1}$, thus it has arithmetic genus
$p_{a}(C^{\prime\prime})=3\cdot5=15>9$. $C^{\prime\prime}$ has no singularity
of multiplicity $3$ or higher as projection from such a point would give a
$g_{d}^{2}$ with $d\leq7.$ It follows that $C^{\prime\prime}$ has exactly
$15-9=6$ double points as only singularities hence projection from one of them
leads to a $g_{8}^{2}$ with a quadruple point and $6$ double points. Thus we
can assume the existence of a plane model $C^{\prime}$ of $C$ with a quadruple
point $p_{0}$ and $6$ double points $p_{1},...,p_{6}$ as only singularities.
After blowing up $\mathbb{P}^{2}$ in these points%
\[
\sigma:S:=\mathbb{\tilde{P}}^{2}(p_{0},...,p_{6})\rightarrow\mathbb{P}^{2}%
\]
with exceptional divisors $E_{0},...,E_{6}$ and hyperplane class $H,$ we can
assume that $C$ is the strict transform of $C^{\prime},$ i.e. $C\sim8H-4E_{0}-%
{\textstyle\sum\nolimits_{i=1}^{6}}
E_{i}$. Then the canonical system is cut out on $C$ by the adjoint series
$\left\vert L^{\prime}\right\vert =\left\vert 5H-3E_{0}-%
{\textstyle\sum\nolimits_{i=1}^{6}}
E_{i}\right\vert .$ With the methods of Chapter 2 we can show that this series
is base point free and the image $S^{\prime}\subset\mathbb{P}^{8}$ of $S$
under the morphism defined by $\left\vert L^{\prime}\right\vert $ has only
isolated singularities and is arithmetically Cohen-Macaulay.
{\normalsize Therefore the projective dimension of its homogenous coordinate
ring $R_{S^{\prime}}$ is $6=8-\dim S^{\prime}$ and
further$\ \operatorname*{reg}R_{S^{\prime}}=2.$ According to Theorem
\ref{hilbpoly} the Hilbertfunction $H_{R_{S^{\prime}}}$ of $R_{S^{\prime}}$ is
given by the Hilbert polynomial $p_{S^{\prime}}$ of $S^{\prime}$ calculated in
\cite{hartshorne} Chapter V, Exercise 1.2.:%
\[
P_{S^{\prime}}(n)=\frac{1}{2}an^{2}+bn+c
\]
with $a=(K_{S}+C)^{2}=(5H-3E_{0}-%
{\textstyle\sum\nolimits_{i=1}^{6}}
E_{i})^{2}=10$, $b=\frac{1}{2}(K_{S}+C)^{2}+1-g(K_{S}+C)=3$ and $c=1.$
Therefore we obtain%
\[
H_{R_{S^{\prime}}}(n)=P_{S^{\prime}}(n)=5n^{2}+3n+1\text{ for all }%
n\in\mathbb{N}%
\]
Now we consider the same approach as in \textbf{1.} to obtain the Betti table
for $C:$ As in the previous consideration the minimal free resolution $F$ of $R_{S^{\prime}}$ over $R$ reduces
modulo to a minimal free resolution of $R_{S^{\prime}%
}^{\prime}:=R_{S^{\prime}}/(y_{1},y_{2},y_{3})R_{S^{\prime}}$ over
$R/(y_{1},y_{2},y_{3})R\cong R^{\prime}:=\Bbbk\lbrack x_{0}^{\prime}%
,...,x_{5}^{\prime}]$. The corresponding Hilbert function of
$R_{S^{\prime}}^{\prime}$ has values $(1,6,3)$ and $H_{R_{S^{\prime
}}^{\prime}}(n)=0$ for $n\geq3.$ Tensoring the Koszul complex of length 6 with
$R_{S^{\prime}}^{\prime}$ a similar argumentation as in \textbf{1. }shows that
$\beta_{m,m}=0$, $m=1,...,6$ and $\beta_{m,m+2}=0$, $m=0,...,3$ for the Betti
numbers of }$S^{\prime}.$ {\normalsize The linear strand of the minimal free
resolution of }$R_{S^{\prime}}$ is a subcomplex of the minimal free resolution
of $R_{C},$ thus we must have {\normalsize $\beta_{6,7}=0$ as }%
$\operatorname*{Cliff}(C)=2.$ {\normalsize The Betti table of the minimal free
resolution of $R_{S^{\prime}}$ then takes the following form%
\[
\begin{tabular}{c|ccccccccc}
& 0 & 1 & 2 & 3 & 4 & 5 & 6 \cr \hline
0 & 1 & - & - & - & - & - & - \cr
1 & - & 18 & 52 & 60 & 24 & 5 & - \cr
2 & - & - & - & - & 15 & 12 & 3 \cr
\end{tabular} 
\]
and from the exact sequence
\[
0\rightarrow\mathcal{O}_{S^{\prime}}(-C)\rightarrow\mathcal{O}_{S^{\prime}%
}\rightarrow\mathcal{O}_{C}\rightarrow0
\]
the minimal free resolution of $\mathcal{O}_{C}$ as $\mathcal{O}%
_{\mathbb{P}^{8}}$-module is given as the mapping cone of the minimal free
resolutions of $\mathcal{O}_{S^{\prime}}$ and $\omega_{S^{\prime}}$:%
\[
\begin{tabular}{c|cccccccccc}
& 0 & 1 & 2 & 3 & 4 & 5 & 6 & 7 \cr \hline
0 & 1 & - & - & - & - & - & - & - \cr
1 & - & 21 & 64 & 70 & 24 & 5 & - & - \cr
2 & - & - & 5 & 24  & 70  & 64  & 21 & -\cr
3 & - & - & - & - & - & - & - & 1 \cr
\end{tabular} 
\]
Therefore we have shown that for a tetragonal curve $C$ of genus $9$ that admits a
$g_{4}^{1}$ and no $g_{6}^{2}$, $g_{8}^{3}$ or a further }$g_{d}^{r}$ of
Clifford index $3${\normalsize , its Betti table looks as above, which is
exactly the case for $(b_{1},b_{2})=(2,2)$. }
\\\\
{\normalsize \noindent\textbf{3. }$(b_{1},b_{2})=(3,1):$ In this situation
$C$ has Betti table as follows%
\[
\begin{tabular}{c|cccccccccc}
& 0 & 1 & 2 & 3 & 4 & 5 & 6 & 7 \cr \hline
0 & 1 & - & - & - & - & - & - & - \cr
1 & - & 21 & 64 & 75 & 44 & 5 & - & - \cr
2 & - & - & 5 & 44  & 75  & 64  & 21 & -\cr
3 & - & - & - & - & - & - & - & 1 \cr
\end{tabular} 
\]
According to our results in \textbf{1. }}and \textbf{2.}, a curve $C$ that is
given as complete intersection of divisors of type $Y\sim2H-3R$ \ ,
\ $Z\sim2H-R$ on the scroll $X$ admits a base point free $g_{d}^{r}$ of Clifford index $3.$ It remains to show that from the existence of such a linear system it follows that $(b_{1},b_{2})=(3,1)$. {\normalsize From a base point free linear system of Clifford index }$3$ we
always get a $g_{5}^{1}$ or a $g_{7}^{2}$ (The Brill Noether dual of a
$g_{9}^{3}$ or $g_{11}^{4}$ is a $g_{7}^{2}$ or $g_{5}^{1}$ respectively!). A
$g_{7}^{2}$ cannot have a singularity with higher multiplicity than 3 as in
this case projection from this singular point would lead to $g_{d}^{1}$ with
$d\leq3$. The $g_{7}^{2}$ cannot have exactly two triple points as then from the existence of two different $g_{4}^{1}$$^{\prime}s$ we would deduce that $(b_{1},b_{2})=(4,0)$ as in \textbf{1.}, hence the existence of a $g_{6}^{2}$ or a $g_{8}^{3}$. {\normalsize Thus the }$g_{7}^{2}$ {\normalsize has exactly 6 double
points or one triple point and 3 double points. Except in the case where three double points are infinitely near a triple point, projection from one of the double points leads to a $g_{5}^{1}.$ Therefore $C$ admits a $g_{4}^{1}\times g_{5}^{1}$ with the mentioned exception. In the case where $C$ admits a $g_{4}^{1}\times g_{5}^{1}$ but no $g_{6}^{2}$ or $g_{8}^{3}$, we consider the image $C^{\prime
}\subset\mathbb{P}^{1}\times\mathbb{P}^{1}\subset\mathbb{P}^{3}$ of degree $9$
under the morphism obtained from the $g_{4}^{1}\times g_{5}^{1}.$ Then
$C^{\prime}$ is a divisor of class $(4,5)$ on $\mathbb{P}^{1}\times
\mathbb{P}^{1}$ with arithmetic genus $p_{a}(C^{\prime})=12$. Further
$C^{\prime}$ has exactly three double points as only singularities, as
projection from a singular point of higher multiplicity would lead to a
$g_{d}^{2}$ with $d\leq6.$ Then projection from one of the double points gives
a plane model $C^{\prime\prime}\subset\mathbb{P}^{2}$ of degree 7 with one
triple point $p_{0}$ and 3 double points $p_{1},p_{2}$ and $p_{3}.$ Thus in general we can assume that we have a $g_{7}^{2}$ with one triple point $p_{0}$ and 3 (possibly infinitely near) double points $p_{1},...,p_{3}$. Blowing up
$\mathbb{P}^{2}$ in the these points%
\[
\sigma:S=\mathbb{\tilde{P}}^{2}(p_{0},...,p_{3})\rightarrow\mathbb{P}^{2}%
\]
with exceptional divisors $E_{0},...,E_{3},$ we can consider $C\sim7H-3E_{0}-%
{\textstyle\sum\nolimits_{i=1}^{3}}
2E_{i}$ to be the strict transform of $C^{\prime\prime}$, $H$ denoting the
hyperplane class on $\mathbb{P}^{2}$ and by abuse of notation also on $S$. The
} {\normalsize canonical series is then cut out by the adjoint series
$\left\vert L^{\prime}\right\vert :=\left\vert 4H-2E_{0}-%
{\textstyle\sum\nolimits_{i=1}^{3}}
2E_{i}\right\vert $, which is base point free on $S$ (cf. Corollary
\ref{ALSP2}). The image $S^{\prime}$ of $S$ under the morphism given by
$\left\vert L^{\prime}\right\vert $ is then a Castelnuovo surface and it is
even arithmetically Cohen-Macaulay. Therefore the projective dimension
$R_{S^{\prime}}$ is $6=8-\dim S^{\prime}$ and further$\ \operatorname*{reg}%
R_{S^{\prime}}=2.$ According to Theorem \ref{hilbpoly} the Hilbertfunction
$H_{R_{S^{\prime}}}$ of $R_{S^{\prime}}$ is given by the Hilbert polynomial
$p_{S^{\prime}}$ of $S^{\prime}$ calculated in \cite{hartshorne} Chapter V, Ex. 1.2.:%
\[
P_{S^{\prime}}(n)=\frac{1}{2}an^{2}+bn+c
\]
with $a=(K_{S}+C)^{2}=(4H-2E_{0}-%
{\textstyle\sum\nolimits_{i=1}^{3}}
2E_{i})^{2}=9$, $b=\frac{1}{2}(K_{S}+C)^{2}+1-g(K_{S}+C)=5$ and $c=1.$
Therefore we obtain%
\[
H_{R_{S^{\prime}}}(n)=P_{S^{\prime}}(n)=4n^{2}+4n+1\text{ for all }%
n\in\mathbb{N}%
\]
Now we consider the same approach as in \textbf{1.} and \textbf{2. }to obtain
the Betti table for $C:$ The minimal free resolution $F$ of $R_{S^{\prime}}$
over $R$ reduces to a minimal free resolution of
$R_{S^{\prime}}^{\prime}:=R_{S^{\prime}}/(y_{1},y_{2},y_{3})R_{S^{\prime}}$
over $R/(y_{1},y_{2},y_{3})R\cong R^{\prime}:=\Bbbk\lbrack x_{0}^{\prime
},...,x_{5}^{\prime}]$. The Hilbertfunction of
$R_{S^{\prime}}^{\prime}$ has values $(1,6,2)$ and $H_{R_{S^{\prime
}}^{\prime}}(n)=0$ for $n\geq3.$ Tensoring the Koszul complex of length 6 with
$R_{S^{\prime}}^{\prime}$ a similar argumentation as in \textbf{1. }shows that
$\beta_{m,m}=0$, $m=1,...,6$ and $\beta_{m,m+2}=0$, $m=0,...,4$ for the Betti
numbers of }$S^{\prime}.$ {\normalsize The linear strand of the minimal free
resolution of }$R_{S^{\prime}}$ is a subcomplex of the minimal free resolution
of $R_{C},$ thus we must have {\normalsize $\beta_{6,7}=0$ as }%
$\operatorname*{Cliff}(C)=2.$ {\normalsize The Betti table of the minimal free
resolution of $R_{S^{\prime}}$ then takes the following form%
\[
\begin{tabular}{c|cccccccccc}
& 0 & 1 & 2 & 3 & 4 & 5 & 6  \cr \hline
0 & 1 & - & - & - & - & - & - \cr
1 & - & 19 & 58 & 75 & 44 & 5 & - \cr
2 & - & - & - & -  & - & 6  & 2 \cr
\end{tabular} 
\]
and from the exact sequence
\[
0\rightarrow\mathcal{O}_{S^{\prime}}(-C)\rightarrow\mathcal{O}_{S^{\prime}%
}\rightarrow\mathcal{O}_{C}\rightarrow0
\]
the minimal free resolution of $\mathcal{O}_{C}$ as $\mathcal{O}%
_{\mathbb{P}^{8}}$-module is given as the mapping cone of the minimal free
resolutions of $\mathcal{O}_{S^{\prime}}$ and $\omega_{S^{\prime}}$:%
\[
\begin{tabular}{c|cccccccccc}
& 0 & 1 & 2 & 3 & 4 & 5 & 6 & 7 \cr \hline
0 & 1 & - & - & - & - & - & - & - \cr
1 & - & 21 & 64 & 75 & 44 & 5 & - & - \cr
2 & - & - & 5 & 44  & 75  & 64  & 21 & -\cr
3 & - & - & - & - & - & - & - & 1 \cr
\end{tabular} 
\]
Therefore we have shown that a tetragonal curve $C$ of genus $9$ admits a
$g_{4}^{1}$, a base pont free linear system of Clifford index 3 and no $g_{6}^{2}$ or $g_{8}^{3}$ exactly if its
Betti table looks as above, which is exactly the case for $(b_{1}%
,b_{2})=(3,1)$. }
\end{proof}

{\normalsize \newpage}
°
\newpage
\chapter{{\protect\normalsize Curves C with Clifford index Cliff(C)=3}}

\section{Results}
{\normalsize Let $C\subset\mathbb{P}^{8}$ be a canonical curve of genus 9 and
Clifford index 3, then there exists a $g_{2r+3}^{r}$ with }$r\in\{1,...,6\}.$
Its Brill Noether dual is of type $g_{13-2r}^{5-r}$, thus for a curve of genus 9 it remains to consider the case where $C$ has a $g_{5}^{1}$ or a $g_{7}^{2}$. A plane curve of degree 7 has arithmetic genus
$\binom{6}{2}=15$, therefore it must have singular points. If it has a
singular point of multiplicity $d\geq3$ then projection from this point leads
to a $g_{7-d}^{1}$ which is a linear system of Clifford index less or equal to
2, a contradiction. Therefore the $g_{7}^{2}$ has exactly six double points as
only singularities. Projection from one of them gives a $g_{5}^{1}.$ We conclude:

\begin{corollary}
{\normalsize Let $C\subset\mathbb{P}^{8}$ be a canonical curve of genus 9 with} Clifford index {\normalsize $\operatorname*{Cliff}(C)=3$. Then} $C$ is pentagonal.
\end{corollary}

{\normalsize In the following section we first want to give some examples
of pentagonal curves that exist. Afterwards we will consider the special
case where }$C$ admits a $g_{7}^{2}. $ {\normalsize  It turns out that such}%
$C$ {\normalsize  has the following Betti table (cf. Theorem \ref{g27}):%

\[
\begin{tabular}{c|cccccccccc}
& 0 & 1 & 2 & 3 & 4 & 5 & 6 & 7 \cr \hline
0 & 1 & - & - & - & - & - & - & - \cr
1 & - & 21 & 64 & 70 & 24 & - & - & - \cr
2 & - & - & - & 24  & 70  & 64  & 21 & -\cr
3 & - & - & - & - & - & - & - & 1 \cr
\end{tabular} 
\]
\\\\
The corresponding plane model $C^{\prime}\subset\mathbb{P}^{2}$ of degree 7
has exactly 6 (possibly infinitely near)\ double points. Hence we deduce 6
(possible infinitely near) $g_{5}^{1}$$^{\prime}s$ by projection from these double
points. It turns out that these are exactly all linear series of this type
(cf. Theorem \ref{g27a}). }

{\normalsize In the main part of this thesis, we focus on pentagonal curves
$C,$ where none of the appearing $g_{5}^{1}$$^{\prime}s$ can be obtained from a $g_{7}^{2}.$
Starting with $D$ an effective divisor on $C,$ such that $\left\vert
D\right\vert $ is of type $g_{5}^{1},$ we construct the corresponding }%
$4-$dimensional {\normalsize scroll }$X$ of type $S(e_{1},...,e_{4}%
)$ {\normalsize  from the complete linear system }$\left\vert D\right\vert $
(cf. Section 1.4)$.$ {\normalsize  Applying Riemann-Roch one easily determines
the following possibilities for $D$ that can occur (cf. Theorem
\ref{scrolltypes}):%
\begin{align*}
(1)~~h^{0}(C,\mathcal{O}_{C}(2D)) &  =3\text{ and }h^{0}(C,\mathcal{O}%
_{C}(3D))=7\\
(2)~~h^{0}(C,\mathcal{O}_{C}(2D)) &  =4\text{ and }h^{0}(C,\mathcal{O}%
_{C}(3D))=7\\
(3)~~h^{0}(C,\mathcal{O}_{C}(2D)) &  =4\text{ and }h^{0}(C,\mathcal{O}%
_{C}(3D))=8
\end{align*}
Then the scroll }$X$ is of type $S(2,1,1,1),$ $S(2,2,1,0)$ or $S(3,1,1,0)$ respectively.

{\normalsize \noindent In this situation the results of Sections
\ref{g15multi1}-\ref{deformation} justify the definition \index{$m_{|D|}$ multiplicity od a $\left\vert D\right\vert :$ of type $g_{5}^{1}$} of the
\textit{multiplicity} $m_{|D|}$ of $\left\vert D\right\vert :$ }

\begin{definition}
{\normalsize Let $C$ be an irreducible, smooth, canonical curve of genus $9$
and $\operatorname*{Cliff}(C)=3.$ Further assume that $C$ has no $g_{7}^{2}$.
Then for each linear system $\left\vert D\right\vert $ of type $g_{5}^{1}$ on
$C$ we define its multiplicity $m_{\left\vert D\right\vert }$ to be equal to
$i$ if and only if $D$ fullfills $(i).$ }
\end{definition}

{\normalsize \noindent Let $k=%
{\textstyle\sum\nolimits_{\left\vert D\right\vert =g_{5}^{1}}}
m_{\left\vert D\right\vert }$ be the number of all appearing $g_{5}^{1}$,
counted with multiplicities, then one of the main results of this thesis is
that $k\leq3$ if $\operatorname*{Cliff}(C)=3$ and $C$ has no $g_{7}^{2}$. This exactly occurs if and only if }$C$ {\normalsize  has Betti table as
follows:%
\[
\begin{tabular}{c|cccccccccc}
& 0 & 1 & 2 & 3 & 4 & 5 & 6 & 7 \cr \hline
0 & 1 & - & - & - & - & - & - & - \cr
1 & - & 21 & 64 & 70 & 4k & - & - & - \cr
2 & - & - & -  & 4k  & 70  & 64  & 21 & -\cr
3 & - & - & - & - & - & - & - & 1 \cr
\end{tabular} 
\]
For }$k=1$ there exists a plane model for $C$ of degree $8$ with exactly one
triple point and $9$ double points as only singularities. {\normalsize Furthermore
in the case $k\geq2$, we obtain a space model of $C$ on }a quadric surface
$Y\subset\mathbb{P}^{3}${\normalsize . In the case where }$C$ has only
$g_{5}^{1\prime}s$ of multiplicity one, the quadric $Y$ is smooth otherwise it is a
cone over a conic. {\normalsize  We will also show that in all cases
where $C$ has a $g_{5}^{1}$ with multiplicity 2 or 3, there exists a local one-parameter family of curves $C_{t}$ with $C_{0}=C$ and $C_{t}$ having the
correspondent number of distinct $g_{5}^{1}$ for $t\neq0$, such that this case
can be seen as a specialization of those where all $g_{5}^{1}$ are different
(cf. Section \ref{deformation}).}

\section{Examples of pentagonal curves}

{\normalsize We have remarked above that for an irreducible, canonical curve
$C$ of genus $9$ and Clifford index $3$ it is possible to have a certain
number $k\in\{1,2,3\}$ of $g_{5}^{1\prime}s$ (counted with multiplicities) or
even a $g_{7}^{2}$.\ We give examples of such curves over a field of finite characteristic $p.$ Applying the main results of
this work and taking into account the semicontinuity of the Betti numbers, we get that
each of these curves has the right number of $g_{5}^{1\prime}s.$ 
\begin{example}{\normalsize (Canonical Curve of genus 9 with exactly one $g_{5}^{1}$) }{\textnormal{Let us start
with the construction of a curve $C$ that has exactly one $g_{5}^{1}.$ From
the Brill Noether number
\[
\varrho(9,8,2)=9-(2+1)\cdot(9-8+2)=0
\]
we deduce the existence of a plane model $C^{\prime}$ of degree 8. If
$C^{\prime}$ has a triple point, projection from this point would lead to a
$g_{5}^{1}.$ 
{\normalsize \noindent We consider a curve $C^{\prime}$ that has exactly one
triple point $q$ and 9 double points $p_{1},...,p_{9}$ in general
position.}}
\textnormal{
{\normalsize \noindent Then we determine its normalisation $C\subset\mathbb{P}^{8}$ given as the image of $C^{\prime}$ under the adjoint series $\left\vert 5H-2q-%
{\textstyle\sum\nolimits_{i=1}^{9}}
p_{i}\right\vert $. Obviously $C$ has at least one $g_{5}^{1}$ which is
obtained from projection from the triple point in the plane model. The Betti
table for }$C$ {\normalsize then takes the following form%
\[
\begin{tabular}{c|cccccccccc}
& 0 & 1 & 2 & 3 & 4 & 5 & 6 & 7 \cr \hline
0 & 1 & - & - & - & - & - & - & - \cr
1 & - & 21 & 64 & 70 & 4 & - & - & - \cr
2 & - & - & - & 4  & 70  & 64  & 21 & -\cr
3 & - & - & - & - & - & - & - & 1 \cr
\end{tabular} 
\]
According to our main results $C$ has exactly one $g_{5}%
^{1}$ of multiplicity one and $\operatorname*{Cliff}(C)=3. $ 
The following theorem says that every pentagonal curve $C$ of genus 9 admits a plane model of degree 8 with a triple point.
\begin{theorem}{\label{triplepoint}}
a) Let $C$ be a canonical curve of genus 9 with $\operatorname{Cliff}(C)=3$ that admits no $g_{7}^{2}$, then there exists a
$g_{8}^{2}$ on $C$ with at least one triple point. 
\\b) If $C$ has exactly one $g_{5}^{1}$ (of multiplicity one) then the
$g_{8}^{2}$ in a) has exactly one triple point and $9$ double points.
\end{theorem}
\begin{proof}
We start with the second claim. Let $C$ admit a plane model $C^{\prime}$ of
degree 8 with two triple points $p_{1},$ $p_{2}$. If $p_{1}$ and $p_{2}$ are
not infinitely near, projection from each of them gives two distinct
$g_{5}^{1},$ so it remains to discuss the case where $p_{1}$ and
$p_{2}$ are infinitely near. Blowing up $\mathbb{P}^{2}$ in the triple
points $p_{1},$ $p_{2}$ and the double points $q_{1},...,q_{6}$:%
\[
\sigma:S\rightarrow\mathbb{P}^{2}%
\]
with exceptional divisors $E_{p_{1}},E_{p_{2}}$, $E_{q_{1}},...,E_{q_{6}}$ and
hyperplane class $H,$ we can assume that $C$ is the strict transform of
$C^{\prime}:$%
\[
C\sim8H-3E_{p_{1}}-3E_{p_{2}}-%
{\textstyle\sum\nolimits_{i=1}^{6}}
2E_{q_{i}}%
\]
The canonical system is cut out by the adjoint \index{$PGL(n)$ projective general linear group} series%
\[
|A_{C}|=|5H-2E_{p_{1}}-2E_{p_{2}}-%
{\textstyle\sum\nolimits_{i=1}^{6}}
E_{q_{i}}|%
\]
As $p_{2}$ lies infinitely near to $p_{1}$ we obtain $E_{p_{1}}|_{C}=E_{p_{2}}|_{C}.$
For the complete linear series $g_{5}^{1}=\left\vert D\right\vert =\left\vert
(H-E_{p_{1}})|_{C}\right\vert =\left\vert (H-E_{p_{2}})|_{C}\right\vert $ it
follows that%
\[
K_{C}-2D=3H|_{C}-E_{p_{1}}|_{C}-E_{p_{2}}|_{C}-%
{\textstyle\sum\nolimits_{i=1}^{6}}
E_{q_{i}}|_{C}%
\]
and because of $h^{0}(S,\mathcal{O}_{S}(3H-E_{p_{1}}-E_{p_{2}}-%
{\textstyle\sum\nolimits_{i=1}^{6}}
E_{q_{i}}))\geq10-8=2$ we get $h^{0}(C,\mathcal{O}_{C}(K_{C}-2D))=2$ (as $\operatorname*{Cliff}%
(C)=3$), thus $m_{\left\vert D\right\vert }\geq2.$ This proves b).
It remains the question if there always exists a $g_{8}^{2}$ with a triple point. We know that there always exists a $g_{8}^{2}$ on $C$. This series is base point free as we have assumed that $\operatorname{Cliff}(C)=3$ and $C$ admits no $g_{7}^{2}$. Let $H$ be an effective divisor of this system, then $h^{0}(C,\mathcal{O}_{C}(H-D))=0$ or $1$. For $h^{0}(C,\mathcal{O}_{C}(H-D))=1$ it follows the existence of points $q_{1},q_{2}$ and $q_{3}$ on $C$ with
\[
H-D\sim q_{1}+q_{2}+q_{3}
\]
Thus $|H|$ has a triple point. In the case $h^{0}(C,\mathcal{O}_{C}(H-D))=0$ we can apply the base point free pencil trick to obtain
\[
h^{0}(C,\mathcal{O}_{C}(H+D))\geq h^{0}(C,\mathcal{O}_{C}(H))\cdot h^{0}(C,\mathcal{O}_{C}(D))-h^{0}(C,\mathcal{O}_{C}(H-D))=6
\]
Riemann-Roch then says that $h^{0}(C,\mathcal{O}_{C}((K_{C}-H)-D))\geq 1$ and equality holds. The $g_{8}^{2}$ given by the divisor $L:=K_{C}-H$ is then base point free and has a triple point. The first property follows from the assumption that $\operatorname{Cliff}(C)=3$ and $C$ admits no $g_{7}^{2}$ and the second from the existence of points $q_{1},q_{2}$ and $q_{3}$ on $C$ with $L-D\sim K_{C}-D\sim q_{1}+q_{2}+q_{3}$
\end{proof}
\\\\Counting dimensions we have $2\cdot10$ possibilities to choose 10 points in $\mathbb{P}^{2}$ and a parameter space of plane curves of degree 8 passing one of the points with multiplicity 3 and the others with multiplicity 2 that has projective dimension $\binom{10}{2}-1-6-3\cdot9$. Taking into account the projective transformations on $\mathbb{P}^{2}$ we obtain the dimension for the subscheme $\mathcal{H}_{(1,5)}\subset \mathcal{M}_{9}$ of all pentagonal curves of genus 9:}
\begin{align*}
\dim \mathcal{H}_{(1,5)} &  =2\cdot10+\binom{10}{2}-1-6-3\cdot9-\dim PGL(3)=23=\\
&  =\dim\mathcal{M}_{9}-1
\end{align*}
}}
\end{example}

\bigskip

\begin{example}
{\normalsize (Canonical curves \index{$\left\vert L-p_{1}-...-p_{m} \right\vert$ effective divisors linear equivalent to $L$, passing the points $p_{1},...,p_{m}$} with exactly two different ordinary $g_{5}^{1}%
$$^{\prime}s$) }{\textnormal{{\normalsize We want to construct a curve $C$ that has at
least two different $g_{5}^{1}$$^{\prime}s$ and no $g_{7}^{2}.$ From these special linear
systems we obtain a space model $C^{\prime}\subset\mathbb{P}^{1}%
\times\mathbb{P}^{1}\subset\mathbb{P}^{3}$ given as the image of the natural
mapping:%
\[
C\overset{g_{5}^{1}\times g_{5}^{1}}{\rightarrow}\mathbb{P}^{1}\times
\mathbb{P}^{1}\subset\mathbb{P}^{3}%
\]
This model in $\mathbb{P}^{3}$ has exactly $7=p_{a}(C^{\prime})-g(C)=4\cdot
4-9$ (possibly infinitely near) double points $p_{1},...,p_{7}$, as
otherwise projection from a singular point of higher multiplicity would lead
to a $g_{d}^{2}$ with $d\leq7.$ For a general curve $C$ having two $g_{5}^{1}$
we would expect that all double points are distinct. Then $C^{\prime}$ is a
divisor of type $(5,5)$ on $\mathbb{P}^{1}\times\mathbb{P}^{1}$ and there
exists two $g_{5}^{1}$$^{\prime}s$ which are cut out by the complete linear systems
$\left\vert (1,0)\right\vert $ and $\left\vert (0,1)\right\vert $ on
$\mathbb{P}^{1}\times\mathbb{P}^{1}.$ We remark that projection from one of
the double points gives a plane model of }$C$ that has exactly two triple points. The following example of a curve in this family is calculated over $\mathbb{Q}$:}
\\\\
{\normalsize \label{C55pic}\textit{Curve of type }$(5,5)$ \textit{on}
$\mathbb{P}^{1}\times\mathbb{P}^{1}$\textit{ with 7 double points as only
singularities}
\\\\\\
\tiny{$I_{C}=ideal(x^{2}+y^{2}-z^{2}-w^{2},w^{3}x^{2}-\frac{70}{3}w^{2}x^{3}-\frac{21}{5}wx^{4}-\frac{49}{9}x^{5}+3w^{4}y-\frac{28}{5}w^{3}y^{2}-\frac{698428568345581244477}{28848956524627500}w^{2}y^{3}+\frac{90215896038289435403}{14424478262313750}wy^{4}-\frac{68220455180173043}{250860491518500}y^{5}-6w^{4}z-7w^{3}xz+14w^{2}x^{2}z-\frac{7}{3}wx^{3}z-\frac{35}{2}x^{4}z+\frac{49}{5}w^{3}yz+\frac{23241926205812880277}{1492187406446250}w^{2}y^{2}z-\frac{28199536240159289641}{1730937391477650}wy^{3}z+\frac{44759596869109840639}{21636717393470625}y^{4}z+\frac{49}{3}w^{3}z^{2}-
\\ \frac{4749324908888289073}{5769791304925500}w^{2}xz^{2}-\frac{83077348288551341}{1923263768308500}wx^{2}z^{2}+\frac{8743506435415727}{83620163839500}x^{3}z^{2}+\frac{94638741549090136267}{22577444236665000}w^{2}yz^{2}+\frac{877466317906635532537}{259640608721647500}wy^{2}z^{2}-\frac{554169997124749656173}{259640608721647500}y^{3}z^{2}-\frac{73873813541874770027}{43273434786941250}w^{2}z^{3}-\frac{487977728485206293}{2404079710385625}wxz^{3}+\frac{114591427898711159}{5769791304925500}x^{2}z^{3}+\frac{23322157426391471302}{7212239131156875}wyz^{3}-\frac{147619621394753090843}{259640608721647500}y^{2}z^{3}-\frac{23712776715338563288}{64910152180411875}wz^{4}-\frac{48249553756301249983}{259640608721647500}xz^{4}+\frac{36743223513209826373}{173093739147765000}yz^{4}+\frac{5859389203919613817}{129820304360823750}z^{5})$
}
\[
{\includegraphics[height=7cm,width=8cm]{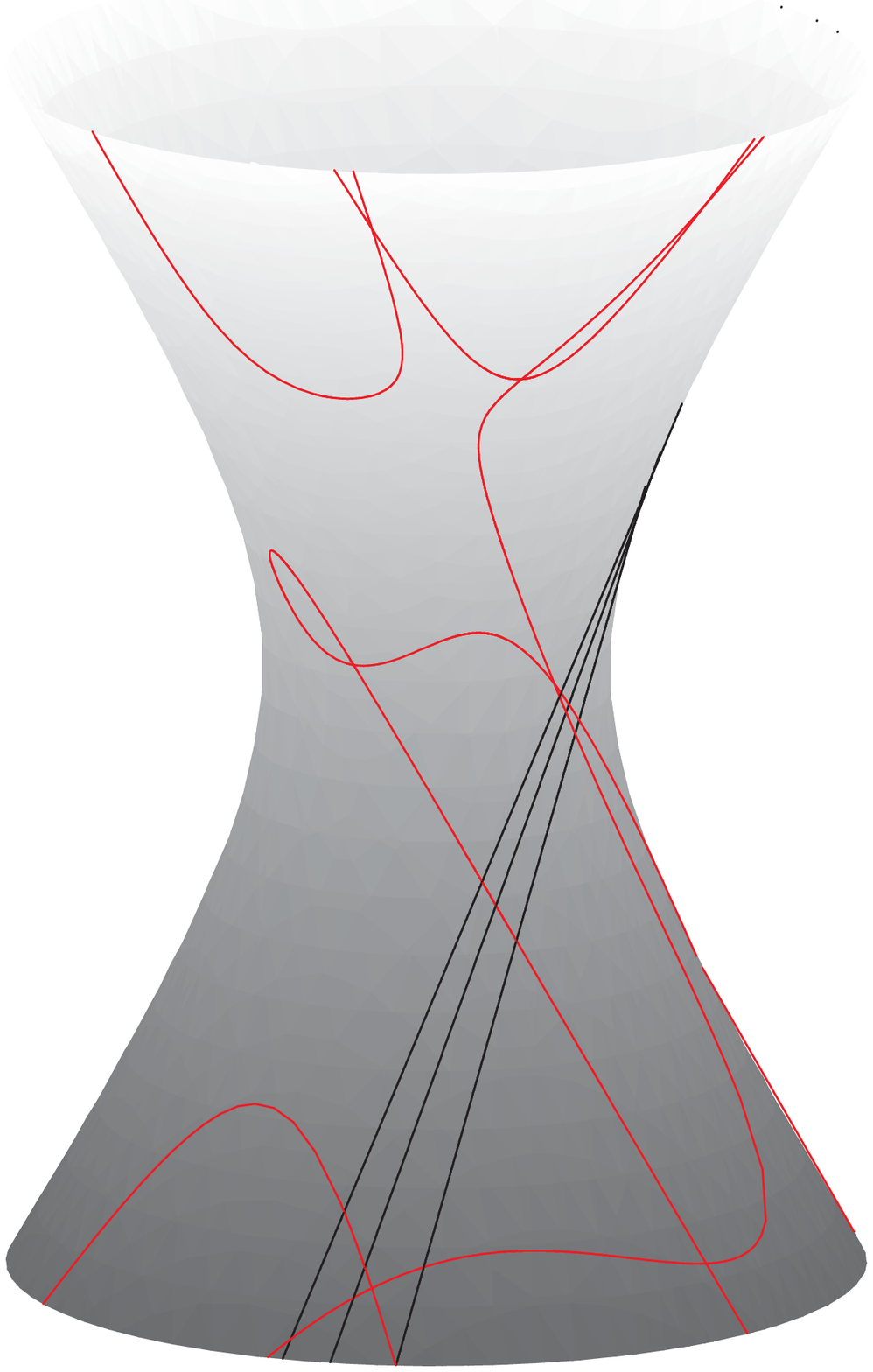}}
\]
}
\textnormal{{\normalsize \noindent We choose 7 points $p_{1},...,p_{7}$ on $\mathbb{P}%
^{1}\times\mathbb{P}^{1}$ in general position and a divisor $C^{\prime}$ of
type $(5,5)$ that passes all these points with multiplicity 2. The canonical
series is then cut out by the linear system $\left\vert (3,3)-%
{\textstyle\sum\nolimits_{i=1}^{7}}
p_{i}\right\vert $ and the normalisation $C$ has
at least two different $g_{5}^{1}.$ Calculating the Betti table for }$C$
gives{\normalsize
\[
\begin{tabular}{c|cccccccccc}
& 0 & 1 & 2 & 3 & 4 & 5 & 6 & 7 \cr \hline
0 & 1 & - & - & - & - & - & - & - \cr
1 & - & 21 & 64 & 70 & 8 & - & - & - \cr
2 & - & - & - & 8  & 70  & 64  & 21 & -\cr
3 & - & - & - & - & - & - & - & 1 \cr
\end{tabular} 
\]
Thus it follows that $C$ has exactly two different $g_{5}^{1}$$^{\prime}s$ of multiplicity
one and $\operatorname*{Cliff}(C)=3$ (cf. Section \ref{2g15res}). }}}
\end{example}
{\normalsize \noindent It will turn out later (cf. Theorem \ref{2g15}) that
for a general curve $C^{\prime}$ of this family there exist exactly these two $g_{5}%
^{1\prime}s$. }
\\To see that the property of having two $g_{5}^{1}$$^{\prime}s$ gives an independent condition in the \index{$\mathcal{M}_{9}$ modular space of curves of genus $g$} subscheme $\mathcal{H}_{(1,5)}\subset \mathcal{M}_{9}$ of all pentagonal curves of genus 9 \index{$\mathcal{M}_{(r,d)}$ curves with a $g_{d}^{r}$}, we count the dimension of the family of plane curves of degree 8 with 2 triple points and 6 double points as only singularities:
\begin{align*}
2\cdot8+\binom{10}{2}-1-2\cdot6-3\cdot6-\dim PGL(3)=22=\dim\mathcal{M}_{9}-2
\end{align*}
\begin{example}
{\normalsize (Canonical curves with exactly three different ordinary
$g_{5}^{1}$$^{\prime}s$) }{\textnormal{{\normalsize The question arises if it is possible for $C$ to have an
additional $g_{5}^{1}$ without admitting a $g_{7}^{2}.$ If we start with a
plane curve of degree }$8$ with $3$ triple points and $3$ double points as
only singularities, then these curves fail to have the right number of
$g_{5}^{1\prime}s$. The reason for this is that the conics passing through the three
triple points cut out a $g_{7}^{2}$ on $C$. However if there exists a curve with
exactly three $g_{5}^{1}$$^{\prime}s$ then we would expect that the condition of having a
third $g_{5}^{1}$ is of codimension one in the variety of all curves with two
$g_{5}^{1}$. {\normalsize Therefore we have examined several curves in this
family, which have been randomly constructed over a finite field
$\mathbb{F}_{p}, p\neq3.$ We observed that in approximately one out of $p$ examples the
Betti table of the minimal free resolution of $S_{C}$ looks like%
\[
\begin{tabular}{c|cccccccccc}
& 0 & 1 & 2 & 3 & 4 & 5 & 6 & 7 \cr \hline
0 & 1 & - & - & - & - & - & - & - \cr
1 & - & 21 & 64 & 70 & 12 & - & - & - \cr
2 & - & - & - & 12 & 70  & 64  & 21 & -\cr
3 & - & - & - & - & - & - & - & 1 \cr
\end{tabular} 
\]
A closer examination in these cases led us to the assumption that they exactly
occur if the following geometric situation on $\mathbb{P}^{1}\times
\mathbb{P}^{1}$ is given: From $\dim\left\vert (2,2)-%
{\textstyle\sum\nolimits_{i=1}^{7}}
p_{i}\right\vert =2,$ it follows the existence of a pencil $(D_{\lambda})_{\lambda}$ of
divisors of type $(2,2)$ on $\mathbb{P}^{1}\times\mathbb{P}^{1}$ passing through the
points $p_{1},...,p_{7}.$ Let $q$ denote the basepoint of the linear system
}$\left\vert (2,2)-%
{\textstyle\sum\nolimits_{i=1}^{7}}
p_{i}\right\vert ${\normalsize :}}

{\normalsize \bigskip}

{\normalsize \label{ConicsPic}\textit{Pencil of divisors of type }$(2,2)$
\textit{on }$\mathbb{P}^{1}\times\mathbb{P}^{1}$ \textit{passing through }%
 $p_{1},...,p_{7}$\textit{ }%
\[{\includegraphics[height=6.5cm,width=6.5cm]{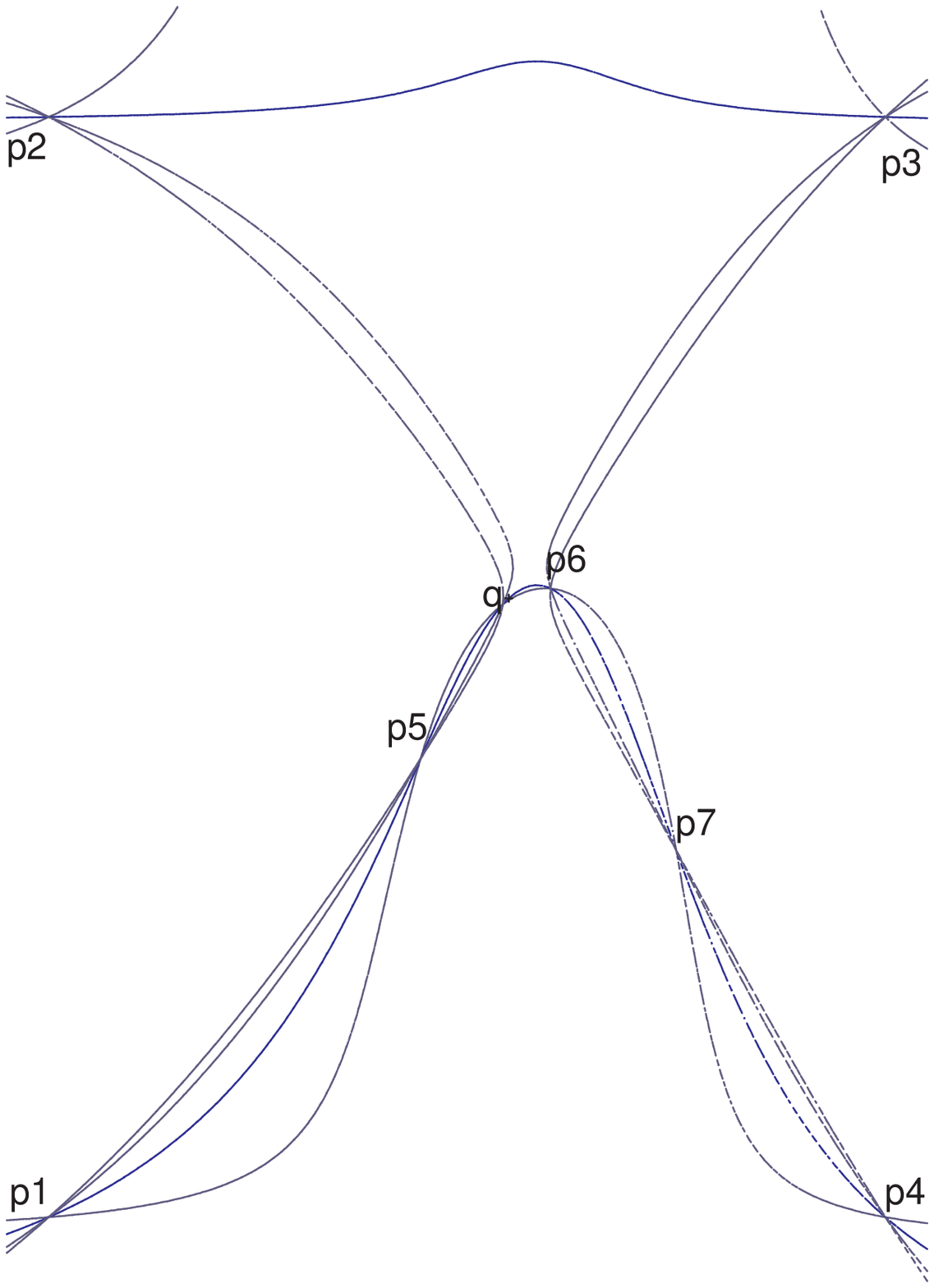}}%
\]
\textnormal{It turns out that $q\in C^{\prime}$ if and only if we get a further $g_{5}^{1}.$ This
$g_{5}^{1}$ is then cut out on $C$ by the pencil $(D_{\lambda})_{\lambda
}:\label{3dg15pic}$%
\[{\includegraphics[height=8.5cm,width=7.5cm]{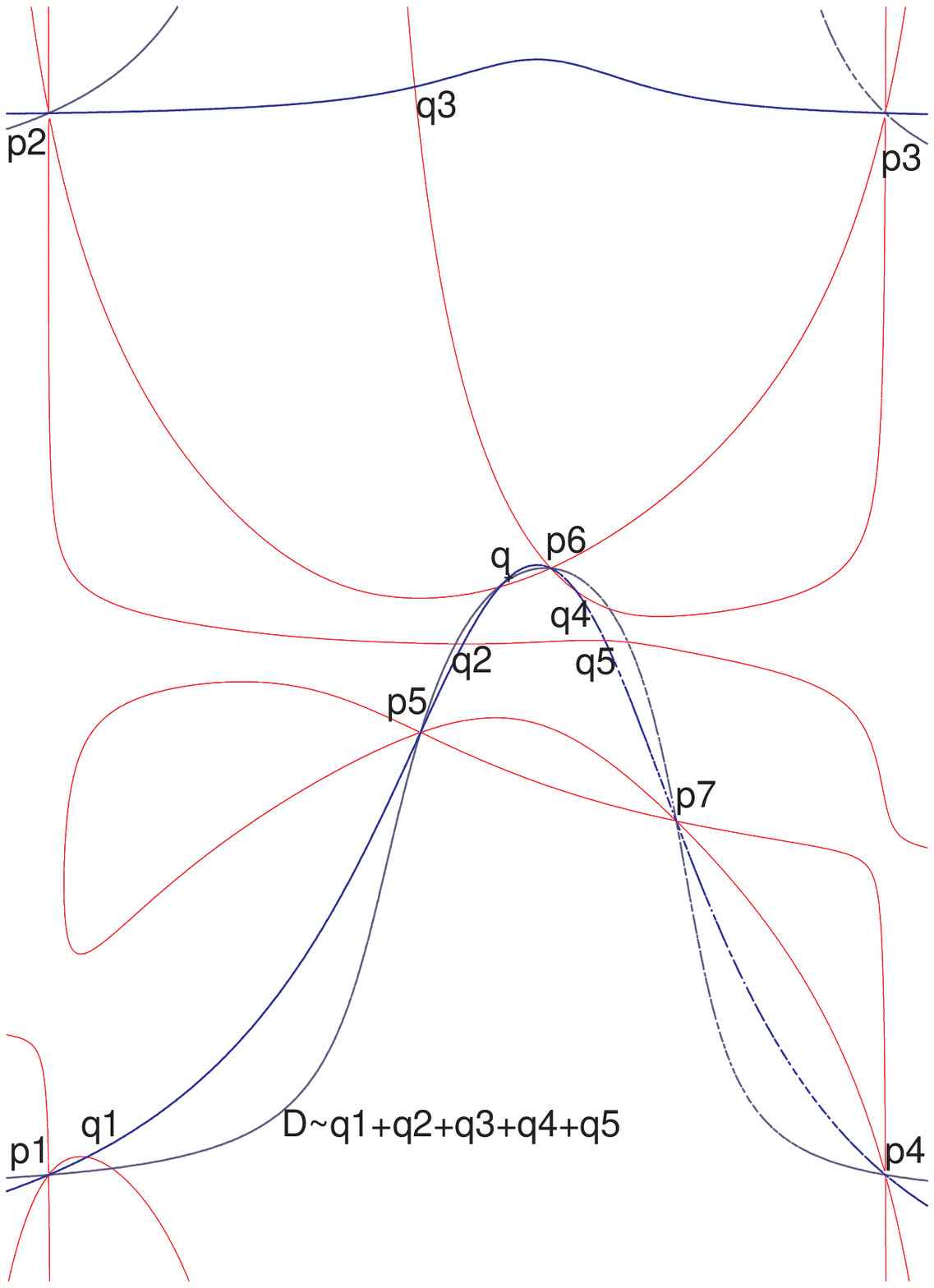}}%
\]
\textnormal{For the plane model }$C^{\prime\prime}\subset\mathbb{P}^{2}$ {\normalsize  of
}$C$ obtained from projection from one of the double points, this third
$g_{5}^{1}$ is cut out by the pencil of cubics passing through the triple points $s_{1},s_{2}$
and the double points $r_{1},...,r_{6}$ of $C^{\prime\prime}$. The base point
$p$ of $\left\vert 3H-s_{1}%
-s_{2}-r_{1}-...-r_{6}\right\vert $ then lies on $C$. Therefore the family of canonical curves of genus 9 that admit exactly three different $g_{5}^{1}$$^{\prime}s$ has dimension $\dim\mathcal{M}_{9}-3=21$.}}} 
\end{example}
{\normalsize From Theorem \ref{2g15} we know that there cannot exist a fourth
$g_{5}^{1}$ on }$C$ if $C$ {\normalsize has no }$g_{7}^{2}$ {\normalsize and
Theorem \ref{2g15b} states that this $g_{5}^{1}$ is different from the two
others if and only if
\[
\dim\left\vert (2,1)-%
{\textstyle\sum\nolimits_{i=1}^{7}}
p_{i}\right\vert ,\dim\left\vert (1,2)-%
{\textstyle\sum\nolimits_{i=1}^{7}}
p_{i}\right\vert =0
\]
i.e. there exists no divisor of class $(2,1)$ or $(1,2)$ passing through the points
$p_{1},...,p_{7}.$ }

{\normalsize \bigskip}

{\normalsize We will see in the following sections that it is possible that
some of the $g_{5}^{1}$$^{\prime}s$ become equal and that in these cases the multiplicity
of these linear systems takes the corresponding value. Furthermore the Betti table
stays the same as in the situation where all $g_{5}^{1}$ are different (for
details see Section \ref{deformation}). }

\bigskip

\begin{example}
{\normalsize (Canonical curve with a $g_{7}^{2}$) }{\textnormal{{\normalsize We complete this section with the most special case for a curve
}$C$ of $\operatorname*{Cliff}(C)=3 $ where {\normalsize  $C$ admits a
$g_{7}^{2}$. Here the construction is straightforward: }The $g_{7}^{2}%
$ {\normalsize  must have exactly 6 double points. \ We start with
$p_{1},...,p_{6}\in\mathbb{P}^{2}$ in general position and $C^{\prime
}\subset\mathbb{P}^{2}$ a curve of degree 7 passing through these points with
multiplicity 2. The canonical series of $C^{\prime}$ is then cut out by quartics
passing through the points $p_{1},...,p_{6}.$ Let $C\subset\mathbb{P}^{8}$ denote its
canonical image$.$ The Betti table }for $C$ {\normalsize  takes the following
form:%
\[
\begin{tabular}{c|cccccccccc}
& 0 & 1 & 2 & 3 & 4 & 5 & 6 & 7 \cr \hline
0 & 1 & - & - & - & - & - & - & - \cr
1 & - & 21 & 64 & 70 & 24 & - & - & - \cr
2 & - & - & - & 24  & 70  & 64  & 21 & -\cr
3 & - & - & - & - & - & - & - & 1 \cr
\end{tabular} 
\]
It follows from our results that $C$ has $\operatorname*{Cliff}(C)=3$. }
{\normalsize \noindent In special cases the double points can also lie
infinitely near. It can even happen that $C^{\prime}$ has an $A_{11}%
-$singularity as only singularity. Using a work of Lossen (\cite{Lossen}
Chapter 1) we constructed the following example $C^{\prime}$ given by the
affine equation$:$%
\[
8y^{2}-16x^{3}y-8x^{2}y^{3}+2xy^{5}-y^{7}+8x^{6}+8x^{5}y^{2}-8x^{3}y^{3}=0
\]
that has such a special singularity at the origin: }}}
\end{example}
\bigskip
{\normalsize \textit{Plane curve }$C$ \textit{of degree 7 with an }$A_{11}%
$-\textit{singularity at the origin} }
{\normalsize
\[{\includegraphics[height=10cm,width=8cm]%
{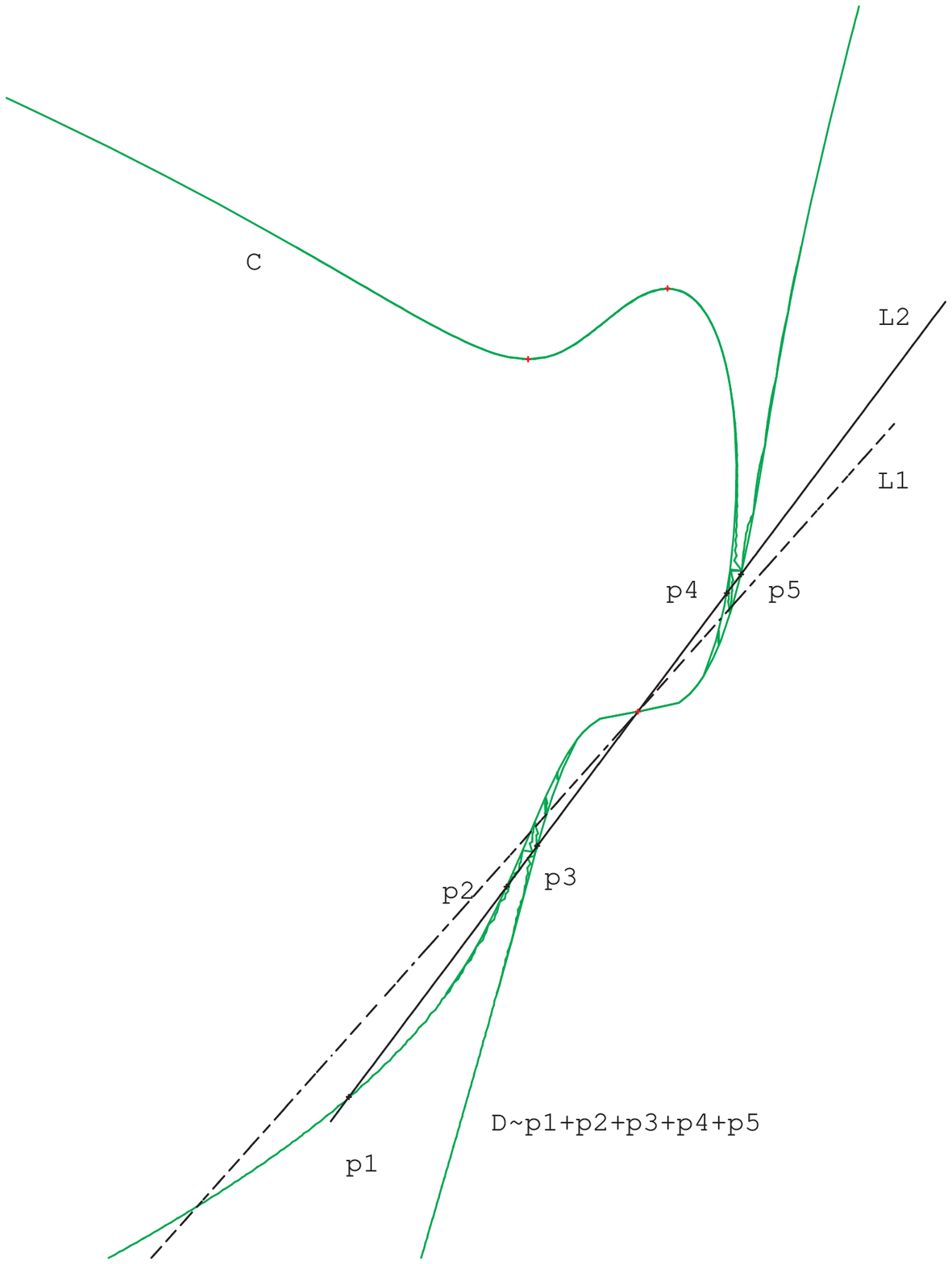}%
}%
\]
}

\section{The $g_{7}^{2}$ case}

{\normalsize Let $C$ be a canonical curve of genus 9 with
$\operatorname*{Cliff}(C)=3$ that admits a plane model $C^{\prime}%
\subset\mathbb{P}^{2}$ of degree 7. Then $C^{\prime}$ has exactly $6$ double
points $p_{1},...,p_{6}$ as singularities. We will show that all $g_{5}^{1}$$^{\prime}s$
on $C^{\prime}$ are given by projection from one of these double points: }

\begin{theorem}
{\normalsize \label{g27a}Let $C$ be a canonical curve with $\operatorname{Cliff}(C)=3$ that admits
a $g_{7}^{2}$. Then every $g_{5}^{1}$ is obtained from projection from one of the double points of the $g_{7}^{2}.$}
\end{theorem}

\begin{proof}
Let $H$ be an effective divisor of the $g_{7}^{2}$ and $|D|$ be a base point free linear series of type $g_{5}^{1}$, then $h^{0}(C,\mathcal{O}_{C}(H-D))=0$ or $1$. From the base point free pencil trick we obtain
\begin{align*}
h^{0}(C,\mathcal{O}_{C}(H+D))\geq h^{0}(C,\mathcal{O}_{C}(H))\cdot h^{0}(C,\mathcal{O}_{C}(D))-h^{0}(C,\mathcal{O}_{C}(H-D))=\\
=6-h^{0}(C,\mathcal{O}_{C}(H-D))
\end{align*}
Hence because of $\operatorname{Cliff}(C)=3$ we must have $h^{0}(C,\mathcal{O}_{C}(H+D))\leq5$, thus $h^{0}(C,\mathcal{O}_{C}(H-D))=1$. It follows the existence of points $q_{1}$ and $q_{2}$ on $C$ with
\[
H-D\sim q_{1}+q_{2}
\]
Therefore $|D|$ is obtained from projection from a double point.
\end{proof}

\bigskip

{\normalsize A dimension count shows that the subscheme $\mathcal{H}_{(2,7)}$
of all $C\in\mathcal{M}_{9}$ with a $g_{7}^{2}$ has the expected codimension
$-\varrho(9,7,2)=3$ in $\mathcal{M}_{9}$: We can choose 6 points in general position and a curve of
degree 7 passing through these points with multiplicity 2. Therefore we obtain a parameter space of dimension%
\[
2\cdot6+\binom{7+2}{2}-1-3\cdot6-\dim PGL(3)=21=\dim\mathcal{M}_{9}-3
\]
}
It is remarkable that the dimension of the stratum of all curves with a $g_{7}^{2}$ equals to that of the pentagonal curves that admit no $g_{7}^{2}$ but exactly three $g_{5}^{1}$$^{\prime}s$.

\begin{theorem}
{\normalsize \label{g27}Let $C\subset\mathbb{P}^{8}$ be a smooth, irreducible,
canonical curve of genus 9 that has a special linear series of type $g_{7}%
^{2}$ and $\operatorname*{Cliff}(C)=3$. Then $C$ has the
following Betti table: 
\textnormal{{\normalsize
\[
\begin{tabular}{c|cccccccccc}
& 0 & 1 & 2 & 3 & 4 & 5 & 6 & 7 \cr \hline
0 & 1 & - & - & - & - & - & - & - \cr
1 & - & 21 & 64 & 70 & 24 & - & - & - \cr
2 & - & - & - & 24  & 70  & 64  & 21 & -\cr
3 & - & - & - & - & - & - & - & 1 \cr
\end{tabular} 
\]}
}}
\end{theorem}

\begin{proof}
{\normalsize Let us consider the blowup of $\mathbb{P}^{2}$ in the singular
points of $C^{\prime}:$%
\[
\sigma:S\rightarrow\mathbb{P}^{2}%
\]
Further }$E_{i}.$ {\normalsize denotes the exceptional divisors and $H$ the class of a
hyperplane on $\mathbb{P}^{2}$ and by abuse of notation also its pullback to
}$S${\normalsize . Then we can assume that $C$ is the strict transform of $C^{\prime}$:%
\[
C\sim7H-%
{\textstyle\sum\nolimits_{i=1}^{6}}
2E_{i}%
\]
The adjoint series can be obtained by using the adjunction formula:%
\[
K_{S}+C\sim-3H+%
{\textstyle\sum\nolimits_{i=1}^{6}}
E_{i}+7H-%
{\textstyle\sum\nolimits_{i=1}^{6}}
2E_{i}\sim4H-%
{\textstyle\sum\nolimits_{i=1}^{6}}
E_{i}%
\]
We first want to focus on the image $S^{\prime}\subset\mathbb{P}^{8}$ of the
map $\varphi$\ defined by the adjoint series%
\[
\varphi:S\rightarrow S^{\prime}\subset\mathbb{P}^{8}%
\]
Let $R=\Bbbk\lbrack x_{0},...,x_{8}]$ and $R_{S^{\prime}}$ denotes the
homogenous coordinate ring of $\mathbb{P}^{8}$ and $S^{\prime}$ respectively.
From Corollary \ref{ALSP2} we already know that $S^{\prime}$ is arithmetically
Cohen-Macaulay, thus the projective dimension $R_{S^{\prime}}$ is $6=8-\dim
S^{\prime}$ and $\ \operatorname*{reg}R_{S^{\prime}}=2.$ Now we
consider the same approach as in the proof of Theorem \ref{cliff2}: The
Hilbert function $H_{R_{S^{\prime}}}$ of $R_{S^{\prime}}$ is given by the
Hilbert polynomial
\[
P_{S^{\prime}}(n)=\frac{1}{2}an^{2}+bn+c
\]
with $a=(K_{S}+C)^{2}=(4\sigma^{\ast}H-%
{\textstyle\sum\nolimits_{i=1}^{6}}
E_{i})^{2}=16-6=10$, $b=\frac{1}{2}(K_{S}+C)^{2}+1-g(K_{S}+C)=5+1-3=3$ and
$c=1.$ Therefore, we obtain%
\[
H_{R_{S^{\prime}}}(n)=P_{S^{\prime}}(n)=5n^{2}+3n+1\text{ for all }%
n\in\mathbb{N}%
\]
The minimal free resolution $F$ of $R_{S^{\prime}}$ over $R$ reduces to a minimal free resolution of $R_{S^{\prime}}^{\prime}:=R_{S^{\prime}}/(y_{1},y_{2},y_{3})R_{S^{\prime}}$ over $R/(y_{1}%
,y_{2},y_{3})R\cong R^{\prime}:=\Bbbk\lbrack x_{0}^{\prime},...,x_{5}^{\prime
}]$ with $(y_{1},y_{2},y_{3})$ being an $R_{S^{\prime}}$ sequence of linear
polynomials in $x_{0},...,x_{8}$. The Hilbertfunction of $R_{S^{\prime}%
}^{\prime}$ has values $(1,6,3)$ and $H_{R_{S^{\prime}}^{\prime}}(n)=0$
for $n\geq3.$ Tensoring $R_{S^{\prime}}^{\prime}$ with the Koszul complex of
length 6$:$%
\[
0\leftarrow\underset{M^{(0)}}{\underbrace{R_{S^{\prime}}^{\prime}}}%
\overset{\varphi_{1}}{\leftarrow}\underset{M^{(1)}}{\underbrace{R^{\prime
6}\otimes R_{S^{\prime}}^{\prime}}}\overset{\varphi_{2}}{\leftarrow}%
\underset{M^{(2)}}{\underbrace{R^{\prime15}\otimes R_{S^{\prime}}^{\prime}}%
}\leftarrow...\overset{\varphi_{6}}{\leftarrow}\underset{M^{(6)}}%
{\underbrace{R_{S^{\prime}}^{\prime}}}\leftarrow0
\]
and taking into acount the graduation we get the following format:%
\[
\text{\begin{xy}
\xymatrix{
M^{(0)} & M^{(1)} \ar[l] & M^{(2)} \ar[l] & M^{(3)} \ar[l] & M^{(4)} \ar[l] & M^{(5)} \ar[l] & M^{(6)} \ar[l] \\
1 & 6 \ar[ld] & 15 \ar[ld] & 20 \ar[ld] & 15 \ar[ld] & 6 \ar[ld] & 1 \ar[ld] \\
6 & 36 \ar[ld] & 90 \ar[ld] & 120 \ar[ld] & 90 \ar[ld] & 36 \ar[ld] & 6 \ar[ld] \\
3 & 18 & 45 & 60 & 45 & 18 & 3 \\
}
\end{xy}}%
\]
where the arrows stand for the maps $\varphi_{k}^{(l)}:M_{l}%
^{(k)}\rightarrow M_{l}^{(k-1)},$ which give a decomposition of $\varphi_{k}$
in the parts of degree $l$, and the numbers in the format are given by
$c_{kl}:=\dim M_{l}^{(k)}. $ As $S^{\prime}$ is not contained in any hyperplane
we get $\beta_{11}=0$ and therefore $\beta_{m,m}=0$ for $m=1,...,6.$ The dual
$F^{\ast}$ of $F$ is a free resolution of $\omega_{S^{\prime}}$ (up to a shift
of degrees) and therefore we obtain from Green's Linear Syzygy Theorem (cf.
\cite{eisenbud2} Theorem 7.1) that the length $n$ of the linear strand of
$F^{\ast}$ satisfies $n\leq\beta_{68}-1=3-1=2,$ hence $\beta_{m,m+2}=0$ for
$m=0,...,3.$ The linear strand of the minimal free resolution of
}$R_{S^{\prime}}$ is a subcomplex of the minimal free resolution of $R_{C},$
hence we must have {\normalsize $\beta_{56}=\beta_{67}=0$ as }%
$\operatorname*{Cliff}(C)=3. $ {\normalsize It follows that the Betti table }for
$S^{\prime}$ takes the following form{\normalsize
\[
\begin{tabular}{c|ccccccccc}
& 0 & 1 & 2 & 3 & 4 & 5 & 6 \cr \hline
0 & 1 & - & - & - & - & - & - \cr
1 & - & 18 & 52 & 60 & 24 & - & - \cr
2 & - & - & - & -  & 10  & 12 & 3 \cr
\end{tabular} 
\]
From the exact sequence
\[
0\rightarrow\mathcal{O}_{S^{\prime}}(-C)\rightarrow\mathcal{O}_{S^{\prime}%
}\rightarrow\mathcal{O}_{C}\rightarrow0
\]
and $\mathcal{O}_{S^{\prime}}(-C)\cong\omega_{S^{\prime}}$ the minimal free
resolution of $\mathcal{O}_{C}$ as $\mathcal{O}_{\mathbb{P}^{8}}$-module is
given as the mapping cone of the minimal free resolutions of $\mathcal{O}%
_{S^{\prime}}$ and $\omega_{S^{\prime}}$ as every map is minimal, thus we
obtain the Betti table for $C$ as claimed. }
\end{proof}

{\normalsize \bigskip}

{\normalsize In this section we have shown that for a canonical curve
$C\subset\mathbb{P}^{8}$ of genus $9$ and Clifford index 3 that admits a
$g_{7}^{2},$ every $g_{5}^{1}$ can be obtained from projection from one of the
doublepoints of the $g_{7}^{2}$. Therefore in general, if none of the doublepoints lies infinitely near to another, then there exist exactly 6 different $g_{5}^{1\prime
}s$ on $C$. The Betti table of $C$ is uniquely determined in this case:%
\[
\begin{tabular}{c|cccccccccc}
& 0 & 1 & 2 & 3 & 4 & 5 & 6 & 7 \cr \hline
0 & 1 & - & - & - & - & - & - & - \cr
1 & - & 21 & 64 & 70 & 24 & - & - & - \cr
2 & - & - & - & 24  & 70  & 64  & 21 & -\cr
3 & - & - & - & - & - & - & - & 1 \cr
\end{tabular} 
\]
In the following considerations we concentrate on those pentagonal canonical
curves that admit no $g_{7}^{2}.$ }

{\normalsize \bigskip}

{\normalsize \newpage}

\section{\label{curvesscrolls}Pentagonal Curves and
Scrolls}

{\normalsize Let $C$ be a pentagonal curve of genus $9$ that admits no
$g_{7}^{2}.$ The variety swept out by the linear spans of these divisors is a
$4$-dimensional rational normal scroll%
\[
X=%
{\displaystyle\bigcup\limits_{D\in g_{5}^{1}}}
\overline{D}\subset\mathbb{P}^{8}%
\]
At first we determine all different possible types $S(e_{1},e_{2},e_{3}%
,e_{4})$ of $X$ which can occur: These are $S(2,1,1,1),$ $S(3,1,1,0)$ and
$S(2,2,1,0).$ Then we resolve $\mathcal{O}_{C}$ as an $\mathcal{O}%
_{\mathbb{P}^{8}}$-module as described in Chapter 1 to get certain conditions
on $(a_{1},...,a_{5})$ in this resolution. The number of different $(a_{1},...,a_{5})$ can be reduced to the remaining three cases
$(a_{1},...,a_{5})=(2,1,1,1,1),~(2,2,1,1,0)$ and $(2,2,2,0,0).$ We will then only focus on the
case, where $X$ is of type $S(2,1,1,1).$ For $(a_{1},...,a_{5})=(2,2,2,0,0)$ or $(2,2,1,1,0)$, it turns out that $C$ has a linear system of type $g_{4}^{1}$ or
a plane model of degree $7$ respectively. The case $(a_{1},...,a_{5})=(2,1,1,1,1)$ will
then be discussed in full detail in Section \ref{highermulti}. We will show
that it is possible for $C$ to have one, two or even three different $g_{5}^{1}$$^{\prime}s$.
In Section \ref{deformation} it will turn out that $X\simeq S(2,2,1,0)$ or
$X\simeq S(3,1,1,0)$ occur as specializations of the general case, when
two or three different linear systems of type $g_{5}^{1}$ become infinitely near. }

{\normalsize \bigskip}

\begin{theorem}
{\normalsize \label{scrolltypes}Let C be a smooth, irreducible, canonical
curve of genus 9 with a base point free complete pencil $g_{5}^{1}=\left\vert
D\right\vert $ and $\operatorname{Cliff}(C)=3$. Then the $4$-dimensional rational normal scroll $X$ swept out
by the linear spans of these divisors is of type $S(2,1,1,1),$ $S(2,2,1,0)$ or
$S(3,1,1,0).$ }
\end{theorem}

\begin{proof}
{\normalsize We mention that $H|_{C}\sim K_{C}$ and $R|_{C}\sim D.$ According to Section \ref{scrolldetermination} the type $S(e_{1},e_{2},e_{3},e_{4})$ of
the scroll $X$ can be determined by considering the following partition of $9:$
\[
d_{i}=h^{0}(C,\mathcal{O}_{C}(K_{C}-iD))-h^{0}(C,\mathcal{O}_{C}%
(K_{C}-(i+1)D)),~i=0,...,3
\]
The numbers $e_{i}$ are given by%
\[
e_{i}=\#\{j|d_{j}\geq i\}-1
\]
Applying Riemann-Roch we get%
\[
h^{0}(C,\mathcal{O}_{C}(K_{C}-D))=8-\deg D+h^{0}(C,\mathcal{O}%
_{C}(D))=5
\]
and%
\[
h^{0}(C,\mathcal{O}_{C}(K_{C}-2D)-h^{0}(C,\mathcal{O}_{C}(2D))=9-1-\deg2D=-2
\]%
\[
\Rightarrow h^{0}(C,\mathcal{O}_{C}(K_{C}-2D))=h^{0}(C,\mathcal{O}%
_{C}(2D))-2\geq1
\]
because of $h^{0}(C,\mathcal{O}_{C}(D))=2.$
If $h^{0}(C,\mathcal{O}_{C}(K_{C}-2D))\geq3$ we would get a $g_{6}^{2}$ as
$\deg(K_{C}-2D)=6$, so that we can reduce to $1\leq h^{0}(C,\mathcal{O}%
_{C}(K_{C}-2D))\leq2.$ As $h^{0}(C,\mathcal{O}_{C}(K_{C}-2D))>h^{0}%
(C,\mathcal{O}_{C}(K_{C}-3D))$ we have $0\leq h^{0}(C,\mathcal{O}_{C}%
(K_{C}-3D))\leq1$ and $h^{0}(C,\mathcal{O}_{C}(K_{C}-kR))=0$ for $k\geq4$
because of $\deg(K_{C}-kD)<0.$ Therefore we conclude that there are exactly
the three different possible types for the scroll as given above. }
\end{proof}

{\normalsize \bigskip}

{\normalsize \noindent In following sections we will only focus on the case
where }$X$ is of tye $S(2,1,1,1)${\normalsize . According to Theorem
\ref{ResC} the resolution of $\mathcal{O}_{C}$ as $\mathcal{O}%
_{\mathbb{P(\mathcal{E})}}$-module is given by%
\[
F_{\ast}:\ \ \ \ \ \ 0\rightarrow\mathcal{O}_{\mathbb{P(\mathcal{E})}%
}(-5H+3R)\rightarrow%
{\displaystyle\sum\limits_{i=1}^{5}}
\mathcal{O}_{\mathbb{P(\mathcal{E})}}(-3H+b_{i}R)\overset{\psi}{\rightarrow}%
\]
}

{\normalsize
\[
\overset{\psi}{\rightarrow}%
{\displaystyle\sum\limits_{i=1}^{5}}
\mathcal{O}_{\mathbb{P(\mathcal{E})}}(-2H+a_{i}R)\rightarrow\mathcal{O}%
_{\mathbb{P(\mathcal{E})}}\rightarrow\mathcal{O}_{C}\rightarrow0
\]
}

{\normalsize \noindent We consider all possible values for $a_{i},b_{i}$ in
$F_{\ast}$: }

\begin{theorem}
{\normalsize \label{psitypes0}If $C$ is an irreducible, nonsingular, canonical curve of
genus 9 and $\operatorname{Cliff}(C)=3$ that admits no $g_{7}^{2}$ then the possible values $a_{1}%
,...,a_{5},$ $b_{1},...,b_{5}\in\mathbb{Z}$ in $F_{\ast}$ are:%
\[
(a_{1},...,a_{5})=(2,1,1,1,1),~(2,2,1,1,0)\text{ or }(2,2,2,0,0)
\]
}
\end{theorem}

\begin{proof}
{\normalsize Without loss of generality we assume $a_{1}\geq...\geq a_{5}.$
The calculation of the numbers $e_{i}$ shows that%
\[
\frac{2g-2}{5}\geq e_{1}\geq e_{2}\geq e_{3}\geq e_{4}\geq0
\]
}

{\normalsize \noindent and $f=e_{1}+e_{2}+e_{3}+e_{4}=5.$ The complex
$F_{\ast}$ is selfdual and $a_{i}+b_{i}=f-2=3$, $a_{1}+...+a_{5}=2g-12=6$ (cf.
Theorem \ref{ResC}). From the structure theorem for Gorenstein ideals in
codimension 3 (see \cite{buchsbaumeisenbud}) we obtain further
information from the complex $F_{\ast}$ above: The matrix $\psi$ is
skew-symmetric and its 5 Pfaffians generate the ideal of $C$ in
$\mathbb{P(\mathcal{E})}$, i.e. form the entries of%
\[%
{\displaystyle\sum\limits_{i=1}^{5}}
\mathcal{O}_{\mathbb{P(\mathcal{E})}}(-2H+a_{i}R)\rightarrow\mathcal{O}%
_{\mathbb{P(\mathcal{E})}}%
\]
Thus $C$ is determined by the entries of $\psi.$ If one of the nondiagonal
entries is zero, say $\psi_{45}=\psi_{54}=0$, then $C$ is contained in the
determinantal surface $Y$ defined by the matrix%
\[
\omega\sim\left(
\begin{array}
[c]{ccc}%
\psi_{14} & \psi_{24} & \psi_{34}\\
\psi_{15} & \psi_{25} & \psi_{35}%
\end{array}
\right)
\]
since in this case the $2\times2$ minors of that matrix are among the
Pfaffians of $\psi$. With $C$ also $Y$ is irreducible. A general fibre of
$Y\subset\mathbb{P(\mathcal{E})}$ over $\mathbb{P}^{1}$ is a twisted cubic.
Furthermore we see that none of the entries of this $2\times3$ matrix above can
be made zero by row and column operations if $C$ is irreducible. We distinguish the three cases $X\simeq S(2,1,1,1),$ $X\simeq
S(3,1,1,0)$ and $X\simeq S(2,2,1,0):$ }

{\normalsize For $X\simeq S(2,1,1,1),$ let us assume that $a_{1}\geq3$, then
the Pfaffian $\psi\in H^{0}(\mathbb{P(\mathcal{E})},\mathcal{O}%
_{\mathbb{P(\mathcal{E})}}(2H-a_{1}R))$ can be written as a sum of products of
global sections in $\mathcal{O}_{\mathbb{P(\mathcal{E})}}(H-R)$ and
$\mathcal{O}_{\mathbb{P(\mathcal{E})}}(H-2R).$ As each summand must contain
a factor in $H^{0}(\mathbb{P(\mathcal{E})},\mathcal{O}%
_{\mathbb{P(\mathcal{E})}}(H-2R))$ and this vector space is generated by
exactly one global section $\varphi_{0}$, $\varphi_{0}$ is a factor of the
Pfaffian, which contradicts that $C$ is irreducible. }

{\normalsize In the case $X\simeq S(3,1,1,0)$, a similar argument shows that
for $a_{1}\geq3$ the unique global section $\varphi_{0}\in H^{0}%
(\mathbb{P(\mathcal{E})},\mathcal{O}_{\mathbb{P(\mathcal{E})}}(H-3R))$ must be
a factor in the Pfaffian of type $H^{0}(\mathbb{P(\mathcal{E})},\mathcal{O}%
_{\mathbb{P(\mathcal{E})}}(2H-a_{1}R)).$ }

{\normalsize It remains to exclude the cases where $a_{5}\leq-1.$ This is
only possible for $(a_{1},...,a_{5})=(2,2,2,1,-1)$ or $(a_{1},...,a_{5}%
)=(2,2,2,2,-2)$. Therefore the last column of $\psi$ only
contains entries in $H^{0}(\mathbb{P(\mathcal{E})},\mathcal{O}%
_{\mathbb{P(\mathcal{E})}}(H-mR))$ with $m\geq2,$ hence at least two of these
entries can be made zero by suitable row and column operations as
$h^{0}(\mathbb{P(\mathcal{E})},\mathcal{O}_{\mathbb{P(\mathcal{E})}%
}(H-2R))\leq2.$ This contradicts that the Pfaffians are irreducible. }

{\normalsize Now we consider the case $X\simeq S(2,2,1,0):$ If $a_{5}\leq-1$
then setting $k=a_{4}-a_{5}\geq0$ we obtain
\begin{align*}
a_{5}+4(a_{5}+k)  &  =a_{5}+4a_{4}\leq%
{\textstyle\sum\limits_{i=1}^{5}}
a_{i}=6<7\leq4-3a_{5}=4(1-a_{5})+a_{5}\\
&  \Rightarrow a_{4}=a_{5}+k<1-a_{5}=(3-a_{5})-2\Rightarrow b_{5}-a_{4}\leq3
\end{align*}
and therefore the entry $\psi_{45}$ of the matrix $\psi$ vanishs. It follows
that the Pfaffians of $\psi$ include the $2\times2-$minors of a $2\times
3-$matrix $\omega$ with entries indicated as below:%

\small
\[
\left(
\begin{array}
[c]{ccc}%
H-(3-(a_{5}+a_{1}))R & H-(3-(a_{5}+a_{2}))R & H-(3-(a_{5}+a_{3}))R\\
H-(3+k-(a_{5}+a_{1}))R & H-(3+k-(a_{5}+a_{2}))R & H-(3+k-(a_{5}+a_{3}))R
\end{array}
\right)
\]

\normalsize
Then $C$ is contained in a determinantal surface $Y\subset\mathbb{P}^{8}$
given by the minors of $\omega$ on $\mathbb{P(\mathcal{E})}$. If we apply
Theorem \ref{ConPk} we see that $Y$ is the image of $P_{k}%
:=\mathbb{P(\mathcal{O}}_{\mathbb{P}^{1}}(k)\oplus\mathbb{\mathcal{O}%
}_{\mathbb{P}^{1}}\mathbb{)}$ under a rational map defined by a subseries of%
\[
H^{0}(P_{k},\mathcal{O}_{P_{k}}(3A+(5-3k-a)B))
\]
with hyperplane class $A$ and ruling $B$ on $P_{k}$ where $a:=3-(a_{5}%
+a_{1})+3-(a_{5}+a_{2})+3-(a_{5}+a_{3})=9-3a_{5}-(6-2a_{5}-k)=3+k-a_{5}.$ Then,
the strict transform $C^{\prime}$ of $C$ in $P_{k}$ is a divisor of class%
\[
C^{\prime}\sim5A+(7-4k-a)B=5A+(4+a_{5}-5k)B
\]
It follows that the hyperplanes on $P_{k}$ cut out a $g_{4+a_{5}}^{k+1}$ on
$C^{\prime}$ which has Clifford index less than 2, a contradiction. For
$a_{4}=a_{5}=0$ the same argument as above leads to the existence of a
$g_{4}^{1}$ on $C.$ Thus, for $(a_{1},...,a_{5})$ with $a_{1}\geq3,$ there is
only the possibility $(a_{1},...,a_{5})=(3,1,1,1,0)$ left. In this case, the
matrix $\psi$ is of type%
\begin{align*}
\psi &  \sim\left(
\begin{array}
[c]{ccccc}%
0 & H+R & H+R & H+R & H\\
& 0 & H-R & H-R & H-2R\\
&  & 0 & H-R & H-2R\\
&  &  & 0 & H-2R\\
&  &  &  & 0
\end{array}
\right)  \sim\\
&  \sim\left(
\begin{array}
[c]{ccccc}%
0 & H+R & H+R & \mathbf{H+R} & \mathbf{H}\\
& 0 & H-R & \mathbf{H-R} & \mathbf{H-2R}\\
&  & 0 & \mathbf{H-R} & \mathbf{H-2R}\\
&  &  & 0 & 0\\
&  &  &  & 0
\end{array}
\right)
\end{align*}
hence $C$ is contained in a determinantal surface $Y$ given by the $2\times2$
minors of
\[
\omega\sim\left(
\begin{array}
[c]{ccc}%
H+R & H-R & H-R\\
H & H-2R & H-2R
\end{array}
\right)
\]
on $\mathbb{P(\mathcal{E})\simeq}S(2,2,1,0).$ $Y$ is the image of
$P_{1}=\mathbb{P(\mathcal{O}}_{\mathbb{P}^{1}}\mathbb{\mathcal{(}%
}1)\mathbb{\mathcal{\oplus O}}_{\mathbb{P}^{1}}\mathbb{)}$ under a rational
map defined by a subseries of%
\[
H^{0}(P_{1},\mathcal{O}_{P_{1}}(3A+B))
\]
and the image $C^{\prime}$ of $C$ in $P_{1}$ is a divisor of class%
\[
C^{\prime}\sim5A+2B
\]
(cf. Theorem \ref{ConPk}) from which we deduce the existence of a $g_{7}^{2}$
cut out by the hyperplane sections. }

{\normalsize Then the given possibilities for $(a_{1},...,a_{5})$ remain. }
\end{proof}

{\normalsize \bigskip}

{\normalsize For the rest of this section we will examine the two cases
$(a_{1},...,a_{5})=(2,2,2,0,0)$ and $(a_{1},...,a_{5})=(2,2,1,1,0),$ where it
will turn out that either we get a linear system of type $g_{4}^{1}$ or a
plane model of $C$ of degree 7 respectively. }

\begin{theorem}
{\normalsize Let C be a curve given by the Pfaffians of a matrix $\psi$ with
entries on a scroll of type $S(2,1,1,1)$ as above. Then }
\\\\
{\noindent \normalsize a) For $(a_{1},...,a_{5})=(2,2,2,0,0)$ the curve C has Clifford
index 2. There exists a special linear series of type $g_{4}^{1}$ and the
Betti table for $C$ has the following form:%
\textnormal{\[
\begin{tabular}{c|cccccccccc}
& 0 & 1 & 2 & 3 & 4 & 5 & 6 & 7 \cr \hline
0 & 1 & - & - & - & - & - & - & - \cr
1 & - & 21 & 64 & 70 & 44 & 5 & - & - \cr
2 & - & - & 5 & 44  & 70  & 64  & 21 & -\cr
3 & - & - & - & - & - & - & - & 1 \cr
\end{tabular} 
\]}}
\\\\
{\noindent \normalsize b) For $(a_{1},...,a_{5})=(2,2,1,1,0)$ the curve $C$ has Clifford
index 3 and there exists a plane model of degree 7. The Betti table for $C$ is
given by:%
\textnormal{\[
\begin{tabular}{c|cccccccccc}
& 0 & 1 & 2 & 3 & 4 & 5 & 6 & 7 \cr \hline
0 & 1 & - & - & - & - & - & - & - \cr
1 & - & 21 & 64 & 70 & 24 & - & - & - \cr
2 & - & - & - & 24  & 70  & 64  & 21 & -\cr
3 & - & - & - & - & - & - & - & 1 \cr
\end{tabular} 
\]}
}
\end{theorem}

\begin{proof}
{\normalsize a) For $(a_{1},...,a_{5})=(2,2,2,0,0)$ the matrix $\psi$ in the
resolution $F_{\ast}$ of $\mathcal{O}_{C}$ on the corresponding $\mathbb{P}%
^{3}-$bundle $\mathbb{P(\mathcal{E})}$ is of type%
\begin{align*}
\psi &  \sim\left(
\begin{array}
[c]{ccccc}%
0 & H+R & H+R & H-R & H-R\\
& 0 & H+R & H-R & H-R\\
&  & 0 & H-R & H-R\\
&  &  & 0 & H-3R\\
&  &  &  & 0
\end{array}
\right)  \sim\\
&  \sim\left(
\begin{array}
[c]{ccccc}%
0 & H+R & H+R & \mathbf{H-R} & \mathbf{H-R}\\
& 0 & H+R & \mathbf{H-R} & \mathbf{H-R}\\
&  & 0 & \mathbf{H-R} & \mathbf{H-R}\\
&  &  & 0 & 0\\
&  &  &  & 0
\end{array}
\right)
\end{align*}
Applying Theorem \ref{ConPk} again, we see that $C$ is contained in a
determinantal surface $Y,$ given by the minors of the $2\times3$ matrix
\[
\omega=\left(
\begin{array}
[c]{ccc}%
\psi_{14} & \psi_{24} & \psi_{34}\\
\psi_{15} & \psi_{25} & \psi_{35}%
\end{array}
\right)  \sim\left(
\begin{array}
[c]{ccc}%
H-R & H-R & H-R\\
H-R & H-R & H-R
\end{array}
\right)
\]
with $(H-R)$-entries on $\mathbb{P(\mathcal{E})}$. Therefore $Y$ is a blowup of
$\mathbb{P}^{1}\times\mathbb{P}^{1}$ in 3 points. The strict transform
$C^{\prime}$ of $C$ in $Y$ is a divisor of type $(5,4)$, thus we get a
$g_{4}^{1}$ from projection onto the second factor of $\mathbb{P}^{1}%
\times\mathbb{P}^{1}.$ From Theorem \ref{cliff2}, we expect for $C$ to admit the following Betti table:%
\[
\begin{tabular}{c|cccccccccc}
& 0 & 1 & 2 & 3 & 4 & 5 & 6 & 7 \cr \hline
0 & 1 & - & - & - & - & - & - & - \cr
1 & - & 21 & 64 & 70 & 44 & 5 & - & - \cr
2 & - & - & 5 & 44  & 70  & 64  & 21 & -\cr
3 & - & - & - & - & - & - & - & 1 \cr
\end{tabular} 
\]
Considering the corresponding mapping cone construction (cf. Theorem
\ref{ResC}) and calculating the ranks of the non minimal maps leads to the
same Betti table, thus it follows that $C$ has Clifford index
$\operatorname*{Cliff}(C)=2,$ and $C$ admits a $g_{4}^{1}\times g_{5}^{1}$,
but no $g_{6}^{2}$ or $g_{8}^{3}.$ }

b) In the case $(a_{1},...,a_{5})=(2,2,1,1,0)$ the matrix $\psi$
has the following form:%

\begin{align*}
\psi &  \sim\left(
\footnotesize{\begin{array}
[c]{ccccc}%
0 & H+R & H & H & H-R\\
& 0 & H & H & H-R\\
&  & 0 & H-R & H-2R\\
&  &  & 0 & H-2R\\
&  &  &  & 0
\end{array}}
\right)  \sim\left(
\footnotesize{\begin{array}
[c]{ccccc}%
0 & H+R & H & \mathbf{H} & \mathbf{H-R}\\
& 0 & H & \mathbf{H} & \mathbf{H-R}\\
&  & 0 & \mathbf{H-R} & \mathbf{H-2R}\\
&  &  & 0 & 0\\
&  &  &  & 0
\end{array}}
\right) \\
\end{align*}

As one of the $H-2R$ entries can be made to zero, we can assume that $\psi
_{45}=0.$ It follows from Theorem \ref{ConPk} that $C$ is contained in a
determinantal surface $Y$ given by the $2\times2$ minors of the matrix%
\[
\left(
\begin{array}
[c]{ccc}%
\psi_{14} & \psi_{24} & \psi_{34}\\
\psi_{15} & \psi_{25} & \psi_{35}%
\end{array}
\right)  \sim\left(
\begin{array}
[c]{ccc}%
H & H & H-R\\
H-R & H-R & H-2R
\end{array}
\right)
\]
on $\mathbb{P(\mathcal{E})}$. $Y$ is a blowup of $\mathbb{P}^{2}$ in 6 points.
The strict transform $C^{\prime}$ of $C$ in $Y$ is a divisor of type
$C^{\prime}\sim5A+2B,$ hence there exists a $g_{7}^{2}$ cut out by the class
of a hyperplane divisor $A.$ Calculating the ranks of the non minimal maps in
the corresponding mapping cone, it turns out that the Betti table of the
minimal free resolution of $\mathcal{O}_{C}$ is the same as determined in
Theorem \ref{g27}, hence $C$ has Clifford index 3. 
\end{proof}

{\normalsize \bigskip}

{\normalsize So far we have seen that in the case, where $C$ is a smooth,
irreducible, canonical curve of Clifford index 3 that is contained in a scroll of type
$S(2,1,1,1),$ for any possible configuration $(a_{1},...,a_{5})\neq
(2,1,1,1,1)$ the curve $C$ has a plane model of degree 7. In the next section we will concentrate
on the "general case" $(a_{1},...,a_{5})=(2,1,1,1,1).$ We will show that it is
possible for $C$ to have exactly one, two or even three different linear series of
type $g_{5}^{1}$ and that they correspond to syzygies in the minimal free
resolution of $\mathcal{O}_{C}.$ }

{\normalsize \newpage}

\section{Curves with ordinary $g_{5}^{1}$}

{\normalsize In this section we will formulate our main theorems. We have seen
that for a pentagonal curve $C\subset X$ with $\operatorname*{Cliff}(C)=3$,
that admits no $g_{7}^{2}$ and is contained in a scroll $X$ of type $S(2,1,1,1)$, there exists
exactly one type for the matrix $\psi,$ whose Pfaffians generate the vanishing
ideal of $C$ on the scroll $X.$ The following theorems give the correspondence
between the possible number of different $g_{5}^{1\prime}s$ and the Betti
numbers for $C\subset\mathbb{P}^{8}.$ }

{\normalsize \bigskip}

\subsection{Types of $\psi$ and non minimal maps}

{\normalsize For $(a_{1},...,a_{5})=(2,1,1,1,1),$ the resolution of $C$ on the
associated $\mathbb{P}^{3}-$bundle $\mathbb{P(\mathcal{E})}$ is of type
\[
F_{\ast}:~~~~0\rightarrow\mathcal{O}_{\mathbb{P(\mathcal{E})}}%
(-5H+3R)\rightarrow\mathcal{O}_{\mathbb{P(\mathcal{E})}}(-3H+R)\oplus
\mathcal{O}_{\mathbb{P(\mathcal{E})}}(-3H+2R)^{\oplus4}\overset{\psi
}{\rightarrow}%
\]%
\[
\overset{\psi}{\rightarrow}\mathcal{O}_{\mathbb{P(\mathcal{E})}}%
(-2H+2R)\oplus\mathcal{O}_{\mathbb{P(\mathcal{E})}}(-3H+R)^{\oplus
4}\rightarrow\mathcal{O}_{\mathbb{P(\mathcal{E})}}\rightarrow\mathcal{O}%
_{C}\rightarrow0
\]
with $\psi$ a skew-symmetric matrix with entries as indicated below: }

{\normalsize \bigskip}

{\normalsize \textbf{(*)}%
\begin{align*}
\psi &  \sim\left(
\begin{array}
[c]{ccccc}%
0 & H & H & H & H\\
& \mathbf{0} & \mathbf{H-R} & \mathbf{H-R} & \mathbf{H-R}\\
&  & \mathbf{0} & \mathbf{H-R} & \mathbf{H-R}\\
&  &  & \mathbf{0} & \mathbf{H-R}\\
&  &  &  & \mathbf{0}
\end{array}
\right) \\
&
\end{align*}
We have already mentioned that the vanishing ideal of $C\subset\mathbb{P(}%
\mathcal{E)}$ is given by the Pfaffians of $\psi.$ As we assumed that $C$ is
irreducible, at most one entry in each row and column of $\psi$ can be made to
zero by suitable row and column operations. Especially if one of the $(H-R)$-entries is zero, we can assume $\psi_{45}=0,$ then none of the remaining ones,
except $\psi_{23}$, can be made zero. In the case where $\psi_{45}=\psi_{23}=0$
the global sections $\psi_{24},\psi_{25},\psi_{34,}\psi_{35}\in H^{0}%
(\mathbb{P}(\mathcal{E}),\mathcal{O}_{\mathbb{P(\mathcal{E})}}(H-R))$ are
linear independent as otherwise $C$ contains a proper one dimensional component,
given by the vanishing locus of these sections,
hence $C$ cannot be irreducible. In the following theorems we will show that
the $H-R$ entries that can be made to zero exactly correspond to different
additional linear systems of type $g_{5}^{1}$. It turns out that up to conjugation there occur 4
different types for the $4\times4-$submatrix $\widetilde{\psi}$ of $\psi$ with
entries in $H^{0}(\mathbb{P}(\mathcal{E}),\mathcal{O}_{\mathbb{P(\mathcal{E}%
)}}(H-R))$.\newpage}

\begin{lemma}
{\normalsize \label{psitypes}Let C be an irreducible smooth pentagonal curve
$C\mathcal{\ }$of genus 9 that is contained in a scroll of type $S(2,1,1,1)$. If $C$ is
given by the Pfaffians of a matrix $\psi$ with entries as in \textbf{(*)} then
$\psi$ is conjugated to one of the following types:\medskip}

{\normalsize A
\[
\psi\sim\left(
\begin{array}
[c]{ccccc}%
0 & H & H & H & H\\
& 0 & f_{1} & f_{2} & f_{3}\\
&  & 0 & f_{4} & f_{5}\\
&  &  & 0 & f_{1}\\
&  &  &  & 0
\end{array}
\right)
\]
\medskip}

{\normalsize B%
\[
\psi\sim\left(
\begin{array}
[c]{ccccc}%
0 & H & H & H & H\\
& 0 & f_{1} & f_{2} & f_{3}\\
&  & 0 & f_{4} & f_{5}\\
&  &  & 0 & 0\\
&  &  &  & 0
\end{array}
\right)
\]
\medskip}

{\normalsize C%
\[
\psi\sim\left(
\begin{array}
[c]{ccccc}%
0 & H & H & H & H\\
& 0 & 0 & f_{2} & f_{3}\\
&  & 0 & f_{4} & f_{5}\\
&  &  & 0 & 0\\
&  &  &  & 0
\end{array}
\right)
\]
}

{\normalsize D%
\[
\psi\sim\left(
\begin{array}
[c]{ccccc}%
0 & H & H & H & H\\
& 0 & f_{1} & f_{2} & f_{3}\\
&  & 0 & f_{4} & f_{2}\\
&  &  & 0 & 0\\
&  &  &  & 0
\end{array}
\right)
\]
with linear independent $f_{1},...,f_{5}\in H^{0}(\mathbb{P}(\mathcal{E}%
),\mathcal{O}_{\mathbb{P(\mathcal{E})}}(H-R))$. }
\end{lemma}

\begin{proof}
{\normalsize As we have already mentioned, the case%
\[
\psi\sim\left(
\begin{array}
[c]{ccccc}%
0 & H & H & H & H\\
& 0 & 0 & f_{2} & f_{3}\\
&  & 0 & f_{4} & f_{2}\\
&  &  & 0 & 0\\
&  &  &  & 0
\end{array}
\right)
\]
cannot occur since in this situation $C$ contains a proper one dimensional
component on $\mathbb{P(}\mathcal{E)}$ given by the vanishing locus of the
sections $f_{2},...,f_{4}\in H^{0}(\mathbb{P}(\mathcal{E}),\mathcal{O}%
_{\mathbb{P(\mathcal{E})}}(H-R))$, so $C$ cannot be irreducible. It remains
to remark that D is a special case of B, where $f_{2},...,f_{5}$ are linear
dependent and $f_{1}\notin\left\langle f_{2},...,f_{5}\right\rangle $.
Therefore, we get%
\[
\psi\sim\left(
\begin{array}
[c]{ccccc}%
0 & H & H & H & H\\
& 0 & f_{1} & f_{2} & f_{3}\\
&  & 0 & f_{4} & f_{5}\\
&  &  & 0 & 0\\
&  &  &  & 0
\end{array}
\right)
\]
with one of the entries $f_{2},...,f_{5}$ being a linear combination of the
others. It follows, that $\psi$ has the form as given above. }
\end{proof}

{\normalsize \bigskip}

{\normalsize To obtain the minimal free resolution of $C\subset\mathbb{P}^{8},$
we have to determine the rank of the only non-minimal map in the
corresponding mapping cone $(MC)$ of $F_{\ast}$ which arise from the
$4\times4-$submatrix $\widetilde{\psi}$ of $\psi$:
\begin{align*}
& \\
&
\begin{array}
[c]{ccc}%
\mathcal{O}_{\mathbb{P(\mathcal{E})}}(-3H+2R)^{\oplus4} & ^{\underrightarrow
{\text{
\ \ \ \ \ \ \ \ \ \ \ \ \ \ \ \ \ \ \ \ \ \ \ \ \ \ \ \ \ \ \ \ \ \ \ \ \ }}}
& \mathcal{O}_{\mathbb{P(\mathcal{E})}}(-2H+R)^{\oplus4}\\
\uparrow &  & \uparrow\\
----|---- & --------- & ----|----\\
&  & \\
S_{2}G(-3)^{\oplus4} & ^{\underrightarrow{\text{
\ \ \ \ \ \ \ \ \ \ \ \ \ \ \ \ \ \ \ \ \ \ \ \ \ \ \ \ \ \ \ \ \ \ \ \ \ }}}
& G(-2)^{\oplus4}\\
\uparrow &  & \uparrow\\
F\otimes G(-4)^{\oplus4} & ^{\underrightarrow{\text{
\ \ \ \ \ \ \ \ \ \ \ \ \ \ \ \ \ \ \ \ \ \ \ \ \ \ \ \ \ \ \ \ \ \ \ \ \ }}}
& F(-3)^{\oplus4}\\
\uparrow &  & \uparrow\\
\wedge^{2}F(-5)^{\oplus4} & ^{\underrightarrow{\text{
\ \ \ \ \ \ \ \ \ \ \ \ \ \ \ \ \ \ }\alpha
\text{\ \ \ \ \ \ \ \ \ \ \ \ \ \ \ \ \ \ \ }}} & \wedge^{3}F(-5)^{\oplus4}\\
\uparrow &  & \uparrow\\
\wedge^{4}F(-7)^{\oplus4} & ^{\underrightarrow{\text{
\ \ \ \ \ \ \ \ \ \ \ \ \ \ \ \ \ \ \ \ \ \ \ \ \ \ \ \ \ \ \ \ \ \ \ \ \ }}}
& \wedge^{4}F\otimes DG^{\ast}(-6)^{\oplus4}\\
\uparrow &  & \uparrow\\
DG^{\ast}(-8)^{\oplus4} & ^{\underrightarrow{\text{
\ \ \ \ \ \ \ \ \ \ \ \ \ \ \ \ \ \ \ \ \ \ \ \ \ \ \ \ \ \ \ \ \ \ \ \ \ }}}
& D_{2}G^{\ast}(-7)^{\oplus4}%
\end{array}
\\
&
\end{align*}
Identifying $\psi_{ij}\in H^{0}(\mathbb{P}(\mathcal{E}),\mathcal{O}%
_{\mathbb{P(\mathcal{E})}}(H-R))$ $\cong H^{0}(\mathbb{P}^{8},F),$
$i,j=2,...,5,$ and $H^{0}(\mathbb{P}(\mathcal{E}),\mathcal{O}%
_{\mathbb{P(\mathcal{E})}}(R))$ $\cong H^{0}(\mathbb{P}^{8},G)$ (cf. page \pageref{ident}),  $\alpha$ is given by the wedge product with the corresponding
section:%
\[
\alpha:\wedge^{2}F^{\oplus4\underrightarrow{\text{ \ \ \ \ }\wedge
\widetilde{\psi}\text{\ \ \ \ \ \ \ \ }}}\wedge^{3}F^{\oplus4}%
\]
The following lemma answers the question of how the rank of $\alpha$ depends on the
type of $\psi$ in Lemma \ref{psitypes}: }

\begin{lemma}
{\normalsize \label{rank}According to the 4 different types of $\psi$ in Lemma
\ref{psitypes}$,$ we get:\medskip}

{\normalsize a) $\operatorname*{rank}~\alpha=40\Leftrightarrow\psi\sim A$ }

{\normalsize b) $\operatorname*{rank}~\alpha=36\Leftrightarrow\psi\sim B$ }

{\normalsize c) $\operatorname*{rank}~\alpha=32\Leftrightarrow\psi\sim C~or~\psi\sim D$ }

\noindent If $C$ is given by the Pfaffians of a matrix $\psi$, then the Betti table for $C$ has the following form
\textnormal{
\[
\begin{tabular}{c|cccccccccc}
& 0 & 1 & 2 & 3 & 4 & 5 & 6 & 7 \cr \hline
0 & 1 & - & - & - & - & - & - & - \cr
1 & - & 21 & 64 & 70 & a & - & - & - \cr
2 & - & - & - & a & 70  & 64  & 21 & -\cr
3 & - & - & - & - & - & - & - & 1 \cr
\end{tabular} 
\]}
with $a=44-\operatorname*{rank}~\alpha$
\end{lemma}

\begin{proof}
{\normalsize a) The matrix $\psi$ has the following form:
\[
\psi\sim\left(
\begin{array}
[c]{ccccc}%
0 & H & H & H & H\\
& 0 & f_{1} & f_{2} & f_{3}\\
&  & 0 & f_{4} & f_{5}\\
&  &  & 0 & f_{1}\\
&  &  &  & 0
\end{array}
\right)
\]
with linear independent entries $f_{1},...,f_{5}\in H^{0}(\mathbb{P}%
(\mathcal{E}),\mathcal{O}_{\mathbb{P(\mathcal{E})}}(H-R))$ $\cong
H^{0}(\mathbb{P}^{8},F)$. An easy calculation (c.f. Appendix 6.1) shows, that $\ker\alpha=0.$ }

{\normalsize \noindent b) A similar calculation (c.f. Appendix 6.1) as in a) shows that the kernel
of $\alpha$ is given by:%
\small{\[
\ker\alpha=\left\langle \left(
\begin{array}
[c]{c}%
0\\
0\\
f_{2}\wedge f_{4}\\
0
\end{array}
\right)  ,\left(
\begin{array}
[c]{c}%
0\\
0\\
0\\
f_{3}\wedge f_{5}%
\end{array}
\right)  ,\left(
\begin{array}
[c]{c}%
0\\
0\\
f_{3}\wedge f_{5}\\
f_{2}\wedge f_{5}+f_{3}\wedge f_{4}%
\end{array}
\right)  ,\left(
\begin{array}
[c]{c}%
0\\
0\\
f_{2}\wedge f_{5}+f_{3}\wedge f_{4}\\
f_{2}\wedge f_{4}%
\end{array}
\right)  \right\rangle
\]}
}
i.e. $\operatorname*{rank}~\alpha=36.$

\noindent {\normalsize c) For $\psi\sim D$, $\ker\alpha$ is generated by 8 elements (cf. Appendix 6.1):%
\[
\ker\alpha=\left\langle
\begin{array}
[c]{c}%
\left(
\begin{array}
[c]{c}%
0\\
0\\
f_{2}\wedge f_{4}\\
0
\end{array}
\right)  ,\left(
\begin{array}
[c]{c}%
0\\
0\\
0\\
f_{2}\wedge f_{3}%
\end{array}
\right)  ,\left(
\begin{array}
[c]{c}%
0\\
0\\
f_{2}\wedge f_{3}\\
f_{3}\wedge f_{4}%
\end{array}
\right)  ,\left(
\begin{array}
[c]{c}%
0\\
0\\
f_{3}\wedge f_{4}\\
f_{2}\wedge f_{4}%
\end{array}
\right)  ,\\\\
\left(
\begin{array}
[c]{c}%
f_{2}\wedge f_{3}\\
0\\
0\\
f_{1}\wedge f_{3}%
\end{array}
\right)  ,\left(
\begin{array}
[c]{c}%
0\\
f_{2}\wedge f_{4}\\
f_{1}\wedge f_{4}\\
0
\end{array}
\right)  ,\left(
\begin{array}
[c]{c}%
f_{2}\wedge f_{4}\\
f_{3}\wedge f_{4}\\
0\\
f_{1}\wedge f_{4}%
\end{array}
\right)  ,\left(
\begin{array}
[c]{c}%
f_{3}\wedge f_{4}\\
f_{2}\wedge f_{3}\\
f_{1}\wedge f_{3}\\
0
\end{array}
\right)
\end{array}
\right\rangle
\]
and in the in the case $\psi\sim C$ we obtain (cf. Appendix 6.1)%
\small{\[
\ker\alpha=\left\langle
\begin{array}
[c]{c}%
\left(
\begin{array}
[c]{c}%
0\\
f_{4}\wedge f_{5}\\
0\\
0
\end{array}
\right)  ,\left(
\begin{array}
[c]{c}%
0\\
0\\
0\\
f_{3}\wedge f_{5}%
\end{array}
\right)  ,\left(
\begin{array}
[c]{c}%
0\\
0\\
f_{3}\wedge f_{5}\\
f_{3}\wedge f_{4}+f_{2}\wedge f_{5}%
\end{array}
\right)  ,\left(
\begin{array}
[c]{c}%
-f_{4}\wedge f_{5}\\
f_{3}\wedge f_{4}-f_{2}\wedge f_{5}\\
0\\
0
\end{array}
\right)  ,\\\\
\left(
\begin{array}
[c]{c}%
0\\
0\\
f_{2}\wedge f_{4}\\
0
\end{array}
\right)  \left(
\begin{array}
[c]{c}%
0\\
0\\
f_{3}\wedge f_{4}+f_{2}\wedge f_{5}\\
f_{2}\wedge f_{4}%
\end{array}
\right)  ,\left(
\begin{array}
[c]{c}%
f_{2}\wedge f_{3}\\
0\\
0\\
0
\end{array}
\right)  ,\left(
\begin{array}
[c]{c}%
f_{2}\wedge f_{5}-f_{3}\wedge f_{4}\\
f_{2}\wedge f_{3}\\
0\\
0
\end{array}
\right)
\end{array}
\right\rangle
\]
}}\end{proof}

\begin{remark}\label{rankremark}
\textnormal{For $\operatorname{char}(\Bbbk)\neq3$ the rank of $\alpha$ is determined by the matrix $\psi$ as in the above lemma. In the case $\operatorname{char}(\Bbbk)=3$ we obtain further elements in $\operatorname{ker}(\alpha)$ for $\psi$ of type $A$ and $B$. To be more precisely in this special situation we get $\dim(\ker{\alpha})=2$ and $6$ for $\psi\sim A$ and $\psi\sim B$ respectively (cf. Appendix 6.1).
\\In the situation where $C$ is given by the pfaffians of a matrix $\psi$ as in $\mathbf{(*)}$ on a scroll $X$ of type $S(2,2,1,0)$ or $S(3,1,1,0)$ the Betti table for $C$ can be calculated in the same way as for $S(2,1,1,1)$, i.e. the only non minimal map has rank depending on the type of $\psi$ as in Lemma \ref{rank}. 
}
\end{remark}

{\normalsize \bigskip}

{\normalsize We have provided the basic information to calculate the
Betti tables for curves $C$, that lie on a scroll of type $S(2,1,1,1)$ and have
one, two or three different $g_{5}^{1\prime}s$. It remains to assign the
different cases to the different possible types for the matrix $\psi.$ This
will be done in the next section. }

\bigskip 

\subsection{\label{g15multi1} {Geometric interpretation}}

{\normalsize Lemma \ref{rank} says that the Betti table of the
minimal free resolution of $\mathcal{O}_{C}$ is determined up to the entries
$\beta_{35}=\beta_{45}$, which can be equal to $4,8$ or $12.$ We have seen
that for $\beta_{35}=8$ or $\beta_{35}=12$ the matrix $\psi$ is of type B, C
or D in Lemma \ref{rank}, i.e. at least one $(H-R)$-entry of $\psi$ can be made
to zero by suitable row and column operations. If we apply\textbf{ }Theorem
\ref{ConPk}, it follows that $C$ is contained in a surface $Y$ given by
the $2\times2$ minors of a matrix
\[
\omega\sim\left(
\begin{array}
[c]{ccc}%
H & H-R & H-R\\
H & H-R & H-R
\end{array}
\right)
\]
on $\mathbb{P(\mathcal{E})}$. $Y$ is a blowup of $\mathbb{P}^{1}%
\times\mathbb{P}^{1}$ in $7$ points. The image $C^{\prime}$ of $C$
in $\mathbb{P}^{1}\times\mathbb{P}^{1}$ is a divisor of type $(5,5)$.
The existence of at least two different $g_{5}^{1}$$^{\prime}s$ that are cut out by the 
factor classes $|(1,0)|$ and $|(0,1)|$ on $\mathbb{P}^{1}\times\mathbb{P}^{1}$ follows. The next
theorem shows that even the converse is true, i.e. $\beta_{35}=4$ if and only if $C$ has
exactly one $g_{5}^{1}:$ }

\begin{theorem}
{\normalsize \label{1g15res}An irreducible, nonsingular canonical curve $C$ with Clifford index 3 that admits no $g_{7}^{2}$ has exactly one ordinary linear system of type $g_{5}^{1}$ if and only if it $C$ is given by the pfaffians of matric $\psi$ is of type A on a scroll $X$ of type $S(2,1,1,1)$. The unique
$g_{5}^{1}$ is cut out by the class of a ruling R on $X.$ The minimal
free resolution of $\mathcal{O}_{C}$ as $\mathcal{O}_{\mathbb{P}^{8}}-$module
has the following Betti table}%
\textnormal{\[
\begin{tabular}{c|cccccccccc}
& 0 & 1 & 2 & 3 & 4 & 5 & 6 & 7 \cr \hline
0 & 1 & - & - & - & - & - & - & - \cr
1 & - & 21 & 64 & 70 & 4 & - & - & - \cr
2 & - & - & - & 4 & 70  & 64  & 21 & -\cr
3 & - & - & - & - & - & - & - & 1 \cr
\end{tabular} 
\]}
\end{theorem}

\begin{proof}
{\normalsize It remains to prove that in the case where $C$ has an additional
$g_{5}^{1}=\left\vert D\right\vert ,$ obtained from an effective divisor $D$
of degree $5$ on $C,$ we get $\beta_{35}=\beta_{45}\neq4.$ Now let us assume
this case, then there is a map%
\[
C\overset{\left\vert R|_{C}\right\vert \times\left\vert D\right\vert
}{\rightarrow}C^{\prime}\subset\mathbb{P}^{1}\times\mathbb{P}^{1}%
\subset\mathbb{P}^{3}%
\]
$C^{\prime}$ is a divisor of type $(5,5)$ on $\mathbb{P}^{1}\times
\mathbb{P}^{1}$, hence $p_{a}(C^{\prime})=4\cdot4=16$. Because of
$p_{a}(C^{\prime})-g(C^{\prime})=16-9=7$ the space model $C^{\prime}$ of $C$
has certain singularities. If $C^{\prime}$ has a singular point with multiplicity 3, then projection from this point leads to a $g_{7}^{2}$ in case
of a triple point or a special linear system of lower Clifford index, which is not possible as in this
situation we would get a different Betti table. It follows that $C^{\prime}$ has exactly $7$ double points $p_{1},...,p_{7}$.
Let $S$ be the surface, which is obtained from $X:=\mathbb{P}^{1}%
\times\mathbb{P}^{1}$ after\ blowing up the singularities of $C^{\prime}$:%
\[
\sigma:S=\tilde{X}(p_{1},...,p_{7})\rightarrow X=\mathbb{P}^{1}\times
\mathbb{P}^{1}%
\]
We denote $A\sim(1,0)$ and $B\sim(0,1)$ the two factor classses of
$P_{1}:=\mathbb{P(}\mathcal{O}_{\mathbb{P}^{1}}\oplus\mathcal{O}%
_{\mathbb{P}^{1}})\simeq\mathbb{P}^{1}\times\mathbb{P}^{1}$, and by abuse of
notation also their pullbacks to $S.$ $E_{i}$ denotes the total transforms of the point $p_{i}$ for $i=1,...,7$. We can assume that $C$ is the strict transform of
$C^{\prime}$ on $S$:%
\[
C\sim5A+5B-\sum_{i=1}^{7}2E_{i}%
\]
We consider the rational map
\[
\varphi:S\rightarrow S^{\prime}\subset\mathbb{P}^{8}%
\]
defined by the adjoint series%
\[
V=H^{0}(S,\omega_{S}(C))=H^{0}(S,\mathcal{O}_{S}(3A+3B-%
{\textstyle\sum_{i=1}^{7}}
E_{i}))
\]
which is base point free according to Corollary \ref{ALSP1P1}. We want to
apply our results on page \pageref{PkonScroll} to show that the variety%
\[
X=%
{\textstyle\bigcup\limits_{B_{\lambda}\in\left\vert B\right\vert }}
\bar{B}_{\lambda}\subset\mathbb{P}^{8}%
\]
is a $4$-dimensional rational normal scroll. Therefore we have to check the
following conditions for $H=3A+3B-%
{\textstyle\sum_{i=1}^{7}}
E_{i}$: }
\\\\
{\normalsize 1. $h^{0}(\mathcal{O}_{S}(H-B))\geq2$ }

\noindent {\normalsize 2. $H^{1}(\mathcal{O}_{S}(kH-B))=0$ for $k\geq1$ and }

\noindent {\normalsize 3. the map $S_{k}H^{0}\mathcal{O}_{S}(H)\rightarrow
H^{0}\mathcal{O}_{S}(kH)$ is surjective }
\\\\
{\normalsize The first condition is trivial because of $h^{0}(\mathcal{O}%
_{S}(H-B))\geq12-7=5$ and 3. follows from Corollary \ref{ALSP1P1} where
$S^{\prime}$ turns out to be arithmetically Cohen Macaulay. It remains to
examine the second condition. Consider the exact sequence%
\[
0\rightarrow\mathcal{O}_{S}(kH-B)\rightarrow\mathcal{O}_{S}(kH)\rightarrow
\mathcal{O}_{B}(kH|_{B})\rightarrow0
\]
and the corresponding long exact sequence of cohomology groups%
\begin{align*}
0  &  \rightarrow H^{0}\mathcal{O}_{S}(kH-B)\rightarrow H^{0}\mathcal{O}%
_{S}(kH)\overset{\delta_{k}}{\rightarrow}H^{0}\mathcal{O}_{B}(kH|_{B}%
)\rightarrow\\
&  \rightarrow H^{1}\mathcal{O}_{S}(kH-B)\rightarrow H^{1}\mathcal{O}%
_{S}(kH)=0
\end{align*}
We first show that $H^{0}\mathcal{O}_{S}(H)\overset{\delta_{1}}{\rightarrow
}H^{0}\mathcal{O}_{B}(H|_{B})\cong H^{0}\mathcal{O}_{\mathbb{P}^{1}%
}(H.B)=H^{0}\mathcal{O}_{\mathbb{P}^{1}}(3)$ is surjective. If $h^{0}%
(\mathcal{O}_{S}(H-B))=d>5$ then the complete linear system $\left\vert
(H-B)|_{C}\right\vert $ would be of type $g_{11}^{d-1}$ and therefore
$\operatorname*{Cliff}(C)\leq2$. It follows that $h^{0}(\mathcal{O}%
_{S}(H-B))=5$ and thus $\dim\delta\geq9-5=4=h^{0}\mathcal{O}_{\mathbb{P}^{1}%
}(3),$ so $\delta_{1}$ is surjective. As further consequence we also obtain
that $H^{0}\mathcal{O}_{S}(kH)\overset{\delta_{k}}{\rightarrow}H^{0}%
\mathcal{O}_{B}(kH|_{B})$ is surjective (the image of $B$ under $\varphi$ is a
rational normal curve), hence $H^{1}\mathcal{O}_{S}(kH-B)=0.$ }

\noindent {\normalsize Let $\pi:\mathbb{P}(\mathcal{E})\rightarrow\mathbb{P}^{1}$ denote
the corresponding $\mathbb{P}^{3}$-bundle and $S^{\prime\prime}$ the strict
transform of $S^{\prime}$ in $\mathbb{P}(\mathcal{E}).$ Then $\mathbb{P}%
(\mathcal{E})$ is of type $S(2,1,1,1)$ an Theorem \ref{Resolutiononscroll}%
\ tells us that $\mathcal{O}_{S^{\prime\prime}}$ has an $\mathcal{O}%
_{\mathbb{P(\mathcal{E})}}$-module resolution of type
\begin{align*}
F_{\ast}  &  :\\
0  &  \rightarrow\mathcal{O}_{\mathbb{P(\mathcal{E})}}(-3H+2R)^{\oplus
2}\overset{\omega}{\rightarrow}\mathcal{O}_{\mathbb{P(\mathcal{E})}%
}(-2H+R)^{\oplus2}\oplus\mathcal{O}_{\mathbb{P(\mathcal{E})}}%
(-2H+2R)\rightarrow
\end{align*}
\hspace{0.8cm} $\rightarrow\mathcal{O}_{\mathbb{P(\mathcal{E})}}\rightarrow
\mathcal{O}_{S^{\prime\prime}}\rightarrow0$
\\\\
\noindent where $\omega$ is given by a matrix 
\[
\omega\sim\left(
\begin{array}
[c]{ccc}%
H-a_{1}R & H-a_{2}R & H-a_{3}R\\
H-a_{1}R & H-a_{2}R & H-a_{3}R
\end{array}
\right)
\]
with entries in $H^{0}(\mathbb{P}(\mathcal{E}),\mathcal{O}%
_{\mathbb{P(\mathcal{E})}}(H-a_{i}R)),$ $a_{i}\in\mathbb{Z}$ for $i=1,2,3.$
From Theorem \ref{ConPk}. we obtain certain conditions on the numbers $a_{i}:$%
\[
5-3\cdot0-(a_{1}+a_{2}+a_{3})=3\Rightarrow
a_{1}+a_{2}+a_{3}=5-3=2
\]
As $S^{\prime}$ is irreducible, we must have $a_{i}\leq1$ for all $i$ and
therefore, assuming $a_{1}\geq a_{2}\geq a_{3}:$%
\[
\omega\sim\left(
\begin{array}
[c]{ccc}%
H & H-R & H-R\\
H & H-R & H-R
\end{array}
\right)
\]
The corresponding mapping cone%
\begin{align*}
& \\
&  \bigskip%
\begin{array}
[c]{cccccc}%
0 & \rightarrow & \mathcal{O}_{\mathbb{P(\mathcal{E})}}(-3H+2R)^{\oplus2} &
^{\underrightarrow{\text{\ \ \ \ \ \ \ \ \ \ \ \ \ \ \ \ \ \ }}} &
\mathcal{O}_{\mathbb{P(\mathcal{E})}}(-2H+R)^{\oplus2} & \oplus\\
&  & \uparrow &  & \uparrow & \\
&  & -----|----- & - & -----|----- & -\\
&  &  &  &  & \\
&  & S_{2}G(-3)^{\oplus2} & ^{\underrightarrow
{\text{\ \ \ \ \ \ \ \ \ \ \ \ \ \ \ \ \ \ }}} & G(-2)^{\oplus2} & \oplus\\
&  & \uparrow &  & \uparrow & \\
&  & F\otimes G(-4)^{\oplus2} & ^{\underrightarrow
{\text{\ \ \ \ \ \ \ \ \ \ \ \ \ \ \ \ \ \ }}} & F(-3)^{\oplus2} & \oplus\\
&  & \uparrow &  & \uparrow & \\
&  & \wedge^{2}F(-5)^{\oplus2} & ^{\underrightarrow{\text{\ \ \ \ \ \ \ \ \ }%
\gamma\text{\ \ \ \ \ \ \ }}} & \wedge^{3}F(-5)^{\oplus2} & \oplus\\
&  & \uparrow &  & \uparrow & \\
&  & \wedge^{4}F(-7)^{\oplus2} & ^{\underrightarrow
{\text{\ \ \ \ \ \ \ \ \ \ \ \ \ \ \ \ \ \ }}} & \wedge^{4}F\otimes DG^{\ast
}(-6)^{\oplus2} & \oplus\\
&  & \uparrow &  & \uparrow & \\
&  & DG^{\ast}(-8)^{\oplus2} & ^{\underrightarrow
{\text{\ \ \ \ \ \ \ \ \ \ \ \ \ \ \ \ \ \ }}} & D_{2}G^{\ast}(-7)^{\oplus2} &
\oplus
\end{array}
\\
&
\end{align*}%
\begin{align*}
&
\begin{array}
[c]{ccccc}%
\mathcal{O}_{\mathbb{P(\mathcal{E})}}(-2H+2R) & ^{\underrightarrow
{\text{\ \ \ \ \ \ \ \ \ \ \ \ \ \ \ \ \ \ }}} & \mathcal{O}%
_{\mathbb{P(\mathcal{E})}} & \rightarrow & 0\\
\uparrow &  & \uparrow &  & \\
----|---- &  & ----|---- &  & \\
&  &  &  & \\
S_{2}G(-2)^{\oplus2} & ^{\underrightarrow
{\text{\ \ \ \ \ \ \ \ \ \ \ \ \ \ \ \ \ \ }}} & \mathcal{O}_{\mathbb{P}^{8}}
&  & \\
\uparrow &  & \uparrow &  & \\
F\otimes G(-3)^{\oplus2} & ^{\underrightarrow
{\text{\ \ \ \ \ \ \ \ \ \ \ \ \ \ \ \ \ \ }}} & \wedge^{2}F(-2) &  & \\
\uparrow &  & \uparrow &  & \\
\wedge^{2}F(-4)^{\oplus2} & ^{\underrightarrow
{\text{\ \ \ \ \ \ \ \ \ \ \ \ \ \ \ \ \ \ }}} & \wedge^{3}F\otimes DG^{\ast
}(-3) &  & \\
\uparrow &  & \uparrow &  & \\
\wedge^{4}F(-6)^{\oplus2} & ^{\underrightarrow
{\text{\ \ \ \ \ \ \ \ \ \ \ \ \ \ \ \ \ \ }}} & \wedge^{4}F\otimes
D_{2}G^{\ast}(-4) &  & \\
\uparrow &  & \uparrow &  & \\
DG^{\ast}(-7)^{\oplus2} & ^{\underrightarrow
{\text{\ \ \ \ \ \ \ \ \ \ \ \ \ \ \ \ \ \ }}} & D_{3}G^{\ast}(-5) &  &
\end{array}
\\
&
\end{align*}
gives us a (not necessarily) minimal free resolution of $\mathcal{O}%
_{S^{\prime}}$. We calculate the rank of the only non minimal map (cf. Appendix 6.2)%
\[
\gamma:\wedge^{2}F(-5)^{\oplus2}\overset{\alpha}{\rightarrow}\wedge
^{3}F(-5)^{\oplus2}%
\]
which is obtained from the matrix $\left(
\begin{array}
[c]{cc}%
\omega_{12} & \omega_{13}\\
\omega_{22} & \omega_{23}%
\end{array}
\right)  \sim\left(
\begin{array}
[c]{cc}%
H-R & H-R\\
H-R & H-R
\end{array}
\right)  :$%
\[
\gamma:\binom{f_{1}\wedge f_{2}}{f_{3}\wedge f_{4}}\rightarrow\binom
{f_{1}\wedge f_{2}\wedge\omega_{12}+f_{3}\wedge f_{4}\wedge\omega_{22}}%
{f_{1}\wedge f_{2}\wedge\omega_{13}+f_{3}\wedge f_{4}\wedge\omega_{23}}%
\]
with $f_{1},...,f_{3}\in H^{0}F.$ The matrix $\left(
\begin{array}
[c]{cc}%
\omega_{12} & \omega_{13}\\
\omega_{22} & \omega_{23}%
\end{array}
\right)  $ has full rank since $S$ is irreduciblew. A
calculation shows that $\dim\ker\gamma=4$:%
\[
\ker\gamma=\langle\binom{\omega_{12}\wedge\omega_{13}}{0},\binom{0}%
{\omega_{22}\wedge\omega_{23}},\binom{\omega_{12}\wedge\omega_{23}+\omega
_{22}\wedge\omega_{13}}{\omega_{12}\wedge\omega_{13}+\omega_{22}\wedge
\omega_{23}},
\]%
\[
\binom{\omega_{12}\wedge\omega_{13}+\omega_{22}\wedge\omega_{23}}{\omega
_{12}\wedge\omega_{23}+\omega_{22}\wedge\omega_{13}}\rangle
\]
Therefore, the Betti table of the minimal free resolution of $\mathcal{O}%
_{S^{\prime}}$ is given by%
\[
\begin{tabular}{c|cccccccc}
& 0 & 1 & 2 & 3 & 4 & 5 & 6 \cr \hline
0 & 1 & - & - & - & - & - & -  \cr
1 & - & 17 & 46 & 45 & 8 & - & - \cr
2 & - & - & - & 4 & 25  & 18 & 4 & \cr
\end{tabular} 
\]
$C$ is contained in $S^{\prime}$, hence the Betti number $\beta_{45}$ in the
minimal free resolution of $\mathcal{O}_{C}$ has to be greater or equal than
8, which contradicts $\beta_{34}=\beta_{45}=4$. }
\end{proof}

\bigskip
\begin{remark}\textnormal{
In the proof of Theorem \ref{1g15res} it turns out that for a pentagonal curve $C$
on a scroll $X$ of type $S(2,1,1,1)$ that has two different $g_{5}^{1}$, there
exists a determinantal surface $S^{\prime}\subset X$ that contains $C.$ This
surface has Betti table as follows%
\[%
\begin{array}
[c]{ccccccc}%
1 &  &  &  &  &  & \\
& 17 & 46 & 45 & 8 &  & \\
&  &  & 4 & 25 & 18 & 4
\end{array}
\]
Considering the exact sequence
\[
0\rightarrow\mathcal{O}_{S^{\prime}}(-C)\rightarrow\mathcal{O}_{S^{\prime}%
}\rightarrow\mathcal{O}_{C}\rightarrow0
\]
a mapping cone construction gives a (not necessarily minimal) free resolution
for $C:$ As {\normalsize $\mathcal{O}_{S^{\prime}}(-C)\cong\omega_{S^{\prime}%
}$ the Betti table for the minimal free resolution of }$\omega_{S^{\prime}}$ as
$\mathcal{O}_{\mathbb{P}^{8}}-$module is given by%
\[%
\begin{array}
[c]{ccccccc}%
4 & 18 & 25 & 4 &  &  & \\
&  & 8 & 45 & 46 & 17 & \\
&  &  &  &  &  & 1
\end{array}
\]
Then the Betti table for $C$ can be obtained as sum of these two Betti tables,
except the values for $\beta_{35}=\beta_{45}$ which are determined by the rank
of the non minimal map "`in the middle"':%
\begin{align*}
&
\begin{array}
[c]{ccccccc}%
1 &  &  &  &  &  & \\
& 17 & 46 & 45 & 8 &  & \\
&  &  & 4 & 25 & 18 & 4\\
&  &  &  &  &  &
\end{array}
\text{ \ \ }+\text{ \ \ }%
\begin{array}
[c]{ccccccc}
&  &  &  &  &  & \\
4 & 18 & 25 & 4 &  &  & \\
&  & 8 & 45 & 46 & 17 & \\
&  &  &  &  &  & 1
\end{array}
\\
& =%
\begin{array}
[c]{cccccccc}%
1 &  &  &  &  &  &  & \\
& 21 & 64 & 70 & \beta_{45} &  &  & \\
&  &  & \beta_{35} & 70 & 64 & 21 & \\
&  &  &  &  &  &  & 1
\end{array}
\end{align*}
In this situation, we know that $8\leq\beta_{35}=\beta_{45}\leq12$. 
In the case of a trigonal and tetragonal curves the Betti table for $C$ could
always be obtained from that of a determinantal surface $S^{\prime
}$ that contains $C$ (cf. Theorems \ref{cliff1} and \ref{cliff2}): It is given as the direct sum of the Betti tables for
{\normalsize $\mathcal{O}_{S^{\prime}}$ and }$\omega_{S^{\prime}}.$ However, in Theorem \ref{2g15res} we will see, that for a pentagonal curve with exactly two different $g_{5}^{1}$ this property fails for pentagonal curves. In this
situation we get $\beta_{35}=\beta_{45}=8.$
In Theorem \ref{triplepoint} we considered pentagonal curves $C$ with exactly one $g_{5}^{1}$
of multiplicity one. It turns out that $C$ admits a $g_{8}^{2}$ with exactly
one triple point and $9$ double points as only singularities. Blowing up
$\mathbb{P}^{2}$ in the singular points we get a surface $S$ and with the
methods of Chapter 2 and Theorem \ref{ConPk} we can see that its image $S^{\prime}\subset
X\subset\mathbb{P}^{8}$ under the adjoint series on $S$ has a determinantal
representation, i.e. on the corresponding {\normalsize $\mathbb{P}^{3}$-bundle
$\mathbb{P}(\mathcal{E})$ the surface }$S^{\prime}$ is given by the $2\times2$
minors of a matrix%
\[
{\normalsize \omega\sim\left(
\begin{array}
[c]{ccc}%
H-a_{1}R & H-a_{2}R & H-a_{3}R\\
H-(a_{1}+1)R & H-(a_{2}+1)R & H-(a_{3}+1)R
\end{array}
\right)  }%
\]
with $a_{1}+a_{2}+a_{3}=0.$}
\end{remark}
\bigskip

{\normalsize \label{conp1p1}The next step is to examine the case where $C$ has
at least two different special linear series of type $g_{5}^{1}.$ As in the
proof of Theorem \ref{1g15res}, we consider the space model $C^{\prime}%
\subset\mathbb{P}^{1}\times\mathbb{P}^{1}\subset\mathbb{P}^{3}$ of $C$ given
by $\left\vert D_{1}\right\vert $ and $\left\vert D_{2}\right\vert :$%
\[
C\overset{\left\vert D_{1}\right\vert \times\left\vert D_{2}\right\vert
}{\rightarrow}C^{\prime}\subset\mathbb{P}^{1}\times\mathbb{P}^{1}%
\subset\mathbb{P}^{3}%
\]
With the notations as above $C^{\prime}$ is a divisor of type $(5,5)\sim5A+5B$
on $\mathbb{P}^{1}\times\mathbb{P}^{1}\cong P_{1}=\mathbb{P(}\mathcal{O}%
_{\mathbb{P}^{1}}\oplus\mathcal{O}_{\mathbb{P}^{1}})$ and we can assume that
$C^{\prime}$ has exactly 7 (possibly infinitely near) double points since
otherwise $C$ has a $g_{7}^{2}$ or a special linear series of lower Clifford
index: }

{\normalsize Blowing up the singularities of $C^{\prime}$ we can assume that $C$ is given as
strict transform of $C^{\prime}$:
\[
C\sim5A+5B-%
{\textstyle\sum\nolimits_{i=1}^{7}}
2E_{i}%
\]
and $K_{C}\sim(3A+3B-\sum_{i=1}^{7}E_{i})|_{C}.$ Then the two
different linear systems of type $g_{5}^{1}$ are given by the divisors $D_{1}\sim
A|_{C}$ and $D_{2}\sim B|_{C}.$ In the next
step, we will examine whether it is possible for $C$ to have a third $g_{5}%
^{1}$ and how such a special linear series can be obtained: }

{\normalsize \bigskip}

\begin{theorem}
{\normalsize \label{2g15}Let $S$ be the iterated blowup in 7 (possibly
infinitely near) points $p_{1},...,p_{7}$ on $\mathbb{P}^{1}\times
\mathbb{P}^{1}$ and $C\sim5A+5B-%
{\textstyle\sum\nolimits_{i=1}^{7}}
2E_{i}$ an irreducible, nonsingular curve with $\operatorname*{Cliff}(C)=3,$
that admits no $g_{7}^{2}.$ Then $C$ has two different $g_{5}^{1}$ given by the divisors
$D_{1}\sim A|_{C}$ and $D_{2}\sim B|_{C}.$ If $C$ admits a third $g_{5}%
^{1}=\left\vert D\right\vert ,$ different from $\left\vert D_{1}\right\vert $
and $\left\vert D_{2}\right\vert ,$ then there exists a point $p\in C$ with
$D\sim K_{C}-D_{1}-D_{2}-p.$ }
\end{theorem}
\begin{proof}
{\normalsize Let us assume that there exists a further $g_{5}^{1}=\left\vert
D\right\vert ,$ $D$ an effective divisor on $C.$ From the base point free
pencil trick we get a long exact sequence:%
\begin{align*}
0 &  \rightarrow H^{0}(C,\mathcal{O}_{C}\mathcal{(}D_{1}-D_{2}))\rightarrow
H^{0}(C,\mathcal{O}_{C}\mathcal{(}D_{1}))\otimes H^{0}(C,\mathcal{O}%
_{C}\mathcal{(}D_{2}))\rightarrow\\
&  \rightarrow H^{0}(C,\mathcal{O}_{C}\mathcal{(}D_{1}+D_{2}))\rightarrow
H^{1}(C,\mathcal{O}_{C}\mathcal{(}D_{1}-D_{2}))\rightarrow...
\end{align*}
As $D_{1}\nsim D_{2}$ it follows that $h^{0}(C,\mathcal{O}_{C}\mathcal{(}%
D_{1}-D_{2}))=0$ and therefore
\[
h^{0}(C,\mathcal{O}_{C}\mathcal{(}D_{1}+D_{2}))\geq h^{0}(C,\mathcal{O}%
_{C}\mathcal{(}D_{1}))\cdot h^{0}(C,\mathcal{O}_{C}\mathcal{(}D_{2}))=4
\]
Hence $h^{0}(C,\mathcal{O}_{C}\mathcal{(}D_{1}+D_{2}))=4$ as for
$h^{0}(C,\mathcal{O}_{C}\mathcal{(}D_{1}+D_{2}))>4$ the linear system
$|D_{1}+D_{2}|$ would have Clifford index
\[
i=\deg(D_{1}+D_{2})-2(h^{0}(C,\mathcal{O}_{C}\mathcal{(}D_{1}+D_{2}))-1)\leq2
\]
Applying the base point free pencil trick again, we obtain the exact sequence%
\begin{align*}
0 &  \rightarrow H^{0}(C,\mathcal{O}_{C}\mathcal{(}D_{1}+D_{2}-D))\rightarrow
H^{0}(C,\mathcal{O}_{C}\mathcal{(}D_{1}+D_{2}))\otimes H^{0}(C,\mathcal{O}%
_{C}\mathcal{(}D))\rightarrow\\
&  \rightarrow H^{0}(C,\mathcal{O}_{C}\mathcal{(}D_{1}+D_{2}+D))\rightarrow
H^{0}(C,\mathcal{O}_{C}\mathcal{(}D_{1}+D_{2}-D))\rightarrow...
\end{align*}
from which we deduce%
\[
h^{0}(C,\mathcal{O}_{C}\mathcal{(}D_{1}+D_{2}+D))\geq4\cdot2-h^{0}%
(C,\mathcal{O}_{C}\mathcal{(}D_{1}+D_{2}-D))=8-h^{0}(C,\mathcal{O}%
_{C}\mathcal{(}D_{1}+D_{2}-D))
\]
The first step will be to exclude $h^{0}(C,\mathcal{O}_{C}\mathcal{(}%
D_{1}+D_{2}-D))\geq2,$ where we must have at least two distinct effective
divisors
\begin{align*}
D^{\ast},D^{\ast\ast} &  \sim D_{1}+D_{2}-D\\
&  \Rightarrow(A+B)|_{C}\sim D_{1}+D_{2}\sim D+D^{\ast},\text{ }D+D^{\ast\ast
}>0
\end{align*}
The divisors $D+D^{\ast}$ and $D+D^{\ast\ast}$ are effective divisors and linear
equivalent to $D_{1}+D_{2}.$ We already know, that $h^{0}(C,\mathcal{O}%
_{C}\mathcal{(}D_{1}+D_{2}))=4$ and therefore every effective divisor, linear
equivalent to $D_{1}+D_{2}$ can be written as $(h|_{C})$ with an element $h\in
H^{0}(S,\mathcal{O}_{S}(A+B))$. Let $h^{\ast},h^{\ast\ast}\in H^{0}%
(S,\mathcal{O}_{S}(A+B))$ with $D^{\ast}+D=(h^{\ast}|_{C})$ and $D^{\ast\ast
}+D=(h^{\ast\ast}|_{C}).$ It follows that $H_{1}:=(h^{\ast})$ and
$H_{2}:=(h^{\ast\ast})$ must have at least the 5 points of $D$ in common. As
$(A+B).(A+B)=2<5,$ this is only possible, if they contain a common divisor
}$\Gamma\sim aA+bB-%
{\textstyle\sum\nolimits_{i=1}^{7}}
\lambda_{i}E_{i}${\normalsize . We assume that }$\Gamma$ is maximal in the
sense that $H_{1}-\Gamma$ and $H_{2}-\Gamma$ have no further component. Then
we must have $0\leq(H_{1}-\Gamma).(H_{2}-\Gamma)=2(1-a)(1-b)-%
{\textstyle\sum\nolimits_{i=1}^{7}}
\lambda_{i}^{2}$. For $(a,b)=(0,0)$ this is only possible if $\Gamma\sim
E_{i}+E_{j},$ $\Gamma\sim E_{i}-E_{j}$ or $\Gamma\sim E_{i}$ for distinct
$i,j\in\{1,...,7\}$. If $\Gamma\sim E_{i}+E_{j}$ or $\Gamma\sim E_{i}-E_{j}$
then we get $(H_{1}-\Gamma).(H_{2}-\Gamma)=0$ and therefore $(H_{1}%
-\Gamma)|_{C}$ and $(H_{2}-\Gamma)|_{C}$ have no common points. Then
{\normalsize $D^{\ast}-D^{\ast\ast}=$}$(H_{1}-\Gamma)|_{C}-(H_{2}-\Gamma
)|_{C}$ is not possible as $(H_{1}-\Gamma).C=(H_{2}-\Gamma).C\geq6$. For
$\Gamma\sim E_{i}$ the correspondance {\normalsize $D^{\ast}-D^{\ast\ast}=$%
}$(H_{1}-\Gamma)|_{C}-(H_{2}-\Gamma)|_{C}$ says that $(H_{1}-\Gamma)|_{C}$ and
$(H_{2}-\Gamma)|_{C}$ must have at least $(H_{1}-\Gamma).C-5=3$ points in
common, which contradicts to $(H_{1}-\Gamma).(H_{2}-\Gamma)=1.$ Thus
{\normalsize we deduce the existence of a common component $l\in
H^{0}(S,\mathcal{O}_{S}(A))$ or $l\in H^{0}(S,\mathcal{O}_{S}(B)).$ We
consider the first case where $h^{\ast}=l\cdot l^{\ast}$ and $h^{\ast\ast
}=l\cdot l^{\ast\ast}$ with $l\in H^{0}(S,\mathcal{O}_{S}(A))$ and $l^{\ast
},l^{\ast\ast}\in H^{0}(S,\mathcal{O}_{S}(B)).$ Now we conclude that
$(l^{\ast}|_{C})=D^{\ast}$ and $(l^{\ast\ast}|_{C})=D^{\ast\ast}$ and
therefore $D_{1}+D_{2}-D\sim D_{2}\Rightarrow D\sim D_{1}$, a contradiction.
For $l\in H^{0}(S,\mathcal{O}_{S}(B))$ the same argument holds. It remains to
show that it is not possible for $\mathcal{O}_{C}\mathcal{(}D_{1}+D_{2}-D)$ to
have exactly one global section. From $\dim|D|=1$ we obtain two different,
effective divisors $D^{\ast},D^{\ast\ast}\sim D.$ As%
\[
h^{0}(C,\mathcal{O}_{C}\mathcal{(}D_{1}+D_{2}-D))=h^{0}(C,\mathcal{O}%
_{C}\mathcal{(}D_{1}+D_{2}-D^{\ast}))=h^{0}(C,\mathcal{O}_{C}\mathcal{(}%
D_{1}+D_{2}-D^{\ast\ast}))=1
\]
there exist effective divisors $\bar{D}^{\ast},\bar{D}^{\ast\ast}$ on $C$ with
$D_{1}+D_{2}-D^{\ast}\sim\bar{D}^{\ast}$ and $D_{1}+D_{2}-D^{\ast\ast}\sim
\bar{D}^{\ast\ast}.$ We get two effective divisors%
\begin{align*}
\tilde{D}^{\ast} &  :=D^{\ast}+\bar{D}^{\ast}\sim D_{1}+D_{2}\\
\tilde{D}^{\ast\ast} &  :=D^{\ast\ast}+\bar{D}^{\ast\ast}\sim D_{1}+D_{2}%
\end{align*}
linear equivalent to $D_{1}+D_{2}$ and in consequence $\tilde{D}^{\ast
}=(h^{\ast}|_{C}),~\tilde{D}^{\ast\ast}=(h^{\ast\ast}|_{C})$ with $h^{\ast
},h^{\ast\ast}\in H^{0}(S,\mathcal{O}_{S}(A+B))$. Because of $h^{0}%
(C,\mathcal{O}_{C}\mathcal{(}D_{1}+D_{2}-D))=1,$ the effective divisors
$\bar{D}^{\ast},\bar{D}^{\ast\ast}$ have to be equal. Therefore%
\[
(h^{\ast}|_{C})-(h^{\ast\ast}|_{C})=D^{\ast}-D^{\ast\ast}%
\]
and the same argument as above leads to $D\sim D_{1}$ or $D\sim D_{2},$ hence
we get a contradiction. The inequality%
\[
h^{0}(C,\mathcal{O}_{C}\mathcal{(}D_{1}+D_{2}+D))\geq8-h^{0}(C,\mathcal{O}%
_{C}\mathcal{(}D_{1}+D_{2}-D))=8
\]
becomes an equality since from Riemann-Roch we obtain%
\begin{align*}
h^{0}(C,\mathcal{O}_{C}\mathcal{(}D_{1}+D_{2}+D)) &  =h^{0}(C,\mathcal{O}%
_{C}\mathcal{(}K_{C}-D_{1}-D_{2}-D)+7\\
&  \leq7+\deg(K_{C}-D_{1}-D_{2}-D)=8
\end{align*}
(cf. \cite{hartshorne} page 298 Ex. 1.5.). It follows that $h^{0}%
(C,\mathcal{O}_{C}\mathcal{(}K_{C}-D_{1}-D_{2}-D)=1,$ hence the existence of a
point $p\in C$ with%
\[
p\sim K_{C}-D_{1}-D_{2}-D
\]%
so
\[
D\sim K_{C}-D_{1}-D_{2}-p\sim(2A+2B)|_{C}-\sum_{i=1}^{7}%
E_{i}|_{C}-p.
\]
}
\end{proof}

{\normalsize \bigskip}

Given {\normalsize an irreducible pentagonal curve $C\sim5A+5B-%
{\textstyle\sum\nolimits_{i=1}^{7}}
2E_{i}$ on a blowup of $\mathbb{P}^{1}\times\mathbb{P}^{1},$ having two
different $g_{5}^{1},$ given by $D_{1}\sim A|_{C}$ and $D_{2}\sim B|_{C},$ and
no }$g_{7}^{2},$ {\normalsize in the case where $C$ has a third $g_{5}%
^{1}=\left\vert D\right\vert $ this linear system is uniquely obtained from a divisor
\[
D\sim K_{C}-D_{1}-D_{2}-p\sim(2A+2B)|_{C}-\sum_{i=1}^{7}E_{i}|_{C}-p
\]
with $p\in C.$ In the proof of the above theorem, we have seen that
$h^{0}(C,K_{C}-D_{1}-D_{2})=2,$ hence there exist two global generators
$g_{1},g_{2}\in H^{0}(S,\mathcal{O}_{S}(2A+2B-\sum_{i=1}^{7}E_{i})).$
Geometrically, the complete linear system $|2A+2B-\sum_{i=1}^{7}E_{i}|$ is given by a pencil $(Q_{\lambda
})_{\lambda}$ of conics passing through the points $p_{1},...,p_{7}$ (cf. picture on
page \pageref{C55pic}). Conversely if $D$ is an effective divisor linear
equivalent to $(2A+2B)|_{C}-\sum_{i=1}^{7}E_{i}|_{C}-p,$ we get $(2A+2B)|_{C}%
-\sum_{i=1}^{7}E_{i}|_{C}\sim D+p.$ Thus $\left\vert
D\right\vert $ is cut out by the pencil $(Q_{\lambda})_{\lambda}.$ As
$(2A+2B).(2A+2B)=8$ we get $h^{0}(C,\mathcal{O}_{C^{\prime}}\mathcal{(}%
K_{C}-D_{1}-D_{2}-p))=2,$ hence $\dim\left\vert D\right\vert =1$ if and only if $p$ is
the basepoint of the linear system $\left\vert 2A+2B-\sum_{i=1}^{7}%
E_{i}\right\vert .$ In general we would expect for $D$ to give a different
linear system of type $g_{5}^{1}$ from $\left\vert D_{1}\right\vert $ and
$\left\vert D_{2}\right\vert $. In the following theorem, we will see that
$D\sim D_{1}$ or $D\sim D_{2}$ can occur in special configurations. For this
purpose let us denote the global sections of $\mathcal{O}_{S}\mathcal{(}A)$
and $\mathcal{O}_{S}\mathcal{(}B)$ by $\lambda,\mu$ and $s,t$ respectively and
w.l.o.g. $(\lambda|_{C})=D_{1}$ and $(s|_{C})=D_{2}$,
then $H^{0}(C,\mathcal{O}_{C}\mathcal{(}D_{1}))=\left\langle \lambda|_{C}%
,\mu|_{C}\right\rangle $ and $H^{0}(C,\mathcal{O}_{C}\mathcal{(}%
D_{2}))=\left\langle s|_{C},t|_{C}\right\rangle .$ We consider the natural
maps%
\[
H^{0}(S,\mathcal{O}_{S}(A))\otimes H^{0}(S,\mathcal{O}_{S}(2A+2B-\sum
_{i=1}^{7}E_{i}))\overset{\delta_{1}}{\rightarrow}H^{0}(S,\mathcal{O}%
_{S}(3A+2B-\sum_{i=1}^{7}E_{i}))
\]
and%
\[
H^{0}(S,\mathcal{O}_{S}(B))\otimes H^{0}(S,\mathcal{O}_{S}(2A+2B-\sum
_{i=1}^{7}E_{i}))\overset{\delta_{2}}{\rightarrow}H^{0}(S,\mathcal{O}%
_{S}(2A+3B-\sum_{i=1}^{7}E_{i}))
\]
Then it turns out that $\ker~\delta_{i}=0$ if and only if $D$ is not linear equivalent to
$D_{i}$ for $i=1,2:$ }

\begin{theorem}
{\normalsize \label{2g15b}Let C be given as in Theorem \ref{2g15} and
$D\sim(2A+2B)|_{C}-\sum_{i=1}^{7}E_{i}|_{C}-p$ with $p$ the unique base point
of $\left\vert 2A+2B-\sum_{i=1}^{7}E_{i}\right\vert $ on $C$. Then, for $i=1,2,$ we
get $m_{\left\vert D\right\vert }\in\{1,2\}$ and
\[
\ker~\delta_{i}\not =0\Leftrightarrow D\sim D_{i}\Leftrightarrow m_{\left\vert
D\right\vert }=2
\]
}
\end{theorem}

\begin{proof}
{\normalsize If $D\sim D_{1}\Leftrightarrow$ $(2A+2B)|_{C}-\sum_{i=1}^{7}%
E_{i}|_{C}-p\sim A|_{C}\Leftrightarrow(A+2B)|_{C}\sim\sum_{i=1}^{7}E_{i}%
|_{C}+p$ then it follows the existence of an element $\gamma\in H^{0}%
(C,\mathcal{O}_{C}((A+2B)|_{C}-%
{\textstyle\sum\nolimits_{i=1}^{7}}
E_{i}|_{C}))$ and therefore $H^{0}(C,\mathcal{O}_{C}((2A+2B)|_{C}-%
{\textstyle\sum\nolimits_{i=1}^{7}}
E_{i}|_{C}))=\left\langle \lambda|_{C}\cdot\gamma,\mu|_{C}\cdot\gamma
\right\rangle .$ We can assume that $(g_{1}|_{C})=(\lambda|_{C}\cdot\gamma)$
and $(g_{2}|_{C})=(\mu|_{C}\cdot\gamma).$ Hence $(\mu\cdot g_{1}-\lambda\cdot
g_{2})|_{C}=0$ and taking into account that%
\begin{align*}
H^{0}(C,\mathcal{O}_{C}((2A+3B)|_{C}-%
{\textstyle\sum\nolimits_{i=1}^{7}}
E_{i}|_{C})) &  =H^{0}(C,\mathcal{O}_{C}(K_{C}-A|_{C}))=\\
&  =H^{0}(S,\mathcal{O}_{S}(2A+3B-%
{\textstyle\sum\nolimits_{i=1}^{7}}
E_{i}))|_{C}%
\end{align*}
we conclude that $\mu\cdot g_{2}-\lambda\cdot g_{1}=0\Rightarrow\ker
~\delta_{1}=\left\langle \mu\otimes g_{2}-\lambda\otimes g_{1}\right\rangle
\neq0.$ Moreover the divisors $Q_{1},Q_{2}\sim2A+2B-\sum_{i=1}^{7}E_{i}$ given
by $g_{1}$ and $g_{2}$ must have a common component $\Gamma$ of type $A+2B-\sum_{i=1}^{7}E_{i}$. Then $\Gamma$ passes through the point $p.$ Conversely from
$\ker~\delta_{1}\neq0$ we deduce the existence of $\tilde{\lambda},\tilde{\mu
}\in$ $H^{0}(S,\mathcal{O}_{S}(A))$ with $\tilde{\lambda}\cdot g_{1}%
=\tilde{\mu}\cdot g_{2},$ hence $(g_{2}|_{C})-(g_{1}|_{C})=(\tilde{\mu}%
|_{C})-(\tilde{\lambda}|_{C})\Rightarrow D\sim D_{2}.$ For $i=2$ the statement
follows in analogous way. }

{\normalsize For the second equivalence we assume that $D\sim D_{1},$ then according to the
approach above, we must have a curve $\Gamma$ of type $A+2B-\sum_{i=1}%
^{7}E_{i}$ that passes through the point $p$. Thus $h^{0}(C,\mathcal{O}_{C}%
\mathcal{(}K_{C}-2D_{1}))=h^{0}(C,\mathcal{O}_{C}\mathcal{(}(A+3B)|_{C}%
-\sum_{i=1}^{7}E_{i}|_{C}))=2$ and applying Riemann Roch, it follows
$h^{0}(C,\mathcal{O}_{C}\mathcal{(}2D_{1})=4.$ The same holds for $D_{2}.$ For
the other direction, let us assume that $h^{0}(C,\mathcal{O}_{C}%
\mathcal{(}2D_{1})=4$ equivalently $h^{0}(C,\mathcal{O}_{C}\mathcal{(}%
(A+3B)|_{C}-\sum_{i=1}^{7}E_{i}|_{C})=h^{0}(C,\mathcal{O}_{C}\mathcal{(}%
K_{C}-2D_{1}))=2.$ Hence the existence of two different effective
divisors $K_{1}$ and $K_{2}$ of type $(A+3B-\sum_{i=1}^{7}E_{i})|_{C}$ on $C$ follows.
The following argument shows that they are cut out by a pencil of divisors of
class $A+3B-\sum_{i=1}^{7}E_{i}$ as above: Without loss of generalisation we
can assume that $(\lambda)$ is irreducible and that $K_{1}$ is cut out by an
effective divisor $\Sigma_{1}\in H^{0}(S,\mathcal{O}_{S}(A+3B-%
{\textstyle\sum\nolimits_{i=1}^{7}}
E_{i}))\neq0$, then
\[
K_{1}+(\lambda^{2}|_{C})=\Sigma_{1}|_{C}+2(\lambda|_{C})\sim K_{C}%
\]
and
\[
K_{2}+2(\lambda|_{C})\sim K_{C}%
\]
hence we obtain
\[
K_{2}=\Omega|_{C}-2(\lambda|_{C})
\]
with an effective divisor $\Omega\sim3A+3B-%
{\textstyle\sum\nolimits_{i=1}^{7}}
E_{i}.$ Because of $\deg K_{2}=6$ and $\deg\Omega|_{C}=16$, the divisors
$\Omega$ and $\Lambda:=2(\lambda)$ must intersect in $10$ points in the
case where }$(\lambda)$ is not contained in $\Omega${\normalsize . As
$2A.\Omega=6$ this is not possible. Thus $\Gamma=(\lambda)$ or $2(\lambda)$ is a
common divisor of $\Omega$ and $\Lambda$, i.e. $\Omega-\Gamma$ and $\Lambda-\Gamma$ are effective. After substracting $(\lambda)$ once we get
$K_{2}=(\Omega-(\lambda))|_{C}-(\lambda|_{C})$ and as $\deg(\Omega
-(\lambda))|_{C}=11,$ $A.(\Omega-(\lambda))=3$ the curves $(\Omega-(\lambda))$
and $(\lambda)$ must contain a further common component $(\lambda).$ It follows the
existence of an effective divisor $\Sigma_{2}:=\Omega-2(\lambda)\sim A+3B-%
{\textstyle\sum\nolimits_{i=1}^{7}}
E_{i}$, that is distinct from $\Sigma_{1}.$ Now from $\Sigma_{1}.\Sigma
_{2}=6-7=-1<0$ it follows that they contain a maximal common divisor
$\Gamma\sim aA+bB-%
{\textstyle\sum\nolimits_{i=1}^{7}}
\lambda_{i}E_{i}$, $a,b\in\mathbb{N}$ and $\lambda_{i}\in\mathbb{Z}$, in the
sense that }$\Sigma_{1}-\Gamma$ and $\Sigma_{2}-\Gamma$ contain no further
common components{\normalsize . Here we can assume that $\Gamma$ is not a union
of components of type $E_{i}$ or $E_{i}-E_{j}$, as from substracting such a
component the intersection product $(\Sigma_{1}-\Gamma).(\Sigma_{2}%
-\Gamma)<\Sigma_{1}.\Sigma_{2}<0$ would stay negative. We distinguish the
different cases for $\Gamma:$ }
\\\\ \noindent {\normalsize - $a=1,$ $b=0\Rightarrow\Sigma_{i}=\Gamma+\Omega_{i}$ with
$\Omega_{i}\sim3B+%
{\textstyle\sum\nolimits_{i=1}^{7}}
(\lambda_{i}-1)E_{i}$ for $i=1,2$. Thus we obtain $\Omega_{1}.\Omega_{2}%
=-\sum_{i=1}^{7}(\lambda_{i}-1)^{2}\leq0$ and $\Omega_{1}.\Omega
_{2}=0\Leftrightarrow\lambda_{1}=...=\lambda_{7}=1\Leftrightarrow\Gamma\sim A-%
{\textstyle\sum\nolimits_{i=1}^{7}}
E_{i},$ which gives a contradiction as $\Gamma.C<0$ and $C$ was assumed to be
irreducible. }
\\\\
\noindent {\normalsize $a=0,$ $b=1\Rightarrow\Sigma_{i}=\Gamma+\Omega_{i}$ with
$\Omega_{i}\sim A+2B+%
{\textstyle\sum\nolimits_{i=1}^{7}}
(\lambda_{i}-1)E_{i}$ for $i=1,2.$ Because of $\Omega_{1}.\Omega_{2}%
=4-\sum_{i=1}^{7}(\lambda_{i}-1)^{2}\geq0$ at least three of the $\lambda_{i}$
are equal to one, we assume this for $\lambda_{1},...,\lambda_{3}.$ But then
\begin{align*}
\Gamma.C  &  =5-\sum_{i=1}^{7}2\lambda_{i}=5+\sum_{i=1}^{7}(\lambda_{i}%
-1)^{2}-\sum_{i=1}^{7}(1+\lambda_{i}^{2})\leq\\
&  \leq9-\sum_{i=1}^{7}(1+\lambda_{i}^{2})<0
\end{align*}
which contradicts that $C$ is irreducible. }
\\\\
\noindent {\normalsize $a=0,$ $b=2\Rightarrow\Sigma_{i}=\Gamma+\Omega_{i}$ with
$\Omega_{i}\sim A+B+%
{\textstyle\sum\nolimits_{i=1}^{7}}
(\lambda_{i}-1)E_{i}$ for $i=1,2$. Then $0\leq\Omega_{1}.\Omega_{2}%
=2-\sum_{i=1}^{7}(\lambda_{i}-1)^{2}$ and therefore $\lambda_{i}=1$ for at
least five distinct $i$, we assume $i=1,...,5.$ Then consider a divisor $L$ of type $(0,1)$ on $\mathbb{P}^{1}%
\times\mathbb{P}^{1}$ that passes through one of the points $p_{1},...,p_{5}.$
Its strict transform $L^{\ast}:=\sigma^{\ast}L\sim B-\sum_{i=1}^{7}%
\mu_{i}E_{i}$, $\mu_{i}\in\mathbb{N}$ and $\mu_{i_{0}}\geq1$ for at least one
$i_{0}\in\{1,...,5\},$ has to be a component of $\Gamma$ as $\Gamma.L^{\ast}<0$. Therefore $\Gamma$ factors into two components $\Gamma_{1}$ and $\Gamma_{2}$ with one of them, we assume $\Gamma_{}$ being of type $B-\sum_{i=1}^{7}\mu_{i}E_{i},$ $\mu_{i}\in\{0,1\}$ and $\#\{i:\mu
_{i}=1\}\geq3.$ Therefore $\Gamma_{1}.C\leq5-6<0,$ a contradiction. 
\\\\
\noindent {\normalsize $a=1,$ $b=1\Rightarrow\Sigma_{i}=\Gamma+\Omega_{i}$ with
$\Omega_{i}\sim2B+%
{\textstyle\sum\nolimits_{i=1}^{7}}
(\lambda_{i}-1)E_{i}$ for $i=1,2.$From $0\leq\Omega_{1}.\Omega_{2}=-\sum
_{i=1}^{7}(\lambda_{i}-1)^{2}\leq0$ it follows, that $\lambda_{1}%
=...=\lambda_{7}=1,$ which is not possible as we would get $\Gamma.C=10-14<0$.
}
\\\\
\noindent {\normalsize $a=1,$ $b=2\Rightarrow\Sigma_{i}=\Gamma+\Omega_{i}$ with
$\Omega_{i}\sim B+%
{\textstyle\sum\nolimits_{i=1}^{7}}
(\lambda_{i}-1)E_{i}$ for $i=1,2$. As $0\leq\Omega_{1}.\Omega_{2}=-\sum
_{i=1}^{7}(\lambda_{i}-1)^{2}$ we must have $\lambda_{1}=...=\lambda_{7}=1,$
thus the curve $\Gamma$ is of class $A+2B-%
{\textstyle\sum\nolimits_{i=1}^{7}}
E_{i}$. In this situation it follows that $D\sim D_{1}$ as above$.$ }
\\\\
\noindent {\normalsize $a=0,$ $b=3\Rightarrow\Sigma_{i}=\Gamma+\Omega_{i}$ with
$\Omega_{i}\sim A+%
{\textstyle\sum\nolimits_{i=1}^{7}}
(\lambda_{i}-1)E_{i}$ for $i=1,2$, hence we get $\lambda_{1}=...=\lambda
_{7}=1$. Then $\Gamma$ is an effective divisor of type $3B-%
{\textstyle\sum\nolimits_{i=1}^{7}}
E_{i}$ which intersects $C$ in one further point $q.$ Let $Q_{j}:=(g_{j}%
)\sim2A+2B-%
{\textstyle\sum\nolimits_{i=1}^{7}}
E_{i},$ $j=1,2,$ denote two of the conics that span the pencil $(Q_{\lambda
})_{\lambda}$, then $\Gamma.Q_{j}=6-7<0$ and thus it follows the existence of
a maximal common component $\Delta_{j}$ of $\Gamma$ and $Q_{j}.$ For
$\Delta_{1}\sim B-%
{\textstyle\sum\nolimits_{i=1}^{7}}
\tilde{\lambda}_{i}E_{i},$ $\tilde{\lambda}_{i}\in\mathbb{Z},$ we obtain
$0\geq(Q_{1}-\Delta_{1}).(\Gamma-\Delta_{1})=4-\sum_{i=1}^{7}(\tilde{\lambda
}_{i}-1)^{2}$ and therefore we can assume that $\tilde{\lambda}_{1}%
=...=\tilde{\lambda}_{3}=1$ and $\tilde{\lambda}$ non negative for all $i$. However this contradicts to the irreducibility of
$C$ because of $\Delta_{1}.C=5-\sum_{i=1}^{7}2\tilde{\lambda}_{i}\leq
9-\sum_{i=1}^{7}(1+\lambda_{i}^{2})<0.$ The same argument shows that
$\Delta_{1}\sim B-%
{\textstyle\sum\nolimits_{i=1}^{7}}
\tilde{\lambda}_{i}E_{i}$ is not possible. It remains to discuss the case
where $\Delta_{j}\sim2B-%
{\textstyle\sum\nolimits_{i=1}^{7}}
\tilde{\lambda}_{i}E_{i}:$ Here we get $0\geq(Q_{1}-\Delta_{1}).(\Gamma
-\Delta_{1})=2-\sum_{i=1}^{7}(\tilde{\lambda}_{i}-1)^{2}$ and therefore we can
assume $\tilde{\lambda}_{1}=...=\tilde{\lambda}_{5}=1.$ Now from $0\leq
\Delta_{1}.C=10-\sum_{i=1}^{7}2\tilde{\lambda}_{i}\leq12-\sum_{i=1}%
^{7}(1+\lambda_{i}^{2})$ it follows that $\tilde{\lambda}_{5}=\tilde{\lambda
}_{6}=0$, hence $\Delta_{1}$ is an effective divisor of type $2B-%
{\textstyle\sum\nolimits_{i=1}^{5}}
E_{i}.$ In the case $(a,b)=(0,2)$ we have already excluded this possibility.}
\\\\
}
{\normalsize It remains to show that $h^{0}(C,\mathcal{O}_{C}\mathcal{(}%
K_{C}-3D_{1}))\neq1.$ For $h^{0}(C,\mathcal{O}_{C}\mathcal{(}K_{C}-3D_{1}))=1$
we have $h^{0}(C,\mathcal{O}_{C}\mathcal{(}K_{C}-2D_{1}))=2$ and therefore
with the considerations above there exists an effective divisor $\Gamma\sim
A+2B-\sum_{i=1}^{7}E_{i}.$ Let $L\sim K_{C}-3D_{1}$ be an effective divisor on
$C,$ then we get the following representation%
\[
(A+3B-\sum_{i=1}^{7}E_{i})|_{C}\sim L+(\lambda|_{C})=\Gamma|_{C}+(s^{\prime}|_{C})
\]
and%
\[
L+(\mu|_{C})=\Gamma|_{C}+(t^{\prime}|_{C})
\]
with $s^{\prime},t^{\prime}\in H^{0}(S,\mathcal{O}_{S}\mathcal{(}B)).$ We
conlcude $(\lambda|_{C})-(\mu|_{C})=(s^{\prime}|_{C})-(t^{\prime}|_{C})$ and
thus $(\lambda|_{C})\sim(s^{\prime}|_{C})\Rightarrow D_{1}\sim D_{2},$ a
contradiction. }

{\normalsize For $i=2$ the claim can be proven in an analogous way. }
\end{proof}

{\normalsize \bigskip}

{\normalsize It is time for a short overview of what we have done in this
chapter: We have proven that a pentagonal curve $C,$ which is given by the
Pfaffians of a matrix $\psi$ on a scroll of type $S(2,1,1,1)$ as in
\textbf{(*)}, has exactly one special linear series of type $g_{5}^{1}$ if and only if
$\psi$ is of type $A$ as in Lemma \ref{psitypes}. Then the entry }$\beta
_{46}=\beta_{56}$ {\normalsize  of the Betti table for $C$ equals to 4. Furthermore we
have seen that if a pentagonal curve $C$ that admits no }$g_{7}^{2}$ and {\normalsize has two different $g_{5}^{1}$$^{\prime}s$, then at least one of them
has multiplicity one, hence the scroll constructed from $\left\vert
D\right\vert $ is of type $S(2,1,1,1).$ Therefore, in this case we can assume
that $C$ is given by the Pfaffians of a matrix $\psi$ as in \textbf{(*)}.
Moreover, it turns out that if $C$ has a third $g_{5}^{1}$, then this
}linear series {\normalsize  is uniquely determined. In some cases it becomes
equal to one of the two others. Now it remains to give the correspondence
between these cases and the different types of the matrix $\psi$ which gives
us the Betti table for $C.$ }

\bigskip

\begin{theorem}
Let $C$ be an irr{\normalsize educible, smooth, canonical curve of genus 9
with $\operatorname*{Cliff}(C)=3$ that admits no $g_{7}^{2}.$ Then}
\\\\
\noindent {\normalsize a) C has exactly three different $g_{5}^{1}$$^{\prime}s$ if and only if $C$
is given by the Pfaffians of a matrix $\psi$ of type C in Lemma \ref{psitypes}
on a scroll $X\simeq S(2,1,1,1)$. }
\\\\
\noindent {\normalsize b) C has exactly two different linear systems $\left\vert
D_{1}\right\vert ,\left\vert D_{2}\right\vert $ of type $g_{5}^{1}$ with
$m_{\left\vert D_{1}\right\vert }=2$ and $m_{\left\vert D_{2}\right\vert }=1$
if and only if C is given by the Pfaffians of a matrix $\psi$ of type D. }
\\\\
{\normalsize In both cases the minimal free resolution of $\mathcal{O}_{C}$ as
$\mathcal{O}_{\mathbb{P}^{8}}-$module has the following Betti diagram:%
\textnormal{\[
\begin{tabular}{c|cccccccccc}
& 0 & 1 & 2 & 3 & 4 & 5 & 6 & 7 \cr \hline
0 & 1 & - & - & - & - & - & - & - \cr
1 & - & 21 & 64 & 70 & 12 & - & - & - \cr
2 & - & - & - & 12 & 70  & 64  & 21 & -\cr
3 & - & - & - & - & - & - & - & 1 \cr
\end{tabular} 
\]}
}
\end{theorem}

\begin{proof}
{\normalsize a) Let us first assume that $C$ 
 different $g_{5}^{1},$
then we already know that $C$ is given by the Pfaffians of a matrix $\psi$ of
type B, C or D (cf. Lemma \ref{psitypes})\ on a scroll $X\simeq S(2,1,1,1)$.
As we have seen before, there exists a space model $C^{\prime}$ of $C$ on
$\mathbb{P}^{1}\times\mathbb{P}^{1},$ having exactly 7 (possibly infinitely
near) double points $p_{1},...,p_{7}$ as only singularities. We can assume
that $C$ is given as strict transform of $C^{\prime}$ in the blowup in the
points $p_{1},...,p_{7}$ (cf. page \pageref{conp1p1}), i.e. $C\sim5A+5B-%
{\textstyle\sum\nolimits_{i=1}^{7}}
2E_{i}$ and $K_{C}\sim(3A+3B)|_{C}-\sum_{i=1}^{7}E_{i}|_{C}.$ From Theorem
\ref{2g15}, we already know that the three different linear series of type
$g_{5}^{1}$ are given by divisors $D_{1}\sim A|_{C}$, $D_{2}\sim B|_{C}$ and
$D\sim(2A+2B)|_{C}-\sum_{i=1}^{7}E_{i}|_{C}-p$ with $p\in C$ the basepoint of
$\left\vert 2A+2B-\sum_{i=1}^{7}E_{i}\right\vert .$ The global sections of
$\mathcal{O}_{S}(A),$ $\mathcal{O}_{S}(B)$ and $H^{0}(S,\mathcal{O}%
_{S}(2A+2B-\sum_{i=1}^{7}E_{i}))$ are denoted as in the theorem above. $C\subset$ $S$ is given by the Pfaffians of the matrix $\psi$ regarding the
entries of $\psi$ as global sections of vectorbundles on $S.$ As the
$2\times2$ minors of the submatrix%
\[
\omega=\left(
\begin{array}
[c]{ccc}%
\psi_{14} & \psi_{24} & \psi_{34}\\
\psi_{15} & \psi_{25} & \psi_{35}%
\end{array}
\right)
\]
vanish on $S,$ the two rows of $\omega$ are linear dependent on $S$. Therefore
we get%
\[
\psi\sim\left(
\begin{array}
[c]{ccccc}%
0 & h_{1} & h_{2} & \lambda\cdot\varphi & \mu\cdot\varphi\\
& 0 & \gamma & \lambda\cdot g_{1} & \mu\cdot g_{1}\\
&  & 0 & \lambda\cdot g_{2} & \mu\cdot g_{2}\\
&  &  & 0 & 0\\
&  &  &  & 0
\end{array}
\right)
\]
with $\varphi\in H^{0}(S,\mathcal{O}_{S}(2A+3B-\sum_{i=1}^{7}E_{i}))$ and 
$\gamma\in H^{0}(S,\mathcal{O}_{S}(3A+2B-\sum_{i=1}^{7}E_{i}))$. Recall
that $H^{0}(C,\mathcal{O}_{C}(D_{1}))=\left\langle \lambda|_{C},\mu
|_{C}\right\rangle ,$ $H^{0}(C,\mathcal{O}_{C}(D_{2}))=\left\langle
s|_{C},t|_{C}\right\rangle ,$ $H^{0}(C,\mathcal{O}_{C}(D+p))=\left\langle
g_{1}|_{C},g_{2}|_{C}\right\rangle $ The curve $C$ is given by the equation
\[
\delta=h_{1}g_{2}-h_{2}g_{1}+\varphi\gamma\in H^{0}(S,\mathbb{\mathcal{O}}%
_{S}(5A+5B-\sum_{i=1}^{7}2E_{i})
\]
Then $\delta$ also vanishes at $p$ as $p\in C.$ As $h_{1}g_{2}-h_{2}g_{1}$ vanishs
at $p$, too, the same holds for $\varphi\gamma\in H^{0}(S,\mathbb{\mathcal{O}%
}_{S}(5A+5B-\sum_{i=1}^{7}2E_{i})\mathbb{)}$, hence $\gamma$ or $\varphi$ has
to vanish at $p.$ As we assumed for $\left\vert D\right\vert $ to be different
from $\left\vert D_{1}\right\vert ,$ Theorem \ref{2g15b}\textbf{ }tells us,
that $\varphi\in\left\langle s\cdot g_{1},t\cdot g_{1},s\cdot g_{2},t\cdot
g_{2}\right\rangle $ or $\gamma\in\left\langle \lambda\cdot g_{1},\mu\cdot
g_{1},\lambda\cdot g_{2},\mu\cdot g_{2}\right\rangle $. In the case
$\varphi\in\left\langle s\cdot g_{1},t\cdot g_{1},s\cdot g_{2},t\cdot
g_{2}\right\rangle $ the entry $\psi_{14}=\lambda\cdot\varphi$ could be made
to zero by suitable row and column operations, hence we get a contradiction
to our assumption that $C$ is irreducible. Therefore we must have $\psi_{23}=\gamma
\in\left\langle \lambda\cdot g_{1},\mu\cdot g_{1},\lambda\cdot g_{2},\mu\cdot
g_{2}\right\rangle $, so this entry of $\psi$ can be made zero. This proves one direction. }

{\normalsize For the other direction let $C$ be given by the Pfaffians of
\[
\psi\sim\left(
\begin{array}
[c]{ccccc}%
0 & H & H & H & H\\
& 0 & 0 & H-R & H-R\\
&  & 0 & H-R & H-R\\
&  &  & 0 & 0\\
&  &  &  & 0
\end{array}
\right)
\]
on a scroll $\mathbb{P}(\mathcal{E})$ of type $S(2,1,1,1).$ Then the 5
Pfaffians are given by the $2\times2$ minors of%
\[
\omega_{1}=\left(
\begin{array}
[c]{ccc}%
\psi_{14} & \psi_{24} & \psi_{34}\\
\psi_{15} & \psi_{25} & \psi_{35}%
\end{array}
\right)  \sim\left(
\begin{array}
[c]{ccc}%
H & H-R & H-R\\
H & H-R & H-R
\end{array}
\right)
\]
and
\[
\omega_{2}=\left(
\begin{array}
[c]{ccc}%
-\psi_{12} & \psi_{24} & \psi_{25}\\
-\psi_{13} & \psi_{34} & \psi_{35}%
\end{array}
\right)  \sim\left(
\begin{array}
[c]{ccc}%
H & H-R & H-R\\
H & H-R & H-R
\end{array}
\right)
\]
Each of the minors of $\omega_{1}$ and $\omega_{2}$ define a surface
$Y_{1}$ and $Y_{2}$ respectively with $C=Y_{1}\cap Y_{2}.$ Both surfaces are
blowups of $\mathbb{P}^{1}\times\mathbb{P}^{1}$ with $C$ being a divisor of
type $5A+5B$ on each of them (cf. Theorem \ref{ConPk}). Then, the two linear
systems of type $g_{5}^{1}$ cut out by global sections \ of
$\mathbb{\mathcal{O}}_{Y_{1}}(B)$ and $\mathbb{\mathcal{O}}_{Y_{2}}(B)$ are
the same. The other two $g_{5}^{1}$$^{\prime}s$, given by $A|_{C}$ on $Y_{1}$ and $Y_{2}$
respectively, have to be different as the surfaces $Y_{1}$ and $Y_{2}$ are
different. This proves the existence of three different $g_{5}^{1}$$^\prime s$ . }

{\normalsize b) The two different linear series of type $g_{5}^{1}$ are given
by divisors $D_{1}\sim A|_{C}$ and $D_{2}\sim B|_{C}.$ Then according to the
last theorem we can assume that the scroll constructed from $\left\vert
D_{2}\right\vert $ is of type $S(2,1,1,1).$ We have also shown that a third
$g_{5}^{1}$ is obtained from a divisor $D\sim(2A+2B)|_{C}-\sum_{i=1}%
^{7}E_{i}|_{C}-p$ with $p\in C$ the basepoint of $\left\vert 2A+2B-\sum_{i=1}%
^{7}E_{i}\right\vert .$ Furthermore\textbf{ }Theorem \ref{2g15b} tells us that
$D\not \sim D_{2}$\ and $D\sim D_{1}\Leftrightarrow$ $\dim\left\langle
\lambda\cdot g_{1},\mu\cdot g_{1},\lambda\cdot g_{2},\mu\cdot g_{2}%
\right\rangle \leq3.$ The only possibility for $\left\vert D_{1}\right\vert $
to have multiplicity 2 is therefore given exactly in the case, where $D\sim
D_{1}$. This is equivalent to that the four entries $\psi_{24},\psi_{25}%
,\psi_{34}$ and $\psi_{24}$ span a three dimensional vector space, thus
$\psi$ is of type D. The claim according to the Betti table is a direct
consequence of Lemma \ref{rank}. }
\end{proof}

{\normalsize \bigskip}

{\normalsize It remains to consider the case where $C$ has exactly two
different linear systems of type $g_{5}^{1}$, each of them with multiplicity one$.$ The
following theorem is a direct conclusion from the above theorems: }

{\normalsize \bigskip}

\begin{theorem}
{\normalsize \label{2g15res}Let C be an irreducible, smooth, canonical curve
of genus 9 with $\operatorname*{Cliff}(C)=3$, that admits no $g_{7}^{2}$. Then
C has exactly two different $g_{5}^{1}$, each of them with multiplicity one, if and only if it is given by the Pfaffians of a matrix $\psi$ of type B
on a scroll of type $S(2,1,1,1)$. The minimal free resolution of
$\mathcal{O}_{C}$ as $\mathcal{O}_{\mathbb{P}^{8}}-$module has the following
Betti table%
\textnormal{\[
\begin{tabular}{c|cccccccccc}
& 0 & 1 & 2 & 3 & 4 & 5 & 6 & 7 \cr \hline
0 & 1 & - & - & - & - & - & - & - \cr
1 & - & 21 & 64 & 70 & 8 & - & - & - \cr
2 & - & - & - & 8 & 70  & 64  & 21 & -\cr
3 & - & - & - & - & - & - & - & 1 \cr
\end{tabular} 
\]}
}
\end{theorem}

{\normalsize \newpage}

\section{{\label{highermulti}Infinitely near
$g_{5}^{1\prime}s$}}

{\normalsize In the last sections, we have seen that for $C$, an irreducible,
smooth, canonical curve of genus 9, it is possible to have one, two or even
three different special linear series of type $g_{5}^{1},$ such that at least
one of }the $g_{5}^{1}$$^{\prime}s$ {\normalsize  has multiplicity one. We have shown that in
each case $C$ is given by the Pfaffians of a corresponding matrix $\psi$ on a
scroll of type $S(2,1,1,1).$ Moreover we have calculated the Betti tables for
$C\subset$ $\mathbb{P}^{8}.$ In this section we concentrate on the remaining
case where $C$ has Clifford index $3$, admits no $g_{7}^{2},$ and has exactly
one $g_{5}^{1}=\left\vert D\right\vert $ of multiplicity $2$ or $3.$ In this
situation the scroll 
\[
X=%
{\textstyle\bigcup\limits_{D_{\lambda}\in\left\vert D\right\vert }}
\bar{D}_{\lambda}\subset\mathbb{P}^{8}%
\]
constructed from $|D|$ is a 4-dimensional rational scroll of type $S(2,2,1,0)$ or $S(3,1,1,0)$
depending on if $m_{\left\vert D\right\vert }=2$ or $m_{\left\vert
D\right\vert }=3$ respectively (cf. Theorem \ref{scrolltypes})$.$ Then Theorem
\ref{psitypes2} says that also in this situation there exists a unique
representation of }$C$ by the Pfaffians of a skew symmetric matrix $\psi$ on the
scroll $X.$ It is of type%
\[
\psi\sim\left(
\begin{array}
[c]{ccccc}%
0 & H+R & H & H & H-R\\
& 0 & H & H & H-R\\
&  & 0 & H-R & H-2R\\
&  &  & 0 & H-2R\\
&  &  &  & 0
\end{array}
\right)
\]
and the Betti table for $C$ can be obtained by considering the non minimal
maps in the corresponding mapping cone construction. Moreover $C$ is contained
in a determinantal surface $Y\subset X,$ that is obtained as image of a blowup
of a cone in $\mathbb{P}^{3}:\ $We focus on this question first. {\normalsize
A necessary condition for }$m_{\left\vert D\right\vert }\geq2$ {\normalsize  is
$h^{0}(C,\mathcal{O}_{C}(2D))=4$. We consider the space model
$C^{\prime}$ of $C$ obtained from the complete linear series $\left\vert
2D\right\vert .$ }

{\normalsize \bigskip}

\begin{lemma}
{\normalsize \label{cone}Let $X\subset\mathbb{P}^{8}$ be the scroll given by a
special linear series $\left\vert D\right\vert $ of type $g_{5}^{1}$ of a
canonical curve $C$ of genus 9 as above, i.e. $m_{D}\geq2$. Then there exists
a space model $C^{\prime}$ of $C$ lying on the cone $Y\subset\mathbb{P}^{3}$
over the irreducible conic $x_{1}^{2}-x_{0}x_{2}$ in $\mathbb{P}^{3}.$ }
\end{lemma}

\begin{proof}
{\normalsize We have remarked that $X$ is of type $S(2,2,1,0)$ or $S(3,1,1,0)$
if and only if $h^{0}(C,\mathcal{O}_{C}(2D))=4$. Let $r,s\in H^{0}(C,\mathcal{O}_{C}(D))$
denote the global generators of $\mathcal{O}_{C}(D)$, then we get three global
sections $x_{0}$, $x_{1}$, $x_{2}\in H^{0}(C,\mathcal{O}_{C}(2D))$ by $x_{0}%
:=r^{2},~x_{1}:=s\cdot r$ and $x_{2}:=s^{2}$. Further there exists an element
$x_{3}\in H^{0}(C,\mathcal{O}_{C}(2D))\setminus\left\langle x_{0},x_{1}%
,x_{2}\right\rangle .$ Considering the space model $C^{\prime}$ of $C$, which
is given by the complete linear system $\left\vert 2D\right\vert :$
\[
C\overset{\left\vert 2D\right\vert }{\rightarrow}C^{\prime}\subset
\mathbb{P}^{3}%
\]
we see that $C^{\prime}$ is contained in the surface $Y\subset\mathbb{P}^{3}$,
given by the equation $x_{1}^{2}-x_{0}x_{2}=0,$ i.e. $Y$ is a cone over a conic in $\mathbb{P}^{2}$%
\[
{\includegraphics[height=8cm,width=7cm]%
{plot2}%
}%
\]
As $\left\vert D\right\vert $ was assumed to be base point free $C^{\prime}%
$ does not pass the vertex of the cone. }
\end{proof}

\subsection{Resolution and Representation}

{\normalsize In the foregoing lemma we have seen that there exists a space
model $C^{\prime}\subset Y\subset\mathbb{P}^{3}$ of }$C$ that is
contained {\normalsize  in the singular cone $Y\subset\mathbb{P}^{3}$ over the
irreducible conic $x_{1}^{2}-x_{0}x_{2}$ in $\mathbb{P}^{2}.$ Consider%
\[
\pi:P_{2}=\mathbb{P}(\mathcal{O}_{\mathbb{P}^{1}}(2)\oplus\mathcal{O}%
_{\mathbb{P}^{1}})\rightarrow\mathbb{P}^{1}%
\]
the corresponding $\mathbb{P}^{1}-$bundle with hyperplane class $A$ and ruling
$B,$ that represents a desingularisation of the cone$.$ The exceptional
divisor of this blowup is an effective divisor of class }$A-2B.$ {\normalsize
The strict transform $C^{\prime\prime}\subset P_{2}$ of }$C^{\prime}$
{\normalsize is then given as divisor of class $5A$ as $C^{\prime}$ does not
pass the vertex of the cone $Y.$ From $p_{a}(C^{\prime})=(5-1)^{2}=16$ it
follows that $C^{\prime\prime}$ has certain singularities on $P_{2}$ that
cannot be of multiplicity $3$ or higher, as projection from such a point would
lead to a $g_{d}^{2}$ with $d\leq7.$ Thus $C^{\prime\prime}$ has exactly $7$
(possible infinitely near) double points $p_{1},...,p_{7}$ on $P_{2}.$ The
$g_{5}^{1}$ is then cut out by the class of a ruling $B$. We consider the
iterated blowup
\[
\sigma:S=\tilde{P}_{2}(p_{1},...,p_{7})\rightarrow P_{2}%
\]
in the singular points of $C^{\prime\prime}$. We can assume that $C$ is the strict transform of
$C^{\prime\prime}$ under the blowup $\sigma$ and $K_{C}\sim(3A-%
{\textstyle\sum\nolimits_{i=1}^{7}}
E_{i})|_{C}$. Consider the morphism%
\[
\varphi:S\rightarrow S^{\prime}\subset\mathbb{P}^{8}%
\]
defined by the adjoint series%
\[
H^{0}(S,\omega_{S}(\sigma^{\ast}C^{\prime\prime}))=H^{0}(S,\mathcal{O}_{S}(3A-%
{\textstyle\sum\nolimits_{i=1}^{7}}
E_{i}))
\]
According to Corollary \ref{ALSF2} the adjoint linear series is base point
free and the surface $S^{\prime}$ is arithmetically Cohen-Macaulay.
Now we want to apply our results on page \ref{PkonScroll} to
show that the variety%
\[
X=%
{\textstyle\bigcup\limits_{B_{\lambda}\in\left\vert B\right\vert }}
\bar{B}_{\lambda}\subset\mathbb{P}^{8}%
\]
is a $4$-dimensional rational normal scroll. We have to check the
following conditions for $H=3A-%
{\textstyle\sum_{i=1}^{7}}
E_{i}$: }

\bigskip

{\normalsize 1. $h^{0}(\mathcal{O}_{S}(H-B))\geq2$ }

{\normalsize 2. $H^{1}(\mathcal{O}_{S}(kH-B))=0$ for $k\geq1$ and }

{\normalsize 3. the map $S_{k}H^{0}\mathcal{O}_{S}(H)\rightarrow
H^{0}\mathcal{O}_{S}(kH)$ is surjective. }

\bigskip

\noindent{\normalsize The first condition is trivial because of $h^{0}%
(\mathcal{O}_{S}(H-B))\geq12-7=5$. Condition 3. follows as $S^{\prime}$ is
arithmetically Cohen-Macaulay. It remains to examine the second condition.
Consider the exact sequence%
\[
0\rightarrow\mathcal{O}_{S}(kH-B)\rightarrow\mathcal{O}_{S}(kH)\rightarrow
\mathcal{O}_{B}(kH|_{B})\rightarrow0
\]
and the corresponding long exact sequence of cohomology groups%
\begin{align*}
0 &  \rightarrow H^{0}\mathcal{O}_{S}(kH-B)\rightarrow H^{0}\mathcal{O}%
_{S}(kH)\overset{\delta_{k}}{\rightarrow}H^{0}\mathcal{O}_{B}(kH|_{B}%
)\rightarrow\\
&  \rightarrow H^{1}\mathcal{O}_{S}(kH-B)\rightarrow H^{1}\mathcal{O}%
_{S}(kH)=0
\end{align*}
We first show that $H^{0}\mathcal{O}_{S}(H)\overset{\delta_{1}}{\rightarrow
}H^{0}\mathcal{O}_{B}(H|_{B})\cong H^{0}\mathcal{O}_{\mathbb{P}^{1}%
}(H.B)=H^{0}\mathcal{O}_{\mathbb{P}^{1}}(3)$ is surjective. If $h^{0}%
(\mathcal{O}_{S}(H-B))=d>5$ then the complete linear system $\left\vert
(H-B)|_{C}\right\vert $ would be of type $g_{11}^{d-1}$ and therefore
$\operatorname*{Cliff}(C)\leq2$. It follows that $h^{0}(\mathcal{O}%
_{S}(H-B))=5$ and thus $\dim\delta\geq9-5=4=h^{0}\mathcal{O}_{\mathbb{P}^{1}%
}(3),$ so $\delta_{1}$ is surjective. As a further consequence we also obtain
that $H^{0}\mathcal{O}_{S}(kH)\overset{\delta_{k}}{\rightarrow}H^{0}%
\mathcal{O}_{B}(kH|_{B})$ is surjective (the image of $B$ under $\varphi$ is a
rational normal curve), hence $H^{1}\mathcal{O}_{S}(kH-B)=0.$ Let
$\mathbb{P}(\mathcal{E})$ be the corresponding }$\mathbb{P}^{3}-$bundle to the
scroll $X$, then {\normalsize  according to Theorem \ref{Resolutiononscroll}
$S^{\prime}$ can be given by the $2\times2$ minors of a matrix of type%
\[
\omega\sim\left(
\begin{array}
[c]{ccc}%
H-a_{1}R & H-a_{2}R & H-a_{3}R\\
H-(a_{1}+2)R & H-(a_{2}+2)R & H-(a_{3}+2)R
\end{array}
\right)
\]
on $\mathbb{P}(\mathcal{E})$ with $a_{1},a_{2},a_{3}\in\mathbb{Z}$. }

{\normalsize From Theorem \ref{ConPk}, we obtain certain conditions on the
numbers $a_{i}$ $($we denote $f=\deg X=5,$ $d=3$ and $a=a_{1}%
+a_{2}+a_{3}):$%
\[
K_{C}\sim(3A-%
{\textstyle\sum\nolimits_{i=1}^{7}}
E_{i})|_{C}=(3A+(f-d\cdot k-a)B-%
{\textstyle\sum\nolimits_{i=1}^{7}}
E_{i})|_{C}\Rightarrow a=-1
\]
}

{\normalsize If $X$ is of type $S(2,2,1,0)$, it follows that $a_{i}\leq0$ for
all $i=1,2,3$ since otherwise some of the minors of $\omega$ would vanish or
become reducible. In the case, where $X\simeq S(3,1,1,0),$ we must have
$a_{i}\leq1$ for all $i$ for the same reasons. Moreover we cannot have $a_{1}=a_{2}=1$ in the case $X\simeq
S(3,1,1,0),$ since we would get a common factor $\varphi_{0}\in H^{0}%
(\mathbb{P}(\mathcal{E}),\mathcal{O}_{\mathbb{P(\mathcal{E})}}(H-3R))$ of al minors. We
conclude: }

{\normalsize \bigskip}

\begin{corollary}
{\normalsize \label{infinitypes0}Let $X\subset\mathbb{P}^{8}$ be the scroll
constructed from a special linear series $\left\vert D\right\vert $ of type
$g_{5}^{1}$, $m_{\left\vert D\right\vert }\geq2,$ on an irreducible,
nonsingular, canonical curve $C\subset\mathbb{P}^{8}$ of genus 9 with
$\operatorname*{Cliff}(C)=3,$ that admits no $g_{2}^{7}.$ With $\mathbb{P}%
(\mathcal{E})$ denoting the corresponding $\mathbb{P}^{3}-$bundle, $C$ is
contained in a determinantal surface $Y$ given by the minors of a $2\times3$
matrix of type }
\\\\
\noindent {\normalsize a)%
\[
\omega\sim\left(
\begin{array}
[c]{ccc}%
H & H & H+R\\
H-2R & H-2R & H-R
\end{array}
\right)
\]
if $\mathbb{P}(\mathcal{E})\cong \mathbb{P}(\mathcal{O}(2)\oplus\mathcal{O}(2)\oplus\mathcal{O}(1)\oplus \mathcal{O}).\medskip$ }
\\\\
\noindent {\normalsize b)%
\[
\omega\sim\left(
\begin{array}
[c]{ccc}%
H-R & H+R & H+R\\
H-3R & H-R & H-R
\end{array}
\right)
\]
if$\ \mathbb{P}(\mathcal{E})\cong \mathbb{P}(\mathcal{O}(3)\oplus\mathcal{O}(1)\oplus\mathcal{O}(1)\oplus \mathcal{O}).$ }
\end{corollary}

{\normalsize \bigskip}

{\normalsize In Section \ref{curvesscrolls} we mentioned that on
$\mathbb{P}(\mathcal{E})$, the vanishing ideal of $C$ is given by the Pfaffians
of the matrix $\psi$ that occurs in the free resolution of $\mathcal{O}_{C}$
as $\mathcal{O}_{\mathbb{P}(\mathcal{E})}$-module:%
\[
F_{\ast}:\ \ \ \ \ \ \ \ \ \ 0\rightarrow\mathcal{O}_{\mathbb{P(\mathcal{E})}%
}(-5H+3R)\rightarrow%
{\displaystyle\sum\limits_{i=1}^{5}}
\mathcal{O}_{\mathbb{P(\mathcal{E})}}(-3H+b_{i}R)\overset{\psi}{\rightarrow}%
\]
}

{\normalsize
\[
\overset{\psi}{\rightarrow}%
{\displaystyle\sum\limits_{i=1}^{5}}
\mathcal{O}_{\mathbb{P(\mathcal{E})}}(-2H+a_{i}R)\rightarrow\mathcal{O}%
_{\mathbb{P(\mathcal{E})}}\rightarrow\mathcal{O}_{C}\rightarrow0
\]
}

{\normalsize Theorem \ref{psitypes0} then shows that we can restrict to $3$
different possibilities for $(a_{1},...,a_{5})$ and from the proposition
above, we can deduce further information: }

{\normalsize \bigskip}

\begin{theorem}
{\normalsize \label{psitypes2}Let C be a curve as given in Corollary
\ref{infinitypes0}, then the matrix $\psi$ has the following form:%
\[
\psi\sim\left(
\begin{array}
[c]{ccccc}%
0 & H+R & H & H & H-R\\
& 0 & H & H & H-R\\
&  & 0 & H-R & H-2R\\
&  &  & 0 & H-2R\\
&  &  &  & 0
\end{array}
\right)  \text{ }%
\]
}
\end{theorem}

\begin{proof}
{\normalsize In Theorem \ref{psitypes0} we have seen that there are only three
possible types for $\psi$, namely $(a_{1},...,a_{5})=(2,1,1,1,1),$
$(2,2,1,1,0)$ or $(2,2,2,0,0).$ From Corollary \ref{infinitypes0} above,
it turns out that $C$ is contained in a determinantal surface $Y$ given by the
$2\times2$ minors of a $2\times3$ matrix $\omega$ with global sections on
$\mathbb{P}(\mathcal{E})$. Therefore the $2\times2$ minors of $\omega$, which
are elements in $H^{0}(\mathbb{P}(\mathcal{E}),\mathcal{O}%
_{\mathbb{P(\mathcal{E})}}(2H-\tilde{a}_{i}R)),$ $i=1,2,3,$ are contained in
the ideal generated by the Pfaffians $\rho_{i}\in H^{0}(\mathbb{P}%
(\mathcal{E}),\mathcal{O}_{\mathbb{P(\mathcal{E})}}(2H-a_{i}R)),$ $i=1,...,5,$
of $\psi$. Now $(a_{1},...,a_{5})=(2,1,1,1,1)$ and $\mathbb{P}(\mathcal{E}%
)\cong \mathbb{P}(\mathcal{O}(3)\oplus\mathcal{O}(1)\oplus\mathcal{O}(1)\oplus \mathcal{O})$ can be omitted, as in this case the two minors of $\omega$,
which are elements of $H^{0}(\mathbb{P}(\mathcal{E}),\mathcal{O}%
_{\mathbb{P(\mathcal{E})}}(2H-2R)),$ are equal to a multiple of $\rho_{1}\in
H^{0}(\mathbb{P}(\mathcal{E}),\mathcal{O}_{\mathbb{P(\mathcal{E})}}(2H-2R))$,
which is not possible as $Y$ is irreducible and non degenerate. }

The following arguments {\normalsize show that the case $(a_{1},...,a_{5}%
)=(2,1,1,1,1)$ does not occur for $\mathbb{P}(\mathcal{E})\cong \mathbb{P}(\mathcal{O}(2)\oplus\mathcal{O}(2)\oplus\mathcal{O}(1)\oplus \mathcal{O}),$
too: In this case $C$ is given by the Pfaffians of
\[
\psi\sim\left(
\begin{array}
[c]{ccccc}%
0 & H & H & H & H\\
& 0 & H-R & H-R & H-R\\
&  & 0 & H-R & H-R\\
&  &  & 0 & H-R\\
&  &  &  & 0
\end{array}
\right)
\]
We already know that $C$ is contained in a determinantal surface $Y$ given by
the minors of
\[
\omega\sim\left(
\begin{array}
[c]{ccc}%
H & H & H+R\\
H-2R & H-2R & H-R
\end{array}
\right)
\]
on $\mathbb{P}(\mathcal{E}).$ We first calculate the Betti table for
$Y\subset$ $\mathbb{P}^{8}:$ Theorem \ref{Resolutiononscroll}\ tells us that
$\mathcal{O}_{Y}$ has an $\mathcal{O}_{\mathbb{P(\mathcal{E})}}$-module
resolution of type
\begin{align*}
F_{\ast} \text{  }  : & \text{    } 0 & \rightarrow\mathcal{O}_{\mathbb{P(\mathcal{E})}}(-3H+R)\oplus
\mathcal{O}_{\mathbb{P(\mathcal{E})}}(-3H+3R)\overset{\omega}{\rightarrow}%
\end{align*}%
\[
\overset{\omega}{\rightarrow}\mathcal{O}_{\mathbb{P(\mathcal{E})}%
}(-2H+R)^{\oplus2}\oplus\mathcal{O}_{\mathbb{P(\mathcal{E})}}%
(-2H+2R)\rightarrow\mathcal{O}_{\mathbb{P(\mathcal{E})}}\rightarrow
\mathcal{O}_{Y}\rightarrow0
\]
The corresponding mapping cone }

{\footnotesize
\begin{align*}
& \\
&  \bigskip%
\begin{array}
[c]{cccccc}%
0 & \rightarrow & \mathcal{O}_{\mathbb{P(\mathcal{E})}}(-3H+3R) & \oplus &
\mathcal{O}_{\mathbb{P(\mathcal{E})}}(-3H+R) & ^{\underrightarrow
{\text{\ \ \ \ \ \ \ \ \ \ \ \ \ \ \ \ \ \ }}}\\
&  & \uparrow &  & \uparrow & \\
&  & -----|----- &  & -----|----- & -\\
&  &  &  &  & \\
&  & S_{3}G(-3) & \oplus & G(-3) & ^{\underrightarrow
{\text{\ \ \ \ \ \ \ \ \ \ \ \ \ \ \ \ \ \ }}}\\
&  & \uparrow &  & \uparrow & \\
&  & F\otimes S_{2}G(-4) & \oplus & F(-4) & ^{\underrightarrow
{\text{\ \ \ \ \ \ \ \ \ \ \ \ \ \ \ \ \ \ }}}\\
&  & \uparrow &  & \uparrow & \\
&  & \wedge^{2}F\otimes G(-5) & \oplus & \wedge^{3}F(-6) & ^{\underrightarrow
{\text{\ \ \ \ \ \ \ \ \ }\gamma\text{\ \ \ \ \ \ \ }}}\\
&  & \uparrow &  & \uparrow & \\
&  & \wedge^{3}F(-6) & \oplus & DG^{\ast}\otimes\wedge^{4}F(-7) &
^{\underrightarrow{\text{\ \ \ \ \ \ \ \ \ \ \ \ \ \ \ \ \ \ }}}\\
&  & \uparrow &  & \uparrow & \\
&  & \wedge^{5}F(-8) & \oplus & D_{2}G^{\ast}(-8) & ^{\underrightarrow
{\text{\ \ \ \ \ \ \ \ \ \ \ \ \ \ \ \ \ \ }}}%
\end{array}
\\
&
\end{align*}%
\begin{align*}
&
\begin{array}
[c]{ccccccc}%
\mathcal{O}_{\mathbb{P(\mathcal{E})}}(-2H+R)^{\oplus2} & \oplus &
\mathcal{O}_{\mathbb{P(\mathcal{E})}}(-2H+2R) & ^{\underrightarrow
{\text{\ \ \ \ \ \ \ \ \ \ \ \ \ \ \ \ \ \ }}} & \mathcal{O}%
_{\mathbb{P(\mathcal{E})}} & \rightarrow & 0\\
\uparrow &  & \uparrow &  & \uparrow &  & \\
----|---- &  & ----|---- &  & ----|---- &  & \\
&  &  &  &  &  & \\
G(-2)^{\oplus2} & \oplus & S_{2}G(-2)^{\oplus2} & ^{\underrightarrow
{\text{\ \ \ \ \ \ \ \ \ \ \ \ \ \ \ \ \ \ }}} & \mathcal{O}_{\mathbb{P}^{8}}
&  & \\
\uparrow &  & \uparrow &  & \uparrow &  & \\
F(-3)^{\oplus2} & \oplus & F\otimes G(-3)^{\oplus2} & ^{\underrightarrow
{\text{\ \ \ \ \ \ \ \ \ \ \ \ \ \ \ \ \ \ }}} & \wedge^{2}F(-2) &  & \\
\uparrow &  & \uparrow &  & \uparrow &  & \\
\wedge^{3}F(-5)^{\oplus2} & \oplus & \wedge^{2}F(-4)^{\oplus2} &
^{\underrightarrow{\text{\ \ \ \ \ \ \ \ \ \ \ \ \ \ \ \ \ \ }}} & \wedge
^{3}F\otimes DG^{\ast}(-3) &  & \\
\uparrow &  & \uparrow &  & \uparrow &  & \\
\wedge^{4}F\otimes DG^{\ast}(-6)^{\oplus2} & \oplus & \wedge^{4}%
F(-6)^{\oplus2} & ^{\underrightarrow
{\text{\ \ \ \ \ \ \ \ \ \ \ \ \ \ \ \ \ \ }}} & \wedge^{4}F\otimes
D_{2}G^{\ast}(-4) &  & \\
\uparrow &  & \uparrow &  & \uparrow &  & \\
D_{2}G^{\ast}(-7)^{\oplus2} & \oplus & DG^{\ast}(-7)^{\oplus2} &
^{\underrightarrow{\text{\ \ \ \ \ \ \ \ \ \ \ \ \ \ \ \ \ \ }}} &
D_{3}G^{\ast}(-5) &  &
\end{array}
\\
&
\end{align*}}

\noindent gives us a (not necessarily) minimal free resolution of $\mathcal{O}_{Y}$. We
calculate the rank of the non minimal map %
\[
\gamma:\wedge^{2}F\otimes G(-5)\overset{\alpha}{\rightarrow}\wedge
^{3}F(-5)^{\oplus2}%
\]
which is obtained from the submatrix $\left(
\begin{array}
[c]{cc}%
\omega_{22} & \omega_{23}%
\end{array}
\right)  \sim\left(
\begin{array}
[c]{cc}%
H-2R & H-2R
\end{array}
\right)  :$%
\[
\gamma:f_{1}\wedge f_{2}\otimes g\rightarrow\binom{f_{1}\wedge f_{2}\wedge
g\omega_{22}}{f_{1}\wedge f_{2}\wedge g\omega_{23}}%
\]
with $f_{1},...,f_{4}\in H^{0}(\mathbb{P(\mathcal{E})},\mathcal{O}%
_{\mathbb{P(\mathcal{E})}}(H-R))\cong$
$H^{0}(\mathbb{P}^{8},\mathcal{O}_{\mathbb{P}^{8}}^{5})$
{\normalsize  and $g_{1},g_{2}\in H^{0}(\mathbb{P(\mathcal{E})}%
,\mathcal{O}_{\mathbb{P(\mathcal{E})}}(R))\cong H^{0}(\mathbb{P}^{8},\mathcal{O}_{\mathbb{P}^{8}}^{2})$. It turns out that $\dim\ker\gamma=4$:%
\begin{align*}
& \ker\gamma=\langle s\omega_{22}\wedge s\omega_{23}\otimes s,~t\omega
_{22}\wedge t\omega_{23}\otimes t,\\
&  s\omega_{22}\wedge t\omega_{23}\otimes s+t\omega_{22}\wedge s\omega
_{23}\otimes s+s\omega_{22}\wedge s\omega_{23}\otimes t,\\
&  t\omega_{22}\wedge t\omega_{23}\otimes s+s\omega_{22}\wedge t\omega
_{23}\otimes t+t\omega_{22}\wedge s\omega_{23}\otimes t\rangle
\end{align*}
hence the Betti table for $Y\subset\mathbb{P}^{8}$ is given by%
\textnormal{\[
\begin{tabular}{c|ccccccccc}
& 0 & 1 & 2 & 3 & 4 & 5 & 6 \cr \hline
0 & 1 & - & - & - & - & - & - \cr
1 & - & 17 & 46 & 45 & 8 & - & - \cr
2 & - & - & - & 4 & 25  & 18  & 4 \cr
\end{tabular} 
\]}
The linear strand of the minimal free resolution for $Y$ is a subcomplex of
that for $C\subset Y.$ It follows that the Betti number $\beta_{45}(C)$
in the minimal free resolution of $\mathcal{O}_{C}$ has to be at least 8. Now
if $C$ is given by the Pfaffians of a matrix
\[
\psi\sim\left(
\begin{array}
[c]{ccccc}%
0 & H & H & H & H\\
& 0 & H-R & H-R & H-R\\
&  & 0 & H-R & H-R\\
&  &  & 0 & H-R\\
&  &  &  & 0
\end{array}
\right)
\]
then from Lemma \ref{rank} and Remark \ref{rankremark} we already know that $\beta_{45}(C)\geq8$ is only
possible if one of the $(H-R)$-entries of $\psi$ can be made zero by suitable
row and column operations. But then we obtain a second special linear series
of type $g_{5}^{1}$ that is different from the one we have started with. This contradicts to our assumption that $C$ has only one $g_{5}^{1}$. }

\noindent {\normalsize Let us now consider the case $(a_{1},...,a_{5})=(2,2,2,0,0): $ For
$\mathbb{P}(\mathcal{E})\cong \mathbb{P}(\mathcal{O}(2)\oplus\mathcal{O}(2)\oplus\mathcal{O}(1)\oplus \mathcal{O})$, $\psi$ is of type}
{\scriptsize \[
\psi\sim\left(
\begin{array}
[c]{ccccc}%
0 & H+R & H+R & H-R & H-R\\
& 0 & H+R & H-R & H-R\\
&  & 0 & H-R & H-R\\
&  &  & 0 & H-3R\\
&  &  &  & 0
\end{array}
\right)  \sim\left(
\begin{array}
[c]{ccccc}%
0 & H+R & H+R & \mathbf{H-R} & \mathbf{H-R}\\
& 0 & H+R & \mathbf{H-R} & \mathbf{H-R}\\
&  & 0 & \mathbf{H-R} & \mathbf{H-R}\\
&  &  & 0 & 0\\
&  &  &  & 0
\end{array}
\right)
\]}
It follows that $C$ is contained in a determinantal surface $Y$ given by the
minors of
\[
\omega\sim\left(
\begin{array}
[c]{ccc}%
H-R & H-R & H-R\\
H-R & H-R & H-R
\end{array}
\right)
\]
on $\mathbb{P}(\mathcal{E}).$ According to Theorem \ref{ConPk}. $Y$ is a
blowup of $\mathbb{P}^{1}\times\mathbb{P}^{1}$ in 3 points. The direct image
$C^{\prime}$ of $C$ in $\mathbb{P}^{1}\times\mathbb{P}^{1}$ is a divisor of
type $(5,4)$, thus we get a $g_{4}^{1}$ from projection onto the second factor
of $\mathbb{P}^{1}\times\mathbb{P}^{1}$, which we have excluded.

{\normalsize It remains to examine the same possibility for $(a_{1}%
,...,a_{5})$ if $X\simeq S(3,1,1,0):$ After suitable row and column
operations, one of the $(H-R)$-entries in the last column of $\psi$ can be made
to zero. Assuming $\psi_{15}=0,$ our matrix takes the following form%
\[
\psi\sim\left(
\begin{array}
[c]{ccccc}%
0 & H+R & H+R & \varphi^{\prime} & 0\\
& 0 & H+R & \varphi^{\prime\prime} & \mathbf{\varphi}_{2}\\
&  & 0 & \mathbf{\varphi}_{2} & \mathbf{\varphi}_{1}\\
&  &  & 0 & \mathbf{\varphi}_{0}\\
&  &  &  & 0
\end{array}
\right)
\]
with linear independent, global sections $\mathbf{\varphi}_{0}\in H^{0}%
(\mathbb{P(\mathcal{E})},\mathcal{O}_{\mathbb{P(\mathcal{E})}}(H-3R\mathbb{))}%
$, $\mathbf{\varphi}_{1},\mathbf{\varphi}_{2}\in H^{0}(\mathbb{P(\mathcal{E}%
)},\mathcal{O}_{\mathbb{P(\mathcal{E})}}(H-R\mathbb{))}$ and $\varphi^{\prime
},\varphi^{\prime\prime}\in\left\langle \mathbf{\varphi}_{1},\mathbf{\varphi
}_{2}\right\rangle .$\ We apply Theorem \ref{ConPk} again to see that $C$ is
contained in a determinantal surface $Y$ given by the matrix%
\[
\omega:=\left(
\begin{array}
[c]{ccc}%
\psi_{12} & \psi_{13} & \psi_{14}\\
\psi_{25} & \psi_{35} & \psi_{45}%
\end{array}
\right)  \sim\left(
\begin{array}
[c]{ccc}%
H+R & H+R & \mathbf{\varphi}^{\prime}\\
\mathbf{\varphi}_{2} & \mathbf{\varphi}_{1} & \mathbf{\varphi}_{0}%
\end{array}
\right)
\]
on $\mathbb{P(\mathcal{E})}$. The image of this surface in $\mathbb{P}^{8}$ is
given as the image of $P_{2}=\mathbb{P(\mathcal{O}}_{\mathbb{P}^{1}%
}\mathbb{\mathcal{(}}2)+\mathbb{\mathcal{O}}_{\mathbb{P}^{1}}\mathbb{)}$ under
a rational map defined by a subseries $H^{0}(P_{2},\mathbb{\mathcal{O}}%
_{P_{2}}(3A))$ which has 7 assigned base points. The image $C^{\prime}$ of $C$
in $P_{2}$ is given by the Pfaffians of a matrix $\tilde{\psi}$ regarding the
entries of $\psi$ as global sections of vector bundles in $P_{2}.$ As the two
rows of $\omega$ are linear dependend, we get
\[
\tilde{\psi}\sim\left(
\begin{array}
[c]{ccccc}%
0 & \lambda g_{2} & \lambda g_{1} & \varphi^{\prime}=\lambda g_{0} & 0\\
& 0 & 3A+B & \varphi^{\prime\prime} & \mu g_{2}\\
&  & 0 & \mu g_{2} & \mu g_{1}\\
&  &  & 0 & \mu g_{0}\\
&  &  &  & 0
\end{array}
\right)
\]
with elements $\lambda\in H^{0}(P_{2},\mathcal{O}_{P_{2}}(A)\mathbb{)}$,
$\mu\in H^{0}(P_{2},\mathcal{O}_{P_{2}}(A-2B)\mathbb{)}$ and $g_{1},g_{2}\in
H^{0}(P_{2},\mathcal{O}_{P_{2}}(2A+B)\mathbb{)},\ g_{0}\in H^{0}%
(P_{2},\mathcal{O}_{P_{2}}(2A-B)\mathbb{)}$. Because of $\varphi^{\prime
},\varphi^{\prime\prime}\in\left\langle \mathbf{\varphi}_{1},\mathbf{\varphi
}_{2}\right\rangle =\left\langle \mu g_{1},\mu g_{2}\right\rangle ,$ $\mu$ is
a common factor of the last two columns, hence it is a common factor of all Pfaffians of $\tilde{\psi}.$ This contradicts that $C$ is assumed to be
irreducible. It remains only one possibility for $(a_{1},...,a_{5}),$ i.e.
$(a_{1},...,a_{5})=(2,2,1,1,0).$ }
\end{proof}

{\normalsize \bigskip}

\begin{theorem}
{\normalsize \label{2infinig15}Let C be an irreducible, nonsingular, canonical
curve of genus 9 with $\operatorname*{Cliff}(C)=3$ that admits no $g_{7}^{2}%
$. If C has exactly one $g_{5}^{1}$ with multiplicitiy two and no further
$g_{5}^{1}$ then it is given by the Pfaffians of a matrix $\psi$%
\[
\psi\sim\left(
\begin{array}
[c]{ccccc}%
0 & H+R & H & H & H-R\\
& 0 & H & H & H-R\\
&  & 0 & H-R & H-2R\\
&  &  & 0 & H-2R\\
&  &  &  & 0
\end{array}
\right)
\]
on a scroll of type $S(2,2,1,0)$. The minimal free resolution of $\mathcal{O}_{C}$ as $\mathcal{O}_{\mathbb{P}^{8}}-$modul has the following
Betti table:%
\textnormal{\[
\begin{tabular}{c|cccccccccc}
& 0 & 1 & 2 & 3 & 4 & 5 & 6 & 7 \cr \hline
0 & 1 & - & - & - & - & - & - & - \cr
1 & - & 21 & 64 & 70 & 8 & - & - & - \cr
2 & - & - & - & 8 & 70  & 64  & 21 & -\cr
3 & - & - & - & - & - & - & - & 1 \cr
\end{tabular} 
\]}
}
\end{theorem}

\begin{proof}
{\normalsize Let $\mathbb{P}(\mathcal{E})$ be the corresponding $\mathbb{P}%
^{3}$-bundle of type $S(2,2,1,0)$ constructed from a unique $g_{5}%
^{1}=\left\vert D\right\vert $ of multiplicity 2, then $C$ is given by the
Pfaffians of a matrix $\psi$ as in the above theorem. By suitable row and
column operations, $\psi$ has the following form:%
\[
\psi\sim\left(
\begin{array}
[c]{ccccc}%
0 & H+R & H & H & 0\\
& 0 & H & H & \varphi_{2}\\
&  & 0 & \varphi & \varphi_{1}\\
&  &  & 0 & \varphi_{0}\\
&  &  &  & 0
\end{array}
\right)
\]
with independent global sections $\varphi_{0},\varphi_{1}\in H^{0}%
(\mathbb{P}(\mathcal{E}),\mathcal{O}_{\mathbb{P(\mathcal{E})}}(H-2R))$ and
$\varphi_{2}\in H^{0}(\mathbb{P}(\mathcal{E}),\mathcal{O}%
_{\mathbb{P(\mathcal{E})}}(H-R)).$ We can assume that $\varphi=0$ or
$\varphi=\varphi_{2}.$ If $\varphi=0,$ then $C$ is contained in the
determinantal surface $Y$ given by the $2\times2$ minors of
\[
\omega:=\left(
\begin{array}
[c]{ccc}%
\psi_{31} & \psi_{32} & \psi_{35}\\
\psi_{14} & \psi_{24} & \psi_{54}%
\end{array}
\right)  \sim\left(
\begin{array}
[c]{ccc}%
H & H & H-2R\\
H & H & H-2R
\end{array}
\right)
\]}
\noindent We apply Theorem \ref{ConPk} to see that $Y$ is a blow up of $\mathbb{P}%
^{1}\times\mathbb{P}^{1}$ in 7 points and the direct image $C^{\prime}%
\subset\mathbb{P}^{1}\times\mathbb{P}^{1}$ of $C$ is a divisor of type $(5,5)$
on $\mathbb{P}^{1}\times\mathbb{P}^{1}$. 

\noindent \normalsize {Therefore $C$ would have two different special linear series of type $g_{5}^{1}$, which contradicts our
assumptions. It follows that $\varphi=\varphi_{2}.$ Now we want to calculate
the Betti table for the minimal free resolution of $\mathcal{O}_{C}.$
We have to look at the non minimal maps in the corresponding mapping
cone}%

{\footnotesize \begin{align*}
&
\begin{array}
[c]{cccccc}%
\mathcal{O}_{\mathbb{P(\mathcal{E})}}(-3H+R)^{\oplus2} & \oplus &
\mathcal{O}_{\mathbb{P(\mathcal{E})}}(-3H+2R)^{\oplus2} & \oplus &
\mathcal{O}_{\mathbb{P(\mathcal{E})}}(-3H+3R) & ^{\underrightarrow
{\text{\ \ \ \ \ \ \ \ \ \ \ \ \ \ \ \ \ \ }}}\\
\uparrow &  & \uparrow &  & \uparrow & \\
-----|----- & - & ----|---- & - & -----|----- & -----\\
&  &  &  &  & \\
G(-3)^{\oplus2} & \oplus & S_{2}G(-3)^{\oplus2} & \oplus & S_{3}G(-3) &
^{\underrightarrow{\text{\ \ \ \ \ \ \ \ \ \ \ \ \ \ \ \ \ \ }}}\\
\uparrow &  & \uparrow &  & \uparrow & \\
F(-4)^{\oplus2} & \oplus & F\otimes G(-4)^{\oplus2} & \oplus & S_{2}G\otimes
F(-4) & ^{\underrightarrow{\text{\ \ \ \ \ \ \ \ \ \ \ \ \ \ \ \ \ \ }}}\\
\uparrow &  & \uparrow &  & \uparrow & \\
\wedge^{3}F(-6)^{\oplus2} & \oplus & \wedge^{2}F(-5)^{\oplus2} & \oplus &
\wedge^{2}F\otimes G(-5) & ^{\underrightarrow
{\text{\ \ \ \ \ \ \ \ \ \ \ \ \ \ \ \ \ \ }}}\\
\uparrow &  & \uparrow &  & \uparrow & \\
\wedge^{4}F\otimes DG^{\ast}(-7)^{\oplus2} & \oplus & \wedge^{4}%
F(-7)^{\oplus2} & \oplus & \wedge^{3}F(-6) & ^{\underrightarrow
{\text{\ \ \ \ \ \ \ \ \ \ \ \ \ \ \ \ \ \ }}}\\
\uparrow &  & \uparrow &  & \uparrow & \\
D_{2}G^{\ast}(-8)^{\oplus2} & \oplus & DG^{\ast}(-8)^{\oplus2} & \oplus &
\wedge^{5}F(-8) & ^{\underrightarrow
{\text{\ \ \ \ \ \ \ \ \ \ \ \ \ \ \ \ \ \ }}}%
\end{array}
\\
& \\
& \\
& \\
&
\begin{array}
[c]{cccccc}%
^{\underrightarrow{\text{\ \ \ \ \ \ \ \ \ \ \ \ \ \ \ \ \ \ }}} &
\mathcal{O}_{\mathbb{P(\mathcal{E})}}(-2H+2R)^{\oplus2} & \oplus &
\mathcal{O}_{\mathbb{P(\mathcal{E})}}(-3H+R)^{\oplus2} & \oplus &
\mathcal{O}_{\mathbb{P(\mathcal{E})}}(-2H)\\
& \uparrow &  & \uparrow &  & \uparrow\\
----- & -----|----- &  & -----|----- & - & -----|-----\\
&  &  &  &  & \\
^{\underrightarrow{\text{\ \ \ \ \ \ \ \ \ \ \ \ \ \ \ \ \ \ }}} &
S_{2}G(-2)^{\oplus2} & \oplus & G(-2)^{\oplus2} & \oplus & \mathcal{O}(-2)\\
& \uparrow &  & \uparrow &  & \uparrow\\
^{\underrightarrow{\text{\ \ \ \ \ \ \ \ \ \ \ \ \ \ \ \ \ \ }}} & F\otimes
G(-3)^{\oplus2} & \oplus & F(-3)^{\oplus2} & \oplus & \wedge^{2}F(-4)\\
& \uparrow &  & \uparrow &  & \uparrow\\
^{\underrightarrow{\text{\ \ \ \ \ \ \ \ \ \ \ \ \ \ \ \ \ \ }}} & \wedge
^{2}F(-4)^{\oplus2} & \oplus & \wedge^{3}F(-5)^{\oplus2} & \oplus & \wedge
^{3}F\otimes DG^{\ast}(-5)\\
& \uparrow &  & \uparrow &  & \uparrow\\
^{\underrightarrow{\text{\ \ \ \ \ \ \ \ \ \ \ \ \ \ \ \ \ \ }}} & \wedge
^{4}F(-6)^{\oplus2} & \oplus & \wedge^{4}F\otimes DG^{\ast}(-6)^{\oplus2} &
\oplus & \wedge^{4}F\otimes D_{2}G^{\ast}(-6)\\
& \uparrow &  & \uparrow &  & \uparrow\\
^{\underrightarrow{\text{\ \ \ \ \ \ \ \ \ \ \ \ \ \ \ \ \ \ }}} & \wedge
^{5}F\otimes DG^{\ast}(-7)^{\oplus2} & \oplus & D_{2}G^{\ast}(-7)^{\oplus2} &
\oplus & \wedge^{5}F\otimes D_{3}G^{\ast}(-7)
\end{array}
\end{align*}}
\bigskip 
\\
which are given by the following maps:%
\[
\alpha:%
\begin{array}
[c]{ccccc}%
F(-4)^{\oplus2} & \oplus & F\otimes G(-4)^{\oplus2} & \oplus & S_{2}G\otimes
F(-4)
\end{array}
\rightarrow\wedge^{2}F(-4)
\]%
\[
\left(
\begin{array}
[c]{c}%
f_{1}\\
f_{2}\\
g_{1}\otimes f_{3}\\
g_{2}\otimes f_{4}\\
g_{3}^{(2)}\otimes f_{5}%
\end{array}
\right)  \rightarrow f_{2}\wedge\varphi_{2}+g_{1}\varphi_{1}\wedge f_{3}%
+g_{2}\varphi_{0}\wedge f_{4}%
\]
\bigskip its dual%
\[
\alpha^{\ast}:\wedge^{3}F(-6)\rightarrow%
\begin{array}
[c]{ccccc}%
\wedge^{4}F\otimes DG^{\ast}(-6)^{\oplus2} & \oplus & \wedge^{4}F\otimes
D_{2}G^{\ast}(-6) & \oplus & \wedge^{4}F(-6)^{\oplus2}%
\end{array}
\]
and%
\[
\beta:%
\begin{array}
[c]{ccc}%
\wedge^{2}F\otimes G(-5) & \oplus & \wedge^{2}F(-5)^{\oplus2}%
\end{array}
\rightarrow%
\begin{array}
[c]{ccc}%
\wedge^{3}F(-5)^{\oplus2} & \oplus & \wedge^{3}F\otimes DG^{\ast}(-5)
\end{array}
\]%
\[
\left(
\begin{array}
[c]{c}%
f_{1}\wedge f_{2}\\
f_{3}\wedge f_{4}\\
g\otimes f_{5}\wedge f_{6}%
\end{array}
\right)  \rightarrow\left(
\begin{array}
[c]{c}%
-\varphi_{2}\wedge f_{3}\wedge f_{4}-g\varphi_{1}\wedge f_{5}\wedge f_{6}\\
\varphi_{2}\wedge f_{1}\wedge f_{2}-g\varphi_{0}\wedge f_{5}\wedge f_{6}\\
(s\varphi_{1}\wedge f_{1}\wedge f_{2}\otimes s^{\ast}+t\varphi_{1}\wedge
f_{1}\wedge f_{2}\otimes t^{\ast}+\\
s\varphi_{0}\wedge f_{3}\wedge f_{4}\otimes s^{\ast}+t\varphi_{0}\wedge
f_{3}\wedge f_{4}\otimes t^{\ast})
\end{array}
\right)
\]
with $f_{1},..,f_{6}\in H^{0}F,~g,g_{1},g_{2}\in H^{0}G,\ g_{3}^{(2)}\in H^{0}S_{2}G~$and
$G=\left\langle s,t\right\rangle ,~H^{0}DG^{\ast}=\left\langle s^{\ast},t^{\ast
}\right\rangle .$ Trivially $\alpha$ is surjective as every element in $\wedge^{2}F(-4)$ is a sum of elements of type $f\wedge\varphi_{2}$, $g\varphi_{1}\wedge f$ and $g\varphi_{0}\wedge f$ with $f\in F$ and $g\in G$. Thus, if $C$ is given by a matrix $\psi$ as above, then we get $\beta_{24}=0$. A calculation
of $\ker\beta$ gives (cf. Appendix 6.3)%
{\small \[
\ker\beta=\left\langle
\begin{array}
[c]{c}%
\left(
\begin{array}
[c]{c}%
0\\
0\\
t\otimes(t\varphi_{0}\wedge t\varphi_{1})
\end{array}
\right)  ,\left(
\begin{array}
[c]{c}%
0\\
0\\
s\otimes(s\varphi_{0}\wedge s\varphi_{1})
\end{array}
\right)  ,\\\\
\left(
\begin{array}
[c]{c}%
0\\
0\\
s\otimes(t\varphi_{0}\wedge t\varphi_{1})+t\otimes(t\varphi_{0}\wedge
s\varphi_{1}+s\varphi_{0}\wedge t\varphi_{1})
\end{array}
\right)  ,\\\\
\left(
\begin{array}
[c]{c}%
0\\
0\\
s\otimes(t\varphi_{0}\wedge s\varphi_{1}+s\varphi_{0}\wedge t\varphi
_{1})+t\otimes(s\varphi_{0}\wedge s\varphi_{1})
\end{array}
\right)
\end{array}
\right\rangle
\]}
Therefore, the Betti table for $C$ is given as follows{\normalsize
\textnormal{\[
\begin{tabular}{c|cccccccccc}
& 0 & 1 & 2 & 3 & 4 & 5 & 6 & 7 \cr \hline
0 & 1 & - & - & - & - & - & - & - \cr
1 & - & 21 & 64 & 70 & 8 & - & - & - \cr
2 & - & - & - & 8 & 70  & 64  & 21 & -\cr
3 & - & - & - & - & - & - & - & 1 \cr
\end{tabular} 
\]}
}
\end{proof}

{\normalsize \bigskip}

{\normalsize It remains to determine the Betti table }for $C$ {\normalsize  if
$C$ lies on a scroll of type $S(3,1,1,0),$ i.e. we must have $h^{0}%
(C,\mathcal{O(}K_{C}-3D))=1$ for the unique $g_{5}^{1}=\left\vert D\right\vert
$ on $C:$ }

\begin{theorem}
{\normalsize \label{3infinig15}Let C be an irreducible, smooth, canonical
curve of genus 9 with $\operatorname*{Cliff}(C)=3$, that admits no $g_{7}^{2}%
$. If C has exactly one $g_{5}^{1}$ with multiplicitiy $3$ then it is given by the Pfaffians of a matrix $\psi$%
\[
\psi\sim\left(
\begin{array}
[c]{ccccc}%
0 & H+R & H & H & H-R\\
& 0 & H & H & H-R\\
&  & 0 & H-R & H-2R\\
&  &  & 0 & H-2R\\
&  &  &  & 0
\end{array}
\right)
\]
on a scroll of type $S(3,1,1,0)$. The minimal free resolution of $\mathcal{O}_{C}$ as $\mathcal{O}_{\mathbb{P}^{8}}-$module has the following
Betti diagram:%
\textnormal{\[
\begin{tabular}{c|cccccccccc}
& 0 & 1 & 2 & 3 & 4 & 5 & 6 & 7 \cr \hline
0 & 1 & - & - & - & - & - & - & - \cr
1 & - & 21 & 64 & 70 & 12 & - & - & - \cr
2 & - & - & - & 12 & 70  & 64  & 21 & -\cr
3 & - & - & - & - & - & - & - & 1 \cr
\end{tabular} 
\]}
}
\end{theorem}

\begin{proof}
{\normalsize Let }$X$ be the scroll $X$ of {\normalsize  type $S(3,1,1,0)$
constructed from the unique $g_{5}^{1}=\left\vert D\right\vert $ and
$\mathbb{P}(\mathcal{E})$ the corresponding $\mathbb{P}^{3}$-bundle, then $C$
is given by the Pfaffians of a matrix $\psi$ as above. By suitable row and
column operations, $\psi$ gets the following form:%
\[
\psi\sim\left(
\begin{array}
[c]{ccccc}%
0 & H+R & H & H & H-R\\
& 0 & H & H & \varphi^{\prime}\\
&  & 0 & \varphi & t\varphi_{0}\\
&  &  & 0 & s\varphi_{0}\\
&  &  &  & 0
\end{array}
\right)
\]
with global sections $\varphi_{0}\in H^{0}(\mathbb{P}%
(\mathcal{E}),\mathcal{O}_{\mathbb{P(\mathcal{E})}}(H-3R))$ and $\varphi
,\varphi^{\prime}\in H^{0}(\mathbb{P}(\mathcal{E}),\mathcal{O}%
_{\mathbb{P(\mathcal{E})}}(H-R)).$ As in the proof of Theorem \ref{2infinig15} $varphi$ cannot be made zreo by suitable row and column operations. The calculation of the minimal free
resolution of $\mathcal{O}_{C}$ by a mapping cone construction can be done in an
analogous way as in the case where }$X$ is of type $S(2,2,1,0).$ {\normalsize
The only difference in this case is, that the non minimal maps are given by
\[
\alpha:%
\begin{array}
[c]{ccccc}%
F(-4)^{\oplus2} & \oplus & F\otimes G(-4)^{\oplus2} & \oplus & S_{2}G\otimes
F(-4)
\end{array}
\rightarrow\wedge^{2}F(-4)
\]%
\[
\left(
\begin{array}
[c]{c}%
f_{1}\\
f_{2}\\
g_{1}\otimes f_{3}\\
g_{2}\otimes f_{4}\\
g_{3}^{(2)}\otimes f_{5}%
\end{array}
\right)  \rightarrow f_{2}\wedge\varphi^{\prime}+g_{1}t\varphi_{0}\wedge
f_{3}+g_{2}s\varphi_{0}\wedge f_{4}%
\]
\bigskip its dual%
\[
\alpha^{\ast}:\wedge^{3}F(-6)\rightarrow%
\begin{array}
[c]{ccccc}%
\wedge^{4}F\otimes DG^{\ast}(-6)^{\oplus2} & \oplus & \wedge^{4}F\otimes
D_{2}G^{\ast}(-6) & \oplus & \wedge^{4}F(-6)^{\oplus2}%
\end{array}
\]
and%
\[
\beta:%
\begin{array}
[c]{ccc}%
\wedge^{2}F\otimes G(-5) & \oplus & \wedge^{2}F(-5)^{\oplus2}%
\end{array}
\rightarrow%
\begin{array}
[c]{ccc}%
\wedge^{3}F(-5)^{\oplus2} & \oplus & \wedge^{3}F\otimes DG^{\ast}(-5)
\end{array}
\]%
\[
\left(
\begin{array}
[c]{c}%
f_{1}\wedge f_{2}\\
f_{3}\wedge f_{4}\\
g\otimes f_{5}\wedge f_{6}%
\end{array}
\right)  \rightarrow\left(
\begin{array}
[c]{c}%
-\varphi\wedge f_{3}\wedge f_{4}-gt\varphi_{0}\wedge f_{5}\wedge f_{6}\\
\varphi\wedge f_{1}\wedge f_{2}-gs\varphi_{0}\wedge f_{5}\wedge f_{6}\\
(st\varphi_{0}\wedge f_{1}\wedge f_{2}\otimes s^{\ast}+t^{2}\varphi_{0}\wedge
f_{1}\wedge f_{2}\otimes t^{\ast}+\\
s^{2}\varphi_{0}\wedge f_{3}\wedge f_{4}\otimes s^{\ast}+st\varphi_{0}\wedge
f_{3}\wedge f_{4}\otimes t^{\ast})
\end{array}
\right)
\]
with $f_{1},..,f_{6}\in H^{0}F,~g,g_{1},g_{2}\in H^{0}G,\ g_{3}^{(2)}\in H^{0}S_{2}G~$and
$G=\left\langle s,t\right\rangle ,~H^{0}DG^{\ast}=\left\langle s^{\ast},t^{\ast
}\right\rangle .$ $\varphi^{\prime}$ cannot be obtained from a linear
combination of $s^{2}\varphi_{0},st\varphi_{0},t^{2}\varphi_{0}$, as then
$\varphi_{0}$ would be a factor of one of the Pfaffians, thus it is a direct
consequence that $\alpha$ is surjective. For $\beta$ a calculation shows (cf. Appendix 6.4), that}%
{\scriptsize \[
\ker\beta=\left\langle
\begin{array}
[c]{c}%
\left(
\begin{array}
[c]{c}%
st\varphi_{0}\wedge t^{2}\varphi_{0}\\
0\\
t\otimes(t^{2}\varphi_{0}\wedge\varphi)
\end{array}
\right)  ,\left(
\begin{array}
[c]{c}%
0\\
0\\
t\otimes(st\varphi_{0}\wedge t^{2}\varphi_{0})
\end{array}
\right)  ,\left(
\begin{array}
[c]{c}%
s^{2}\varphi_{0}\wedge t^{2}\varphi_{0}\\
st\varphi_{0}\wedge t^{2}\varphi_{0}\\
s\otimes(t^{2}\varphi_{0}\wedge\varphi)+2t\otimes(st\varphi_{0}\wedge\varphi)
\end{array}
\right)  ,\\\\
\left(
\begin{array}
[c]{c}%
0\\
0\\
s\otimes(st\varphi_{0}\wedge t^{2}\varphi_{0})+t\otimes(s^{2}\varphi_{0}\wedge
t^{2}\varphi_{0})
\end{array}
\right)  ,\left(
\begin{array}
[c]{c}%
0\\
s^{2}\varphi_{0}\wedge st\varphi_{0}\\
s\otimes(s^{2}\varphi_{0}\wedge\varphi)
\end{array}
\right)  ,\left(
\begin{array}
[c]{c}%
0\\
0\\
s\otimes(s^{2}\varphi_{0}\wedge st\varphi_{0})
\end{array}
\right)  ,\\\\
\left(
\begin{array}
[c]{c}%
s^{2}\varphi_{0}\wedge st\varphi_{0}\\
s^{2}\varphi_{0}\wedge t^{2}\varphi_{0}\\
2s\otimes(st\varphi_{0}\wedge\varphi)+t\otimes(s^{2}\varphi_{0}\wedge\varphi)
\end{array}
\right)  ,\left(
\begin{array}
[c]{c}%
0\\
0\\
s\otimes(s^{2}\varphi_{0}\wedge t^{2}\varphi_{0})+t\otimes(s^{2}\varphi
_{0}\wedge st\varphi_{0})
\end{array}
\right)
\end{array}
\right\rangle
\]}
Therefore the Betti table for $C$ {\normalsize  is given as follows:}%
\textnormal{\[
\begin{tabular}{c|cccccccccc}
& 0 & 1 & 2 & 3 & 4 & 5 & 6 & 7 \cr \hline
0 & 1 & - & - & - & - & - & - & - \cr
1 & - & 21 & 64 & 70 & 12 & - & - & - \cr
2 & - & - & - & 12 & 70  & 64  & 21 & -\cr
3 & - & - & - & - & - & - & - & 1 \cr
\end{tabular} 
\]}

\end{proof}

\bigskip

\subsection{Deformation\label{deformation}}

{\normalsize In the last two sections we have determined the Betti table for a
canonical curve $C\subset\mathbb{P}^{8}$ of genus 9 that has exactly $k=1,2$
or $3$ special linear series of type $g_{5}^{1}$ (counted with
multiplicities)$:$%
\[%
\begin{array}
[c]{cccccccc}%
1 & \cdot & \cdot & \cdot & \cdot & \cdot & \cdot & \cdot\\
\cdot & 21 & 64 & 70 & 4k & \cdot & \cdot & \cdot\\
\cdot & \cdot & \cdot & 4k & 70 & 64 & 21 & \cdot\\
\cdot & \cdot & \cdot & \cdot & \cdot & \cdot & \cdot & 1
\end{array}
\]
Our definition of the multiplicity of a $g_{5}^{1}$ was rather technically, so
it remains to give a geometric explanation for this procedure. In the situation
where }$C$ {\normalsize has three different ordinary $g_{5}^{1},$ we have seen
that in some cases it is possible for the third one to become equal to one of
the two others. In this sitiuation its multiplicity rises to 2. We want to show that in any case where the curve $C_{0}:=C$
has a $g_{5}^{1}=\left\vert D\right\vert ,$ }$D$ an effective divisor of
degree 5 on $C,${\normalsize \ with higher multiplicity, there exists a local
one parameter family $(C_{\lambda})_{\lambda\in\mathbb{A}^{1}}$ of curves in
$\mathcal{M}_{9}$, such that $C_{\lambda}$\ has the corresponding number of
different $g_{5}^{1}$ with ordinary multiplicity for $\lambda\neq0$. }

{\normalsize According to\textbf{ }Theorem \ref{cone} we can assume that $C$
has a space model $C^{\prime}$ given as the image of }$C$ under
the {\normalsize map defined by the complete linear series $\left\vert
2D\right\vert .$ Then }$C^{\prime}$ {\normalsize is contained in the cone
$Y\subset\mathbb{P}^{3}$ over the irreducible conic $x_{0}x_{2}=x_{1}^{2}$ in
$\mathbb{P}^{2}$ and has exactly 7 (possibly infinitely near)
double points $p_{i},$ $i=1,...,7,$ as only singularities. $C^{\prime}$ does not
pass through the vertex of $Y.$ The corresponding $\mathbb{P}^{1}-$bundle
$P_{2}=\mathcal{\mathbb{P(}O}_{\mathbb{P}^{1}}\mathcal{(}2)\oplus
\mathcal{O}_{\mathbb{P}^{1}})$ gives a resolution of the singularity of $Y.$
Let $A$ denote the class of a hyperplane divisor and $B$ that of a ruling on
$P_{2}.$ The curve $C^{\prime}\subset P_{2}$ is given
by an element $\tilde{F}_{0}\in\left\vert 5A-2p_{1}-...-2p_{7}\right\vert $ on
$P_{2}.$ We denote $Q_{0}$ the defining equation of $Y$ in $\mathbb{P}^{3}$
and $H$ the hyperplane divisor in $\mathbb{P}^{3}.$ Then $A$ is the pullback
of $H|_{Y}$ and thus $C^{\prime}\subset\mathbb{P}^{3}$ is the complete
intersection of a quintic $\Gamma=V(F_{0})$ given by an $F_{0}\in\left\vert
5H\right\vert $ and the cone $Y=V(Q_{0})$.} Blowing up $P_{2}$ in the points $p_{1},...,p_{7}$, we can assume that $C$ is a divisor of type $5A-\sum_{i=1}^{7}2E_{i}$ on $S:=\tilde{P_{2}}(p_{1},...,p_{7})$ where $A$ and $B$ denotes by abuse of notation also the pullback of the class of a hyperplane divisor and a ruling on $P_{2}$ respectively.
\\\\
We give an outline how to construct the local one paramter family with the
desired properties:

\bigskip

\noindent \textbf{1. }In the first step we separate infinitely near double points,
i.e. we give a local one paramter family of curves on the cone $Y$ that have
the same number of $g_{5}^{1\prime}s$ (respecting their multiplicities) as
$C^{\prime},$ but exactly $7$ \textit{distinct} double points. This shows that
we can assume for $C^{\prime}$ to have exactly $7$ distinct double
points.\bigskip

\noindent \textbf{2. }Then we show that we can choose $F_{0}$ in such a way that
{\normalsize $\Gamma=V(F_{0})$} has multiplicity 2 in each of the points $p_{i}$, i.e. the singularities of $C^{\prime}$ in the
points $p_{i}$ come from singularities of $\Gamma$ in $p_{i}$.

\bigskip

\noindent \textbf{3. }In the last step we get a deformation $C_{\lambda}^{\prime
}:=\Gamma\cap Y_{\lambda}$ of $C^{\prime}$ by fixing the quintic $\Gamma$ and
deforming $Y$ into a smooth quadric $Y_{\lambda}\cong\mathbb{P}^{1}%
\times\mathbb{P}^{1}$ that passes through the points $p_{i}.$ Then
$C_{\lambda}^{\prime}$ becomes a divisor of class $(5,5)$ on $\mathbb{P}%
^{1}\times\mathbb{P}^{1}$ with double points $p_{i}.$ This separates the
infinitesimal near $g_{5}^{1\prime}s$ of $C^{\prime}$ that are cut out by the
class of a ruling on the cone $Y$.

\bigskip

We start with the following lemma:

\begin{lemma}
Let $C^{\prime}\in|5A|$ be an irreducible curve on $P_{2}$ that has exactly
$7$ double points $p_{1},...,p_{7}$ as only singularities. Further assume that
$p_{7}$ lies infinitely near to $p_{6}$. Let $\Sigma\subset P_{2}$ be a
rational curve that passes through $p_{7}$ and $(q_{\lambda})_{\lambda
\in\mathbb{A}^{1}}$ be a local parametrisation of $\Sigma$ with $q_{0}=p_{7}%
.$ Then there exists a local one parameter family $(C_{\lambda}^{\prime
})_{\lambda\in\mathbb{A}^{1}}$ of curves in the divisor class $5A$ that have double points
in $p_{1},...,p_{6},q_{\lambda}$ as only singularities.
\end{lemma}

\begin{proof}
We {\normalsize blow up $P_{2}$ in the points $p_{1},...,p_{6},q_{\lambda}:$%
\[
\sigma_{\lambda}:S_{\lambda}=\tilde{P}_{2}(p_{1},...,p_{6,}q_{\lambda})\rightarrow
S=\tilde{P}_{2}(p_{1},...,p_{6})\rightarrow P_{2}%
\]
}$E_{i}$ denotes the total transforms of $p_{i}$ {\normalsize and }$E\sim
A-2B${\normalsize \ the exceptional divisor on }$P_{2}.$ We
{\normalsize consider the complete linear system $\left\vert L\right\vert
=\left\vert 5A-%
{\textstyle\sum\nolimits_{i=1}^{6}}
2E_{i}\right\vert $ of strict transforms of curves passing through the points
$p_{1},...,p_{6}$ with multiplicity $2.$ Due to Corollary \ref{ALSF2} we know
that $\left\vert L\right\vert $ is very ample on }$S\backslash(%
{\textstyle\bigcup\nolimits_{i=1}^{5}}
E_{i}\cup E)$ {\normalsize and therefore, passing through an arbritrary point
$q_{\lambda}$ with multiplicity two gives three independent conditions for
$\left\vert L\right\vert .$ Thus the dimension $d_{\lambda}:=\dim\left\vert
5A-%
{\textstyle\sum\nolimits_{i=1}^{6}}
2E_{i}-2E_{q_{\lambda}}\right\vert $ for the complete linear system of curves
passing through $p_{1},...,p_{6},q_{\lambda}$ with multiplicity 2 is independent
from $\lambda\in\mathbb{A}^{1}.$ Let $\tilde{S}$ be the blowup of
$\mathbb{A}^{1}\times S$ along the subscheme }$\left\{(\lambda,q_{\lambda}):\lambda \in \mathbb{A}^{1}\right\}\subset \mathbb{A}^{1}\times S,$ then we
obtain {\normalsize the following picture:%
\[%
\begin{array}
[c]{ccccc}%
S_{\lambda} & \subset & \tilde{S} & \rightarrow & \mathbb{A}_{1}\times S\\
& \searrow & \text{ }\downarrow\pi & \swarrow & \\
&  & \mathbb{A}_{1} &  &
\end{array}
\]
The fibres of }$\tilde{S}$ {\normalsize are given by the blowups $S_{\lambda}%
$ with exceptional divisor $E_{\lambda}$. For the sheaf $\mathcal{F}%
:=\mathcal{O}_{\tilde{S}}(5A-\sum_{i=1}^{6}2E_{i}-2E_{\lambda})$ on $\tilde{S}$, we
know that $\dim_{\Bbbk(\lambda)}H^{0}(S_{\lambda},\mathcal{F}_{\lambda})$ is
constant on $\mathbb{A}^{1}.$ It follows that $\pi_{\ast}\mathcal{F}$ is a
vector bundle of dimension $d=d_{\lambda}$ on $\mathbb{A}^{1}.$ Then we can choose a
global section of $\pi_{\ast}\mathcal{F}$, that lifts to a one parameter family
of elements }$C_{\lambda}^{\prime}\sim|${\normalsize $5A-%
{\textstyle\sum\nolimits_{i=1}^{6}}
2E_{i}-2q_{\lambda}|$ on $S.$}
\end{proof}

\bigskip

According to the above approach we can successively separate infinitely
near double points. We remark that in the situation where all double points
$p_{i}$ are lying on a rational curve $\Sigma$ (especially in the case where we
have three infinitesimal near $g_{5}^{1\prime}s$) this can be done in such a
way that all double points which have been separated lie again on $\Sigma.$
This shows \textbf{1.}

\bigskip

To handle step \textbf{2. }{\normalsize we have to show that we can choose
}$F_{0}$ in such a way that {\normalsize $p_{i}$ is a point of multiplicity 2
of }$\Gamma=V(F_{0}).$ {\normalsize As $p_{i}$ is a double point of the complete
intersection }$C^{\prime}=${\normalsize $\Gamma$}$\cap${\normalsize $Y$, we
must have%
\[
F_{0}(p_{i})=Q_{0}(p_{i})=0\text{ \ and \ }dF_{0}(p_{i})=\lambda_{i}%
dQ_{0}(p_{i})
\]
with elements $\lambda_{i}\in\Bbbk$. Now we replace $F_{0}$ by $F_{G}%
:=F_{0}+G\cdot Q_{0}\in I_{C^{\prime}}$ where $G$ is a homogeneous, cubic
polynomial. If it is possible to choose $G$ in such a way that $F_{G}%
(p_{i})=0$ and $dF_{G}(p_{i})=dF_{0}(p_{i})+dG(p_{i})\cdot Q_{0}(p_{i}%
)+dQ_{0}(p_{i})\cdot G(p_{i})=dF_{0}(p_{i})+dQ_{0}(p_{i})\cdot G(p_{i})=0,$
then $F_{G}$ has the desired property. This is fulfilled if $G(p_{i}%
)=-\lambda_{i}$ for $i=1,...,7.$ If three of the double points are colinear,
projection from one of them would lead to a plane model for $C$ of degree $8$
that has a point with multiplicity of at least four, hence we would get a
$g_{4}^{1}$ or a special linear series with even lower Clifford index, so that
we can ommit this case. Then we can find hyperplanes, such that each of them contains exactly two of the points $p_{i}$ and no two of them intersect in one of the points $p_{i}$. For the product $G_{j}$ of their defining equations  
we obtain $G_{j}\in\Bbbk\lbrack x_{0},...,x_{3}]_{3},$ $j=1,...,7,$ with
$G_{j}(p_{i})=\delta_{ij}\cdot G_{j}(p_{j})$ where $G_{j}(p_{j})\not =0.$ Then
for $G:=%
{\textstyle\sum\nolimits_{j=1}^{7}}
-\frac{\lambda_{j}}{G_{j}(p_{j})}G_{j}$ we get $G(p_{i}):=%
{\textstyle\sum\nolimits_{j=1}^{7}}
-\frac{\lambda_{j}}{G_{j}(p_{j})}G_{j}(p_{i})=-\lambda_{i}.$}

\bigskip

{\normalsize According to above considerations we can assume that $C^{\prime}$
is given by the complete intersection of the cone $Y=V(Q_{0})$ and a
hypersurface $\Gamma=V(F_{0})$ of degree $5$ that passes through the points
$p_{i}$ with multiplicity 2. As mentioned in \textbf{3.}} {\normalsize we want
to deform $C_{0}=C$ by fixing $\Gamma$ as well as the points $p_{i}$ and
varying the quadratic form $Q_{\lambda}$ in such a way that $Q_{\lambda}$ defines a
smooth, irreducible hypersurface of degree 2 in $\mathbb{P}^{3}$ for $\lambda\neq 0,$ thus $Y=V(Q_{\lambda})\cong\mathbb{P}^{1}%
\times\mathbb{P}^{1}$. The only condition on $Y_{\lambda}:=V(Q_{\lambda})$ is that it has to
contain the points $p_{i}$. From $\dim\Bbbk\lbrack x_{0},...,x_{3}]_{2}%
=\binom{3+2}{2}=10$, it follows the existence of three linear independent quadratic
forms $Q^{\prime},Q^{\prime\prime},Q^{\prime\prime\prime}\in\Bbbk\lbrack
x_{0},...,x_{3}]_{2}$ that vanish at all points $p_{i}.$ Further we can assume
that $Q^{\prime\prime\prime}=Q_{0}.$ We denote $Y_{1}:=V(Q^{\prime})$ and $Y_{2}:=V(Q^{\prime\prime})$. On the blowup $S:=\tilde{P_{2}}(p_{1},...,p_{7})$ they give two effective divisors $\Gamma_{1}$ and $\Gamma_{2}$ of type $2A-\sum_{i=1}^{7}E_{i}$ and the pencil spanned by these two divisors cut out the linear system $|K_{C}-2D|=|(3A-2B-\sum_{i=1}^{7}E_{i})|_{C}|=|(2A-\sum_{i=1}^{7}E_{i})_{C}|$ as $E|_{C}\sim (A-2B)|_{C}\sim 0$. Then it follows     
\\
\begin{lemma} \label{component}
The complete linear system $|2A-\sum_{i=1}^{7}E_{i}|$ has a rational curve $\Sigma\sim A+B-\sum_{i=1}^{7}E_{i}$ as base locus if and only if $m_{|D|}=3$. For $m_{|D|}=2$ this system has exactly one base point $q\in C$ if $C$ admits a further $g_{5}^{1}$ otherwise its base locus is a point $q\in S\subset C$ or a divisor of type $A+B-\sum_{i\in \Delta}E_{i}$ with $|\Delta|=5$.   
\end{lemma}
\begin{proof}
{\normalsize For $m_{\left\vert D\right\vert }=3$ we get $1=h^{0}%
(C,K_{C}-3D)=h^{0}(C,\mathcal{O}_{C}((3A-3B)|_{C}-%
{\textstyle\sum\nolimits_{i=1}^{7}}
E_{i}|_{C}))=h^{0}(C,\mathcal{O}_{C}((2A-B)|_{C}-%
{\textstyle\sum\nolimits_{i=1}^{7}}
E_{i}|_{C}))$ since }$(A-2B)|_{C}\sim0.$ {\normalsize Hence, it follows the
existence of an element $\gamma\in H^{0}(C,\mathcal{O}_{C}((2A-B)|_{C}-%
{\textstyle\sum\nolimits_{i=1}^{7}}
E_{i}|_{C})).$ Denoting the global generators of $H^{0}(S,\mathcal{O}_{S}(B))$
by $s$ and $t,$ we can assume that the rulings }$R_{s},$ $R_{t}\sim
R$ {\normalsize given by }$s$ and $t$ do not pass through any point $p_{i}.$
{\normalsize The complete linear system%
\[
H^{0}(C,\mathcal{O}_{C}(2A|_{C}-%
{\textstyle\sum\nolimits_{i=1}^{7}}
E_{i}|_{C}))=\left\langle s|_{C}\cdot\gamma,t|_{C}\cdot\gamma\right\rangle
\]
is cut out by a pencil of effective divisors of type $2A-%
{\textstyle\sum\nolimits_{i=1}^{7}}
E_{i}$, thus \ we have $q_{1}|_{C}=s|_{C}\cdot\gamma$ and $q_{2}|_{C}%
=t|_{C}\cdot\gamma$ with $\left\langle q_{1},q_{2}\right\rangle =H^{0}%
(S,\mathcal{O}_{S}(2A-%
{\textstyle\sum\nolimits_{i=1}^{7}}
E_{i})).$ Since the rulings }$R_{s}$ and $R_{t}$ {\normalsize cut out the same
divisor on }$C$ {\normalsize as $\Gamma_{1}=(q_{1})$ and $\Gamma_{2}=(q_{2})$
respectively they must have }$5$ in common. We have $R_{s}.\Gamma_{1}=R_{t}%
.\Gamma_{2}=2$, that is only possible if $R_{s}$ and $R_{t}$ are components of
$\Gamma_{1}$ and $\Gamma_{2}$ respectively. Hence we get {\normalsize a common component
$\Sigma\sim2A-B-%
{\textstyle\sum\nolimits_{i=1}^{7}}
E_{i}$ of }$\Gamma_{1}$ and $\Gamma_{2}$ {\normalsize on $S.$ The intersection product
of }$\Sigma$ {\normalsize with the exceptional divisor $E\sim A-2B$ is
negative, hence $E$ is a component of $\Sigma$. Substracting $E$ once we get a
rational curve $\Sigma\sim A+B-%
{\textstyle\sum\nolimits_{i=1}^{7}}
E_{i}.$ For the other direction let us assume the existence of such an
effective divisor }$\Sigma$ on $S.$ Then $h^{0}(C,\mathcal{O}_{C}%
(K_{C}-3D))=h^{0}(C,\mathcal{O}_{C}((3A-3B)|_{C}-%
{\textstyle\sum\nolimits_{i=1}^{7}}
E_{i}|_{C}))=h^{0}(C,\mathcal{O}_{C}((A+B)|_{C}-%
{\textstyle\sum\nolimits_{i=1}^{7}}
E_{i}|_{C}))=1$ from which $m_{|D|}=3$ follows.

It remains to consider the case $m_{|D|}=2:$ If two divisors in $|2A-%
{\textstyle\sum\nolimits_{i=1}^{7}}
E_{i}|$ have no common component they intersect in exactly $(2A-%
{\textstyle\sum\nolimits_{i=1}^{7}}
E_{i})^{2}=8-7=1$ point $q\in S$. In the case where $C$ admits no further
$g_{5}^{1}$ the point $q$ cannot lie on $C$ as otherwise we would get a $g_{5}%
^{1}=|(2A-%
{\textstyle\sum\nolimits_{i=1}^{7}}
E_{i})|_{C}-q|$ that is different from $|D|$ because of $0=h^{0}%
(C,\mathcal{O}_{C}(K_{C}-3D))=h^{0}(C,\mathcal{O}_{C}((2A-B)|_{C}-%
{\textstyle\sum\nolimits_{i=1}^{7}}
E_{i}|_{C})).$ If $C$ has a further $g_{5}^{1}$ then according to our results
in Theorem \ref{2g15} this linear series is given by $|K_{C}-2D-q^{\prime
}|=|(2A-B)|_{C}-q^{\prime}|$ with a $q^{\prime}\in C$, but then $q^{\prime
}=q.$

We consider the situation where $|2A-%
{\textstyle\sum\nolimits_{i=1}^{7}}
E_{i}|$ has a common divisor $\Sigma\sim aA+bB-\sum_{i=1}^{7}\lambda_{i}E_{i}%
$, $a,b,\lambda_{i}\in\mathbb{Z}$, in its base locus. Let $\Gamma_{1},\Gamma_{2}\in|2A-%
{\textstyle\sum\nolimits_{i=1}^{7}}
E_{i}|$ be two different divisors. Then we can assume that $\Sigma$ is maximal
in the sense that $\Gamma_{1}-\Sigma$ and $\Gamma_{2}-\Sigma$ contain no further
component, i.e.
\\\\
\noindent (I) $\ 0\leq(\Gamma_{1}-\Sigma).(\Gamma_{2}-\Sigma)=2(2-a)(2-a-b)-\sum_{i=1}%
^{7}(1-\lambda_{i})^{2}$

$\Leftrightarrow\sum_{i=1}^{7}2\lambda_{i}\geq7-2(2-a)(2-a-b)+\sum_{i=1}%
^{7}\lambda_{i}^{2}$
\\\\
\noindent Moreover $|2A-%
{\textstyle\sum\nolimits_{i=1}^{7}}
E_{i}|$ cuts out on $C$ a linear system of dimension $2$ and degree
\\\\
\noindent (II) $d=C.((2-a)A-bB-\sum_{i=1}^{7}(1-\lambda_{i})E_{i})=10(2-a)-5b-\sum
_{i=1}^{7}2(1-\lambda_{i})\geq5$

$\Leftrightarrow\Sigma.C=10a+5b-\sum_{i=1}^{7}2\lambda_{i}\leq1$ \\\\ and as $C$ is
irreducible \\
\[
0\leq\Sigma.C=10a+5b-\sum_{i=1}^{7}2\lambda_{i}\leq1
\]
\\
\noindent If some of the $\lambda_{i}$ are negative then one of the
exceptional divisors $E_{i}$ is contained in the base locus of our complete
linear series,  hence we must have $\Sigma.C\geq E_{i}.C=2,$ a contradiction.
From (I) it follows that at least $7-2(2-a)(2-a-b)$ coefficiants $\lambda_{i}$
have to be equal to one. Therefore we get%
\[
0\leq10a+5b-\sum_{i=1}^{7}2\lambda_{i}\leq10a+5b-14+4(2-a)(2-a-b)
\]
We distinguish the cases $a=0,1$ and $2:$
\\\\
\noindent $a=0:$ Here we have $2-3b\geq0$, hence $b=0$ and therefore $\Sigma\sim0.$
\\\\
\noindent $a=1:$ In this case it follows that $b\geq0$ and from (I) that $b\leq1$. Then
either $\Sigma\sim A-\sum_{i\in\Delta}E_{i}$ for $\Delta\subset\{1,...,7\}$
and $|\Delta|=5$ or $\Sigma\sim A+B-\sum_{i=1}^{7}E_{i}$. The last possibility
can be excluded since then it follows that $m_{|D|}=3.$
\\\\
\noindent $a=2:$ From (I) we see that $\lambda_{i}=1$ for all $i=1,...,7$ and from (II):
$0\leq\Sigma.C=6+5b\leq1,$ hence $b=-1.$ Then we get the existence of an
effective divisor $\Sigma$ of type $2A-B-\sum_{i=1}^{7}E_{i}$ and thus
$m_{|D|}=3$ as above.
In the case where $\Sigma\sim A-\sum_{i\in\Delta}E_{i}$, we assume $\Sigma\sim A-\sum_{i=1}^{5}E_{i}$, the complete linear system cut out on $C$ by the pencil of divisors spanned by $\Gamma_{1}-\Sigma$ and $\Gamma_{2}-\Sigma\sim A-E_{6}-E_{7}$ cannot be of type $g_{5}^{1}$ as the divisors $\Gamma_{1}-\Sigma$ and $\Gamma_{2}-\Sigma$ do not intersect and $(\Gamma_{i}-\Sigma).C=6$ for $i=1,2$. If $C$ admits a further $g_{5}^{1}$ then it is cut out by the pencil of divisors spanned by $\Gamma_{1}$ and $\Gamma_{2}$, hence in this situation there exists no common divisor in the base locus of $|2A-\sum_{i=1}^{7}E_{i}|$. In this case, i.e. if $\Sigma\sim 0$, $\Gamma_{1}$ and $\Gamma_{2}$ intersect in exactly one point $q\in S$ and we have $q\in C$ if and only if $C$ has a further $g_{5}^{1}$. 
\end{proof}
\\\\
We have collected enough information to formulate our first result:
\begin{theorem}
{\normalsize Let C be an irreducible, smooth, canonical curve of genus 9 with
$\operatorname*{Cliff}(C)=3$ that admits no $g_{7}^{2}$. If }$C$%
{\normalsize admits exactly one $g_{5}^{1}=$ $\left\vert D\right\vert ,$ having
multiplicity $m_{\left\vert D\right\vert }=2$, then there exists a local one
parameter family $(C_{\lambda})_{\lambda\in\mathbb{A}^{1}}$ with $C_{0}=C$ and
$C_{\lambda}\in\mathcal{M}_{9}$ having exactly two different $g_{5}^{1}$$^{\prime}s$ of
ordinary multiplicity. }
\end{theorem}
\begin{proof}
{\normalsize We consider the image }$C^{\prime}$ {\normalsize on the cone
}$Y\subset\mathbb{P}^{3}$ {\normalsize of }$C$ under {\normalsize the map
obtained from the complete linear system }$\left\vert 2D\right\vert .$ We have
already shown that we can assume for $C^{\prime}$ to have exactly 7 distinct
double points $p_{1},...,p_{7}$ as only singularities. Further $C^{\prime}%
$ {\normalsize is given as complete intersection of the singular quadric
$Y=V(Q_{0})$ and a quintic }$\Gamma$ that passes through the points $p_{i}$ with
multiplicity 2. From the above lemma we know that the quadrics
{\normalsize $Y_{1}=V(Q^{\prime})$ and $Y_{2}=V(Q^{\prime\prime})$ either intersect in a one dimensional component $\Sigma^{\prime}$ or in the points $p_{i}$ and one further point which
we denote by $q.$ In the first case we have shown that on the blowup $S$ the strict transform $\Sigma$ of $\Sigma^{\prime}$ is a divisor of type $A-\sum_{i\in \Delta}E_{i}$ with $|\Delta|=5$. This situation occurs if five double points are lying on a hyperplane. In the second case $C$ does not pass through $p$. 
Then setting $Q_{\lambda}:=Q_{0}+\lambda Q^{\prime}$ and
accordingly $Y_{\lambda}=V(Q_{\lambda})$, this gives us an adequate one
parameter family of smooth quadrics and $C_{\lambda}:=Y_{\lambda}\cap
\Gamma$ the corresponding deformation of $C^{\prime}=C_{0}:$ $Y_{\lambda}$
is an irreducible, smooth quadric in $\mathbb{P}^{3}$ which can be considered
as an image of $\mathbb{P}^{1}\times\mathbb{P}^{1}$ via a Segre embedding in
$\mathbb{P}^{3}$ (cf. \cite{harris}\textbf{ }page 285)$.$ The quintic $\Gamma$
cuts out a divisor of class $(5,5)$ on $Y_{\lambda}\cong\mathbb{P}%
^{1}\times\mathbb{P}^{1},$ hence $C_{\lambda}$ is a divisor of class $(5,5)$
on $\mathbb{P}^{1}\times\mathbb{P}^{1}$ having double points $p_{i}$ as only
singularities. Now $C_{\lambda}$ has two distinct special linear series of
type $g_{5}^{1}$ obtained from projection from each factor of $\mathbb{P}%
^{1}\times\mathbb{P}^{1}.$ }

{\normalsize From Theorem \ref{2g15} we already know that $C_{\lambda}$ has a
third $g_{5}^{1}$ if and only if it is cut out by the pencil of quadrics passing through the points $p_{i}$. In the case where $Y,Y_{1}$ and $Y_{2}$ intersect in the common component $\Sigma^{\prime}$ the additional $g_{5}^{1}$ has to be cut out on $C_{\lambda}\subset\mathbb{P}^{1}\times \mathbb{P}^{1}$ by the pencil of divisors $|(1,1)-p_{6}-p_{7}|$ which is not possible as this systems cuts out a $g_{6}^{1}$ on $C_{\lambda}$ and two divisors of this linear system do exactly intersect in the points $p_{1}$ and $p_{2}$.
In the situation where $Y,Y_{1}$ and $Y_{2}$ intersect in exactly the points $p_{1},...,p_{7}$ and $q$, we must have $q\in C_{\lambda}$ if $C_{\lambda}$ has a third $g_{5}^{1}$. Thus it follows that $C_{\lambda}$ has exactly two different $g_{5}^{1}$$^{\prime}s.$  
}
\end{proof}

\bigskip
{\normalsize Now let us consider the case where $C$ has an additional
$g_{5}^{1}=\left\vert D^{\ast}\right\vert ,$ $D^{\ast}\nsim D,$ $m_{\left\vert
D^{\ast}\right\vert }=1$. Then the three quadrics $Y=V(Q_{0}),$ $Y
_{1}:=V(Q^{\prime})$ and $Y_{2}:=V(Q^{\prime\prime})$ intersect in the
points }$p_{i}$ {\normalsize  and one further point $q$, the base point of the
complete linear system $\left\vert 2H-p_{1}-...-p_{7}\right\vert .$ According
to Theorem \ref{2g15b} we have $q\in C^{\prime}$ and on the corresponding
$\mathbb{P}^{1}$-bundle $P_{2}$, the linear system $\left\vert D^{\ast
}\right\vert $ is cut out by the pencil of divisors $\left\vert
2A-p_{1}-...-p_{7}\right\vert $. With $Q_{\lambda}:=Q_{0}+\lambda Q^{\prime}$
and accordingly $Y_{\lambda}=V(Q_{\lambda})$ we get an adequate one
parameter family of quadrics $Y_{\lambda}\cong\mathbb{P}^{1}\times
\mathbb{P}^{1}$ and $C_{\lambda}:=Y_{\lambda}\cap\Gamma$ a divisor of class
$(5,5)$ on $\mathbb{P}^{1}\times\mathbb{P}^{1}$ the corresponding local
deformation of $C^{\prime}=C_{0}:$ $C_{\lambda}$ has two
$g_{5}^{1}$$^{\prime}s$ from projection along each factor of $\mathbb{P}^{1}%
\times\mathbb{P}^{1}$ and as $q\in C^{\prime}$, a third $g_{5}^{1}$ is cut out
by the pencil of divisors of class $(2,2)$ passing through the points
$p_{1},...,p_{7}.$ We conclude: }

{\normalsize \bigskip}

\begin{theorem}
{\normalsize Let C be an irreducible, smooth, canonical curve of genus 9 with
$\operatorname*{Cliff}(C)=3$ that admits no $g_{7}^{2}$. If }$C$%
 {\normalsize admits exactly one $g_{5}^{1}=$ $\left\vert D\right\vert $ with
multiplicity $m_{\left\vert D\right\vert }=2$ and a further one with
multiplicity one then there exists a local one parameter family $(C_{\lambda
})_{\lambda\in\mathbb{A}^{1}}$ with $C_{0}=C$ and $C_{\lambda}\in
\mathcal{M}_{9}$ having exactly three different $g_{5}^{1}$ of ordinary
multiplicity.}
\end{theorem}

{\normalsize \bigskip}

{\normalsize It remains to examine the most special case where $C$ has a
$g_{5}^{1}$ of multiplicity 3. We will see that the same deformation as above
leads to a local one parameter family $(C_{\lambda})_{\lambda\in\mathbb{A}^{1}}$
of curves that have a $g_{5}^{1}$ with multiplicity two and a distinct one of
ordinary multiplicity. Together with the theorem above, where we have dicussed
this case, it follows: }

{\normalsize \bigskip}

\begin{theorem}
{\normalsize Let C be an irreducible, smooth, canonical curve of genus 9 with
$\operatorname*{Cliff}(C)=3,$ that has exactly one linear system $\left\vert
D\right\vert $ of type $g_{5}^{1}$ and admits no $g_{7}^{2}.$ If
$m_{\left\vert D\right\vert }=3,$ then there exists a local one parameter
family $(C_{\lambda})_{\lambda\in\mathbb{A}^{1}}$ with $C_{0}=C$ and
$C_{\lambda}\in\mathcal{M}_{9}$ having exactly three different, ordinary
$g_{5}^{1}$. }
\end{theorem}

\begin{proof}
In Lemma \ref{component} we have shown that $C$ has a $g_{5}^{1}$ of multiplicity 3 if and only if the three quadrics $Y, Y_{1}$ and $Y_{2}$ intersect in a rational curve $\Sigma^{\prime}\subset Y$, such that its strict transform $\Sigma$ is a divisor of type $A+B-\sum_{i=1}^{7}E_{i}$ on the blowup $S$. As above we define $Q_{\lambda}:=Q_{0}+\lambda Q^{\prime}$ and accordingly $Y_{\lambda}=V(Q_{\lambda})$. Then, on $Y_{\lambda}\cong\mathbb{P}^{1}\times
\mathbb{P}^{1}$ the curve $\Sigma^{\prime}$ is a divisor of type $(2,1)$ that passes through the points $p_{1},...,p_{7}$. It follows from Theorem \ref{2g15b} that $C_{\lambda}:=Y_{\lambda}\cap\Gamma$, a divisor of class
$(5,5)$ on $Y_{\lambda}$, has exactly two different $g_{5}^{1}$$^{\prime}s$ given by the divisors $D\sim(1,0)|_{C_{\lambda}}$ and $D^{*}\sim(0,1)|_{C_{\lambda}}$ with $m_{|D|}=1$ and $m_{|D^{*}|}=2$.
Therefore $C_{\lambda}$ gives an adequate one parameter family with the desired properties.
\end{proof}

{\normalsize \newpage}
°
\newpage
\chapter{{\protect\normalsize Summary}}

The result of this thesis is a complete description of irreducible, smooth, canonical curves of genus 9 concerning their Betti tables and the corresponding Clifford index of $C$. For $\operatorname{Cliff}(C)\leq2$ a collection of all Betti tables that occur was already determined in \cite{schreyer1}. In the case for a hyperelliptic or a tetragonal curve there exists a unique Betti table. 
\\For tetragonal curves the Betti table depends on the existence of a further special linear system of Clifford index $2$ or $3$. Here we can distinguish three cases. 
\\The most interesting situation is that of a pentagonal curve: From Theorem \ref{HRV} we know that for odd genus $g=2l-1$ there exist extra syzygies if and only if $C$ carries a pencil of degree $l$, i.e. especially for $g=9$ this is exactly the case if $C$ has a $g_{5}^{1}$. For curves $C$ with Clifford index 3 there exist $4$ different Betti tables that correspond to curves having a $g_{7}^{2}$ or $k$ $g_{5}^{1}$$^{\prime}s$ (counted with multiplicities) with $k=1,2$ or $3$. We calculated the minimal free resolution for the homogenous coordinate ring $S_{C}$ of $C\subset\mathbb{P}^{8}$ applying the structure theorem in codimension $3$ to $C\subset X$ where $X$ is the scroll constructed from a $g_{5}^{1}$. This calculation even shows that the given collection of Betti tables is correct for arbritrary fields $\Bbbk$ with $\operatorname{char}(\Bbbk)\neq3$. For $\operatorname{char}(\Bbbk)=3$ the following Betti tables are different:       

\vspace{5mm}
\footnotesize
\noindent
\begin{tabular}{|c@{\hspace{0.4mm}}|c@{\hspace{0.4mm}}|c@{\hspace{0.4mm}}|}
\hline
general & $\exists\; g^1_5$ & $\exists \hbox{ two } g^1_5$ \cr \hline
\begin{tabular}{l@{\hspace{1.8mm}}|c@{\hspace{1mm}}c@{\hspace{1mm}}c@{\hspace{1mm}}c@{\hspace{1mm}}c@{\hspace{1mm}}c@{\hspace{1mm}}c@{\hspace{1mm}}c@{\hspace{1mm}}c@{\hspace{1mm}}c@{\hspace{1mm}}c@{\hspace{1mm}}}
& 0 & 1 & 2 & 3 & 4 & 5 & 6 & 7 \cr \hline
0 & 1 & - & - & - & - & - & - & - \cr
1 & - & 21 & 64 & 70 & 4 & - & - & - \cr
2 & - & - & -  & 4  & 70  & 64  & 21 & -\cr
3 & - & - & - & - & - & - & - & 1 \cr
\end{tabular} &
\begin{tabular}{l@{\hspace{1.8mm}}|c@{\hspace{1mm}}c@{\hspace{1mm}}c@{\hspace{1mm}}c@{\hspace{1mm}}c@{\hspace{1mm}}c@{\hspace{1mm}}c@{\hspace{1mm}}c@{\hspace{1mm}}c@{\hspace{1mm}}c@{\hspace{1mm}}c@{\hspace{1mm}}}
& 0 & 1 & 2 & 3 & 4 & 5 & 6 & 7 \cr \hline
0 & 1 & - & - & - & - & - & - & - \cr
1 & - & 21 & 64 & 70 & 6 & - & - & - \cr
2 & - & - & - & 6 & 70 & 64  & 21 & -\cr
3 & - & - & - & - & - & - & - & 1 \cr
\end{tabular} &
\begin{tabular}{l@{\hspace{1.8mm}}|c@{\hspace{1mm}}c@{\hspace{1mm}}c@{\hspace{1mm}}c@{\hspace{1mm}}c@{\hspace{1mm}}c@{\hspace{1mm}}c@{\hspace{1mm}}c@{\hspace{1mm}}c@{\hspace{1mm}}c@{\hspace{1mm}}c@{\hspace{1mm}}}
& 0 & 1 & 2 & 3 & 4 & 5 & 6 & 7 \cr \hline
0 & 1 & - & - & - & - & - & - & - \cr
1 & - & 21 & 64 & 70 & 10 & - & - & - \cr
2 & - & - & -   & 10 & 70 & 64  & 21 & -\cr
3 & - & - & - & - & - & - & - & 1 \cr
\end{tabular} \cr \hline
\end{tabular}
\\\\
\normalsize
\\
\textnormal{In \cite{schreyer3}, due to computational evidence, the author also gives conjectural Betti tables for canonical curve of genus $g=10$ and $g=11$. For $g=10$ the methods of this thesis can be applied to check these Betti tables with the exception of that of a curve with Clifford index 4 since an analogue of Theorem \ref{HRV} for even genus is not proven yet. In the case $g=11$ the method to determine the minimal free resolution of $S_{C}$ from a free resolution on a scroll fails for hexagonal curves as we do not have any structure theorem for varieties of codimension 4. The reader is encouraged to continue this work.}         

\newpage
°
\newpage
\chapter{{\protect\normalsize Appendix}}

The following calculations are done in Macaulay2:
\\\\
\bigskip\textbf{1. }Lemma 4.5.1: We denote by $f_{1},...,f_{5}$ the generators
of $H^{0}(X,\mathcal{O}_{X}(H-R)).$ We have to distinguish four cases:\bigskip

restart

kk=QQ

kk=ZZ/3

R=kk[f\_1..f\_5,f\_6,SkewCommutative=%
$>$%
true]

Mpsi=matrix\{\{0,f\_1,f\_2,f\_3\},\{-f\_1,0,f\_4,f\_5\},\{-f\_2,-f\_4,0,f\_6\},\{-f\_3,-f\_5,-f\_6,0\}\}

\bigskip

For $\psi\sim A$ we have:

betti syz substitute(Mpsi,\{f\_1=%
$>$%
f\_6\})

\bigskip

$\psi\sim B:$

betti syz substitute(Mpsi,\{f\_6=%
$>$%
0\})

\bigskip

$\psi\sim C:$

betti syz substitute(Mpsi,\{f\_1=%
$>$%
0,f\_6=%
$>$%
0\})

\bigskip

$\psi\sim D$

betti syz substitute(Mpsi,\{f\_1=%
$>$%
f\_6,f\_3=%
$>$%
f\_4\})

\bigskip

\noindent\textbf{2. }Theorem 4.5.2: 

\bigskip

Momega=matrix\{\{f\_1,f\_2\},\{f\_3,f\_4\}\}

betti syz Momega

betti syz substitute(Momega,\{f\_2=%
$>$%
f\_3\})

\bigskip

\noindent\textbf{3. }Theorem 4.6.4: We denote $f_{1}=s\varphi_{0}$ , $f_{2}%
=t\varphi_{0}$ , $f_{3}=s\varphi_{1}$ , $f_{4}=t\varphi_{1}$ and
$f_{5}=\varphi_{2}:$

\bigskip

restart

kk=QQ

kk=ZZ/3

R=kk[f\_1..f\_5,f\_6,SkewCommutative=%
$>$%
true]

Dach2=mingens (ideal(f\_1..f\_5))\symbol{94}2

Dach3=mingens (ideal(f\_1..f\_5))\symbol{94}3

M1a=substitute(matrix\{\{0,0,0,0,0,0,0,0,0,0\}\},R)

M1b=f\_5*Dach2

M1c=f\_3*Dach2

M1d=f\_4*Dach2

M1=diff(Dach3,transpose (M1a%
$\vert$%
M1b%
$\vert$%
M1c%
$\vert$%
M1d))

M2a=-f\_5*Dach2

M2b=substitute(matrix\{\{0,0,0,0,0,0,0,0,0,0\}\},R)

M2c=f\_1*Dach2

M2d=f\_2*Dach2

M2=diff(Dach3,transpose (M2a%
$\vert$%
M2b%
$\vert$%
M2c%
$\vert$%
M2d))

M3a=-f\_3*Dach2

M3b=-f\_1*Dach2

M3c=substitute(matrix\{\{0,0,0,0,0,0,0,0,0,0\}\},R)

M3d=substitute(matrix\{\{0,0,0,0,0,0,0,0,0,0\}\},R)

M3=diff(Dach3,transpose (M3a%
$\vert$%
M3b%
$\vert$%
M3c%
$\vert$%
M3d))

M4a=-f\_4*Dach2

M4b=-f\_2*Dach2

M4c=substitute(matrix\{\{0,0,0,0,0,0,0,0,0,0\}\},R)

M4d=substitute(matrix\{\{0,0,0,0,0,0,0,0,0,0\}\},R)

M4=diff(Dach3,transpose (M4a%
$\vert$%
M4b%
$\vert$%
M4c%
$\vert$%
M4d))

betti ker (M1%
$\vert$%
M2%
$\vert$%
M3%
$\vert$%
M4)

\bigskip

\noindent\textbf{4. }Theorem 4.6.5: We denote $f_{1}=s^{2}\varphi_{0}$ , $f_{2}%
=st\varphi_{0}$ , $f_{3}=t^{2}\varphi_{0}$ , $f_{4}=\varphi_{1}=\varphi$ and
$f_{5}=\varphi_{2}:$

\bigskip

Dach2=mingens (ideal(f\_1..f\_5))\symbol{94}2

Dach3=mingens (ideal(f\_1..f\_5))\symbol{94}3

M1a=substitute(matrix\{\{0,0,0,0,0,0,0,0,0,0\}\},R)

M1b=f\_4*Dach2

M1c=f\_2*Dach2

M1d=f\_3*Dach2

M1=diff(Dach3,transpose (M1a%
$\vert$%
M1b%
$\vert$%
M1c%
$\vert$%
M1d))

M2a=-f\_4*Dach2

M2b=substitute(matrix\{\{0,0,0,0,0,0,0,0,0,0\}\},R)

M2c=f\_1*Dach2

M2d=f\_2*Dach2

M2=diff(Dach3,transpose (M2a%
$\vert$%
M2b%
$\vert$%
M2c%
$\vert$%
M2d))

M3a=-f\_2*Dach2

M3b=-f\_1*Dach2

M3c=substitute(matrix\{\{0,0,0,0,0,0,0,0,0,0\}\},R)

M3d=substitute(matrix\{\{0,0,0,0,0,0,0,0,0,0\}\},R)

M3=diff(Dach3,transpose (M3a%
$\vert$%
M3b%
$\vert$%
M3c%
$\vert$%
M3d))

M4a=-f\_3*Dach2

M4b=-f\_2*Dach2

M4c=substitute(matrix\{\{0,0,0,0,0,0,0,0,0,0\}\},R)

M4d=substitute(matrix\{\{0,0,0,0,0,0,0,0,0,0\}\},R)

M4=diff(Dach3,transpose (M4a%
$\vert$%
M4b%
$\vert$%
M4c%
$\vert$%
M4d))

betti ker (M1%
$\vert$%
M2%
$\vert$%
M3%
$\vert$%
M4)\bigskip

\bigskip

\newpage
°
\newpage
\printindex
\newpage
°


\newpage
\begin{thebibliography}{999999}                                                                                           
\bibitem[ACGH85]{arbarello}{\normalsize Arbarello, E. ;\ Cornalba, M. ;
Griffiths, P.A. ; Harris, J.: \textit{Geometry of Algebraic Curves.
}Heidelberg:Springer, 1985 (Grundlehren der mathematischen Wissenschaften 129)
}

\bibitem[BE77]{buchsbaumeisenbud}{\normalsize Buchsbaum, D.A. ; Eisenbud, D.:
Algebra Structures for finite free Resolutions and some Structure Theorems for
Ideals of Codimension 3. In:\ \textit{Amer. J. Math. }99 (1977), Nr. 3, p.
447-485 }

\bibitem[BFS89]{BFS}{Belrametti M., Francia P., Sommese A.J.: On Reider's method and higher order embeddings. In:\ \textit{Duke M.J.} Vol. 58 No.2 (1989), p. 435-439}

\bibitem[E95]{eisenbud1}{\normalsize Eisenbud, D.: \textit{Commutative Algebra
with a View Toward Algebraic Geometry. }Springer, 1995 (Graduate Texts in
Mathematics 150) }

\bibitem[E05]{eisenbud2}{\normalsize Eisenbud, D.: \textit{Geometry of
Syzygies. Heidelberg: Springer }2005 (Graduate Texts in Mathematics
229)\textit{ } }

\bibitem[GL84]{GreenLazarsfeld}{\normalsize Green, M. ; Lazarsfeld, R.: The
non-vanishing of certain Koszul cohomology groups. In: \textit{J. Differential
Geometry }19 (1984), p. 168-170 }

\bibitem[Ha92]{harris}{\normalsize Harris, J.: \textit{Algebraic Geometry}.
Heidelberg:\ Springer 1992 (Graduate Texts in Math. 133) }

\bibitem[Hs77]{hartshorne}{\normalsize Hartshorne, R.: \textit{Algebraic
Geometry}. Heidelberg: Springer, 1977 (Graduate Texts in Math. 52) }

\bibitem[HR98]{hirschowitz}{\normalsize Hirschowitz, A., Ramanan, S.: New
evidence for Green's conjecture on syzygies of canonical curves. \textit{Ann.
Sci. \'{E}cole Norm. Sup.} 31 (1998), p. 145-152 }

\bibitem[L99]{Lossen}{\normalsize Lossen, C.: New asymptotics for the
existence of plane curves with prescribed singularities. In.: \textit{Comm. in
Alg}. 27 (1999), p. 3263-3282 }

\bibitem[M95]{mukai}{\normalsize Mukai, S.: Curves and symmetric spaces. In:
\textit{Am. J. Math. }117 (1995), p. 1627-1644 }

\bibitem[N1880]{noether}{\normalsize Noether, M.: \"{U}ber die invariante
Darstellung algebraischer Funktionen. In: \textit{Math. Ann. }17 (1880), p.
263-284 }

\bibitem[P23]{petri}{\normalsize Petri, K.: \"{U}ber die invariante
Darstellung algebraischer Funktionen einer Ver\"{a}nderlichen. In:
\textit{Math. Ann. }88 (1923), p. 242-289 }

\bibitem[R88]{Reider}{\normalsize Reider, I.: Vector bundles of rank 2 and
linear systems on algebraic surfaces. In.: \textit{Ann. of Math. }(2) 127
(1988), no. 2, p. 309-316 }

\bibitem[S-D73]{saint-donat}{\normalsize Saint-Donat, B.: On Petri's analysis
of the linear system of quadrics through a canonical curve. In: \textit{Math.
Ann. }206 (1973), p. 157-175 }

\bibitem[S86]{schreyer1}{\normalsize Schreyer, F.-O.: Syzygies of Canonical
Curves and Special Linear Series. In: \textit{Math. Ann. }275 (1986), p.
105-137 }

\bibitem[S91]{schreyer2}{\normalsize Schreyer, F.-O.:\ A standard basis
approach to syzygies of canonical curves. In.: \textit{J. Reine Angew. Math.
}421 (1991), p. 83-123 }

\bibitem[S03]{schreyer3}{\normalsize Schreyer, F.-O.: Some Topics in
Computational Algebraic Geometry. In.: \textit{Advances in Algebra and
Geometry}. Hindustan Book Agency (2003), p. 263-278. }

\bibitem[V88]{voisin1}{\normalsize Voisin, C.: Courbes t\'{e}tragonales et
cohomologie de Koszul. In: \textit{Reine Angew. Math. }387 (1988), p. 111-121
}

\bibitem[V01]{voisin2}{\normalsize Voisin, C.: Green's generic syzygy
conjecture for curves of even genus lying on a K3 surface. Preprint, Paris,
France, 2001 }

\bibitem[V05]{voisin3}{\normalsize Voisin, C.: Green's Canonical syzygy
conjecture for curves of odd genus. In.: \textit{Compositio Math. (2005)} in
press }

\bibitem[W]{weinfurtner}{Weinfurtner R.: Singul$\ddot{\text{a}}$re Bordigafl$\ddot{\text{a}}$chen, Dissertation Universit$\ddot{\text{a}}$t Bayreuth (1991).}
\end{thebibliography}
\end{document}